%% file: barton-mayboroda.tex
\documentclass[draft]{memo-l}

\usepackage{amssymb,esint,stmaryrd}

\usepackage[pdfborder={0 0 0},draft=false]{hyperref}

\newtheorem{thm}[equation]{Theorem}
\newtheorem{lem}[equation]{Lemma}
\newtheorem{cor}[equation]{Corollary}
\newtheorem{prp}[equation]{Proposition}

\theoremstyle{definition}
\newtheorem{dfn}[equation]{Definition}
\theoremstyle{remark}
\newtheorem{rmk}[equation]{Remark}

\numberwithin{section}{chapter}
\numberwithin{equation}{chapter}
\numberwithin{figure}{chapter}

\newcommand\abs[2][empty]{\csname#1\endcsname \lvert{#2}\csname#1\endcsname\rvert}
\newcommand\doublebar[2][empty]{\csname#1\endcsname \lVert{#2}\csname#1\endcsname\rVert}

\newcommand\dist{\mathop{\mathrm{dist}}}
\newcommand\Div{\mathop{\mathrm{div}}}
\newcommand\Tr{\mathop{\mathrm{Tr}}}
\newcommand\Ext{\mathop{\mathrm{Ext}}}
\newcommand\supp{\mathop{\mathrm{supp}}}
\newcommand\diam{\mathop{\mathrm{diam}}}
\newcommand\esssup{\mathop{\mathrm{ess\,sup}}}
\newcommand\re{\mathop{\mathrm{Re}}}
\newcommand\im{\mathop{\mathrm{Im}}}
\newcommand\RR{\mathbb{R}} \let\R\RR
 \let\C\CC
 \let\Z\ZZ
\newcommand\NN{\mathbb{N}} 
\newcommand\1{\mathbf{1}}
\newcommand\e{\vec{e}}
\newcommand\D{\mathcal{D}}
\newcommand\s{\mathcal{S}}

\newcommand\F{\vec{F}}
\newcommand\G{\mathcal{G}}
\newcommand\T{\mathcal{T}}

\usepackage{tikz}
\usetikzlibrary{intersections}

\tikzset{ill-posed/.style={color=white!85!black,draw=black}}
\tikzset{bounded/.style={color=white!60!black}}
\tikzset{boundary bounded/.style={color=white!50!black, ultra thick}}
\tikzset{well-posed/.style={color=white!30!black}}
\tikzset{boundary well-posed/.style={color=white!20!black, ultra thick}}
\newcommand\drawunitsquare{
	\fill [white!90!black] (0,1)--(1,1)--(1,0)--(0,0)--cycle;
	}

\def\figureaxes{%
	\drawunitsquare
	\node [left] at (0,1) {$1$};
	\node [below] at (1,0) {$1$};
	\plainfigureaxes
	}
\def\plainfigureaxes{
	\draw[line width=0.5pt, ->] (-0.3,0)--(1.3,0) node [below] {$\theta$};
	\draw[line width=0.5pt, ->] (0,-0.3)--(0,1.3) node [left] {$1/p$};
}
\def\boundedhexagon{
	\fill [bounded] (0,0)--
	(\alph,0)--
	(1,1/2-\eps)--
	(1,1+\alph/\enn)--
	(1-\alph,1)--
	(0,1/2+\eps)--
	cycle;
	\draw [boundary bounded] (0,0)--(\alph,0);
	}
\def\eps{0.1}
\def\alph{0.4}
\def\alphasharp{0.25}
\def\enn{3}
\def\pnought{0.6}
\def\pone{0.7}

\begin{document}

\frontmatter

\title[Layer Potentials and BVPs in Besov Spaces]{Layer Potentials and Boundary-Value Problems for Second Order Elliptic Operators with Data in Besov Spaces}


\author{Ariel Barton}
\address{Ariel Barton, 202 Math Sciences Bldg., University of Missouri, Columbia, MO 65211}
\email{bartonae@missouri.edu}

\author{Svitlana Mayboroda}
\address{Svitlana Mayboroda, School of Mathematics, University of Minnesota, 127 Vincent Hall, 206 Church St.\ SE, Minneapolis, Minnesota 55455}
\email{svitlana@math.umn.edu}

\date{}




\begin{abstract} This monograph presents a comprehensive treatment of second order divergence form elliptic operators with bounded measurable $t$-independent coefficients in spaces of fractional smoothness, in Besov and weighted $L^p$ classes. We establish:
\begin{enumerate}
\item Mapping properties for the double and single layer potentials, as well as the Newton potential;
\item Extrapolation-type solvability results: the fact that solvability of the Dirichlet or Neumann boundary value problem at any given $L^p$ space automatically assures their solvability in an extended range of  Besov spaces;
\item Well-posedness for the non-homogeneous boundary value problems. 
\end{enumerate}  
In particular, we prove well-posedness of the non-homogeneous Dirichlet problem with data in Besov spaces for operators with real, not necessarily symmetric coefficients.
\end{abstract}

\maketitle

\tableofcontents


\addtocounter{tocdepth}{1}

\mainmatter

\input sec-1-introduction

\input sec-2-dfn
\input sec-3-main

\input sec-4-background

\input sec-5-bounded

\input sec-6-trace
\input sec-7-sobolev

\input sec-8-green

\input sec-9-invertible

\input sec-10-besov


\backmatter

\bibliographystyle{amsalpha}\bibliography{barton-mayboroda}

\end{document}

%% file: sec-1-introduction.tex
\chapter{Introduction}

In this monograph we will discuss boundary-value problems and layer potential operators associated to the elliptic differential operator
$-\Div A\nabla u$. Here $A$ is an elliptic matrix; that is, there are some numbers $\Lambda\geq\lambda>0$ such that, if $\eta$, $\xi\in\C^{n+1}$ and if $x\in\R^n$, $t\in\R$, then
\begin{equation*}\lambda\abs{\eta}^2\leq \re \overline\eta\cdot A(x,t)\eta, \qquad \abs{\xi\cdot A(x,t)\eta}\leq \Lambda\abs{\xi}\abs{\eta}.\end{equation*}
Specifically, we will be concerned with the Dirichlet problem
\begin{equation}\label{eqn:dirichlet:introduction}
\Div A\nabla u=0 \text{ in }\R^{n+1}_+,\qquad\quad \Tr u=f \text{ on } \R^n\end{equation}
and the Neumann problem
\begin{equation}\label{eqn:neumann:introduction}\Div A\nabla u=0 \text{ in }\R^{n+1}_+,\qquad \nu\cdot A\nabla u=f \text{ on } \R^n.\end{equation}
Here we identify $\R^n$ with $\partial\R^{n+1}_+$. 

Stimulated, in part, by the celebrated resolution of the Kato conjecture in \cite{AusHLMT02}, recent years have witnessed a surge of activity devoted to the problems \eqref{eqn:dirichlet:introduction} and~\eqref{eqn:neumann:introduction} with data in $L^p$ (and Sobolev, $\dot W^{p}_1$), spaces. The present monograph concentrates on boundary data in ``intermediate'' spaces $B_\theta^{p,p}(\RR^n)$, $0<\theta<1$,  and establishes well-posedness of the corresponding boundary-value problems, and associated properties of layer potentials. An important new aspect is a comprehensive treatment of the non-homogeneous boundary-value problems, which have not been addressed in this context before, for any type of boundary data.

Unless otherwise specified we will assume that the coefficients $A(x,t)=A(x)$ are independent of the $n+1$st coordinate, often called the $t$-coordinate. This is a natural starting point in this context. First of all, it is known that some smoothness in the transversal direction to the boundary is necessary, for otherwise the corresponding elliptic measure may be mutually singular with respect to the Lebesgue measure; see \cite{CafFK81}. 
Furthermore, such coefficients are suggested by considering a change of variables that straightens the boundary of a Lipschitz domain. We will return to this point momentarily.

\section{History of the problem: $L^p$ setting}

Recent results have brought a complete understanding of the Dirichlet boundary problem in $L^p$ for elliptic operators with real (possibly non-symmetric) $t$-independent coefficients. It has been established in \cite{KenKPT00, HofKMP12} (see \cite{JerK81a} for the symmetric case) that given any such matrix $A$ there exists a $p$ such that the solutions~$u$ to the Dirichlet problem~\eqref{eqn:dirichlet:introduction} with $L^p$ boundary data exist and satisfy the bound $N_+u\in L^p(\R^n)$, where the nontangential maximal function $N_+ u$ is given by $N_+ u(x)=\sup_{\abs{x-y}<t} \abs{u(y,t)}$. Moreover, such solutions are unique among functions that satisfy $N_+u\in L^p(\R^n)$. The result is sharp, in the sense that for any given $p<\infty$ there exists an elliptic operator with real non-symmetric coefficients such that well-posedness in $L^p$ is violated; see \cite{KenKPT00}. 
 
These advances were followed by results for the Dirichlet problem with boundary data $f$ whose \emph{gradient} lies in $L^p(\R^n)$, that is, data $f$ in the Sobolev space $\dot W^{p}_1(\R^n)$. We will refer to this as the $\dot W^p_1$-Dirichlet problem; in the literature it is usually called the $L^p$-regularity problem. In this case we expect solutions~$u$ to satisfy $\widetilde N_+(\nabla u)\in L^p(\R^n)$, where $\widetilde N_+$ is the modified nontangential maximal operator introduced in \cite{KenP93}. Existence and uniqueness of solutions for $1<p\leq 2$ was established in \cite{KenP93} for real symmetric coefficients, and for $p>1$ small enough in \cite{KenR09, Rul07,HofKMP13} for real nonsymmetric coefficients. The proof exploited  a certain duality between the $L^{p'}$-Dirichlet and $\dot W^p_1$-Dirichlet problems,  $1/p+1/p'=1$; see also \cite{AusM} and \cite[Proposition~2.52]{AusAH08}.

The Neumann problem for real non-symmetric coefficients remains beyond reach, at least in higher dimensions. Well-posedness with boundary data in~$L^p(\R^n)$ and solutions $u$ satisfying $\widetilde N_+(\nabla u)\in L^p(\R^n)$ is only known in the case of real symmetric coefficients (see \cite{KenP93}) or real nonsymmetric coefficients in two dimensions (see \cite{KenR09,Rul07}). In addition, inspired by the aforementioned duality between the Dirichlet and regularity problems, one can consider  Neumann problems with boundary data in the negative Sobolev space $\dot W^p_{-1}(\R^n)$, the dual of $\dot W^{p'}_1(\R^{n})$, in which case $N_+ u\in L^p(\R^n)$; see~\cite{AusM}.

Finally, some (still narrowly specialized) results are available in the case of complex coefficients. It is well known that the $L^2$-Dirichlet, $L^2$-Neumann, and $\dot W^2_1$-Dirichlet problems are well-posed for constant coefficients. These problems are also well-posed for complex self-adjoint coefficients~$A$ (see \cite{AusAM10}); this generalizes the corresponding results of \cite{JerK81a,KenP93} for the real symmetric case. Furthermore, the resolution of the Kato problem \cite{AusHLMT02} allows us to treat such problems for complex $t$-independent matrices in a block form; see \cite[Remark~2.5.6]{Ken94}, as well as \cite{AusAH08,May10a}. 

One of the most recent advances in this direction, and an important part of the background of the present monograph, is a certain self-improvement, or extrapolation, property for the well-posedness results. 
In the case of real coefficients, interpolation with the maximum principle allows us to pass from the $L^{p_0}$-Dirichlet problem to the $L^p$-Dirichlet problem for any $p>p_0$. By the aforementioned duality considerations this also means that well-posedness of the $\dot W^{p_0}_1$-Dirichlet problem automatically extends to all $1<p<p_0$. In the much more general context of complex coefficient $t$-independent operators whose solutions satisfy a certain version of the De Giorgi-Nash-Moser (H\"older continuity) bounds, such an  ``extrapolation" was established in \cite{AusM}, for  the $W^p_{-1}$, $L^p$-Neumann and $L^p$, $\dot W^p_1$-Dirichlet problems. 

Most of the results above extend, as appropriate, to the Hardy spaces when $p<1$ (in the $\dot W^p_1$-Dirichlet and $L^p$-Neumann case) and H\"older spaces (in the $L^p$-Dirichlet case).  
Let us also point out for the sake of completeness that a number of results concerning boundary-value problems under perturbation of the coefficients~$A$ are also known;  see \cite{FabJK84, AusAH08, AusAM10, AlfAAHK11, Bar13} in the case of $L^\infty$ perturbation, and \cite{Dah86a, Fef89, FefKP91, Fef93, KenP93, KenP95, DinPP07, DinR10, AusA11, AusR12, HofMayMou} in the case of $t$-dependent Carleson measure perturbations.

\begin{figure}
\begin{center}
\includegraphics[draft=false]{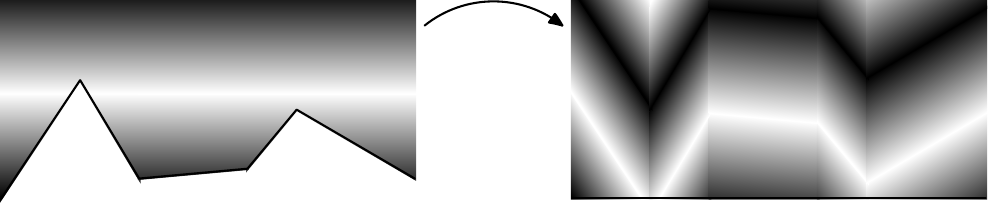}
\caption{If $u$ is harmonic, then $\tilde u(x,t)=u(x,t+\varphi(x))$ satisfies $\Div \widetilde A\nabla \tilde u=0$ for some real symmetric $t$-independent matrix~$\widetilde A$.}
\label{fig:change-of-variables}
\end{center}
\end{figure}

To clarify the role of $t$-independent coefficients, consider a bilipschitz change of variables~$\rho$. If $u$ is harmonic in some domain~$\Omega$, or more generally if $\Div A\nabla u=0$ in~$\Omega$, then $\tilde u=u\circ\rho^{-1}$ satisfies $\Div \widetilde A \nabla\tilde u=0$ in~$\rho(\Omega)$, where $(\det J_\rho) \widetilde A\circ\rho=J_\rho\, A\, J_\rho^T$ for $J_\rho$ the Jacobean matrix. In particular, let $\Omega=\{(x,t):t>\varphi(x)\}$ be the domain above a Lipschitz graph and let $\rho(x,t)=(x,t-\varphi(x))$. See  Figure~\ref{fig:change-of-variables}. For this choice of~$\rho$, if $A$ is elliptic, $t$-independent, real, or symmetric, then so is~$\widetilde A$. 
Notice that this change of variables transforms $\Omega$ to the upper half-space, and so the results of the present monograph (to be proven only in the half-space), \emph{may immediately be generalized to Lipschitz domains.}

A detailed statement of our main theorems requires an extensive notation and terminology discussion and will be presented in Chapter~\ref{chap:main}. In this introduction, we will only outline the principal results with figures elucidating their general scope. 

\section{The nature of the problem and our main results}


The principal goal of this monograph is a comprehensive treatment of boundary-value problems with boundary data in {intermediate} smoothness spaces; that is, we will consider problems that are in some sense between the $L^p$-Dirichlet and $\dot W^p_1$-Dirichlet problems, or between the $\dot W^p_{-1}$-Neumann and $L^p$-Neumann problems. Specifically, we will consider boundary data in the Besov spaces $\dot B^{p,p}_\theta(\R^n)$, for $0<p\leq \infty$ and $0<\theta<1$ (the Dirichlet problem) or $-1<\theta<0$ (the Neumann problem). The parameter $\theta$ measures smoothness; the spaces $L^p(\R^n)$, $\dot W^p_1(\R^n)$ and $\dot W^p_{-1}(\R^n)$ have the same order of smoothness as the spaces $\dot B^{p,p}_0(\R^n)$, $\dot B^{p,p}_1(\R^n)$ and $\dot B^{p,p}_{-1}(\R^n)$, and if $p=2$ then they are in fact the same spaces.

Such boundary-value problems are thoroughly understood in the case of the Laplacian, that is, for the equation $\Delta u=\Div \nabla u=0$, in a Lipschitz domain; see \cite{JerK95,FabMM98,Zan00,May05,MayMit04a}. Both problems have been investigated in the case of (possibly $t$-dependent) $C^1$ coefficients in \cite{Agr07,Agr09}, and the Dirichlet problem has been investigated in \cite{MazMS10} in the case of  coefficients~$A$ having vanishing mean oscillation. \emph{In this monograph, we impose no regularity or oscillation control on the coefficients of the underlying operator.} We will only assume that the coefficient matrix $A$ is elliptic, $t$-independent and satisfies the De Giorgi-Nash-Moser condition, that is, that solutions to $\Div A\nabla u=0$ and $\Div A^*\nabla u=0$ are locally H\"older continuous. The De Giorgi-Nash-Moser condition is always true, e.g., for real coefficients. 
We remark that the results mentioned above have heavily employed Calder\'on-Zygmund theory and local regularity of solutions, which are not available in our case of rough coefficients.

An important consideration in formulating boundary value problems is the nature of the sharp estimates on the solutions. For example, recall that if $\Tr u$ lies in $L^p(\R^n)$ then we expect solutions $u$ to satisfy $N_+ u\in L^p(\R^n)$. If $\Tr u$ lies in the Besov space $\dot B^{p,p}_\theta(\R^n)$, what bounds should we expect $u$ to satisfy?

In \cite{JerK95,FabMM98,Zan00,May05,MayMit04a}, three different types of estimates were sought for harmonic functions $u$ with $\Tr u\in \dot B^{p,p}_\theta(\partial\Omega)$ or $\nu\cdot\nabla u\in \dot B^{p,p}_{\theta-1}(\partial\Omega)$, $0<\theta<1$. Solutions $u$ were expected to lie in the Besov spaces $\dot B^{p,p}_{\theta+1/p}(\Omega)$, the Bessel potential spaces $L^p_{\theta+1/p}(\Omega)$, or to satisfy the estimate
\begin{equation}\label{eqn:space:harmonic}
\int_{\Omega} \abs{\nabla u(X)}^p\,\dist(X,\partial\Omega)^{p-1-p\theta}\,dX<\infty.
\end{equation}
The requirement $u\in \dot B^{p,p}_{\theta+1/p}(\Omega)$ is suggested by the well-known trace and extension theorems for the Besov spaces, that is, by theorems that state that the trace operator $\Tr$ is bounded $\dot B^{p,p}_{\theta+1/p}(\Omega)\mapsto \dot B^{p,p}_{\theta}(\partial\Omega)$ and that there is a bounded extension operator $\Ext: \dot B^{p,p}_{\theta}(\partial\Omega)\mapsto \dot B^{p,p}_{\theta+1/p}(\Omega)$ such that $\Tr \Ext f=f$. See, for example, \cite[Section~2.7.2]{Tri83} in the case of the half-space, and \cite[Chapter~V]{JonW84} in more general domains.
For harmonic functions, the bound \eqref{eqn:space:harmonic} was shown to be equivalent to the requirement $u\in \dot B^{p,p}_{\theta+1/p}(\Omega)$, and (if $1<p<\infty$) to the requirement $u\in L^p_{\theta+1/p}(\Omega)$.

For our purposes, if $\theta+1/p>1$ then the estimate $u\in \dot B^{p,p}_{\theta+1/p}(\R^{n+1}_+)$ is unreasonable. Recall that harmonic functions are smooth, and so locally lie in $\dot B^{p,p}_{\theta+1/p}$ for any $p$ and~$\theta$. However, solutions to $\Div A\nabla u=0$, for general $t$-independent matrices $A$, are \emph{not} smooth. For example, consider a harmonic function $u$ after a bilipschitz change of variables, such as that illustrated in Figure~\ref{fig:change-of-variables}. Such a function is a solution to an elliptic equation but is not smooth; its gradient $\nabla u$ is discontinuous on sets of codimension~$1$. Thus, if $\theta+1/p>1$ then we do not expect solutions to lie in $\dot B^{p,p}_{\theta+1/p}(\R^{n+1}_+)$. (If $\theta+1/p<1$ then we do; see Chapter~\ref{chap:besov}. We will not consider the Bessel potential spaces~$L^p_{\theta+1/p}$.)

Let us turn to estimate~\eqref{eqn:space:harmonic}. We remark that it was this bound that was used in \cite{MazMS10} to formulate the Dirichlet problem for coefficients~$A$ in~$VMO$. This estimate is also a natural choice coinciding, for $p=2$, with the classical square function bounds satisfied by solutions (see, e.g., \cite{DahJK84,DahKPV97,KenKPT00,AusAH08,AlfAAHK11,AusAM10,DinKP11,HofKMP12}) and supported by the long and celebrated history of boundary value problems in domains with isolated singularities, traditionally stated in weighted $L^p$ and Sobolev spaces \cite{KonO83,Gri85,KozMR01}. 

However, we intend to study rougher coefficients~$A$, and thus the bound~\eqref{eqn:space:harmonic} requires one final modification. For general coefficients~$A$, the best we may expect of solutions $u$ to $\Div A\nabla u=0$ is that their gradients are locally in $L^{2+\varepsilon}$ for some (possibly small) $\varepsilon>0$. See {\cite[Theorem~2]{Mey63}}, reproduced as Lemma~\ref{lem:PDE2} below.
However, we wish to consider boundary data $f\in \dot B^{p,p}_\theta(\R^n)$ for $p$ potentially very large, and thus the requirement~\eqref{eqn:space:harmonic} that $\nabla u$ be locally in $L^p(\R^{n+1}_+)$ is again too strong. We will seek solutions $u$ that instead satisfy the averaged bound
\begin{equation}\label{eqn:space:average}
\int_{\R^{n+1}_+} \biggl(\fint_{B((x,t),t/2)}\abs{\nabla u}^2\biggr)^{p/2}t^{p-1-p\theta}\,dx\,dt<\infty.
\end{equation}
The idea of taking $L^2$ averages of the gradient over Whitney balls $B((x,t),t/2)$ is not a new one; such averages are used in defining the modified nontangential maximal function $\widetilde N_+(\nabla u)$ of \cite{KenP93}, and for much the same reason. We will be able to prove appropriate trace theorems for functions $u$ satisfying this estimate; see Chapter~\ref{chap:trace}.


\subsection{Main result I: Well-posedness of the homogeneous problems in Besov spaces and extrapolation}
\label{sec:introduction:homogeneous}
Taken all together, we will study solutions to the Dirichlet problem
\begin{equation}
\label{eqn:dirichlet:besov}
\left\{
\begin{gathered}
\Div A\nabla u=0 \text{ in }\R^{n+1}_+,
\qquad
\Tr u = f\in \dot B^{p,p}_\theta(\R^n),
\\
\int_{\R^{n+1}_+} \biggl(\fint_{B((x,t),t/2)}\abs{\nabla u}^2\biggr)^{p/2}t^{p-1-p\theta}\,dx\,dt<\infty
\end{gathered}
\right.
\end{equation}
and the Neumann problem
\begin{equation}
\label{eqn:neumann:besov}
\left\{
\begin{gathered}
\Div A\nabla u=0 \text{ in }\R^{n+1}_+,
\qquad
\nu\cdot A\nabla u = f\in \dot B^{p,p}_{\theta-1}(\R^n),
\\
\int_{\R^{n+1}_+} \biggl(\fint_{B((x,t),t/2)}\abs{\nabla u}^2\biggr)^{p/2}t^{p-1-p\theta}\,dx\,dt<\infty
.\end{gathered}
\right.
\end{equation}
\emph{One of our main results concerning these problems (and, even more importantly, their non-homogeneous analogues discussed below) are the following ``extrapolation'' theorems. }

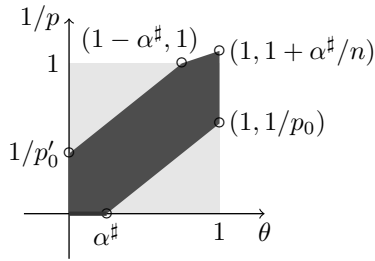
\begin{figure}
\begin{center}
\begin{tikzpicture}[scale=2]
\figureaxes
\fill[well-posed] 
	(0,1-\pnought) node [black] {$\circ$} node [black, left] {$1/p_0'$}--
	(0,0)--
	(\alphasharp,0) node [black] {$\circ$} node [black, below] {$\alpha^\sharp$} --
	(1,\pnought) node [black] {$\circ$} node [black, right] {$(1,1/p_0)$}--
	(1,1) --
	(1,1+\alphasharp/\enn) node [black] {$\circ$} node [black, right] {$(1,1+\alpha^\sharp/n)$}--
	(1-\alphasharp,1) node [black] {$\circ$} node [black, above left, at = {(1-0.3*\alphasharp,1)}] {$(1-\alpha^\sharp,1)$}--
	cycle;
\draw[boundary well-posed] (0,0)--(\alphasharp,0);
\end{tikzpicture}
\caption{Well-posedness of the Dirichlet problem: if $(D)^A_{p_0,0}$ and $(D)^{A^*}_{p_0,0}$ are well-posed, then under certain assumptions on~$A$, we have that $(D)^A_{p,\theta}$ is well-posed for the values of $(\theta,1/p)$ shown here.}\label{fig:well-posed:intro}
\end{center}
\end{figure}

Given a $p_0$ with $1<p_0\leq 2$ such that either the $L^{p_0'}$-Dirichlet or the $L^{p_0}$-Neumann problem is well-posed, for both $A$ and~$A^*$, we have that if $A$ is elliptic, $t$-independent and satisfies certain De Giorgi-Nash-Moser bounds (see Corollary~\ref{cor:well-posed:AusM}), the Dirichlet problem~\eqref{eqn:dirichlet:besov} or the Neumann problem~\eqref{eqn:neumann:besov} is well-posed whenever the numbers $p$ and $\theta$ are such that the point $(\theta,1/p)$ lies in the hexagonal region shown in Figure~\ref{fig:well-posed:intro}. In other words, \emph{well posedness at one single point $(0,1/p_0')$ (or at the dual one), for the Dirichlet or for Neumann problem, automatically implies well-posedness in the entire region depicted in Figure~\ref{fig:well-posed:intro}.}

In particular, given any operator with real, possibly not symmetric, $t$-inde\-pen\-dent coefficients, there exists a $p_0<\infty$, such that the corresponding Dirichlet boundary problem is well-posed in the entire region in Figure~\ref{fig:well-posed:intro}.

For any \emph{given} $p$ and $\theta$ with $0<\theta<1$ and $0\leq 1/p<1+\theta/n$, and with $(\theta,1/p)$ not lying on the diagonal $\theta=1/p$, one can construct a counterexample, that is, a real $t$-independent matrix $A$ such that the Dirichlet problem~\eqref{eqn:dirichlet:besov} is ill-posed. (For some values of $\theta$, $p$, these boundary-value problems will be ill-posed in the sense that solutions fail to exist; for other values these problems will be ill-posed in the sense that solutions fail to be unique.) See Section~\ref{sec:sharp}.

We remark that the underlying special De Giorgi-Nash-Moser assumptions on~$A$ are necessarily valid if $A$ is real and elliptic, or if $A$ is elliptic and the ambient dimension $n+1=2$. Under weaker assumptions, the  appropriate results are detailed in Corollaries~\ref{cor:well-posed:open} and~\ref{cor:well-posed:1}.

\subsection{Main result II: Inhomogeneous problem and extrapolation}
\label{sec:introduction:inhomogeneous}

In fact, we establish a considerably stronger result than the one stated in the previous section. We prove that the well-posedness of the \emph{homogeneous} problem in any $L^{p_0}$ (that is, at  one single point $(0,1/p_0')$), for the Dirichlet or for Neumann problem, automatically implies well-posedness of the corresponding \emph{non-homogeneous} problem in the entire region depicted in Figure~\ref{fig:well-posed:intro}. 

Indeed, Besov spaces offer a natural environment for investigation of the \emph{inhomogeneous} Dirichlet boundary value problem,
\begin{equation}
\label{eqn:dirichlet:inhomogeneous}
\Delta u=\Div \F \text{ in }\Omega,\qquad \Tr u=f \text{ on } \partial\Omega,
\end{equation}
and for its Neumann analogue. The appropriate well-posedness results for Poisson's equation in Lipschitz domains were obtained in  \cite{JerK95,FabMM98,Zan00,MayMit04a}; see also \cite{Mit08}. Corresponding results were obtained for operators with $C^1$ and $VMO$ coefficients in \cite{Agr07,Agr09,MazMS10}). Some related results for operators with constant coefficients were established in \cite{MitM11}. 
However, to the authors' knowledge, the present monograph is the first investigation of the inhomogeneous problem for general elliptic operators with no additional smoothness or oscillation restrictions on coefficients.

We will consider the inhomogeneous problems for the equation $\Div A\nabla u=\Div \F$, where $\F$ satisfies the same estimates as $\nabla u$; that is, where 
\begin{equation}\label{eqn:space:average:F}
\int_{\R^{n+1}_+} \biggl(\fint_{B((x,t),t/2)}\abs{\F}^q\biggr)^{p/q}t^{p-1-p\theta}\,dx\,dt<\infty
\end{equation}
for $q$ sufficiently close to $2$. We will show (see Theorem~\ref{thm:inhomogeneous}) that \emph{well-posedness of the inhomogeneous Dirichlet or Neumann problems is equivalent to well-posedness of the homogeneous problems~\eqref{eqn:dirichlet:besov} and~\eqref{eqn:neumann:besov};} under the assumptions discussed above, this implies well-posedness of the inhomogeneous problems whenever the point $(\theta,1/p)$ lies in the hexagonal region in Figure~\ref{fig:well-posed:intro}. In particular, we prove well-posedness of the inhomogeneous Dirichlet problem for all elliptic operators with real non-symmetric coefficients. (We note, in passing, that even in the case of real and symmetric coefficients, our results for both homogeneous and non-homogeneous Dirichlet and Neumann problems in Besov spaces are new). 

Observe that the results for the inhomogeneous problem naturally yield   sharp estimates for the Green potential in the corresponding weighted Sobolev spaces, as well as an array of new estimates for the underlying Green's function. We shall develop this subject in the next publication.

\subsection[Main result III: Mapping properties of layer potentials]{Main result III: Mapping properties of the single and double layer potentials and the Newton potential}

One of the leading methods for constructing solutions to boundary problems, which remains amenable to an extremely rough context of elliptic operators with non-smooth coefficients, is the method of layer potentials. Layer potentials have been employed, in particular, in \cite{DahK87, KenR09, Rul07, Mit08, Agr09, Bar13, BarM13B, HofKMP13, HofMayMou} and, for relatively nice operators in Besov spaces, in \cite{FabMM98,Zan00,MayMit04a,Agr09,Mit08,MitM11}. Quite recently an alternative approach via the functional calculus of first order Dirac-type operators was proven to be equivalent to layer potentials as well; see \cite{Ros12A}. Most closely related to the subject of this monograph are the bounds for layer potentials associated to general elliptic $t$-independent operators in $L^p$, Hardy, and H\"older spaces in \cite{AlfAAHK11,HofMitMor}.  

\begin{figure}[tbp]
\begin{tikzpicture}[scale=2]
\figureaxes
\fill [bounded] (0,0)--
(\alph,0) node [black] {$\circ$} node [black, below] {$(\alpha,0)$}--
(1,1/2-\eps) node [black] {$\circ$} node [black, right] {$(1,1/p^+)$}--
(1,1+\alph/\enn) node [black] {$\circ$} node [black, right] {$(1,1+\alpha/n)$}--
(1-\alph,1) node [black] {$\circ$} node [black, above, at = {(1-0.2-\alph,1)}] {$(1-\alpha,1)$}--
(0,1/2+\eps) node [black] {$\circ$} node [black, left] {$(0,1-1/p^+)$}--
cycle;
\draw [boundary bounded] (0,0)--(\alph,0);
\draw [line width=0.5pt, dotted] (0,1/2)--(1,1/2);
\end{tikzpicture}
\caption{Values of $(\theta,1/p)$ such that the double layer potential is appropriately bounded on $\dot B^{p,p}_\theta(\R^n)$, the single layer potential is bounded on $\dot B^{p,p}_{\theta-1}(\R^n)$, and the Newton potential is bounded on spaces of functions $\F$ that satisfy the bound~\eqref{eqn:space:average:F}. Here $p^+$ is a parameter depending on~$A$ that necessarily satisfies the inequality $1/p^+<1/2$, and $\alpha$ is the De Giorgi-Nash-Moser exponent.}\label{fig:bounded:introduction}
\end{figure}
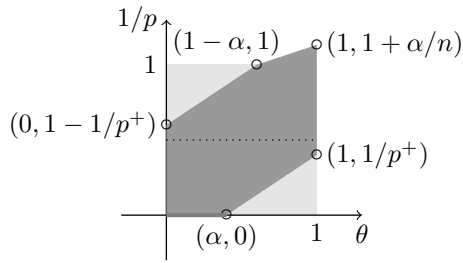

The corresponding mapping properties for layer potentials constitute the technical core of the present work and underpin the new well-posedness results listed in Sections~\ref{sec:introduction:homogeneous} and~\ref{sec:introduction:inhomogeneous}. \emph{We establish sharp bounds on the single and double layer potentials in Besov spaces $\dot B^{p,p}_\theta(\R^n)$, as well as bounds for the Newton potential needed for inhomogeneous problems, in the full range of $(\theta, 1/p)$ depicted in Figure~\ref{fig:bounded:introduction}; see Theorem~\ref{thm:bounded}.} We note that, as usual, the boundedness range exceeds the well-posedness one, for the well-posedness requires invertibility of the boundary potentials and normally introduces further restrictions. Furthermore, our boundedness results automatically give analytic perturbation in the spirit of \cite{AlfAAHK11,HofMitMor}.

As pointed out above, similar boundedness results for Lebesgue and H\"older spaces were established in \cite{HofMitMor} and ``framed" our range (they correspond, in some sense, to the lines $\theta=0$, $\theta=1$, and $p=\infty$ in Figure~\ref{fig:bounded:introduction}). However, our mapping properties will not be a direct consequence of the results for Lebesgue spaces, and in fact cannot be proven using interpolation of known results; see Remark~\ref{rmk:not-interpolation}. 

\section{Outline of the monograph}

The outline of this monograph is as follows. We will establish our terminology in Chapter~\ref{chap:dfn}; having done so, we will state our main results more precisely in Chapter~\ref{chap:main}. In Chapter~\ref{chap:background}, we will review some known results concerning interpolation functors, function spaces, and solutions to elliptic equations. 

We will establish boundedness of the Newton potential, and of the double and single layer potentials acting on fractional smoothness spaces, in Chapter~\ref{chap:bounded}. In Chapter~\ref{chap:trace}, we will prove trace theorems; combined with the results of Chapter~\ref{chap:bounded}, this will establish boundedness of the boundary layer potential operators.

Our results concerning well-posedness of the Dirichlet and Neumann problems with boundary data in Besov spaces will be proven in Chapter~\ref{chap:invertibility}.
Chapters~\ref{chap:sobolev} and~\ref{chap:green} contain important preliminary results.

More precisely, in Chapter~\ref{chap:invertibility}, we will show that invertibility of layer potentials on Besov spaces is equivalent to well-posedness of boundary-value problems, and will use interpolation and functional analysis to prove extrapolation-type results. We will need a Green's formula representation for solutions; we will prove this representation formula in Chapter~\ref{chap:green}. We will want to extrapolate from well-posedness of boundary-value problems with data in Lebesgue or Sobolev spaces, and will need good behavior of layer potentials on such spaces; thus, Chapter~\ref{chap:sobolev} will be devoted to reviewing known results in such spaces.
In the case where the matrix~$A$ has real coefficients, the known results of \cite{JerK81a,KenP93,KenR09,HofKMP12,HofKMP13,AusM} combine with our results to give a particularly complete and satisfactory understanding of the Dirichlet problem in fractional smoothness spaces; we describe these results in Section~\ref{sec:real}.

Finally, in Chapter~\ref{chap:besov}, we will return to the notion of Besov-space estimates $u\in \dot B^{p,p}_{\theta+1/p}(\R^{n+1}_+)$ and show that for appropriate $\theta$ and~$p$, such estimates are equivalent to the estimate~\eqref{eqn:space:average}.

\section*{Acknowledgements} 

Svitlana Mayboroda is partially supported by the Alfred P. Sloan Fellowship, the NSF CAREER Award DMS 1056004, and the NSF Materials Research Science and Engineering Center Seed Grant. We would like to thank Steve Hofmann for making the unpublished work \cite{HofMitMor} available to us.

%% file: sec-2-dfn.tex
\chapter{Definitions} \label{chap:dfn}

In this chapter we define the notation used throughout this monograph.

We work in the upper half-space $\R^{n+1}_+=\R^n\times(0,\infty)$ and the lower half-space $\R^{n+1}_-=\R^n\times(-\infty,0)$.
We identify $\partial\R^{n+1}_\pm$ with $\R^n$. We let $\nu_\pm$ denote the outward unit normal to $\R^{n+1}_\pm$ and let $\nu=\nu_+$. Observe that $\nu_+=-\e_{n+1}$ and that $\nu_-=\e_{n+1}$ where $\e_{n+1}$ is the unit vector in the $n+1$st direction. We will reserve the letter $t$ to denote the $(n+1)$st coordinate in~$\R^{n+1}$. 

We let $B(X,r)$ denote balls in $\R^{n+1}$ and let $\Delta(x,r)$ denote ``surface balls'' on $\partial\R^{n+1}_\pm$, that is, balls in~$\R^n$. 
We will let $\Omega(x,t)$ denote the Whitney ball $B((x,t),t/2)\subset\R^{n+1}_+$.

If $Q\subset\R^n$ or $Q\subset\R^{n+1}$ is a cube, we let $\ell(Q)$ denote its side-length, and let $r Q$ denote the concentric cube with side-length $r\ell(Q)$. If $\mu$ is a measure and $E$ is a set with $\mu(E)<\infty$, we let $\fint$ denote the average integral $\fint_E f\,d\mu=\frac{1}{\mu(E)} \int_E f\,d\mu$.

If $u$ is defined in $U$ for some open set $U\subset\R^{n+1}$, we let $\nabla_\parallel u$ denote the gradient of $u$ in the first $n$ variables, that is, $\nabla_\parallel u=(\partial_1 u,\partial_2 u,\dots,\partial_n u)$. We will also use $\nabla_\parallel$ to denote the \emph{full} gradient of a function defined on~$\R^n$.

We say that $Q\subset\R^n$ is a \emph{dyadic cube} if $\ell(Q)=2^j$ for some integer $-\infty<j<\infty$, and if each vertex $x$ of $Q$ may be written $x=(k_12^j,k_22^j,\dots,k_n2^j)$ for some integers~$k_i$.

We let $\G$ be the grid of \emph{dyadic Whitney cubes} given by
\begin{equation}\label{eqn:whitney:grid}
\G=\{Q\times(\ell(Q),2\ell(Q)):Q\text{ dyadic}\}.
\end{equation}
Then $\overline{\R^{n+1}_+}=\cup_{Q\in\G} \overline{Q}$ and if $Q$, $R\in\G$, then the interiors of $Q$ and~$R$ are disjoint.

\section{Function spaces}
\label{sec:dfn:function}

Let $0<p<\infty$. If $(U,\mu)$ is a measure space and $B$ is a Banach space, we denote the standard Lebesgue space by
\begin{equation*}
L^p(U\mapsto B,d\mu)=\{u:\doublebar{u}_{L^p(U)}<\infty\},
\quad\text{where}\quad \doublebar{u}_{L^p(U)}= \biggl(\int_U \doublebar{u}_B^p\,d\mu\biggr)^{1/p}.\end{equation*}
As usual, we let $\doublebar{u}_{L^\infty(U)}$ be the essential supremum of $\doublebar{u}_B$ in~$U$. If not otherwise specified, $d\mu$ will be the Lebesgue measure and the Banach space $B$ will be the complex numbers $\C$ or the vector space $\C^k$.

If $U\subseteq\R^{k}$ for some~$k$, we denote the homogeneous Sobolev space by
\begin{equation*}\dot W^p_1(U)=\{u:\nabla u\in L^p(U)\}\end{equation*}
with the norm $\doublebar{u}_{\dot W^p_1(U)}=\doublebar{\nabla u}_{L^p(U)}$. (Elements of $\dot W^p_1(U)$ are then defined only up to an additive constant.) If $1< p \leq \infty$ and $1/p+1/p'=1$, we let $\dot W^p_{-1}(\R^n)$ denote the dual space to $\dot W^{p'}_1(\R^n)$. We say that $u\in L^p_{loc}(U)$ or $u\in W^p_{1,loc}(U)$ if $u\in L^p(V)$ or $u\in \dot W^p_1(V)$ for every $V$ compactly contained in~$U$.

If $0<p\leq 1$, it is often more appropriate to consider Hardy spaces rather than Lebesgue or Sobolev spaces. There are several equivalent characterizations; see \cite[Chapter III]{Ste93} for an extensive discussion of these spaces. We will only consider Hardy spaces with $p>n/(n+1)$; we will characterize these spaces by the following atomic decomposition. If $n/(n+1)<p\leq 1$, then let
\begin{equation*}H^p_{at}(\R^n)=\Bigl\{\sum_Q \lambda_Q a_Q: \Bigl(\sum_Q\abs{\lambda_Q}^p\Bigr)^{1/p}<\infty\Bigr\}\end{equation*}
where $a_Q$ is a $H^p(\R^n)$ atom; that is, for some cube~$Q\subset\R^n$,
\begin{itemize}
\item $\supp  a_Q\subseteq Q$,
\item $\int_Q a_Q=0$,
\item $\doublebar{a_Q}_{L^\infty(Q)}\leq \abs{Q}^{-1/p}$.
\end{itemize}
$H^p(\R^n)$ is then the completion of $H^p_{at}(\R^n)$ with respect to the norm 
\[\doublebar{u}_{H^p(\R^n)}=\inf \Bigl\{\Bigl(\sum_Q\abs{\lambda_Q}^p\Bigr)^{1/p}:u=\sum_Q \lambda_Q a_Q
\text{ for some atoms }a_Q
\Bigr\}
.\]

We let $\dot H^p_1(\R^n)$ be the set of functions whose gradients lie in $H^p(\R^n)$; if $p>n/(n+1)$ then a $\dot H^p_1(\R^n)$ atom is a function $a_Q$ that satisfies
\begin{itemize}
\item $\supp  a_Q\subseteq Q$,
\item $\doublebar{\nabla a_Q}_{L^\infty(Q)}\leq \abs{Q}^{-1/p}$
\end{itemize}
for some cube~$Q$.

If $p\leq 1$ then $H^p(\R^n)\subsetneq L^p(\R^n)$. If $1<p<\infty$, then by  \cite[Chapter III]{Ste93}, $H^p(\R^n)=L^p(\R^n)$ and so $\dot H^p_1(\R^n)=\dot W^p_1(\R^n)$; when dealing with a broad range of $p$ we will use this fact to simplify our notation.

If $U$ is a metric space, then the homogeneous H\"older spaces $\dot C^\theta(U)$ are given by 
\begin{equation*}
\dot C^\theta(U)=\{u:\doublebar{u}_{\dot C^\theta(U)}<\infty\},
\quad\text{where}\quad \doublebar{u}_{\dot C^\theta(U)}= \esssup_{x,y\in U, \>x\neq y}\frac{\abs{u(x)-u(y)}}{\abs{x-y}^\theta}.\end{equation*}


If $B$ is a space, we let $B'$ denote the dual space.
If $p$ is an extended real number with $1\leq p\leq \infty$, we will let $p'$ be the extended real number that satisfies $1/p+1/p'=1$; as is well known, if $1\leq p<\infty$ then $(L^p(U))'=L^{p'}(U)$. If $0<p<1$, then we let $p'=\infty$.

The purpose of this monograph is to study boundary-value problems with boundary data in fractional smoothness spaces, that is, in spaces that in some sense are between $L^p(\R^n)$ and $\dot W^p_1(\R^n)$ (or $\dot W^p_{-1}(\R^n)$). Specifically, we will be concerned with the Besov spaces $\dot B^{p,p}_\theta(\R^n)$, to be defined momentarily. In some arguments we will also use the more general Besov and Triebel-Lizorkin spaces $\dot B^{p,r}_\theta(\R^n)$ and~$\dot F^{p,r}_\theta(\R^n)$.

The classical Littlewood-Paley definition of homogeneous
Triebel-Lizorkin and Besov spaces (see, for example, \cite[Section~5.1.3]{Tri83} or 
\cite[Section~2.6]{RunS96}) is as follows. 
Let $\mathcal{F}$ denote the Fourier transform in $\RR^n$.
Let $S$ denote the space of Schwartz functions defined on~$\R^n$ and let \[Z=\{\varphi\in S: \partial^\beta \mathcal{F} \varphi(0)=0 \text{ for every multiindex }\beta\}.\]
Let $\dot \Xi$ be the collection of all systems
$\{\zeta_j\}_{j=-\infty}^{\infty}\subset S$ with the properties
\begin{itemize}
\item $\supp \zeta_j\subset \{x:2^{j-1}<\abs{x}<2^{j+1}\}$,
\item for every multiindex $\beta$ there exists a positive number
$c_{\beta}$ such that 
\[2^{j|\beta|}|\partial^{\beta}\zeta_j(x)| \leq c_{\beta}\]
for all integers $j$ and all $x\in\R^n$,
\item $\sum_{j=-\infty}^{\infty} \zeta_j(x)=1$   for every  $ x\in
\R^n\setminus\{0\}$.
\end{itemize}

Let $\theta\in \RR$ and $0<r \leq \infty$ and fix some family
$\{\zeta_j\}_{j=-\infty}^{\infty} \in \dot\Xi$. If $0<p<\infty$ then the
Triebel-Lizorkin spaces are defined as
\[\dot F_{\theta}^{p,r}(\RR^n):=\bigl\{f\in Z'(\RR^n):
\|f\|_{\dot F^{p,r}_{\theta}(\RR^n)} <\infty \bigr\}\]
and if $0<p\leq\infty$ then the Besov spaces are defined as
\begin{equation*}
\dot B^{p,r}_{\theta}(\RR^n):=\bigl\{f\in Z'(\RR^n):
\|f\|_{\dot B^{p,r}_s(\RR^n)}<\infty \bigr\}
\end{equation*}
where the norms are given by
\begin{gather}
\label{eqn:triebel:norm}
\|f\|_{\dot F^{p,r}_{\theta}(\RR^n)}:= \Bigl\|\Bigl(\sum_{j=-\infty}^{\infty}
|2^{j\theta}{\mathcal{F}}^{-1}(\zeta_j{\mathcal{F}} f)|^r\Bigr)^{1/r}
\Bigr\|_{L^p(\RR^n)}
,\\
\label{eqn:besov:norm}
\|f\|_{\dot B^{p,r}_\theta(\RR^n)}:= \Bigl(\sum_{j=-\infty}^{\infty}
\|2^{j\theta}{\mathcal{F}}^{-1}(\zeta_j {\mathcal{F}}
f)\|_{L^p(\RR^n)}^r\Bigr)^{1/r}
.\end{gather}

A different choice of the system
$\{\zeta_j\}_{j=-\infty}^{\infty}\in\dot\Xi$ yields the same spaces
\eqref{eqn:triebel:norm} and \eqref{eqn:besov:norm}, albeit equipped with equivalent norms. 

We observe that if $p=r$ then $\dot B^{p,p}_\theta(\R^n)=\dot F^{p,p}_\theta(\R^n)$; in this case we will usually denote these spaces as $\dot B^{p,p}_\theta(\R^n)$ rather than~$\dot F^{p,p}_\theta(\R^n)$.

\begin{rmk}
In the literature it is also common to consider the inhomogeneous Besov spaces 
\begin{equation*}\dot B^{p,r}_\theta(\R^n)\cap L^p(\R^n) = \{u:\doublebar{u}_{L^p(\R^n)}+\doublebar{u}_{\dot B^{p,p}_\theta(\R^n)}<\infty\}.\end{equation*}
This space is usually referred to as $B^{p,r}_\theta(\R^n)$. However, we will work almost exclusively with the homogeneous spaces of \eqref{eqn:besov:norm}; to avoid confusion we will consistently use the terminology $\dot B^{p,r}_\theta(\R^n)\cap L^p(\R^n)$ to refer to inhomogeneous spaces. 
\end{rmk}

The main results of this monograph (to be discussed in Chapter~\ref{chap:main}) concern solutions $u$ to partial differential equations in $\R^{n+1}_+$. We will use the following terminology to discuss the behavior of functions in~$\R^{n+1}_+$.

Recall that $\Omega(x,t)=B((x,t),t/2)\subset\R^{n+1}_+$.
We define the space $L(p,\theta,q)$  by
\begin{equation}\label{eqn:weighted:norm}
\doublebar{u}_{L(p,\theta,q)}
= 
\biggl(\int_{\R^{n+1}_+} \biggl(\fint_{\Omega(x,t)} \abs{u}^q \biggr)^{p/q}  t^{p-1-p\theta}\,dt\,dx\biggl)^{1/p}.
\end{equation}
Observe that if $0<a<1$, then we may replace $\Omega(x,t)$ by $B((x,t),at)$ in formula~\eqref{eqn:weighted:norm} and produce an equivalent norm.
We let $\dot W(p,\theta,q)$ be the set of functions~$u$ defined on $\R^{n+1}_+$ with $\nabla u\in L(p,\theta,q)$. Notice that 
\[\doublebar{u}_{L(\infty,\theta,q)} = \sup_{(x,t)\in\R^{n+1}_+} t^{1-\theta}\biggl(\fint_{\Omega(x,t)} \abs{u}^q \biggr)^{1/q}  .\]

We will occasionally use the standard nontangential maximal function~$N$, as well as the modified nontangential maximal function $\widetilde N$ introduced in \cite{KenP93}, to describe functions defined in~$\R^{n+1}_+$. These functions are defined as follows.
If $a>0$ is a constant and $x\in\R^n$, then the \emph{nontangential cone} $\gamma_\pm(x)$ is given by
\begin{align}
\label{eqn:cone}
\gamma_{\pm}(x)=\bigl\{(y,s)\in\R^{n+1}_\pm:\abs{x-y}< \abs{s}\bigr\}.\end{align}

The nontangential maximal function and modified nontangential maximal function are given by 
\begin{align}
\label{eqn:NTM}
 N_{\pm} F(x) &= \sup\bigl\{\abs{F(y,s)}:(y,s)\in \gamma_\pm(x)\bigr\},
\\
\label{eqn:NTM:averaged}
\widetilde N_{\pm} F(x) &= \sup
\biggl\{\biggl(\fint_{B((y,s),\abs{s}/2)} \abs{F}^2\biggr)^{1/2}:
(y,s)\in \gamma_\pm(x)
\biggr\}.
\end{align}
We remark that if $\widetilde N_+ F\in L^{p_0}(\R^n)$, then $F\in L(p,\theta,2)$ whenever $\theta<1$, $0<p\leq\infty$ and $1/p-\theta/n=1/p_0-1/n$; see Theorem~\ref{thm:NTM:embedding}. Thus, spaces of nontangentially bounded functions are a natural $\theta=1$ endpoint for our weighted, averaged Sobolev spaces $L(p,\theta,2)$.

\section{Elliptic equations}
\label{sec:dfn:elliptic}

Suppose that $A:\R^{n+1}\mapsto\C^{(n+1)\times(n+1)}$ is a bounded measurable matrix-valued function defined on $\R^{n+1}$. We let $A^T$ denote the transpose matrix and let $A^*$ denote the adjoint matrix $\overline {A^T}$.
We say that $A$ is \emph{elliptic} if there exist constants $\Lambda>\lambda>0$ such that
\begin{equation}
\label{eqn:elliptic}
\lambda \abs{\eta}^2\leq \re\bar\eta\cdot A(x,t)\eta,
\quad
\abs{\xi\cdot A(x,t)\eta}\leq \Lambda\abs{\eta}\abs{\xi}
\end{equation}
for all $(x,t)\in\R^{n+1}$ and all vectors $\eta$, $\xi\in\C^{n+1}$. We refer to $\lambda$ and $\Lambda$ as the \emph{ellipticity constants of~$A$}.

We say that $A$ is \emph{$t$-independent} if
$A(x,t)=A(x,s)$ for all $x\in\R^n$ and all $s$, $t\in\R$.

If $u\in W^1_{1,loc}(U)$ and $\F\in L^1_{loc}( U\mapsto\C^{n+1})$ for some open set~$U\subset\R^{n+1}$, we say that $\Div A\nabla u=\Div \F$ if
\begin{equation}
\label{eqn:solution}
\int_ U \nabla\varphi\cdot A\nabla u=\int_ U \nabla\varphi\cdot \F
\quad\text{for all }\varphi\in C^\infty_0( U)
.\end{equation}

If $u\in W^1_{1,loc}(\overline U)$ and $\Div A\nabla u=0$ in $U$ for some reasonably well-behaved open set~$ U$, we may define the conormal derivative $\nu\cdot A\nabla u\big\vert_{\partial U}$ weakly by the formula
\begin{equation}\label{eqn:conormal:1}
\int_{\partial U} \varphi\,\nu\cdot A\nabla u\,d\sigma
= \int_ U\nabla\varphi\cdot A\nabla u
\quad\text{for all }\varphi\in C^\infty_0(\R^{n+1}).
\end{equation}
We observe that the value of $\langle \varphi,\nu\cdot A\nabla u\rangle$ given by formula~\ref{eqn:conormal:1} does not depend on the choice of extension of~$\varphi$. 
Recall that if $\Div A\nabla u=0$ in $U$, then $\int_U \nabla\psi\cdot A\nabla u=0$ whenever $\psi$ is smooth and compactly supported in~$U$. If $\nabla u\in L^p_{loc}(\overline U)$ for any $p>1$, then by density this is still true whenever $\psi$ is compactly supported, $\nabla \psi$ is bounded and $\Tr \psi=0$, and so if $\Phi_1$ and $\Phi_2$ are two different smooth extensions of~$\varphi$ then $\int_U \nabla(\Phi_1-\Phi_2)\cdot A\nabla u=0$.

In this monograph we will only consider the conormal derivatives in the case $U=\R^{n+1}_+$ and $U=\R^{n+1}_-$. 
Recall that $\nu_\pm$ is the unit outward normal to $\R^{n+1}_\pm$ and so $\nu_+\cdot A\nabla u\big\vert_{\partial\R^{n+1}_+}$ and $\nu_- \cdot A\nabla u\big\vert_{\partial\R^{n+1}_-}$ are given by formula~\eqref{eqn:conormal:1}. 
We will occasionally use $\nu\cdot A\nabla u$ as shorthand for $\nu_+\cdot A\nabla u\big\vert_{\partial\R^{n+1}_+}$, and will let $\nu_+ \cdot A\nabla u\big\vert_{\partial\R^{n+1}_-}$ denote $-\nu_- \cdot A\nabla u\big\vert_{\partial\R^{n+1}_-}$.

The following theorem concerning solutions to $\Div A\nabla u=\Div \F$ has been known for some time.

\begin{lem}[{\cite[Theorem~2]{Mey63}}]\label{lem:PDE2} 
	Let $A$ be elliptic. Then there exists a number $p^+>2$, depending only on the constants $\lambda$, $\Lambda$ in \eqref{eqn:elliptic} and the dimension~$n+1$, such that if 
	\begin{equation}
	\label{eqn:PDE2:exponents}
	\frac{1}{p^-}+\frac{1}{p^+}=1\quad\text{and}\quad p^-<p<q<p^+,\end{equation}
	then there is a number $C$ (depending only on $p$, $q$, $\lambda$, $\Lambda$ and~$n$) such that if $u$ satisfies $\Div A\nabla u=\Div \F$ in $B(X,2r)$, then
	\begin{equation}
	\label{eqn:PDE2}
	\biggl(\fint_{B(X,r)}\abs{\nabla u}^q\biggr)^{1/q}
	\leq
	C\biggl(\fint_{B(X,2r)}\abs{\nabla u}^p\biggr)^{1/p}
	+
	C\biggl(\fint_{B(X,2r)}\abs{\F}^q\biggr)^{1/q}
	\end{equation}	
	whenever the right-hand side is finite.
\end{lem}

Throughout this monograph we will let $p^+$ and $p^-$ be as in Lemma~\ref{lem:PDE2}. More precisely, if an elliptic matrix $A$ is mentioned in a particular theorem, then $p^+$ is the largest number such that if formula~\eqref{eqn:PDE2:exponents} is valid then the bound~\eqref{eqn:PDE2} is valid for solutions $u$ to either $\Div A\nabla u=0$ or $\Div A^*\nabla u=0$.

\begin{dfn}\label{dfn:DGNM}
If $A$ is an elliptic matrix, we say that \emph{$A$ satisfies the De Giorgi-Nash-Moser condition with exponent~$\alpha$}
if there is some constant~$H$ such that, for every ball $B=B(X_0,2r)$ and every function $u$ such that $\Div A\nabla u=0$ in~$B$, we have that 
\begin{equation}
\label{eqn:holder}
\abs{u(X)-u(X')}\leq H\biggl(\frac{\abs{X-X'}}{r}\biggr)^\alpha
\biggl(\fint_{B(X_0,2r)} \abs{u}^2\biggr)^{1/2}
\end{equation}
for all $X$, $X'\in B(X_0,r)$.
\end{dfn}
We say that $A$ satisfies the De Giorgi-Nash-Moser condition if there is some $\alpha>0$ such that $A$ satisfies the De Giorgi-Nash-Moser condition with exponent~$\alpha$. Notice that if the bound~\eqref{eqn:holder} is true, then it remains true if we replace $\alpha$ by $\beta$ for any $0<\beta<\alpha$.

We will only consider coefficient matrices $A$ such that both $A$ and its adjoint $A^*$ satisfy the De Giorgi-Nash-Moser condition. Throughout we will let $\alpha$ be such that the estimate~\eqref{eqn:holder} is valid for solutions to either $\Div A\nabla u=0$ or $\Div A^*\nabla u=0$. It will sometimes be convenient to let $\alpha_0$ be the supremum of all such~$\alpha$; note that $A$ may or may not satisfy the De Giorgi-Nash-Moser condition with exponent~$\alpha_0$.

All real elliptic coefficients~$A$ satisfy the De Giorgi-Nash-Moser condition; this was proven independently by De Giorgi and Nash in \cite{DeG57,Nas58} in the case of symmetric coefficients and extended to real nonsymmetric coefficients by Morrey in \cite{Mor66}. Notice that if a solution $u$ satisfies the condition~\eqref{eqn:holder}, then $u$ also satisfies Moser's local boundedness estimate
\begin{equation}\label{eqn:local-bound:2}
\abs{u(X)}\leq C\biggl(\fint_{B(X,r)} \abs{u}^2\biggr)^{1/2}
\end{equation}
for some constant $C$ depending only on $H$ and the dimension $n+1$.

We remark that throughout this monograph, we will let $C$ denote a constant whose value may change from line to line, but that depends only on the dimension $n+1$ and the values of the ellipticity constants $\lambda$, $\Lambda$ in formula~\eqref{eqn:elliptic}, the exponents $p^+$ and $p^-$ in Lemma~\ref{lem:PDE2}, and the constants $H$ and~$\alpha$ in formula~\eqref{eqn:holder}. We will refer to these numbers as the standard constants. We will let the symbol $\approx$ indicate that two quantities are comparable, that is, $a\approx b$ if $\frac{1}{C}a\leq b\leq Ca$ for some $C$ depending only on the standard constants.

\section{Layer potentials}

Let $\Gamma_{(x,t)}^A$ be the fundamental solution to $-\Div A\nabla$, that is, $-\Div A\nabla\Gamma_{(x,t)}^A=\delta_{(x,t)}$. It was constructed in \cite{HofK07} for $n+1\geq 3$, and in \cite{Ros12A} in the cases $n+1=2$. Many of the important properties are listed in Section~\ref{sec:fundamental} below. One important fact is the adjoint relation
\begin{equation}
\label{eqn:FS:switch}
\Gamma_{(x,t)}^A(y,s)=\overline{\Gamma_{(y,s)}^{A^*}(x,t)}.
\end{equation}
We will let $\nabla \Gamma_{(x,t)}^A(y,s)$ denote the gradient in the variables~$y$, $s$; we will indicate derivatives with respect to $x$ and~$t$ by subscripts.

Let $f:\R^n\mapsto\C$ be a bounded, compactly supported function. The double layer potential is given by the formula
\begin{equation}\label{eqn:D}
\D f(x,t)=\D^A f(x,t)=\int_{\R^n} \overline{\nu\cdot A^*\nabla\Gamma_{(x,t)}^{A^*}(y,0)}  f(y)\,dy.\end{equation}
The single layer potential is given by 
\begin{equation}\label{eqn:S}
\s f(x,t)=\s^A f(x,t)=\int_{\R^n} \overline{\Gamma_{(x,t)}^{A^*}(y,0)}  f(y)\,dy=\int_{\R^n} \Gamma_{(y,0)}^{A}(x,t) f(y)\,dy.\end{equation}
We will establish boundedness of $\D$ and $\s$ on Besov spaces $\dot B^{p,p}_\theta(\R^n)$ or $\dot B^{p,p}_{\theta-1}(\R^n)$. If $p<\infty$ this will allow us to extend these operators by continuity. If $p=\infty$ then, for appropriate~$\theta$, we may define $\D$ and $\s$ up to an additive constant; see formula~\eqref{eqn:T:difference} and the following discussion.

The Newton potential is traditionally given by 
\[\int_{\R^{n+1}}\Gamma_{(y,s)}^{A}(x,t) F(y,s)\,dy\,ds.\]
However, we will find it more convenient to define
\begin{equation}
\label{eqn:newton}
\Pi\F(x,t)=
\Pi^A \F(x,t)=\int_{\R^{n+1}_+}
\overline{\nabla \Gamma_{(x,t)}^{A^*}(y,s)} \cdot \F(y,s)\,dy\,ds
\end{equation}
for $\F:\R^{n+1}_+\mapsto\C^{n+1}$ bounded and compactly supported; again, we will be able to extend $\Pi$ to larger spaces of functions by continuity or (in the case of, for example, $L^\infty$ functions) by duality.

We remark that in $\R^{n+1}_+$,
\begin{equation}
\label{eqn:Div:plus}
\Div A\nabla \D^A f=\Div A\nabla\s^A f=0,\quad \Div A\nabla \Pi^A\F = \Div \F
\end{equation}
and that in $\R^{n+1}_-$,
\begin{equation}
\label{eqn:Div:minus}
\Div A\nabla \D^A f=\Div A\nabla\s^A f=\Div A\nabla \Pi^A\F = 0
.\end{equation}

We will show that for appropriate $p$, $\theta$ and~$q$, $\D^A:\dot B^{p,p}_\theta(\R^n)\mapsto \dot W(p,\theta,2)$, $\s^A:\dot B^{p,p}_{\theta-1}(\R^n)\mapsto \dot W(p,\theta,2)$, and $\Pi^A:\dot W(p,\theta,q)\mapsto \dot W(p,\theta,q)$ are bounded; see Section~\ref{sec:layers:bounded}. We will also show that the trace operator is bounded $\dot W(p,\theta,q)\mapsto \dot B^{p,p}_\theta(\R^n)$ for such $p$, $\theta$,~$q$. We will define the boundary operators $\D^A_\pm$, $\s^A_\pm$, $\Pi^A_\pm$ by
\[\D^A_\pm f = \Tr \bigl(\D^A f\big\vert_{\R^{n+1}_\pm}\bigr),
\qquad
\s^A_\pm f =  \Tr \bigl(\s^A f\big\vert_{\R^{n+1}_\pm}\bigr),
\qquad
\Pi^A_\pm \F =  \Tr \bigl(\Pi^A \F\big\vert_{\R^{n+1}_\pm}\bigr)\]
whenever $f$ or $\F$ lies in a space such that the right-hand side is a composition of two bounded operators.
It is known (see \cite[Lemma~4.18]{AlfAAHK11}) that $\s^A_+=\s^A_-$. 

Similarly, if $f$ lies in a space such that $\nabla\D f$ or $\nabla \s f$ lies in $L^1_{loc}(\R^{n+1}_\pm)$ (see Section~\ref{sec:layers:bounded} and Theorem~\ref{thm:whitney:embedding}), then we will denote the boundary conormal derivatives by
\[
(\partial_\nu^A\D^A)_\pm f = \nu_\pm\cdot A\nabla \D^A f\big\vert_{\partial\R^{n+1}_\pm},
\quad
(\partial_\nu^A\s^A)_\pm f = \nu_\pm\cdot A\nabla \s^A f\big\vert_{\partial\R^{n+1}_\pm}
\]
where the conormal derivative is in the sense of formula~\eqref{eqn:conormal:1}. In Proposition~\ref{prp:jump} we will see that $(\partial_\nu^A\D^A)_+=(\partial_\nu^A\D^A)_-$. If $\F$ is compactly supported in $\R^{n+1}_+$, then we may let
\[
(\partial_\nu^A\Pi^A)_\pm \F = \nu_\pm\cdot A\nabla \Pi^A \F\big\vert_{\partial\R^{n+1}_\pm}
\]
in the sense of formula~\eqref{eqn:conormal:2}; we will see that we may extend this operator by density or duality to all of $L(p,\theta,q)$ for appropriate $p$, $\theta$,~$q$.

If $\vec f:\R^{n}\mapsto\C^{n+1}$ is a vector-valued function, it is by now standard to define
\begin{equation*}(\s^A\nabla_{\text{full}})\cdot\vec f(x,t)=\int_{\R^n} \overline{\nabla\Gamma_{(x,t)}^{A^*}(y,0)} \cdot \vec f(y)\,dy=\int_{\R^n} \nabla_{y,s}\Gamma_{(y,s)}^{A}(x,t)\bigr\vert_{s=0} \cdot\vec f(y)\,dy.\end{equation*}
We will generally use the special case
\begin{equation*}(\s^A\nabla_\|)\cdot\vec f(x,t)=\int_{\R^n} \overline{\nabla_\parallel\Gamma_{(x,t)}^{A^*}(y,0)} \cdot \vec f(y)\,dy=\int_{\R^n} \nabla_{y}\Gamma_{(y,0)}^{A}(x,t) \cdot\vec f(y)\,dy\end{equation*}
where $\vec f:\R^n\mapsto\C^n$, not $\vec f:\R^n\mapsto\C^{n+1}$. Observe that if $f\in \dot B^{p,p}_{\theta-1}(\R^n)$, then by \cite[Section~5.2.3]{Tri83} we have that $f=\Div \F$ for some $\F\in \dot B^{p,p}_\theta(\R^n \mapsto \C^n)$; thus, to establish (for example) that $\s$ is bounded $\dot B^{p,p}_{\theta-1}(\R^n)\mapsto \dot W(p,\theta,2)$, it suffices to show that $\s\nabla_\|$ is bounded $\dot B^{p,p}_{\theta}(\R^n)\mapsto \dot W(p,\theta,2)$.

Many results concerning $\D$ and $\s\nabla_\|$ can be proven using identical techniques. We will define the operator $\T$ in order to simplify our notation. Let $\vec f:\R^n\mapsto \C^{n+1}$ be a vector-valued function. Then
\begin{align}\label{eqn:T}
\T^A \vec f(x,t) 
&=
	\int_{\R^n} \overline{\bigl(\nabla_\parallel \Gamma_{(x,t)}^{A^*}(y,0), -\nu\cdot A^*\nabla\Gamma_{(x,t)}^{A^*}(y,0)\bigr)} \cdot \vec f(y)\,dy
\\\nonumber
&=
	\int_{\R^n} \overline{\nabla_{A^*} \Gamma_{(x,t)}^{A^*}(y,0)} \cdot \vec f(y)\,dy
\\\nonumber
&=
	(\s\nabla_\parallel)\cdot \vec f_\parallel(x,t) - \D f_\perp(x,t)
.\end{align}
Here the conormal gradient $\nabla_A$, defined by $\nabla_{A} u = \bigl(\nabla_\parallel u, \e_{n+1}\cdot A\nabla u\bigr)$, is as in, e.g., \cite{AusA11}. Also, $\vec f=(\vec f_\parallel, f_\perp)$ where $\vec f_\parallel(x)\in \C^{n}$ and $f_\perp(x)$ is a scalar.
Observe that if $1 < p < \infty$, then boundedness of $\T$ on $L^p(\R^n)$ is equivalent to boundedness of $\s\nabla_{\text{full}}$, but that this need not be true for higher smoothness spaces, or in the Hardy spaces $H^p(\R^n)$ for $p\leq 1$.

\section{Boundary-value problems}

In this section we will define well-posedness of the homogeneous and inhomogeneous Dirichlet and Neumann problems, with boundary data in fractional smoothness spaces or Lebesgue or Sobolev spaces. We will provide an example of ill-posedness in Section~\ref{sec:sharp}; we will also establish some terminology that will let us discuss the nature of ill-posedness.

\begin{dfn}\label{dfn:Dirichlet}
Let $0<\theta<1$ and let $0<p\leq\infty$, $0<q\leq\infty$. 

We say that the Dirichlet problem $(D)^A_{p,\theta,q}$ has the \emph{existence property} in $\R^{n+1}_+$ if, 
for every $f\in \dot B^{p,p}_\theta(\R^n)$, there exists a function $u$ that satisfies
\begin{equation}
\label{eqn:D:existence}
\left\{\begin{aligned} \Div A\nabla u&=0 &&\text{in }\R^{n+1}_+,\\
\Tr u &=f && \text{on } \R^n=\partial\R^{n+1}_+,\\
u & \in \dot W(p,\theta,q)
.\hskip-20pt
\end{aligned}\right.
\end{equation}
We say that $(D)^A_{p,\theta,q}$ is \emph{solvable} in $\R^{n+1}_+$ if it has the existence property and there is some constant $C$ such that, for every $f\in \dot B^{p,p}_\theta(\R^n)$, one of the functions $u$ that satisfies~\eqref{eqn:D:existence} also satisfies $\doublebar{u}_{\dot W(p,\theta,q)}\leq C \doublebar{f}_{\dot B^{p,p}_\theta(\R^n)}$.

We say that $(D)^A_{p,\theta,q}$ has the \emph{uniqueness property} in $\R^{n+1}_+$ if the only solution to the problem~\eqref{eqn:D:existence} with $f\equiv 0$ is the function $u\equiv 0$. Notice that $(D)^A_{p,\theta,q}$ need not have the existence property to have the uniqueness property.

If $(D)^A_{p,\theta,q}$ is solvable and has the uniqueness property, we say that $(D)^A_{p,\theta,q}$ is \emph{well-posed}.

Notice that by Lemma~\ref{lem:PDE2}, if $(D)^A_{p,\theta,q}$ is well-posed for some $q$ with $p^-<q<p^+$, then $(D)^A_{p,\theta,q}$ is well-posed for all such~$q$. We let $(D)^A_{p,\theta}$ denote any of these equivalent boundary-value problems.

We say that the \emph{inhomogeneous version of $(D)^A_{p,\theta,q}$} has the existence property in $\R^{n+1}$ if for every $f\in \dot B^{p,p}_\theta(\R^n)$ and every $\F\in L(p,\theta,q)$, there exists a unique function $u$ that satisfies
\begin{equation} \label{eqn:D:existence:inhomogeneous}
\left\{\begin{aligned} \Div A\nabla u&=\Div\F &&\text{in }\R^{n+1}_+,\\
\Tr u &=f && \text{on } \R^n=\partial\R^{n+1}_+,\\
u&\in \dot W(p,\theta,q)
.
\end{aligned}\right.
\end{equation}
If there is a constant $C$, independent of~$f$ and~$\F$, such that some solution to the problem~\eqref{eqn:D:existence:inhomogeneous} also satisfies
$\doublebar{u}_{\dot W(p,\theta,q)}
\leq C \doublebar{f}_{\dot B^{p,p}_\theta(\R^n)}
+ C\doublebar{\F}_{L(p,\theta,q)}$,
then we say that the inhomogeneous version of $(D)^A_{p,\theta,q}$ is solvable. If $(D)^A_{p,\theta,q}$ has the uniqueness property and the inhomogeneous version is solvable then we say that the inhomogeneous version is well-posed.
\end{dfn}

\begin{rmk}\label{rmk:inhomogeneous}
The existence property or solvability of the inhomogeneous problem implies the corresponding property of the homogeneous problem; simply take $\F\equiv 0$. 
Conversely, notice that if $\F\equiv 0$ and $f\equiv 0$ then the problems~\eqref{eqn:D:existence} and~\eqref{eqn:D:existence:inhomogeneous} are the same; thus, we need only define one uniqueness property.
We will see that for appropriate $p$, $\theta$ and~$q$, the existence property or solvability for the homogeneous problem also imply the corresponding properties for the inhomogeneous problem; see Theorem~\ref{thm:inhomogeneous} below.
\end{rmk}

\begin{rmk} Throughout this section, if we define both solvability and the uniqueness property of a boundary-value problem, we let well-posedness denote the combination of the two properties.
\end{rmk}

We now consider the Neumann problem.

\begin{dfn}
Let $0<\theta<1$ and let $0<p\leq\infty$, $0<q\leq\infty$.  

We say that the Neumann problem $(N)^A_{p,\theta,q}$ has the existence property in $\R^{n+1}$ if, for every $g\in \dot B^{p,p}_{\theta-1}(\R^n)$, there exists a function $u$ that satisfies
\begin{equation}\label{eqn:N:existence}
\left\{\begin{aligned} \Div A\nabla u&=0 &&\text{in }\R^{n+1}_+,\\
\nu\cdot A\nabla u &=g && \text{on } \R^n=\partial\R^{n+1}_+,\\
u&\in \dot W(p,\theta,q)
.
\end{aligned}\right.
\end{equation}
If the only solutions $u$ to the problem~\eqref{eqn:N:existence} with $g\equiv 0$ are $u\equiv c$ for some constant~$c$, then we say that $(N)^A_{p,\theta,q}$ has the uniqueness property.
If there is a constant $C$, independent of~$g$, such that some solution $u$ to the problem~\eqref{eqn:N:existence} satisfies $\doublebar{u}_{\dot W(p,\theta,q)}
\leq C \doublebar{g}_{\dot B^{p,p}_{\theta-1}(\R^n)}$, then we say that $(N)^A_{p,\theta,q}$ is solvable. 
\end{dfn}
Again, if $p^-<q<p^+$ then $(N)^A_{p,\theta,q}$ is equivalent to $(N)^A_{p,\theta,2}$; we denote these equivalent boundary-value problems by $(N)^A_{p,\theta}$.

The following theorem is a straightforward and by now well known consequence of the Lax-Milgram lemma and certain trace and extension theorems relating the space $\dot W^2_1(\R^{n+1}_+)$ and $\dot B^{2,2}_{1/2}(\R^n)$. See, for example, {\cite[Section~3]{AusMM13}}.
\begin{thm}\label{thm:2,1/2}
If $A$ is elliptic, then the problems $(D)^A_{2,1/2}$ and $(N)^A_{2,1/2}$ are well-posed in $\R^{n+1}_+$ and $\R^{n+1}_-$.
\end{thm}
This theorem is true even if $A$ is not $t$-inde\-pen\-dent or fails to satisfy the De Giorgi-Nash-Moser condition.

We would like to define an inhomogeneous Neumann problem. We need some notion of conormal derivative for functions~$u\in \dot W(p,\theta,q)$ with $\Div A\nabla u\neq 0$.

If $\Div A\nabla u=\Div \F$ for some $\F$ compactly supported in~$\Omega$ (and so $\F\equiv 0$ in some neighborhood of~$\partial\Omega$), we may still define $\nu\cdot A\nabla u$ either by restricting the support of $\varphi$ in formula~\eqref{eqn:conormal:1} to that neighborhood, or equivalently by writing
\begin{equation}\label{eqn:conormal:2}
\int_{\partial\Omega} \varphi\,\nu\cdot A\nabla u\,d\sigma
= \int_\Omega\nabla\varphi\cdot A\nabla u
-\int_\Omega \nabla \varphi\cdot \F
\quad\text{for all }\varphi\in C^\infty_0(\R^{n+1}).
\end{equation}
Now, observe that if $p<\infty$ and $q<\infty$, then smooth and compactly supported functions~$\F$ are dense in $L(p,\theta,q)$. In Chapter~\ref{chap:trace} (see in particular Remark~\ref{rmk:newton:trace}), we will see that for appropriate $p$, $\theta$ and~$q$, the Newton potential $\Pi^A$ of formula~\eqref{eqn:newton} extends to a bounded operator $L(p,\theta,q)\mapsto \dot W(p,\theta,q)$, and that $(\partial_\nu^A\Pi^A)_+$ extends to a bounded operator $L(p,\theta,q)\mapsto \dot B^{p,p}_{\theta-1}(\R^n)$.

Thus, if $u\in\dot W(p,\theta,q)$, we may define its conormal derivative by
\begin{equation}\label{eqn:conormal:3}
\nu\cdot A\nabla u\big\vert_{\partial\R^{n+1}_+}=\nu\cdot A\nabla(u-\Pi^A(A\nabla u))\big\vert_{\partial\R^{n+1}_+}+(\partial_\nu^A\Pi^A)_+(A\nabla u)
.\end{equation}
Because $\Div A\nabla (u-\Pi^A(A\nabla u))=0$, we may use 
formula~\eqref{eqn:conormal:1} to define the term $\nu\cdot A\nabla(u-\Pi(A\nabla u))$.

\begin{dfn}
Let $A$ be an elliptic matrix and let $p^-<q<p^+$. Suppose that $0<p\leq \infty$, $0<\theta<1$ are such that $\Pi^A$ extends to a bounded operator $L(p,\theta,q)\mapsto \dot W(p,\theta,q)$, and such that $\nu\cdot A\nabla\Pi^A$ extends to a bounded operator $L(p,\theta,q)\mapsto \dot B^{p,p}_{\theta-1}(\R^n)$.

We say that the inhomogeneous version of the Neumann problem $(N)^A_{p,\theta,q}$ has the existence property in $\R^{n+1}$ if for every $g\in \dot B^{p,p}_{\theta-1}(\R^n)$ and every $\F\in L(p,\theta,q)$, there exists a function $u$ that satisfies
\begin{equation}\label{eqn:N:existence:inhomogeneous}
\left\{\begin{aligned} \Div A\nabla u&=\Div\F &&\text{in }\R^{n+1}_+,\\
\nu\cdot A\nabla u &=g && \text{on } \R^n=\partial\R^{n+1}_+,\\
u&\in \dot W(p,\theta,q)
\end{aligned}\right.
\end{equation}
in the sense of formula~\eqref{eqn:conormal:3}. If there is some constant $C$ such that some solution to the problem~\eqref{eqn:N:existence:inhomogeneous} satisfies $\doublebar{u}_{\dot W(p,\theta,q)}
\leq C \doublebar{g}_{\dot B^{p,p}_{\theta-1}(\R^n)}
+C\doublebar{\F}_{L(p,\theta,q)}$, then we say that the problem is solvable.
\end{dfn}
Again the homogeneous and inhomogeneous uniqueness properties are the same.

There is an extensive literature concerning boundary-value problems with data in integer smoothness spaces, that is, in the spaces $L^p(\R^n)$, $\dot W^p_1(\R^n)$, $\dot W^p_{-1}(\R^n)$, $\dot H^p(\R^n)$, or~$\dot H^p_1(\R^n)$. In order to discuss and use such results, we establish the following terminology.

\begin{dfn}
We say that the Dirichlet problem $(D)^A_{p,0}$ is solvable in $\R^{n+1}_+$ if there is a constant $C$ such that, for every $f\in L^p(\R^n)$, there exists a function $u$ that satisfies
\begin{equation*}
\left\{\begin{aligned} \Div A\nabla u&=0 &&\text{in }\R^{n+1}_+,\\
\Tr u &=f && \text{on } \R^n=\partial\R^{n+1}_+,\\
\doublebar{N_+u}_{L^p(\R^n)}
&\leq C \doublebar{f}_{H^p(\R^n)}
.\hskip-20pt
\end{aligned}\right.
\end{equation*}

We say that the regularity problem $(D)^A_{p,1}$ is solvable in $\R^{n+1}_+$ if there is a constant $C$ such that, for every $f\in\dot H^p_1(\R^n)$, there exists a function $u$ that satisfies
\begin{equation*}
\left\{\begin{aligned} \Div A\nabla u&=0 &&\text{in }\R^{n+1}_+,\\
\Tr u &=f && \text{on } \R^n=\partial\R^{n+1}_+,\\
\doublebar{\widetilde N_+(\nabla u)}_{L^p(\R^n)}
&\leq C \doublebar{\nabla_\parallel f}_{H^p(\R^n)}
.\hskip-20pt
\end{aligned}\right.
\end{equation*}

We say that the Neumann problem $(N)^A_{p,1}$ is solvable in $\R^{n+1}_+$ if there is a constant $C$ such that, for every $f\in H^p(\R^n)$, there exists a function $u$ that satisfies
\begin{equation*}
\left\{\begin{aligned} \Div A\nabla u&=0 &&\text{in }\R^{n+1}_+,\\
\nu\cdot A\nabla u &=f && \text{on } \R^n=\partial\R^{n+1}_+,\\
\doublebar{\widetilde N_+(\nabla u)}_{L^p(\R^n)}
&\leq C \doublebar{f}_{H^p(\R^n)}
.\hskip-20pt
\end{aligned}\right.
\end{equation*}

We define the existence property for each problem by relaxing the requirement in the last line to $N_+ u\in L^p(\R^n)$ or $\widetilde N_+(\nabla u)\in L^p(\R^n)$. We define the uniqueness property by requiring that the only solution with $\Tr u\equiv 0$ or $\nu\cdot A\nabla u\equiv 0$, and $N_+ u\in L^p(\R^n)$ or $\widetilde N_+(\nabla u)\in L^p(\R^n)$ (as appropriate), be $u\equiv 0$.
\end{dfn}

We refer to the regularity problem as $(D)^A_{p,1}$ for consistency with the Dirichlet problems $(D)^A_{p,\theta}$; however, it is somewhat more common to refer to this problem as $(R)^A_p$. In some papers (such as \cite{AlfAAHK11, AusM}), solutions to Dirichlet problem $(D)^A_{p,0}$ are not required to satisfy $N_+u\in L^p(\R^n)$; instead they are required to satisfy the square-function estimate
\begin{equation}\label{eqn:square-function}
\int_{\R^n} \biggl(\int_0^\infty \int_{\abs{x-y}<t} \abs{\nabla u(y,t)}^2 \frac{dy\,dt}{t^{n-1}}\biggr)^{p/2}\,dx
\leq C \doublebar{f}_{H^p(\R^n)}.\end{equation}
In the case of real coefficients these two estimates are equivalent (see \cite{DahJK84,KenKPT00,HofKMP12}). See also the related results in \cite{DahKPV97,AusAH08,AlfAAHK11,AusAM10,DinKP11}. 

We now remark on a certain \emph{compatibility} issue. The paper \cite{Axe10} considers solvability of boundary-value problems for the (two-dimensional) coefficient matrix
\[A_k(x,t)=A_k(x)=\begin{pmatrix}1&k\mathop{\mathrm{sgn}}(x)\\ -k\mathop{\mathrm{sgn}}&1\end{pmatrix}\]
where $k$ is a real number. (Ill-posedness of the Dirichlet, Neumann and regularity problems for this matrix were considered earlier in \cite{KenKPT00} and \cite{KenR09}.) It turns out that for certain values of $k$ and~$p$, the Dirichlet problem $(D)^{A_k}_{p,0}$ is solvable (and, combined with a uniqueness result of \cite{HofMayMou}, is even well-posed) in the sense given above, but that solutions are not compatible with solutions to $(D)^{A_k}_{2,1/2}$. That is, there is some $f\in L^p(\R^1)\cap \dot B^{2,2}_{1/2}(\R^1)$ such that there are two \emph{different} functions $u$ that satisfy
\[\Div A_k\nabla u=0\text{ in }\R^{2}_+,\qquad\Tr u=f,\]
one of which satisfies $N_+u\in L^p(\R^1)$ and the other of which satisfies $\nabla u\in L^2(\R^2_+)$.
Similar statements are valid for the Neumann and regularity problems $(N)^{A_k}_{p,1}$ and~$(D)^{A_k}_{p,1}$.

We intend to apply interpolation methods, and thus we need the solution operators in various spaces to be compatible; that is, we wish to disallow coefficients displaying this type of behavior. We say that the boundary-value problem $(D)^A_{p,\theta}$ is \emph{compatibly solvable} or \emph{compatibly well-posed} if it is solvable (or well-posed) and if the solution operator is compatible with that of $(D)^A_{2,1/2}$; that is, if whenever $f\in \dot B^{2,2}_{1/2}(\R^n)\cap \dot B^{p,p}_\theta(\R^n)$, then the (unique) solution to 
\[\Div A\nabla u=0 \text{ in }\R^{n+1}_+,
\quad \Tr u=f,
\quad \nabla u\in L^2(\R^{n+1}_+)\]
satisfies $u\in \dot W(p,\theta,2)$.
We define compatible solvability or compatible well-posedness for $(D)^A_{p,0}$, $(D)^A_{p,1}$, $(N)^A_{p,\theta}$, or $(N)^A_{p,1}$ similarly. In Remark~\ref{rmk:compatible}, we will see that under certain conditions on $\theta$ and~$p$, compatible solvability is equivalent to compatible well-posedness.

The rough Neumann problem $(N)^A_{p,0}$, as defined in \cite{AusM}, is compatibly solvable if, for every $f\in \dot W^p_{-1}(\R^n)\cap\dot B^{2,2}_{-1/2}(\R^n)$, the solution $u$ to the Neumann problem $(N)^A_{2,1/2}$ with boundary data~$f$ satisfies
\begin{equation*}
\left\{\begin{aligned} \Div A\nabla u&=0 &&\text{in }\R^{n+1}_+,\\
\nu\cdot A\nabla u &=f && \text{on } \R^n=\partial\R^{n+1}_+,\\
\int_{\R^n} \biggl(\int_0^\infty \int_{\abs{x-y}<t} \abs{\nabla u(y,t)}^2 \frac{dy\,dt}{t^{n-1}}\biggr)^{p/2}\,dx
&\leq C\doublebar{f}_{\dot W^p_{-1}(\R^n)}
.\hskip-20pt
\end{aligned}\right.
\end{equation*}
Notice that the space $\dot W^p_{-1}(\R^n)$ is defined only for $p>1$; however, our results and the results of \cite{AusM} do not involve $(N)^A_{p,0}$ or $(D)^A_{p,0}$ in the case $p\leq 1$.

We remark that many papers consider compatible solvability of boundary-value problems; see, for example, \cite{KenR09,AusMM13, AusM}.

%% file: sec-3-main.tex
\chapter{The Main Theorems}
\label{chap:main}

In this chapter we will state the main theorems of this monograph.

\begin{figure}[tbp]
\begin{tikzpicture}[scale=2]
\figureaxes
\fill [bounded] (0,0)--
(\alph,0) node [black] {$\circ$} node [black, below] {$(\alpha,0)$}--
(1,1/2-\eps) node [black] {$\circ$} node [black, right] {$(1,1/p^+)$}--
(1,1+\alph/\enn) node [black] {$\circ$} node [black, right] {$(1,1+\alpha/n)$}--
(1-\alph,1) node [black] {$\circ$} node [black, above, at = {(1-0.2-\alph,1)}] {$(1-\alpha,1)$}--
(0,1/2+\eps) node [black] {$\circ$} node [black, left] {$(0,1/p^-)$}--
cycle;
\draw [boundary bounded] (0,0)--(\alph,0);
\draw [line width=0.5pt, dotted] (0,1/2)--(1,1/2);
\end{tikzpicture}
\caption{Values of $p$, $\theta$ such that one of the inequalities \eqref{eqn:hexagon:1}--\eqref{eqn:hexagon:4} are valid. Here $\alpha$ is the De Giorgi-Nash-Moser exponent, and $p^-$, $p^+$ are the exponents in Lemma~\ref{lem:PDE2}; then $1/p^-+1/p^+=1$. The dotted line indicates $1/p=1/2$.}\label{fig:bounded}
\end{figure}
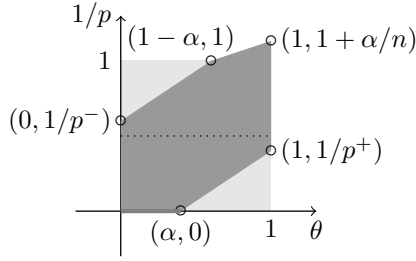

\begin{thm}\label{thm:bounded}
Let $A:\R^{n+1}\mapsto \C^{(n+1)\times(n+1)}$ be elliptic and $t$-independent, such that $A$ and $A^*$ satisfy the De Giorgi-Nash-Moser condition. Let $\alpha$, $p^-$ and~$p^+$ be as in Section~\ref{sec:dfn:elliptic}. We require $n+1\geq 2$.

Suppose that $p^-<q<p^+$ and that $p$, $\theta$ satisfy one of the following pairs of inequalities:
\begin{align}
\label{eqn:hexagon:1}
\frac{n}{n+\alpha}&<p\leq 1, & n(1/p - 1) +(1 -\alpha )&< \theta<1,\\
1&\leq p\leq p^-,& (1-\alpha)\frac{1/p-1/p^-}{1-1/p^-}&<\theta<1,\\
p^-&\leq p\leq p^+,& 0&<\theta<1, \quad\text{or}\\
p^+&\leq p\leq \infty,& 0&<\theta<\alpha+(1-\alpha)p^+/p.
\label{eqn:hexagon:4}
\end{align}
We then have the following bounds.
\begin{align}
\label{eqn:space:D}
\doublebar{\D^A f}_{\dot W(p,\theta,q)}
&\leq C \doublebar{f}_{\dot B^{p,p}_\theta(\R^n)}
,\\
\label{eqn:space:S}
\doublebar{\s^A f}_{\dot W(p,\theta,q)}
&\leq C \doublebar{f}_{\dot B^{p,p}_{\theta-1}(\R^n)}
,\quad\text{and}\\
\label{eqn:space:newton}
\doublebar{\Pi^A\F}_{\dot W(p,\theta,q)}
&\leq C \doublebar{\F}_{L(p,\theta,q)}
.\end{align}
\end{thm}
We will establish the bound \eqref{eqn:space:newton} in Section~\ref{sec:newton:bounded} and will establish the bounds \eqref{eqn:space:D} and~\eqref{eqn:space:S} in Section~\ref{sec:layers:bounded}.
The inequalities~\eqref{eqn:hexagon:1}--\eqref{eqn:hexagon:4} are valid for the values of $\theta$ and $p$ indicated in Figure~\ref{fig:bounded}. We remark that if $\alpha_0$ is as in Section~\ref{sec:dfn:elliptic}, then the theorem remains true if we replace $\alpha$ by~$\alpha_0$. By Lemma~\ref{lem:morrey} below, we have that $\alpha_0\geq 1-n/p^+$; it is simple to establish that if $\alpha\geq 1-n/p^+$ then the region in Figure~\ref{fig:bounded} is convex. The bounds \eqref{eqn:space:D} and~\eqref{eqn:space:S} do not follow from known results concerning the behavior of $\D$ and $\s$ on Lebesgue and Sobolev spaces; see Remark~\ref{rmk:not-interpolation}.

\begin{thm} \label{thm:trace}
Suppose that $n+1\geq 2$, $1\leq q\leq \infty$, $0<\theta<1$, and $0< p\leq \infty$ with $1/p<1+\theta/n$.
Then the trace operator extends to a bounded operator 
\begin{equation*}\Tr:\dot W(p,\theta,q)\mapsto \dot B^{p,p}_\theta(\R^n).\end{equation*}
Furthermore, if $\Div A\nabla u=0$ in~$\R^{n+1}_+$, then the conormal derivative $\nu\cdot A\nabla u\big\vert_{\partial\R^{n+1}_+}$ given by formula~\eqref{eqn:conormal:1} satisfies
\begin{equation*}
\doublebar{\nu\cdot A\nabla u}_{\dot B^{p,p}_{\theta-1}(\R^n)}
\leq C\doublebar{u}_{\dot W(p,\theta,q)}.
\end{equation*}
If in addition $A$, $p$, $q$ and $\theta$ satisfy the conditions of Theorem~\ref{thm:bounded}, then whenever $\F\in L(p,\theta,q)$ is smooth and compactly supported we have that $(\partial_\nu^A\Pi^A)_+\F=\nu\cdot A\nabla\Pi^A\F\big\vert_{\partial\R^{n+1}_+}$ exists in the sense of formula~\eqref{eqn:conormal:2}. We may extend the operator $(\partial_\nu^A\Pi^A)_+$ to an operator that satisfies the bound
\begin{align}
\label{eqn:normal:newton}
\doublebar{(\partial_\nu^A\Pi^A)_+\F}_{\dot B^{p,p}_{\theta-1}(\R^n)}
&\leq C \doublebar{\F}_{L(p,\theta,q)}.
\end{align}
\end{thm}
We will prove Theorem~\ref{thm:trace} in Chapter~\ref{chap:trace}; in that chapter we will also discuss the meaning of these extensions in the case $p=\infty$. Observe that if the conditions of Theorems~\ref{thm:bounded} and~\ref{thm:trace} both hold, we then have the bounds
\begin{align}
\label{eqn:trace:D}
\doublebar{\D_+^A f}_{\dot B^{p,p}_\theta(\R^n)}
&=\doublebar{\Tr\D^A f}_{\dot B^{p,p}_\theta(\R^n)}
\leq C \doublebar{f}_{\dot B^{p,p}_\theta(\R^n)}
,\\
\label{eqn:trace:S}
\doublebar{\s^A_+ f}_{\dot B^{p,p}_\theta(\R^n)}
&=\doublebar{\Tr\s^A f}_{\dot B^{p,p}_\theta(\R^n)}
\leq C \doublebar{f}_{\dot B^{p,p}_{\theta-1}(\R^n)}
,\\
\label{eqn:trace:newton}
\doublebar{\Pi^A_+ \F}_{\dot B^{p,p}_\theta(\R^n)}
&=\doublebar{\Tr\Pi^A \F}_{\dot B^{p,p}_\theta(\R^n)}
\leq C \doublebar{\F}_{L(p,\theta,q)}
\end{align}
and
\begin{align}
\label{eqn:conormal:D}
\doublebar{(\partial_\nu^A\D^A)_+ f}_{\dot B^{p,p}_{\theta-1}(\R^n)}
&=\doublebar{\nu\cdot A\nabla\D^A f}_{\dot B^{p,p}_{\theta-1}(\R^n)}
\leq C \doublebar{f}_{\dot B^{p,p}_\theta(\R^n)},
\\
\label{eqn:conormal:S}
\doublebar{(\partial_\nu^A\s^A)_+ f}_{\dot B^{p,p}_{\theta-1}(\R^n)}
&=\doublebar{\nu\cdot A\nabla\s^A f}_{\dot B^{p,p}_{\theta-1}(\R^n)}
\leq C \doublebar{f}_{\dot B^{p,p}_{\theta-1}(\R^n)}.
\end{align}

The remaining main results of this monograph concern well-posedness of the Dirichlet and Neumann problems. We begin with the homogeneous problems; we will state a series of results that show that, given well-posedness of certain boundary-value problems with data in Lebesgue or Sobolev spaces, we have well-posedness of boundary-value problems with data in Besov spaces as well.

\begin{thm}\label{thm:well-posed:invertible}
Let $A$ be an elliptic, $t$-independent matrix defined on $\R^{n+1}$, $n+1\geq 2$, such that $A$ and $A^*$ satisfy the De Giorgi-Nash-Moser condition. Let $\theta$ and $p$ satisfy the conditions of Theorem~\ref{thm:bounded}.

The boundary-value problem $(D)^A_{p,\theta}$ is well-posed in both $\R^{n+1}_+$ and $\R^{n+1}_-$ if and only if $\s^A_+:\dot B^{p,p}_{\theta-1}(\R^n)\mapsto \dot B^{p,p}_{\theta}(\R^n)$ is invertible with bounded inverse. The solution $u$ to $(D)^A_{p,\theta}$ with boundary data $f$ is given by $u=\s^A((\s^A_+)^{-1} f)$.

Furthermore, $\s^A_+$ is onto if and only if $(D)^A_{p,\theta}$ has the existence property in both $\R^{n+1}_+$ and $\R^{n+1}_-$, and $\s^A_+$ is one-to-one if and only if $(D)^A_{p,\theta}$ has the uniqueness property in both $\R^{n+1}_+$ and~$\R^{n+1}_-$. 

Similarly, $(N)^A_{p,\theta}$ has the existence property, the uniqueness property, or is well-posed in both $\R^{n+1}_+$ and $\R^{n+1}_-$, if and only if $(\partial_\nu^A \D^A)_+:\dot B^{p,p}_{\theta}(\R^n)\mapsto \dot B^{p,p}_{\theta-1}(\R^n)$ is onto, one-to-one, or invertible with bounded inverse, respectively. If $(N)^A_{p,\theta}$ is well-posed then the solution $u$ with boundary data $f$ is given by $u=\D^A((\partial_\nu^A \D^A)_+^{-1} f)$.
\end{thm}
This theorem will be proven in Section~\ref{sec:invertible:well-posed}. We remark that surjectivity of $\s^A_+$ is enough to guarantee existence of solutions $u\in \dot W(p,\theta,2)$ to $(D)^A_{p,\theta}$, but is not enough to give the uniform bound $\doublebar{u}_{\dot W(p,\theta,2)}\leq C(p,\theta)\doublebar{f}_{\dot B^{p,p}_\theta(\R^n)}$; that is, surjectivity gives the existence property but not solvability.

The $\theta=1$ endpoint of Theorem~\ref{thm:well-posed:invertible} is already known (see \cite{Ver84,AlfAAHK11,HofKMP13,HofMitMor,BarM13B}) and may be formulated as follows.
\begin{thm}\label{thm:well-posed:invertible:1}
Let $A$ be an elliptic, $t$-independent matrix defined on $\R^{n+1}$ such that $A$ and $A^*$ satisfy the De Giorgi-Nash-Moser condition. Then there is some $\varepsilon>0$ such that, if $1/2-\varepsilon<1/p<1+\alpha/n$, then the following statements are true.

First, $(D)^A_{p,1}$ has the existence property, the uniqueness property, or is well-posed in both $\R^{n+1}_+$ and $\R^{n+1}_-$, if and only if the operator $\s^A_+:L^p(\R^n) \mapsto \dot W^p_1(\R^n)$ is onto, one-to-one, or invertible, respectively.

Similarly, $(N)^A_{p,1}$ has the existence property, the uniqueness property, or is well-posed in both $\R^{n+1}_+$ and $\R^{n+1}_-$, if and only if $(\partial_\nu^A \D^A)_+:\dot W^p_1(\R^n)\mapsto L^p(\R^n)$ is onto, one-to-one, or invertible, respectively.
\end{thm}

If $(D)^A_{p,\theta}$ is well-posed, then $\s^A_+$ is invertible $\dot B^{p,p}_{\theta-1}(\R^n)\mapsto \dot B^{p,p}_\theta(\R^n)$; by Theorem~\ref{thm:2,1/2}, $\s^A_+$ is also invertible $\dot B^{2,2}_{-1/2}(\R^n)\mapsto \dot B^{2,2}_{1/2}(\R^n)$. However, if $(D)^A_{p,\theta}$ is well-posed but \emph{not} compatible, then the inverses to the operators 
$\s^A_+:\dot B^{p,p}_{\theta-1}(\R^n)\mapsto \dot B^{p,p}_{\theta}(\R^n)$ and $\s^A_+:\dot B^{2,2}_{-1/2}(\R^n)\mapsto \dot B^{2,2}_{1/2}(\R^n)$ are necessarily different operators; that is, there is a function  $f\in \dot B^{2,2}_{1/2}(\R^n) \cap \dot B^{p,p}_{\theta}(\R^n)$ such that $f=\s^A_+g =\s^A_+ h$ for some $g\neq h$, with $g\in \dot B^{p,p}_{\theta-1}(\R^n)$ and $h\in \dot B^{2,2}_{-1/2}(\R^n)$. Thus, in this case, the operator $(\s^A_+)^{-1}$ is not well-defined on the space $\dot B^{p,p}_{\theta-1}(\R^n) \cap \dot B^{2,2}_{-1/2}(\R^n)$ and so bounding $(\s^A_+)^{-1}$ on intermediate spaces by interpolation is not possible.

However, if $(D)^A_{p,\theta}$ is compatibly well-posed, then this issue does not arise.
\begin{thm}\label{thm:well-posed:compatible}
Let $A$ be an elliptic, $t$-independent matrix defined on $\R^{n+1}$ such that $A$ and $A^*$ satisfy the De Giorgi-Nash-Moser condition. 

Let $p$, $\theta$ satisfy the conditions of Theorem~\ref{thm:bounded}.
Suppose that $(D)^A_{p,\theta}$ is compatibly well-posed in both $\R^{n+1}_+$ and~$\R^{n+1}_-$. By Theorem~\ref{thm:well-posed:invertible}, $\s^A_+$ is invertible $\dot B^{p,p}_{\theta-1}(\R^n)\mapsto \dot B^{p,p}_\theta(\R^n)$. We have that the inverse is compatible in the sense that if $f\in \dot B^{p,p}_\theta(\R^n)\cap \dot B^{2,2}_{1/2}(\R^n)$, then $(\s^A_+)^{-1} f$ is the same whether we consider $f$ to be a $\dot B^{p,p}_\theta(\R^n)$ function or a $\dot B^{2,2}_{1/2}(\R^n)$ function.

Similarly, if $(D)^A_{p,1}$ is compatibly well-posed for some $p$ satisfying the conditions of Theorem~\ref{thm:well-posed:invertible:1}, then $(\s^A_+)^{-1}:\dot W^p_1(\R^n)\mapsto L^p(\R^n)$ is compatible with $(\s^A_+)^{-1}:\dot B^{2,2}_{1/2}(\R^n)\mapsto \dot B^{2,2}_{-1/2}(\R^n)$.

Finally, the same statements are true of the problems
$(N)^A_{p,\theta}$ and $(N)^A_{p,1}$ and the operator $(\partial_\nu^A\D^A)_\pm$ in appropriate function spaces.
\end{thm}

Because these theorems reduce well-posedness to boundedness and invertibility of linear operators, we may use standard interpolation and duality results to prove a number of useful corollaries.

\begin{cor}\label{cor:well-posed:open}
Let $A$ be an elliptic, $t$-independent matrix defined on $\R^{n+1}$ such that $A$ and $A^*$ satisfy the De Giorgi-Nash-Moser condition. Let $WP(D)$ and $WP(N)$ be the set of all points $(\theta,1/p)$ that satisfy the conditions of Theorem~\ref{thm:bounded} such that $p<\infty$ and such that $(D)^A_{p,\theta}$ or $(N)^A_{p,\theta}$, respectively, are compatibly well-posed in $\R^{n+1}_+$ and~$\R^{n+1}_-$.

Then $WP(D)$ and $WP(N)$ are open, 
convex, and contain the point $(\theta,1/p)=(1/2,1/2)$.
\end{cor}

\begin{cor}\label{cor:well-posed:1}
Let $A$ be an elliptic, $t$-independent matrix defined on $\R^{n+1}$ such that $A$ and $A^*$ satisfy the De Giorgi-Nash-Moser condition. Then there is some $\varepsilon>0$ such that the following is true. 

Suppose that there is some $p_0$ with $1/2-\varepsilon < 1/p_0<1+\alpha/n$ such that
$(D)^A_{p_0,1}$ is compatibly well-posed in $\R^{n+1}_+$ and $\R^{n+1}_-$.
Then $(D)^A_{p,\theta}$ is compatibly well-posed in $\R^{n+1}_+$ provided $1/2\leq \theta<1$ and $1/p=\theta/p_0+(1-\theta)/p_0'=1-\theta+(2\theta-1)/p_0$. 

This is also true of the Neumann problems~$(N)^A_{p,\theta}$.
\end{cor}

\begin{figure}[tbp]
\begin{tikzpicture}[scale=2]
\figureaxes
\boundedhexagon
\draw [well-posed, ultra thick] 
	(1,0.65) node [black, right] {$(1,1/p_0)$} node [black] {$\circ$} --
	(1/2,1/2) node [black] {$\circ$} node [left,black] {$(1/2,1/2)$}
	-- cycle;	
\node at (0.3,1.3) [black, right] {$(D)^A_{p,\theta}$ or $(N)^A_{p,\theta}$};
\end{tikzpicture}
\qquad
\begin{tikzpicture}[scale=2]
\figureaxes
\boundedhexagon
\draw [well-posed, ultra thick] 
	(0,1-0.65) node [black, left] {$(1,1/p_0')$} node [black] {$\circ$} --
	(1/2,1/2) node [black] {$\circ$} node [right,black] {$(1/2,1/2)$}
	-- cycle;	
\node at (0.3,1.3) [black, right] {$(D)^{A^*}_{p,\theta}$ or $(N)^{A^*}_{p,\theta}$};
\end{tikzpicture}
\caption{Corollaries~\ref{cor:well-posed:1} and~\ref{cor:well-posed:duality}: If $D^A_{p_0,1}$ or $(N)^A_{p_0,1}$ are compatibly well-posed in $\R^{n+1}_\pm$, then $(D)^A_{p,\theta}$ and $(D)^{A^*}_{p,\theta}$ or $(N)^A_{p,\theta}$ and $(N)^{A^*}_{p,\theta}$ are well-posed in $\R^{n+1}_\pm$ for the indicated values of $p$ and~$\theta$.}\label{fig:well-posed:1}
\end{figure}
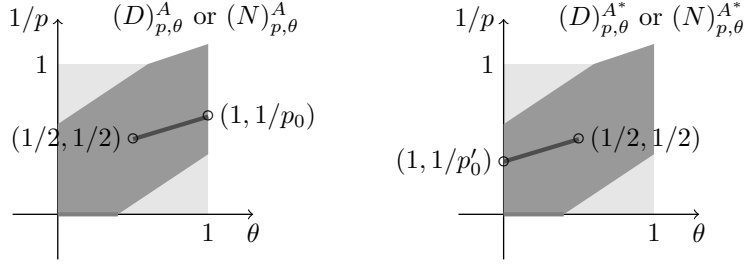

The acceptable values of $\theta$ and $1/p$ in Corollary~\ref{cor:well-posed:1} are indicated in Figure~\ref{fig:well-posed:1}.

Using certain duality results for layer potentials, we may relate boundary-value problems for $A$ to those for~$A^*$.
\begin{cor}\label{cor:well-posed:duality}
Let $A$ be an elliptic, $t$-independent matrix defined on $\R^{n+1}$ such that $A$ and $A^*$ satisfy the De Giorgi-Nash-Moser condition. Let $p$, $\theta$ satisfy the conditions of Theorem~\ref{thm:bounded}, and in addition suppose $p<\infty$.

Suppose $(D)^A_{p,\theta}$ is well-posed in $\R^{n+1}_+$ and~$\R^{n+1}_-$. If $1\leq p<\infty$, then $(D)^{A^*}_{p',1-\theta}$ is well-posed in~$\R^{n+1}_+$, and if $n/(n+\alpha)<p<1$ then $(D)^{A^*}_{\infty,1-\theta+n/p-n}$ is well-posed in~$\R^{n+1}_+$. If $(D)^A_{p,\theta}$ is compatibly well-posed, then so is $(D)^{A^*}_{p',1-\theta}$ or $(D)^{A^*}_{\infty,1-\theta+n/p-n}$.

If $1<p<\infty$, then $(D)^A_{p,\theta}$ has the uniqueness property in both $\R^{n+1}_+$ and $\R^{n+1}_-$ if and only if $(D)^{A^*}_{p',1-\theta}$ has the existence property in both $\R^{n+1}_+$ and $\R^{n+1}_-$.

The theorem remains true if we instead consider the Neumann problem $(N)^A_{p,\theta}$.
\end{cor}

This theorem follows from the fact that the adjoint operator to $\s^A_+$ is $\s^{A^*}_+$, and the adjoint to $(\partial_\nu^A\D^A)_+$ is~$(\partial_\nu^{A^*}\D^{A^*})_+$. See formulas~\eqref{eqn:FS:switch} and~\eqref{eqn:D:conormal:adjoint}. 

The main result of \cite{HofKMP13} is that if the Dirichlet problem $(D)^A_{p,0}$ is well-posed and the square-function estimate \eqref{eqn:square-function} is valid, for some $p$ with $0<1/p<1/2+\epsilon$, then the regularity problem $(D)^{A^*}_{p',1}$ is also well-posed. Thus, the conclusion of Corollary~\ref{cor:well-posed:1} or any other theorem invoking well-posedness of $(D)^A_{p_0,1}$ also follows from well-posedness of~$(D)^{A^*}_{p_0',0}$ with the additional assumption of square-function estimates.

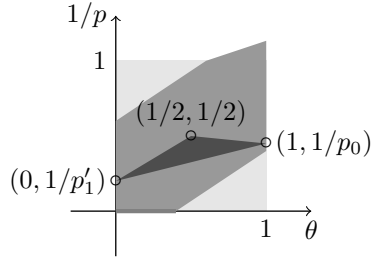
\begin{figure}[tbp]
\begin{tikzpicture}[scale=2]
\figureaxes
\boundedhexagon
\fill [well-posed] 
	(0,0.2) node [black, left] {$(0,1/p_1')$} node [black] {$\circ$} --
	(1,0.45) node [black, right] {$(1,1/p_0)$} node [black] {$\circ$} --
	(1/2,1/2) node [black, above] {$(1/2,1/2)$} node [black] {$\circ$}-- cycle;	
\end{tikzpicture}
\caption{If $(D)^A_{p_0,1}$ and $(D)^{A^*}_{p_1,1}$ or $(N)^A_{p_0,1}$ and $(N)^{A^*}_{p_1,1}$ are compatibly well-posed in $\R^{n+1}_\pm$, then $(D)^A_{p,\theta}$ or $(N)^A_{p,\theta}$ is well-posed in $\R^{n+1}_\pm$ whenever $(\theta,1/p)$ lies in the closure of the indicated triangle.}\label{fig:well-posed:2}
\end{figure}

Notice that if $A$ is self-adjoint, then we may combine the well-posedness regions indicated on the two sides of Figure~\ref{fig:well-posed:1}. More generally, if the regularity or Neumann problems for both $A$ and $A^*$ are compatibly well-posed, then we may combine Corollaries~\ref{cor:well-posed:1} and~\ref{cor:well-posed:duality} and use interpolation to derive well-posedness of $(D)^A_{p,\theta}$ or $(N)^A_{p,\theta}$ for $(\theta,1/p)$ as shown in Figure~\ref{fig:well-posed:2}.

The recent work of Auscher and Mourgoglou \cite{AusM} allows for certain extrapolation theorems; that is, there is some $\alpha^\sharp$ such that if $(D)^A_{p_0,1}$ is well-posed for some $p_0$ with $1<p_0\leq 2$, then $(D)^A_{p,1}$ is also well-posed for any $n/(n+\alpha^\sharp)<p<p_0$. However, their results require a certain technical assumption, closely related to the issue of boundary H\"older continuity. 

Specifically, let the (elliptic but $t$-dependent) matrix $A^\sharp$ be given by
\begin{equation}\label{eqn:AusM}
A^\sharp(x,t)=\begin{cases}A(x), & t>0\\
JA(x)J, &t<0\end{cases}\end{equation}
where $J$ is the (symmetric) Jacobean matrix of the change of variables $(x,t)\mapsto (x,-t)$. The results of \cite{AusM} require that $A^\sharp$ and its adjoint $(A^\sharp)^*$ satisfy the De Giorgi-Nash-Moser condition.

\begin{figure}[tbp]
\begin{tikzpicture}[scale=2]
\figureaxes
\boundedhexagon
\fill [well-posed] 
	(1/2,1/2) node [black] {$\circ$} node [below left,black] {$(1/2,1/2)$}--
	(1-\alphasharp,1) 
		node [black] {$\circ$} 
		node [above left,black, at = {(1-0.3*\alphasharp,1)}] {$(1-\alpha^\sharp,1)$} 
		-- 
	(1,1+\alphasharp/\enn) node [right,black] {$(1,1+\alpha^\sharp/n)$} node [black] {$\circ$} --
	(1,\pnought) node [right,black] {$(1,1/p_0)$} node [black] {$\circ$} --
	cycle;
\node at (0.3,1.5) [black, right] {$(D)^A_{p,\theta}$ or $(N)^A_{p,\theta}$};
\end{tikzpicture}
\qquad
\begin{tikzpicture}[scale=2]
\figureaxes
\boundedhexagon
\fill [well-posed] 
	(\alphasharp,0) node [below,black] {$\alpha^\sharp$} node [black] {$\circ$} --
	(0,0)--
	(0,1-\pnought) node [left,black] {$1/p_0'$} node [black] {$\circ$} --
	(1/2,1/2) node [black] {$\circ$} node [above right,black] {$(1/2,1/2)$}--
	cycle;
\draw [boundary well-posed] (0,0)--(\alphasharp,0);
\node at (0.3,1.5) [black, right] {$(D)^{A^*}_{p,\theta}$ or $(N)^{A^*}_{p,\theta}$};
\end{tikzpicture}
\caption{Corollary~\ref{cor:well-posed:AusM}: If $A$, $A^*$, $A^\sharp$ and $(A^\sharp)^*$ satisfy the De Giorgi-Nash-Moser condition, and if $D^A_{p_0,1}$ or $(N)^A_{p_0,1}$ is well-posed in $\R^{n+1}_\pm$, then $(D)^A_{p,\theta}$ and $(D)^{A^*}_{p,\theta}$ or $(N)^A_{p,\theta}$ and $(N)^{A^*}_{p,\theta}$ are well-posed in $\R^{n+1}_\pm$ for $(\theta,1/p)$ in the indicated quadrilaterals.}\label{fig:well-posed:AusM:1}
\end{figure}
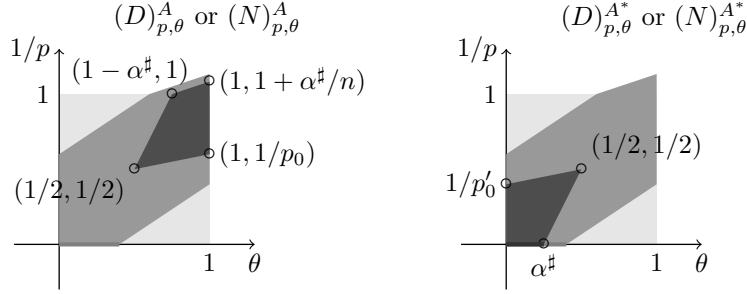

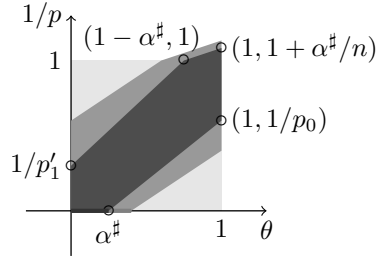
\begin{figure}[tbp]
\begin{tikzpicture}[scale=2]
\figureaxes
\boundedhexagon
\fill [well-posed] 
	(\alphasharp,0) node [below,black] {$\alpha^\sharp$} node [black] {$\circ$} --
	(0,0)--
	(0,1-\pone) node [left,black] {$1/p_1'$} node [black] {$\circ$} --
		(1-\alphasharp,1) 
		node [black] {$\circ$} 
		node [above left,black, at = {(1-0.3*\alphasharp,1)}] {$(1-\alpha^\sharp,1)$} 
		-- 
	(1,1+\alphasharp/\enn) node [right,black] {$(1,1+\alpha^\sharp/n)$} node [black] {$\circ$} --
	(1,\pnought) node [right,black] {$(1,1/p_0)$} node [black] {$\circ$} --
	cycle;
\draw [boundary well-posed] (0,0)--(\alphasharp,0);
\end{tikzpicture}
\caption{Corollary~\ref{cor:well-posed:AusM}: 
If $A$, $A^*$, $A^\sharp$ and $(A^\sharp)^*$ satisfy the De Giorgi-Nash-Moser condition, and if $(D)^A_{p_0,1}$ and $(D)^{A^*}_{p_1,1}$ or $(N)^A_{p_0,1}$ and $(N)^{A^*}_{p_1,1}$ are well-posed in $\R^{n+1}_\pm$, then $(D)^A_{p,\theta}$ or $(N)^A_{p,\theta}$ is well-posed in $\R^{n+1}_\pm$ whenever $(\theta,1/p)$ lies in the indicated hexagon.}\label{fig:well-posed:AusM:2}
\end{figure}

\begin{cor}\label{cor:well-posed:AusM}
Let $A$ be an elliptic, $t$-independent matrix defined on $\R^{n+1}$ such that $A$, $A^*$, $A^\sharp$ and $(A^\sharp)^*$ satisfy the De Giorgi-Nash-Moser condition with exponent~$\alpha^\sharp$. 

Suppose that there is some $p_0$ with $1/2 < 1/p_0<1$ such that
$(D)^A_{p_0,1}$ is compatibly well-posed in $\R^{n+1}_+$ and $\R^{n+1}_-$.
Then $(D)^A_{p,\theta}$ or $(D)^{A^*}_{p,\theta}$ is well-posed in $\R^{n+1}_+$ whenever the point $(\theta,1/p)$ lies in the quadrilateral region shown on the appropriate side of Figure~\ref{fig:well-posed:AusM:1}. 

If there is also some $p_1$ with $1/2 < 1/p_1<1$ such that
$(D)^{A^*}_{p_1,1}$ is compatibly well-posed in $\R^{n+1}_+$ and $\R^{n+1}_-$, then $(D)^A_{p,\theta}$ is well-posed in $\R^{n+1}_+$ whenever the point $(\theta,1/p)$ lies in the hexagonal region shown in Figure~\ref{fig:well-posed:AusM:2}. 

This statement is also true of the Neumann problems~$(N)^A_{p,\theta}$.
\end{cor}

We will prove Corollary~\ref{cor:well-posed:AusM} in Section~\ref{sec:real}.

Notice that if $A$ is real-valued, then $A$ and $A^\sharp$ satisfy the De Giorgi-Nash-Moser condition. Furthermore, in this case 
\cite{KenR09}, \cite{HofKMP13} and \cite{KenP93} give useful well-posedness results; by Corollary~\ref{cor:well-posed:AusM} we have the following result.

\begin{cor}\label{cor:well-posed:real} Suppose that $A$ is real-valued, $t$-independent and elliptic. 

Then there exist some numbers $p_0$, $p_1$ with $1<p_j<\infty$ such that 
$(D)^A_{p,\theta}$ is well-posed for all $(\theta,1/p)$ as in Figure~\ref{fig:well-posed:AusM:2}. If $A$ is real symmetric, then by \cite{KenP93}, we have that $p_0=p_1> 2$ and furthermore the same result is true for the Neumann problem.
\end{cor}

Our last main result concerns well-posedness of the inhomogeneous problems.
\begin{thm}\label{thm:inhomogeneous}
Suppose that $A$, $p$, $\theta$ and $q$ satisfy the conditions of Theorem~\ref{thm:bounded}.

If the homogeneous Dirichlet problem $(D)^A_{p,\theta,q}$ or Neumann problem  $(N)^A_{p,\theta,q}$ is well-posed in $\R^{n+1}_+$, then the same is true of the corresponding inhomogeneous problem. In particular, if $A$ is real then the inhomogeneous version of $(D)^A_{p,\theta,q}$ is well-posed whenever $p^-<q<p^+$ and the point $(\theta,1/p)$ is as in Figure~\ref{fig:well-posed:AusM:2}.

Furthermore, the solution $u$ to 
\begin{equation*}\Div A\nabla u=\Div \F \text{ in }\R^{n+1}_+,\qquad \Tr u=f
\end{equation*}
is given by $u=\Pi^A\F + \s^A((\s^A_+)^{-1} (f-\Pi^A_+\F))$.

The solution $u$ to 
\begin{equation*}\Div A\nabla u=\Div \F \text{ in }\R^{n+1}_+,\qquad \nu_+\cdot A\nabla u\big\vert_{\partial\R^{n+1}_+}=f
\end{equation*}
is given by $u=\Pi^A\F + \D^A((\partial_\nu^A\D^A)_+^{-1} (f-(\partial_\nu^A\Pi^A)_+\F))$.
\end{thm}

\begin{proof}
Let $u$ be as given in the theorem statement; if $(D)^A_{p,\theta,q}$ or $(N)^A_{p,\theta,q}$ is well-posed, then by Theorem~\ref{thm:well-posed:invertible}, the appropriate layer potential is invertible.
By formula~\eqref{eqn:Div:plus}, we have that $\Div A\nabla u=\Div \F$, and by construction, we have that $\Tr u=f$ or $\nu_+\cdot A\nabla u\big\vert_{\partial\R^{n+1}_+}=f$.

We are left with the estimate 
\begin{equation}\label{eqn:inhomogeneous:D}
\doublebar{u}_{\dot W(p,\theta,q)}\leq C\doublebar{f}_{\dot B^{p,p}_\theta(\R^n)}+ C\doublebar{\F}_{L(p,\theta,q)}
\quad\text{(the Dirichlet problem)}
\end{equation}
or 
\begin{equation}\label{eqn:inhomogeneous:N}
\doublebar{u}_{\dot W(p,\theta,q)}\leq C\doublebar{f}_{\dot B^{p,p}_{\theta-1}(\R^n)}+ C\doublebar{\F}_{L(p,\theta,q)}
\quad\text{(the Neumann problem)}
.\end{equation}

Consider the inhomogeneous Dirichlet problem, so $u=\Pi^A\F + \s^A((\s^A_+)^{-1} (f-\Pi^A_+\F))$.
By the bound~\eqref{eqn:space:newton}, $\doublebar{\Pi^A\F}_{\dot W(p,\theta,q)}\leq C\doublebar{\F}_{L(p,\theta,q)}$.
By the bound \eqref{eqn:space:S} and the bound on~$(\s^A_+)^{-1}$, we have that
\begin{align*}
\doublebar{\s^A((\s^A_+)^{-1} (f-\Pi^A_+\F))}_{\dot W(p,\theta,q)}
&\leq C\doublebar{(\s^A_+)^{-1} (f-\Pi^A_+\F)}_{\dot B^{p,p}_{\theta-1}(\R^n)}
\\&\leq C\doublebar{f-\Pi^A_+\F}_{\dot B^{p,p}_{\theta}(\R^n)}
.\end{align*}
By the bound \eqref{eqn:trace:newton}, $\doublebar{\Pi^A_+ \F}_{\dot B^{p,p}_\theta(\R^n)}\leq C\doublebar{\F}_{L(p,\theta,q)}$. Thus, the bound~\eqref{eqn:inhomogeneous:D} is valid.

Similarly, if the inhomogeneous Neumann problem is well-posed, then by the bounds \eqref{eqn:space:D} and~\eqref{eqn:normal:newton} and the bound on $(\partial_\nu^A\D^A)_+^{-1}$, the bound \eqref{eqn:inhomogeneous:N} is valid.

As in Remark~\ref{rmk:inhomogeneous}, the uniqueness property is the same for the homogeneous and inhomogeneous versions.
\end{proof}

\begin{rmk} If the homogeneous version of $(D)^A_{p,\theta,q}$ or $(N)^A_{p,\theta,q}$ is merely solvable or has the existence property, then by writing $u=\Pi^A\F+v$ where $v$ satisfies
\[\Div A\nabla v=0\text{ in }\R^{n+1}_+,\quad \Tr v=f-\Pi^A_+\F\]
or
\[\Div A\nabla v=0\text{ in }\R^{n+1}_+,\quad \nu\cdot A\nabla v\big\vert_{\partial\R^{n+1}_+}=f-(\partial_\nu^A\Pi^A)_+\F,\]
we see that the inhomogeneous problem is also solvable or has the existence property.
\end{rmk}

\section{Sharpness of these results}
\label{sec:sharp}

We remark that Corollary~\ref{cor:well-posed:real} is sharp in the sense that, for any $\alpha^\sharp>0$ and any $p_0$, $p_1$ with $1<p_j<\infty$, there is some matrix $A$ such that either $(D)^A_{p,\theta}$ or  $(D)^{A^*}_{p,\theta}$ is ill-posed whenever $\theta$, $p$ satisfy the conditions of Theorem~\ref{thm:bounded} but not of Corollary~\ref{cor:well-posed:real}.

We return to the example of poor regularity constructed in \cite{KenKPT00}.
Recall that this example concerns the matrix
\begin{equation*}A_k(x,t)=A_k(x)=\begin{pmatrix}1&k\mathop{\mathrm{sgn}}(x)\\ -k\mathop{\mathrm{sgn}}&1\end{pmatrix}\end{equation*}
where $k$ is a real number.
In \cite[Theorem~3.2.1]{KenKPT00}, the authors showed that if $u_k(x,t)=\im(\abs{x}+it)^\beta$, and if $k=\tan\bigl(\frac{\pi}{2}(1-\beta)\bigr)$, then $\Div A_k\nabla u_k=0$ in $\R^2_+$.
This $u_k$ (with $0<\beta<1$) was used in that paper to show that $(D)^{A_k}_{p,0}$ is ill-posed for $p\leq 1/\beta$, and in the appendix to \cite{KenR09} to show that $(D)^{A_k}_{p,1}$ and $(N)^{A_{-k}}_{p,1}$ are ill-posed  for $p>1/(1-\beta)$. The nature of this ill-posedness was investigated in \cite{Axe10}.

It is elementary to show that $\abs{\nabla u_k(x,t)}=\beta (x^2+t^2)^{(\beta-1)/2}$. Let $\Omega_\varepsilon$ be the domain above the level set $u_k(x,t)=\varepsilon^\beta$ for some $\varepsilon>0$. If $0<\beta<1$ then $\Omega_\varepsilon=\{(x,t):t>\varphi_\varepsilon(x)\}$ for some Lipschitz function~$\varphi_\varepsilon$, and so as discussed in the introduction (see Figure~\ref{fig:change-of-variables}), we may consider well-posedness of $(D)^A_{p,\theta}$ in $\Omega_\varepsilon$ rather than the upper half-plane.

Notice that $B(0,\varepsilon)$ lies below the level set $u_k(x,t)=\varepsilon^\beta$. By direct computation 
\begin{equation*}\int_{\Omega_\varepsilon}\abs{\nabla u(X)}^p \,\dist(X,\partial\Omega_\varepsilon)^{p-1-p\theta} \,dX
\leq
\int_{\R^2_+\setminus B(0,\varepsilon)}\abs{\nabla u(x,t)}^p \,t^{p-1-p\theta} \,dx\,dt
<\infty\end{equation*}
whenever $\beta<\theta-1/p$ and $\varepsilon>0$.  Furthermore, $\nabla u$ is locally bounded and so we may pass to averages of $\nabla u$ over Whitney balls; thus, if $1/p<\theta-\beta$ then $u$ is a solution to $(D)^{A_k}_{p,\theta}$ in $\Omega_\varepsilon$. But by construction $u$ is constant along $\partial\Omega_\varepsilon$, and so $(D)^{A_k}_{p,\theta}$ does not have the uniqueness property. (This parallels a counterexample to uniqueness for $(D)^{A_k}_{p,1}$ as presented in \cite{Bar13}. The counterexample of \cite{KenR09} is in fact a counterexample to \emph{compatible} well-posedness.)

By Corollary~\ref{cor:well-posed:duality}, this means that if $\theta+\beta<1/p<1$, then $(D)^{A_{k}^*}_{p,\theta}$ does not have the existence property; for some $f\in\dot B^{p,p}_\theta(\R^n)$, there is no solution $u\in \dot W(p,\theta,2)$ with $\Tr u=f$. 

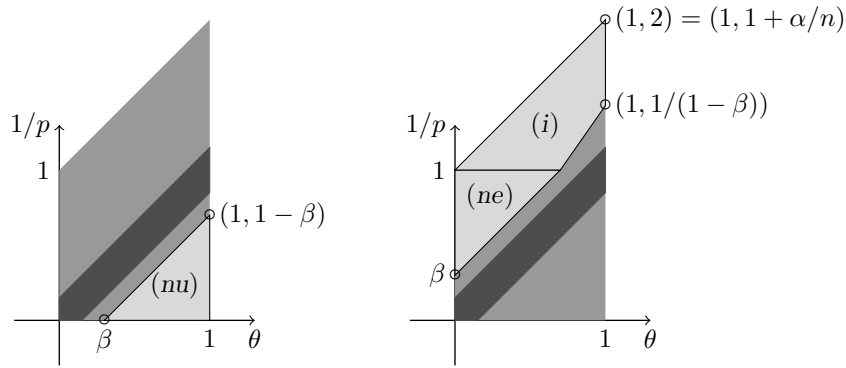
\begin{figure}[tbp]
\def\oneminusbeta{0.7}
\def\mysmaller{0.15}

\begin{tikzpicture}[scale=2]
\plainfigureaxes
\node at (0,1) [black, left] {$1$};
\node at (1,0) [black, below] {$1$};
\fill [bounded] (0,0)--(1,0)--(1,1+1) --(0,1)--cycle;

\filldraw [well-posed]
	(0,0) -- (0,\mysmaller) -- (1,1+\mysmaller) -- (1,1-\mysmaller) -- (\mysmaller,0) -- cycle;

\filldraw [ill-posed]
	(1-\oneminusbeta,0) node [black] {$\circ$} node [black,below] {$\beta$} --
	(1,0) --
	(1,\oneminusbeta) node [black] {$\circ$} node [black,right] {$(1,1-\beta)$} --
	cycle;
\node at (1-\oneminusbeta/3,\oneminusbeta/3) {(\textsl{nu})};

\end{tikzpicture}
\qquad
\begin{tikzpicture}[scale=2]
\plainfigureaxes
\node at (0,1) [black, left] {$1$};
\node at (1,0) [black, below] {$1$};
\fill [bounded] (0,0)--(1,0)--(1,1+1) --(0,1)--cycle;
		
\filldraw [ill-posed]
	(0,1-\oneminusbeta) node [black] {$\circ$} node [black,left] {$\beta$}--
	(\oneminusbeta,1) --
	(0,1)--cycle;
\node at (\oneminusbeta/3,1-\oneminusbeta/4) {(\textsl{ne})};
	
\filldraw [ill-posed]
	(1, 1/\oneminusbeta)
	node [black] {$\circ$}
	node [black,right] {$(1,1/(1-\beta))$}
	--
	(1, 1+ 1)
	node [black] {$\circ$}
	node [black,right] {$(1,2)=(1,1+\alpha/n)$}
	--
	(0,1) -- (\oneminusbeta,1) --
	cycle;
\node at (3/5,1+3/10) {(\textsl{i})};

\filldraw [well-posed]
	(0,0) -- (0,\mysmaller) -- (1,1+\mysmaller) -- (1,1-\mysmaller) -- (\mysmaller,0) -- cycle;

\end{tikzpicture}

\caption{For any $\beta$ with $0<\beta<1$, there is some real, $t$-independent elliptic $2\times 2$ matrix $A$ such that $(D)^A_{p,\theta}$ does not have the uniqueness property in region~(\textsl{nu}), and such that $(D)^{A^*}_{p,\theta}$ does not have the existence property in region~(\textsl{ne}) and is ill-posed in region~(\textsl{i}).}\label{fig:ill-posed}
\end{figure}

That is, $(D)^{A_k}_{p,\theta}$ does not have the uniqueness property whenever $(\theta,1/p)$ lies in region~(\textsl{nu}) of Figure~\ref{fig:ill-posed}, and $(D)^{A_k^*}_{p,\theta}$ does not have the existence property whenever $(\theta,1/p)$ lies in region~(\textsl{ne}) of Figure~\ref{fig:ill-posed}. But by Corollary~\ref{cor:well-posed:real}, $(D)^{A_k^*}_{p,\theta}$ is well-posed whenever $(\theta,1/p)$ lies in the dark region of Figure~\ref{fig:ill-posed}, and in particular whenever $\theta=1/p$. By Corollary~\ref{cor:well-posed:open}, the region $WP(D)$ of well-posedness is convex; thus, $(D)^{A_k^*}_{p,\theta}$ cannot be well-posed in region~(\textsl{i}) of Figure~\ref{fig:ill-posed}.

Thus, $(D)^{A_k}_{p,\theta}$ and $(D)^{A_{k}^*}_{p,\theta}$ are ill-posed in the regions shown in Figure~\ref{fig:ill-posed}. The trapezoidal background in Figure~\ref{fig:ill-posed} is different from the region shown in Figure~\ref{fig:bounded} for two reasons. First, Figure~\ref{fig:bounded} was drawn with $n=\enn$ rather than $n=1$. Second, if $A$ is an elliptic $t$-independent matrix in two dimensions, then $A$ satisfies the De Giorgi-Nash-Moser condition with exponent~$\alpha=1$; see \cite[Th\'eor\`eme~II.2]{AusT95}. Thus, for $t$-independent $2\times 2$ coefficient matrices, the region of boundedness illustrated in Figure~\ref{fig:bounded} expands to the trapezoid with vertices at $(0,0)$, $(1,0)$, $(1,2)$ and $(0,1)$.

%% file: sec-4-background.tex
\chapter{Interpolation, Function Spaces and Elliptic Equations}
\label{chap:background}

In this chapter we will collect some known results and a few extensions; these results will be useful in later chapters. Specifically, in Section~\ref{sec:interpolation} we will discuss the real and complex methods of interpolation, in Section~\ref{sec:function} we will explore properties of the function spaces defined in Section~\ref{sec:dfn:function}, and in Section~\ref{sec:elliptic} we will investigate properties of solutions to elliptic equations. In Chapter~\ref{chap:invertibility}, we will need some known results concerning solutions to boundary-value problems; we will delay reviewing these results until Chapter~\ref{chap:sobolev}.

\section{Interpolation functors}
\label{sec:interpolation}





In the proofs of Theorem~\ref{thm:bounded} and Corollaries~\ref{cor:well-posed:open}, \ref{cor:well-posed:1}, and~\ref{cor:well-posed:AusM}, we will make use of interpolation. That is, we will establish boundedness of $\Pi$, $\D$, or $\s$ for values of $p$, $\theta$ only in some small regions of the hexagon in Figure~\ref{fig:bounded}, and will use interpolation to extend to the entire region. To that end, in this section, we will review some known facts concerning real and complex interpolation; we refer the interested reader to the classic reference \cite{BerL76} for the details of real interpolation and of complex interpolation in Banach spaces, and to \cite{KalM98,KalMM07} for a treatment of complex interpolation in quasi-Banach spaces.

Following \cite{BerL76}, we say that two quasi-normed vector spaces $A_0$, $A_1$ are \emph{compatible} if there is a Hausdorff topological vector space $\mathfrak{A}$ such that $A_0\subset\mathfrak{A}$, $A_1\subset\mathfrak{A}$. Then $A_0\cap A_1$, $A_0+A_1$ may be defined in the natural way.

We will use two interpolation functors, the real method of Lions and Peetre, and the complex method of Lions, Calder\'on and Krejn, to identify special subspaces of $A_0+A_1$.

If $0< r\leq \infty$ and $0<\sigma<1$, then we let the real interpolation space $(A_0,A_1)_{\sigma,r}$ denote $\{a\in A_0+A_1:\doublebar{a}_{(A_0,A_1)_{\sigma,r}}<\infty\}$, where
\begin{equation*}
\doublebar{a}_{(A_0,A_1)_{\sigma,r}}
=\biggl(\int_0^\infty \bigl(s^{-\sigma} \inf\{\doublebar{a_0}_{A_0}+s\doublebar{a_1}_{A_1}:a_0+a_1=a\}\bigr)^r \,ds/s\biggr)^{1/r}.\end{equation*}

If $0<\sigma<1$, and $A_0$, $A_1$ are Banach spaces, then the complex interpolation space $[A_0,A_1]_\sigma$ is defined as follows. Let $S$ be the strip in the complex plane $S=\{z:0<\re z<1\}$ and let $\overline{S}$ be its closure. Let $\mathcal{F}$ denote the set of all continuous functions $f:\overline{S}\mapsto A_0+A_1$ that are analytic $S\mapsto A_0+A_1$.  We then let $[A_0,A_1]_\sigma$  denote $\{a\in A_0+A_1:\doublebar{a}_{[A_0,A_1]_{\sigma}}<\infty\}$, where
\begin{equation*}
\doublebar{a}_{[A_0,A_1]_{\sigma}}
=\inf\bigl\{\max(\sup_{s\in\R} \doublebar{f(is)}_{A_0}, \sup_{s\in\R} \doublebar{f(1+is)}_{A_1}):f\in\mathcal{F},f(\sigma)=a\bigr\}
.\end{equation*}
Observe that the Besov spaces $\dot B^{p,r}_\theta(\R^n)$, with $p<1$ or $r<1$, are only quasi-Banach spaces. The real interpolation method applies to quasi-Banach spaces without modification. Difficulties arise in applying the complex interpolation method to general quasi-Banach spaces. However, the Besov spaces $\dot B^{p,r}_\theta(\R^n)$ have an important property (the property of being \emph{analytically convex}), which means that the complex interpolation method may still be used. We refer the interested reader to \cite{KalM98,MenM00,KalMM07} for the details; here we will only need the properties discussed below, specifically the bound \eqref{eqn:operator:interpolation:complex}, Theorem~\ref{thm:interpolation:extrapolation} and the interpolation formulas \eqref{eqn:F-F:interpolation:complex} and~\eqref{eqn:B-B:interpolation:complex}.

These two interpolation methods have the following interpolation property.
For any compatible pairs $A_0$, $A_1$ and $B_0$, $B_1$, we have that if $T:A_0+A_1\mapsto B_0+B_1$ is a linear operator such that $T(A_0)\subseteq B_0$ and $T(A_1)\subseteq B_1$, then $T$ is bounded on appropriate interpolation spaces, with
\begin{align}
\label{eqn:operator:interpolation:real}
\doublebar{T}_{(A_0,A_1)_{\sigma,r}\mapsto(B_0,B_1)_{\sigma,r}} &\leq 
\doublebar{T}_{A_0\mapsto B_0}^{1-\sigma}\doublebar{T}_{A_1\mapsto B_1}^{\sigma},
\\
\label{eqn:operator:interpolation:complex}
\doublebar{T}_{[A_0,A_1]_{\sigma}\mapsto[B_0,B_1]_{\sigma}} &\leq 
\doublebar{T}_{A_0\mapsto B_0}^{1-\sigma}\doublebar{T}_{A_1\mapsto B_1}^{\sigma}
\end{align}
for $0<\sigma<1$ and $0<r\leq\infty$.
See \cite[Theorems~3.11.2 and~4.1.2]{BerL76}.

We have the following theorem concerning operators that are invertible on an interpolation space. In the case of quasi-Banach spaces, see {\cite[Theorem~2.7]{KalM98}}. In the case of Banach spaces the result was known earlier; see \cite[Lemma~7]{Sne74}. 

\begin{thm}
\label{thm:interpolation:extrapolation}
Let $(A_0,A_1)$ and $(B_0,B_1)$ be two compatible couples, and suppose that $T:A_0+A_1\mapsto B_0+B_1$ is a linear operator with $T(A_0)\subseteq B_0$ and $T(A_1)\subseteq B_1$; by the bound~\eqref{eqn:operator:interpolation:complex}, $T$ is bounded $A_\theta\mapsto B_\theta$ for any $0<\theta<1$.

Suppose that for some $\theta_0$ with $0<\theta_0<1$, we have that $T:A_{\theta_0}\mapsto B_{\theta_0}$ is invertible. Then  there is some $\varepsilon>0$ such that $T:A_\theta\mapsto B_\theta$ is invertible for all $\theta$ with $\theta_0-\varepsilon<\theta<\theta_0+\varepsilon$.
\end{thm}
This theorem has been applied in the context of boundary-value problems for harmonic functions; see \cite{KalM98} and~\cite{KalMM07}.

We will need a few general properties of interpolation spaces. Notice that $(A_0,A_1)_{\sigma,r}=(A_1,A_0)_{1-\sigma,r}$ and that $[A_0,A_1]_{\sigma}=[A_1,A_0]_{1-\sigma}$.
By \cite[Theorem~3.4.2 and remarks on p.~66]{BerL76}, if $0<\sigma<1$ and $0<r<\infty$ then $A_0\cap A_1$ is dense in $(A_0,A_1)_{\sigma,r}$.
Furthermore, if $a\in A_0\cap A_1$, then (for $0<r\leq \infty$) it is straightforward to show that
\begin{equation}\label{eqn:interpolation:intersection}
\doublebar{a}_{(A_0,A_1)_{\sigma,r}}\leq C(\sigma)\doublebar{a}_{A_0}^{1-\sigma} \doublebar{a}_{A_1}^{\sigma}
,\qquad
\doublebar{a}_{[A_0,A_1]_{\sigma}}\leq \doublebar{a}_{A_0}^{1-\sigma} \doublebar{a}_{A_1}^{\sigma}
.\end{equation}

\begin{figure}[tbp]
\begin{tikzpicture}[scale=2]
\plainfigureaxes
\draw [boundary bounded] 
(1/3,1/3) node [black] {$\bullet$} node [black, below] {$(\theta_0,1/p_0)$}
--
(3/4,2/3) node [black] {$\bullet$} node [black, above] {$(\theta_1,1/p_1)$};
\end{tikzpicture}
\caption{If $0<\sigma<1$, then the number $(\theta_\sigma,1/p_\sigma)$ lies on the line segment connecting $(\theta_0,1/p_0)$ and $(\theta_1,1/p_1)$.
}\label{fig:interpolation}
\end{figure}
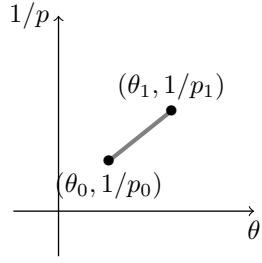

In the remainder of this section we will investigate the effects of interpolation on Besov spaces $\dot B^{p,r}_\theta(\R^n)$ and on the weighted Lebesgue spaces $L(p,\theta,q)$. We will note several persistent patterns in the parameters $p$, $r$ and~$\theta$.
Throughout this section, let $\sigma$, $\theta_0$, $\theta_1$, $p_0$, $p_1$, $p$, $r_0$, $r_1$ and $r$ satisfy the bounds
\begin{equation}\label{eqn:interpolation:numbers}
0<\sigma<1,\quad -\infty<\theta_0,\theta_1<\infty,\quad
0<p_0, p_1, p, r_0, r_1, r\leq \infty.
\end{equation}
Define
\begin{equation}\label{eqn:interpolation:sigma}
\theta_\sigma=(1-\sigma)\theta_0+\sigma\theta_1,\quad
\frac{1}{p_\sigma}=\frac{1-\sigma}{p_0}+\frac{\sigma}{p_1},\quad
\frac{1}{r_\sigma}=\frac{1-\sigma}{r_0}+\frac{\sigma}{r_1}.
\end{equation}
Observe that if we plot $(\theta_0,1/p_0)$ and $(\theta_1,1/p_1)$ in the $xy$-plane, then for any $0<\sigma<1$ the point $(\theta_\sigma,1/p_\sigma)$ lies on the line segment connecting $(\theta_0,1/p_0)$ and $(\theta_1,1/p_1)$. See Figure~\ref{fig:interpolation}.

We begin by quoting the following result concerning real interpolation on Besov spaces.

\begin{thm}[{\cite[Theorem~2.4.2]{Tri83}}]
\label{thm:B-B:interpolation:real}
Let $p$, $r$, $r_0$, $r_1$, $\sigma$ and $\theta_\sigma$ be as in formulas \eqref{eqn:interpolation:numbers} and~\eqref{eqn:interpolation:sigma}, and suppose in addition that $\theta_0\neq\theta_1$. Then
\begin{equation}
\label{eqn:B-B:interpolation:real}
(\dot B^{p,r_0}_{\theta_0}(\R^n),\dot B^{p,r_1}_{\theta_1}(\R^n))_{\sigma,r}
= \dot B^{p,r}_{\theta_\sigma}(\R^n)
.\end{equation}
If in addition $p<\infty$ then
\begin{equation}
\label{eqn:F-F:interpolation:real}
(\dot F^{p,r_0}_{\theta_0}(\R^n),\dot F^{p,r_1}_{\theta_1}(\R^n))_{\sigma,r}
= \dot B^{p,r}_{\theta_\sigma}(\R^n)
.\end{equation}
\end{thm}
Notice that the values of $r_0$, $r_1$ on the left-hand side are irrelevant to the parameters on the right-hand side.

We also have the following result concerning complex interpolation on these spaces.
\begin{thm}[{\cite[Theorem~11]{MenM00}}]
\label{thm:B-B:interpolation:complex}
Let $\sigma$, $\theta_\sigma$, $p_\sigma$ and $r_\sigma$ be as in formulas~\eqref{eqn:interpolation:numbers} and~\eqref{eqn:interpolation:sigma}, and suppose in addition that either $p_0+r_0<\infty$ or $p_1+r_1<\infty$. Then
\begin{equation}
\label{eqn:F-F:interpolation:complex}
[\dot F^{p_0,r_0}_{\theta_0}(\R^n),\dot F^{p_1,r_1}_{\theta_1}(\R^n)]_\sigma
= \dot F^{p_\sigma,r_\sigma}_{\theta_\sigma}(\R^n)
.\end{equation}
If in addition $\theta_0\neq \theta_1$, then
\begin{equation}
\label{eqn:B-B:interpolation:complex}
[\dot B^{p_0,r_0}_{\theta_0}(\R^n),\dot B^{p_1,r_1}_{\theta_1}(\R^n)]_\sigma
= \dot B^{p_\sigma,r_\sigma}_{\theta_\sigma}(\R^n)
.\end{equation}
\end{thm}
Recall that $\dot B^{p,p}_\theta(\R^n)=\dot F^{p,p}_\theta(\R^n)$, and so if $p_0=r_0$ and $p_1=r_1$ then we may drop the requirement that $\theta_0\neq\theta_1$ in formula~\eqref{eqn:B-B:interpolation:complex}. We remark that in the case of Banach spaces (that is, where $p_0,p_1,r_0,r_1\geq 1$), Theorem~\ref{thm:B-B:interpolation:complex} was already known; see \cite[Theorem~6.2.4]{BerL76}, \cite[p.~182--185]{Tri78}, and \cite[Corollary~8.3]{FraJ90}.

We will want to apply interpolation in the spaces $L(p,\theta,q)$. 
The following theorem will let us do this.
\begin{thm}
\label{thm:interpolation:L}
Suppose that $1\leq q\leq\infty$. Let $\sigma$, $\theta_\sigma$, $p_\sigma$ be as in formulas~\eqref{eqn:interpolation:numbers} and~\eqref{eqn:interpolation:sigma}. Suppose further that $0<p_0\leq\infty$ and that $0<p_1<\infty$. Then
\begin{equation}\relax 
(L(p_0,\theta_0,q),L(p_1,\theta_1,q))_{\sigma,p_\sigma} =  
L(p_\sigma,\theta_\sigma,q)
.\end{equation}
If in addition $1\leq p_0$ and $1\leq p_1$ then
\begin{equation}\label{eqn:L-L:interpolation:complex} 
[L(p_0,\theta_0,q),L(p_1,\theta_1,q)]_{\sigma} = 
L(p_\sigma,\theta_\sigma,q)
.\end{equation}
\end{thm}

We conjecture that the complex interpolation formula \eqref{eqn:L-L:interpolation:complex} is also valid for $0<p_0<1$ and $0<p_1<1$; however, the proof does not seem to follow easily from the existing literature, and because we will not need to apply interpolation to the spaces $L(p,\theta,q)$ for $p<1$, we will leave the theorem as stated above.

In order to prove Theorem~\ref{thm:interpolation:L}, we will use the following lemma. This is essentially \cite[Theorem~6.4.2]{BerL76} and \cite[Theorem~1.2.4]{Tri78}. See also \cite[Section~3]{Pee71}.

\begin{lem}\label{lem:retract} Suppose that $(A_0,A_1)$ and $(B_0,B_1)$ are two compatible couples, and that there are linear operators $\mathcal{I}:A_0+A_1\mapsto B_0+B_1$ and $\mathcal{P}:B_0+B_1\mapsto A_0+A_1$ such that $\mathcal{P}\circ\mathcal{I}$ is the identity operator on $A_0+A_1$, and such that $\mathcal{I}:A_0\mapsto B_0$, $\mathcal{I}:A_1\mapsto B_1$, $\mathcal{P}:A_0\mapsto B_0$, and $\mathcal{P}:A_1\mapsto B_1$ are all bounded operators.

If $0<\sigma<1$ and $0<r\leq\infty$, then
\begin{equation*}(A_0,A_1)_{\sigma,r}=\mathcal{P}((B_0,B_1)_{\sigma,r}),
\qquad
[A_0,A_1]_{\sigma}=\mathcal{P}([B_0,B_1]_{\sigma}).\end{equation*}
\end{lem}

\begin{proof}
By boundedness of $\mathcal{I}$ and the interpolation property \eqref{eqn:operator:interpolation:complex}, we have that
$\mathcal{I}([A_0,A_1]_{\sigma})
\subseteq
[B_0,B_1]_{\sigma}$
and so
\begin{equation*}
\relax [A_0,A_1]_{\sigma}
=\mathcal{P}( \mathcal{I} ([A_0,A_1]_{\sigma}))
\subseteq \mathcal{P}([B_0,B_1]_{\sigma})
.\end{equation*}
By boundedness of $\mathcal{P}$, and again by formula~\eqref{eqn:operator:interpolation:complex}, $\mathcal{P}([B_0,B_1]_{\sigma})
\subseteq [A_0,A_1]_{\sigma}$; combining these two inclusions yields that
$\mathcal{P}([B_0,B_1]_{\sigma})
= [A_0,A_1]_{\sigma}$.
The same argument with formula~\eqref{eqn:operator:interpolation:complex} replaced by formula~\eqref{eqn:operator:interpolation:real} establishes the result for the real interpolation method.
\end{proof}

\begin{proof}[Proof of Theorem~\ref{thm:interpolation:L}] 
Let $\G$ be a grid of dyadic Whitney cubes and observe that 
\[\doublebar{\F}_{L(p,\theta,q)}^p\approx\sum_{Q\in\G} \biggl(\fint_Q \abs{\F}^q\biggr)^{p/q}\ell(Q)^{n+p-p\theta}.\]
Let $Q_0$ be the unit cube in $\R^{n+1}$.
For any $1\leq q\leq\infty$ and any $0<p\leq\infty$, define the space of sequences of functions $\ell(p,\theta,q)$ by
\[\ell(p,\theta,q)=\Bigl\{\{\F_Q\}_{Q\in\G}:\F_Q\in L^q(Q_0),\Bigl(\sum_{Q\in\G} \doublebar{\F_Q}_{L^q(Q_0)}^p\ell(Q)^{n+p-p\theta}\Bigr)^{1/p}<\infty \Bigr\}.\]
By examining the proofs of \cite[Theorems~5.6.1 and~5.6.3]{BerL76}, we see that 
\[(\ell(p_0,\theta_0,q),\ell(p_1,\theta_1,q))_{\sigma,p_\sigma} = \ell(p_\sigma,\theta_\sigma,q)\]
and if $1\leq p_0\leq\infty$ and $1\leq p_1\leq \infty$ then
\[\relax[\ell(p_0,\theta_0,q),\ell(p_1,\theta_1,q)]_{\sigma} = \ell(p_\sigma,\theta_\sigma,q).\]
We then apply Lemma~\ref{lem:retract} with $(\mathcal{I}\F)_Q(x,t)=\F(x_Q+\ell(Q)x,t_Q+\ell(Q)t)$ for appropriately chosen $x_Q\in\R^n$ and $t_Q\in \R$. This completes the proof.
\end{proof}

\section{Function spaces}
\label{sec:function}

In this section we collect some useful information regarding the Besov spaces $\dot B^{p,r}_\theta(\R^n)$ and the Triebel-Lizorkin spaces~$\dot F^{p,r}_\theta(\R^n)$.

We identify the following well-known spaces as Triebel-Lizorkin spaces; see \cite[Sections 2.3.5, 5.2.3, and~5.2.4]{Tri83}, as well as \cite[Proposition~3]{MenM00}.
\begin{align}
\label{eqn:H-F}
H^p(\R^n)&=\dot F^{p,2}_0(\R^n), && \frac{n}{n+1}<p<\infty
,\\
\label{eqn:H1-F}
\dot H^p_1(\R^n)&=\dot F^{p,2}_1(\R^n), && \frac{n}{n+1}<p<\infty
,\\
\label{eqn:C-B}
\dot C^\theta(\R^n)&=\dot B^{\infty,\infty}_\theta(\R^n), && 0<\theta<1
.\end{align}
In particular,
\begin{align}
\label{eqn:L-F}
L^p(\R^n)&=H^p(\R^n)=\dot F^{p,2}_0(\R^n), &&1<p<\infty
,\\
\label{eqn:W-F}
\dot W^p_1(\R^n)&=\dot H^p_1(\R^n)=\dot F^{p,2}_1\R^n), && 1<p<\infty
.\end{align}

Next, we provide some comments on the Besov spaces.

First, observe that if $\theta>0$, then by formula~\eqref{eqn:besov:norm} the constant functions have $\dot B^{p,r}_\theta(\R^n)$ norm zero, and so elements of $\dot B^{p,r}_\theta(\R^n)$ are equivalence classes of functions modulo additive constants. Next, observe that by formula~\eqref{eqn:F-F:interpolation:real} and by formulas~\eqref{eqn:H-F}--\eqref{eqn:W-F}, we have that
\begin{align}
\label{eqn:H-H1:interpolation}
\dot B^{p,p}_{\theta}(\R^n)
&=(H^p(\R^n),\dot H^p_1(\R^n))_{\theta,p} &&\text{if $0<\theta<1$ and $\frac{n}{n+1}<p<\infty$}
,\\
\label{eqn:L-W:interpolation}
\dot B^{p,p}_{\theta}(\R^n)
&=(L^p(\R^n),\dot W^p_1(\R^n))_{\theta,p} &&\text{if $0<\theta<1$ and $1<p<\infty$}
.\end{align}

Recall from Section~\ref{sec:interpolation} that if $0<p<\infty$ then $H^p(\R^n)\cap \dot H^p_1(\R^n)$ is dense in $\dot B^{p,p}_{\theta}(\R^n)=(H^p(\R^n),\dot H^p_1(\R^n))_{\theta,p}$. Thus, smooth, compactly supported functions are also dense in~$\dot B^{p,p}_{\theta}(\R^n)$.

Next, we identify the dual spaces to the Besov spaces. Recall that if $1\leq p\leq\infty$ then $p'$ is the extended real number that satisfies $1/p+1/p'=1$, and if $0<p<1$ then $p'=\infty$. By \cite[Section~2.11]{Tri83},
\begin{align}
\label{eqn:dual:banach}
(\dot B^{p,r}_\theta(\R^n))'&=\dot B^{p',r'}_{-\theta}(\R^n)
&& \text{provided } 1\leq p<\infty
,\\
\label{eqn:dual:quasi-banach}
(\dot B^{p,r}_\theta(\R^n))'&=\dot B^{\infty,r'}_{-\theta+n(1/p-1)}(\R^n)
&& \text{provided } 0<p<1
.\end{align}

We have the following embedding theorems for the homogeneous Besov and Triebel-Lizorkin spaces. First, by the norms~\eqref{eqn:triebel:norm} and\eqref{eqn:besov:norm}, the spaces $\dot F^{p,r}_\theta(\R^n)$ and  $\dot B^{p,r}_\theta(\R^n)$ are increasing in~$r$; that is, if $0<r_0<r_1\leq\infty$ then $\dot F^{p,r_0}_\theta(\R^n)\subset \dot F^{p,r_1}_\theta(\R^n)$ and $\dot B^{p,r_0}_\theta(\R^n)\subset \dot B^{p,r_1}_\theta(\R^n)$. 

Another embedding theorem is also possible.
\begin{lem}[{\cite[Theorem~2.1]{Jaw77}}]\label{lem:besov:embedding}
Let $0<p_0< p_1<\infty$ and suppose that the real numbers ${\theta_1}$, ${\theta_0}$ satisfy ${\theta_1}-n/p_1={\theta_0}-n/p_0$.

Then for any $r$ with $0<r\leq\infty$, we have that
\[\dot B^{p_0,r}_{\theta_0}(\R^n)\subset \dot B^{p_1,r}_{\theta_1}(\R^n),
\qquad
\dot F^{p_0,r}_{\theta_0}(\R^n)\subset\dot B^{p_1,p_0}_{\theta_1}
.\]
\end{lem}
We may embed Besov spaces into local $L^p$ spaces as well.
\begin{lem}\label{lem:besov:poincare}
Let $f\in \dot B^{p,p}_\theta(\R^n)$ for some $1\leq p< \infty$ and some $0<\theta<1$. If $\Delta=\Delta(x_0,R)\subset \R^n$ is a ball of radius~$R$ and $c_R=\fint_{\Delta} f$, then
\begin{equation}
\label{eqn:besov-poincare}
\doublebar{f-c_R}_{L^p(\Delta)} \leq C(p,\theta) R^\theta \doublebar{f}_{\dot B^{p,p}_\theta(\R^n)}.\end{equation}
\end{lem}
\begin{proof}
Observe that 
\begin{equation*}\doublebar{f-c_R}_{L^p(\Delta)} \leq 2\doublebar{f}_{L^p(\R^n)},
\qquad
\doublebar{f-c_R}_{L^p(\Delta)} \leq C(p) R\doublebar{f}_{\dot W^p_1(\R^n)}.\end{equation*}
The first inequality follows by direct computation and the second follows from the Poincar\'e inequality. Applying the interpolation property \eqref{eqn:operator:interpolation:real} and formula~\eqref{eqn:L-W:interpolation}, we immediately derive the bound~\eqref{eqn:besov-poincare}.
\end{proof}

We now provide some equivalent characterizations of $\dot B^{p,p}_\theta(\R^n)$; we will provide separate characterizations in the cases $1\leq p\leq\infty$ and $0<p\leq 1$.

Suppose first that $0<\theta<1$ and that $1\leq p\leq \infty$. Then by \cite[Section~5.2.3]{Tri83},
\begin{align}\label{eqn:besov:norm:differences}
\doublebar{f}_{\dot B^{p,p}_\theta(\R^n)}
&\approx \biggl(\int_{\R^n}\int_{\R^n} \frac{\abs{f(y)-f(x)}^p}{\abs{x-y}^{n+p\theta}}\,dx\,dy\biggr)^{1/p}
.\end{align}

In the special case $p=\infty$, this reduces to
\begin{align}\label{eqn:besov:infty:norm}
\doublebar{f}_{\dot B^{\infty,\infty}_\theta(\R^n)}
&\approx\doublebar{f}_{\dot C^\theta(\R^n)}
= \esssup_{(x,y)\in \R^n\times\R^n} \frac{\abs{f(x)-f(y)}}{\abs{x-y}^\theta}.
\end{align}

We can decompose Besov spaces into atoms, as in \cite{FraJ85}. Suppose that $0<p\leq \infty$ and that $\theta\in\R$. Let $ a_Q$ be a function, and suppose that there is some dyadic cube $Q$ such that
\begin{itemize}
\item $\supp a_Q \subseteq 3Q$,
\item $\doublebar{a_Q}_{L^\infty(3Q)} \leq \abs{Q}^{\theta/n-1/p}$,
\item $\doublebar{\partial^\beta a_Q}_{L^\infty(3Q)} \leq \abs{Q}^{\theta/n-1/p-\abs{\beta}/n}$ whenever $\abs{\beta}\leq \theta+1$,
\item $\int x^\gamma \,a_Q(x)\,dx=0$ whenever $\abs{\gamma}\leq \max(n(1/p-1),0)-\theta $.
\end{itemize}
Here $\beta$, $\gamma$ are multiindices.
Then we say that $a_Q$ is a $(\theta,p)$ atom. (Notice that if $\theta>0$ or $\theta<\max(0,n(1/p-1))$, then only one of the last two conditions need be considered.)

If $f\in \dot B^{p,p}_\theta(\R^n)$, then there are constants $\lambda_Q$ and atoms $a_Q$ such that 
\begin{equation*}f=\sum_{Q \>\text{dyadic}} \lambda_Q \, a_Q,\qquad \sum_{Q \>\text{dyadic}} \abs{\lambda_Q}^p\approx \doublebar{f}_{\dot B^{p,p}_\theta(\R^n)}^p.\end{equation*}

\begin{rmk}\label{rmk:atom:strange-norm} If $a_Q$ is a $(\theta_0,p_0)$ atom and $0<p_0\leq \infty$,  $0<p_1\leq \infty$, $0< \theta_0< 1$ and $0< \theta_1< 1$, then
\begin{equation*}
\doublebar{a_Q}_{\dot B^{p_1,p_1}_{\theta_1}(\R^n)}
\leq C \abs{Q}^{\theta_0/n-1/p_0-\theta_1/n+1/p_1}.\end{equation*}

\end{rmk}

\begin{rmk}\label{rmk:atoms}

Suppose that $1\leq q\leq \infty$ and $0<\theta<1$.
If $p\geq 1$ then $L(p,\theta,q)$ and $\dot B^{p,p}_\theta(\R^n)$ are Banach spaces. If $p<1$ then $L^p(U)$ and $\dot  B^{p,p}_\theta(\R^n)$ are only quasi-Banach spaces; however, the $p$-norm bounds
\begin{equation*}
\doublebar{f+g}_{L(p,\theta,q)}^p\leq \doublebar{f}_{L(p,\theta,q)}^p+\doublebar{g}_{L(p,\theta,q)}^p
,\quad\text{if }0<p\leq 1\end{equation*}
and
\begin{equation*}
\doublebar{f+g}_{\dot B^{p,p}_\theta(\R^n)}^p\leq \doublebar{f}_{\dot B^{p,p}_\theta(\R^n)}^p + \doublebar{g}_{\dot B^{p,p}_\theta(\R^n)}^p
,\quad\text{if }0<p\leq 1\end{equation*}
are valid.

The space $\dot B^{p,p}_\theta(\R^n)$ has an atomic decomposition. The space $L(p,\theta,q)$ has a natural decomposition as well, given by
\begin{equation*}\F=\sum_{Q\in\G} \F_Q=\sum_{Q\in\G} \F\1_Q\end{equation*}
where $\G$ is the grid of dyadic Whitney cubes. If $T$ is a linear operator and $0<p\leq 1$, then we may show that $T:L(p,\theta,q)\mapsto\dot B^{p,p}_\theta(\R^n)$ or $T:\dot B^{p,p}_\theta(\R^n)\mapsto L(p,\theta,q)$ is bounded simply by showing that $\doublebar{T\F_Q}_{\dot B^{p,p}_\theta(\R^n)}\leq C\doublebar{\F_Q}_{L(p,\theta,q)}$ or $\doublebar{T\vec a_Q}_{L(p,\theta,q)}\leq C\doublebar{\vec a_Q}_{\dot B^{p,p}_\theta(\R^n)}$.
\end{rmk}

Another useful tool when dealing with Besov spaces with $p\leq 1$ is the notion of a \emph{molecule}. Suppose that $0<\theta<1$ and that $0<p\leq 1$. As defined in \cite{FraJ85}, a function $m$ is a $(\theta,p)$ molecule if there is some $x_0\in\R^n$, some $r>0$ and some integer $M$ large enough such that
\begin{align*}
\abs{m(x)}&\leq r^{\theta-n/p} (1+\abs{x-x_0}/r)^{-M} \quad\text{for all }x\in\R^n
,\\
\abs{\nabla m(x)}&\leq r^{\theta-n/p-1} (1+\abs{x-x_0}/r)^{-M-1} \quad\text{for all }x\in\R^n
.\end{align*}
If $m$ is a $(\theta,p)$ molecule then $m\in \dot B^{p,p}_{\theta}(\R^n)$; see \cite[Section~3]{FraJ85}.

Such molecules are often of interest because many linear operators map atoms to molecules. However, our operators will map atoms to somewhat rougher functions; thus we will need rough molecules. Such objects have been investigated extensively in the case of Hardy spaces; see \cite{TaiW80}, \cite{ChaGT81}, \cite{GarR85}. We will use Lemma~\ref{lem:besov:embedding} and the formula~\eqref{eqn:H1-F} relating Hardy and Triebel-Lizorkin spaces to generalize to Besov spaces.

\begin{lem}\label{lem:molecule:plus}
Let $0<\theta<1$ and $n/(n+1)<p\leq 1$ with $1/p<1+\theta/n$.
Suppose that $1<q\leq\infty$. 
Fix some cube $Q\subset\R^n$, and let $A_0=Q$, $A_j=2^{j}Q\setminus 2^{j-1} Q$ for $j\geq 1$.
Suppose that there is a number $\varepsilon>0$ such that
\begin{align}
\label{eqn:molecule:3}
\biggl(\int_{A_j}
\abs{\nabla m(x)}^q\,dx\biggr)^{1/q}
&\leq 
2^{-j\varepsilon}
(2^j\ell(Q))^{\theta-1+n/q-n/p}  \quad\text{for all }j\geq 0
.\end{align}
Then $m\in \dot B^{p,p}_{\theta}(\R^n)$ with norm at most~$C(\varepsilon,q)$. 
\end{lem}

\begin{proof}
By \cite[Section~2]{TaiW80}, $\nabla m\in H^{p_0}(\R^n)$ and so $m\in \dot H^{p_0}_1(\R^n)$, where $1/p_0=1/p+(1-\theta)/n$. Notice that $n/(n+1)<p_0<p\leq 1$.
By Lemma~\ref{lem:besov:embedding} and formula~\eqref{eqn:H1-F}, we have that $\dot H^{p_0}_1(\R^n)\subset \dot B^{p,p}_{\theta}(\R^n)$, and so the lemma is proven.
\end{proof}

Finally, we prove the following lemma; this lemma will allow us to localize functions in Besov spaces by multiplying by smooth cutoff functions.
\begin{lem}\label{lem:besov:localize}
Let $f\in \dot B^{p,p}_\theta(\R^n)$ for some $1\leq p< \infty$ and some $0<\theta<1$. Suppose that $\xi_R$ is equal to 1 in $\Delta(x_0,R)$, supported in $\Delta(x_0,2R)$, and satisfies $\abs{\nabla\xi_R}\leq C/R$. Let $\tilde c_R=\fint_{\Delta(x_0,2R)\setminus \Delta(x_0,R)} f$. Then
\begin{equation}
\label{eqn:besov-lebesgue-dense}
\lim_{R\to\infty} \doublebar{f-(f-\tilde c_R)\xi_R}_{\dot B^{p,p}_\theta(\R^n)}=0\end{equation}
provided $1\leq p<\infty$ and $0<\theta<1$.
\end{lem}

\begin{proof}

Because constants have $\dot B^{p,p}_\theta(\R^n)$ norm zero, we may replace formula~\eqref{eqn:besov-lebesgue-dense} by
\begin{equation*}
\lim_{R\to\infty} \doublebar{(f-\tilde c_R)(1-\xi_R)}_{\dot B^{p,p}_\theta(\R^n)}=0.
\end{equation*}
By the definition of real interpolation, we have that
\begin{equation*}\doublebar{f}_{\dot B^{p,p}_\theta(\R^n)}^p
\approx
\int_0^\infty \inf \{t^{-p\theta-1} \doublebar{g}_{L^p(\R^n)}^p+t^{p-1-p\theta} \doublebar{h}_{\dot W^p_1(\R^n)}^p:g+h=f\}\,dt.\end{equation*}
Observe that we may choose $g$ and $h$ consistently over dyadic intervals; thus,
\begin{equation*}\doublebar{f}_{\dot B^{p,p}_\theta(\R^n)}^p
\approx
\sum_{j=-\infty}^\infty
\inf\bigl\{2^{-jp\theta} \doublebar{g_j}_{L^p(\R^n)}^p+2^{j(p-p\theta)} \doublebar{h_j}_{\dot W^p_1(\R^n)}^p:g_j+h_j=f \bigr\}.\end{equation*}

Fix a choice of functions $g_j$ and~$h_j$. Let $A(R)={\Delta(x_0,2R)\setminus \Delta(x_0,R)}$, and let $g_{j,R}=\fint_{A(R)} g_j$, $h_{j,R}=\fint_{A(R)} h_j$ so that $\tilde c_R=g_{j,R}+h_{j,R}$. Notice that
\[(f-\tilde c_R)(1- \xi_R)= 
g_j(1- \xi_R)+ g_{j,R} \xi_R
+(h_j-h_{j,R})(1- \xi_R)
-g_{j,R}\]
and so
\begin{multline*}
\doublebar{(f-\tilde c_R)(1- \xi_R)}_{\dot B^{p,p}_\theta(\R^n)}^p
\leq C
	\sum_{j=-\infty}^\infty
	2^{-jp\theta} \doublebar{g_j(1-\xi_R) + g_{j,R}\xi_R}_{L^p(\R^n)}^p
\\
	+ 2^{j(p-p\theta)} \doublebar{(h_j- h_{j,R})(1-\xi_R)-g_{j,R}}_{\dot W^p_1(\R^n)}^p
.\end{multline*}
$g_{j,R}$ is a constant and so has ${\dot W^p_1(\R^n)}$ norm zero. By direct computation,
\[\doublebar{g_j(1-\xi_R) + g_{j,R}\xi_R}_{L^p(\R^n)}
\leq 2\doublebar{g_j}_{L^p(\R^n\setminus \Delta(x_0,R))}.\]
Furthermore, $\nabla ((h_j- h_{j,R})(1-\xi_R))= (1-\xi_R)\nabla h_j -(h_j-h_{j,R})\nabla\xi_R$, and so
\begin{align*}
\doublebar{(h_j- h_{j,R})(1-\xi_R)}_{\dot W^p_1(\R^n)}
&\leq
	\doublebar{\nabla h_j}_{L^p(\R^n\setminus \Delta(x_0,R))}
	+
	\frac{1}{R}\doublebar{h_j-h_{j,R}}_{L^p(A(R))}
.\end{align*}
By the Poincar\'e inequality the second term is at most $C(p)\doublebar{\nabla h_j}_{L^p(A(R))}$.
Because $A(R)\subset \R^n\setminus \Delta(x_0,R)$, we have that
\begin{multline*}
\doublebar{(f-\tilde c_R)(1- \xi_R)}_{\dot B^{p,p}_\theta(\R^n)}^p
\leq C
	\sum_{j=-\infty}^\infty
	2^{-jp\theta} \doublebar{g_j}_{L^p(\R^n\setminus\Delta(x_0,R))}^p
\\
	+ 2^{j(p-p\theta)} \doublebar{\nabla h_j}_{L^p(\R^n\setminus\Delta(x_0,R))}^p
.\end{multline*}

Choose some some $\varepsilon>0$. Because $f\in \dot B^{p,p}_\theta(\R^n)$, there is some $J$ large enough that
\begin{equation*}\sum_{\abs{j}>J}
2^{-jp\theta} \doublebar{g_j}_{L^p(\R^n)}^p+2^{j(p-p\theta)} \doublebar{\nabla h_j}_{L^p(\R^n)}^p<\varepsilon.\end{equation*}
Because $p<\infty$ and there are only finitely many integers~$j$ with $\abs{j}\leq J$, there is some $R>0$ large enough that
\[2^{-jp\theta} \doublebar{g_j}_{L^p(\R^n\setminus\Delta(x_0,R))}^p
	+ 2^{j(p-p\theta)} \doublebar{\nabla h_j}_{L^p(\R^n\setminus\Delta(x_0,R))}^p\leq \varepsilon/J\]
for all $j$ with $\abs{j}\leq J$.
Summing, we see that if $R$ is large enough then 
\begin{equation*}\doublebar{(f-\tilde c_R)(1- \xi_R)}_{\dot B^{p,p}_\theta(\R^n)}^p
\leq C\varepsilon\end{equation*}
as desired.
\end{proof}

\section{Solutions to elliptic equations}
\label{sec:elliptic}

In establishing bounds on potential operators and well-posedness of boundary-value problems, we will need some properties of solutions to the equation $\Div A\nabla u=0$.
We begin with the following well-known bound.
\begin{lem}[The Caccioppoli inequality]\label{lem:Caccioppoli}
Suppose that $A$ is elliptic. Let $X\in\R^{n+1}$ and let $r>0$. Suppose that $\Div A\nabla u=0$ in $B(X,2r)$, for some $u\in W^2_1(B(X,2r))$. Then
\begin{equation*}\int_{B(X,r)}\abs{\nabla u}^2 \leq \frac{C}{r^2}
\int_{B(X,2r)\setminus B(X,r)} \abs{u}^2\end{equation*}
for some constant $C$ depending only on the ellipticity constants of $A$ and the dimension $n+1$.
\end{lem}

Next, we need local bounds on $u$ and its gradient.
Recall from formula~\ref{eqn:local-bound:2} that if $A$ satisfies the De Giorgi-Nash-Moser condition, then any function~$u$ that satisfies $\Div A\nabla u=0$ in $B(X,2r)$ also satisfies
\begin{equation*}
\abs{u(X)}\leq C\biggl(\fint_{B(X,r)} \abs{u}^2\biggr)^{1/2}.
\end{equation*}
We wish to show that we can replace the exponent $2$ by any exponent $p>0$. In at least the case of real coefficients, this result is fairly well known. We will also need a similar result for the gradient $\nabla u$ of a solution; this result may exist in the literature but we provide a proof for the sake of completeness.

We begin with the following lemma; we will follow the method of proof given in \cite[Section~9, Lemma~2]{FefS72}, where such inequalities were established in the case of harmonic functions. 
\begin{lem}\label{lem:subaverage} 
Let $0<p_0<q\leq\infty$.
Suppose that $ U \subset\R^{n+1}$ is a domain, and that $u$ is a function defined on $ U $ with the property that, whenever $B(X,2r)\subset  U $, we have the bound
\begin{equation*}
r^{-(n+1)/q}\doublebar{u}_{L^q(B(X,r))}
\leq C_0 r^{-(n+1)/p_0}\doublebar{u}_{L^{p_0}(B(X,2r))}\end{equation*}
for some constant~$C_0$ depending only on $u$ (not on the specific ball~$B(X,r)$).

Then for every $p$ with $0<p\leq \infty$, there is some constant $C(p)$, depending only on $p$, $p_0$, $q$ and $C_0$, such that
\begin{equation}
\label{eqn:local-bound:step1}
r^{-(n+1)/q}\doublebar{u}_{L^q(B(X,r))}\leq C(p)r^{-(n+1)/p}\doublebar{u}_{L^{p}(B(X,2r))}.
\end{equation}
\end{lem}

\begin{proof}
Fix some $X$, $r$ with $B(X,2r)\subset U $.
By assumption we have that the bound \eqref{eqn:local-bound:step1} is valid for $p=p_0$. Observe that if \eqref{eqn:local-bound:step1} is valid for some~$p<q$ and for arbitrary balls $B\subset U$, then by considering small balls of radius $\rho'-\rho$, we have that 
\begin{equation}\label{eqn:local-bound:step2}
\rho^{-(n+1)/q}\doublebar{u}_{L^q(B(X,\rho))}\leq C(p)\rho^{-(n+1)/p}\biggl(\frac{\rho}{\rho'-\rho}\biggr)^\beta
\doublebar{u}_{L^{p}(B(X,\rho'))}
\end{equation}
whenever $0<\rho<\rho'\leq 2r$ and $p\leq q$, where $\beta={(n+1)(1/p-1/q)}>0$.

Let $p_0$ be as given in the theorem statement and let $p_j=p_{j-1}/2=2^{-j}p_0$.
We will show inductively that if the bound \eqref{eqn:local-bound:step1} (and thus \eqref{eqn:local-bound:step2}) is valid for $p=p_j$, then these bounds are also valid for $p=p_{j+1}$ (at a cost of increasing the constant~$C(p)$). Since $p_j\to 0$ as $j\to\infty$, H\"older's inequality then implies that \eqref{eqn:local-bound:step1} is valid for all $p>0$.

Let $r=\rho_0<\rho_1<\rho_2<\dots<2r$ for some $\rho_k$ to be chosen momentarily, and let $B_k=B(X,\rho_k)$. If $0<\tau<1$, then
\begin{align*}
\doublebar{u}_{L^{p_j}({B_k})}
&=
\biggl(\int_{B_k} \abs{u}^{p_j}\biggr)^{1/p_j}
=
\biggl(\int_{B_k} \abs{u}^{\tau p_j} \abs{u}^{(1-\tau) p_j}\biggr)^{1/p_j}
\end{align*}
and if $0<\tau\leq 1/2$, then $p_{j+1}/\tau p_j\geq 1$ and so we may apply H\"older's inequality to see that
\begin{align*}
\doublebar{u}_{L^{p_j}({B_k})}
&\leq
\doublebar{u}_{L^{p_{j+1}}({B_k})}^\tau
\doublebar{u}_{L^\gamma({B_k})}^{1-\tau}
\end{align*}
where $\gamma$ satisfies $1/p_j=\tau/p_{j+1}+(1-\tau)/\gamma$.

Choose $\tau=\tau(j)$ so that $\gamma=q$; this means that $\tau=(q-p_j)/(2q-p_j)$, and because $p_j\leq p_0<q$, we do have that $0<\tau\leq (q-p_0)/(2q-p_0)<1/2$. By the bound~\eqref{eqn:local-bound:step2},
\begin{align*}
\frac{\doublebar{u}_{L^{p_j}({B_k})}}{r^{(n+1)/p_j}}
&\leq
	\frac{\doublebar{u}_{L^{p_{j+1}}({B_k})}^\tau} {r^{(n+1)\tau/p_{j+1}}}
	\frac{\doublebar{u}_{L^q({B_k})}^{1-\tau}}
	{r^{(n+1)(1-\tau)/q}}
\\&\leq
	\frac{\doublebar{u}_{L^{p_{j+1}}({B_k})}^\tau} {r^{(n+1)\tau/p_{j+1}}}
	\biggl(C(p_j) \biggl(\frac{r}{\rho_{k+1}-\rho_k}\biggr)^\beta
	\frac{\doublebar{u}_{L^{p_j}(B_{k+1})}}{r^{(n+1)/p_j}} 
	\biggr)^{1-\tau}
.\end{align*}
Recall that $\rho_0=r$. Let $\rho_{k+1}=\rho_k + r(1-\sigma) \sigma^k$ for some constant $0<\sigma<1$ to be chosen momentarily. Notice that $\lim_{k\to\infty} \rho_k=2r$. We then have that
\begin{align*}
\frac{\doublebar{u}_{L^{p_j}({B_k})}}{r^{(n+1)/p_j}}
&\leq
	\biggl(\frac{\doublebar{u}_{L^{p_{j+1}}({B_k})}} {r^{(n+1)/p_{j+1}}}\biggr)^\tau
	\biggl(\frac{C(p_j)}{(1-\sigma)^\beta}\sigma^{-k\beta}
	\frac{\doublebar{u}_{L^{p_j}(B_{k+1})}}{r^{(n+1)/p_j}} 
	\biggr)^{1-\tau}
\\&\leq
	\biggl(C(p_j,\sigma)\sigma^{-k\beta(1-\tau)/\tau} \frac{\doublebar{u}_{L^{p_{j+1}}({B_k})}} {r^{(n+1)/p_{j+1}}}
	\biggr)^\tau
	\biggl(
	\frac{\doublebar{u}_{L^{p_j}(B_{k+1})}}{r^{(n+1)/p_j}} 
	\biggr)^{1-\tau}
.\end{align*}
Young's inequality implies that
\begin{align*}
\frac{\doublebar{u}_{L^{p_j}({B_k})}}{r^{(n+1)/p_j}}
&\leq
	\tau C(p_j,\sigma)\sigma^{-k\beta(1-\tau)/\tau} \frac{\doublebar{u}_{L^{p_{j+1}}({B_k})}} {r^{(n+1)/p_{j+1}}}
	+
	(1-\tau)
	\frac{\doublebar{u}_{L^{p_j}(B_{k+1})}}{r^{(n+1)/p_j}} 
.\end{align*}
Applying this bound to $k=0$ and iterating, we have that
\begin{align*}
\frac{\doublebar{u}_{L^{p_j}({B_0})}}{r^{(n+1)/p_j}}
&\leq
	\sum_{k=1}^{m+1} 
	C(p_j,\sigma)\tau(1-\tau)^{k-1}\sigma^{-k\beta(1-\tau)/\tau} \frac{\doublebar{u}_{L^{p_{j+1}}({B_k})}} {r^{(n+1)/p_{j+1}}}
	\\&\qquad+
	(1-\tau)^{m+1}
	\frac{\doublebar{u}_{L^{p_j}(B_{m+1})}}{r^{(n+1)/p_j}} 
.\end{align*}
We want to take the limit as $m\to \infty$. Choose $\sigma$ so that $(1-\tau)<\sigma^{\beta(1-\tau)/\tau}<1$; then the sum converges and we have that
\begin{align*}
\frac{\doublebar{u}_{L^{p_j}({B(X,r)})}}{r^{(n+1)/p_j}}
&\leq
	C(p_j)\frac{\doublebar{u}_{L^{p_{j+1}}({B(X,2r)})}} {r^{(n+1)/p_{j+1}}}
.\end{align*}
This completes the proof.
\end{proof}

We now apply Lemma~\ref{lem:subaverage} to solutions and their gradients.

\begin{cor}\label{cor:local-bound}
Suppose that $A$ is elliptic and satisfies the De Giorgi-Nash-Moser condition. If $\Div A\nabla u=0$ in $B(X,r)$ and $\nabla u\in L^2(B(X,r))$, then
\begin{equation*}\abs{u(X)}\leq C(p)\biggl(\fint_{B(X,r)}\abs{u}^p\biggr)^{1/p}
\quad\text{for any }0<p<\infty.\end{equation*}
\end{cor}
\begin{proof}
By formula~\eqref{eqn:local-bound:2}, we may apply Lemma~\ref{lem:subaverage} with $q=\infty$ and $p_0=2$.
\end{proof}

\begin{lem}\label{lem:local-bound:gradient}
Suppose that $A$ is elliptic. If $\Div A\nabla u=0$ in $B(X,2r)$ and $\nabla u\in L^2(B(X,2r))$, then
\begin{equation*}\biggl(\fint_{B(X,r)}\abs{\nabla u}^q\biggr)^{1/q}\leq C(p,q)\biggl(\fint_{B(X,r)}\abs{\nabla u}^p\biggr)^{1/p}
\end{equation*}
for any $0<p\leq\infty$ and any $0<q<p^+$.
\end{lem}

\begin{proof} 
By H\"older's inequality (if $q<2$) or by Lemma~\ref{lem:PDE2} (if $2<q<p^+$), we need only prove this lemma for $q=2$. But again by Lemma~\ref{lem:PDE2}, we may apply Lemma~\ref{lem:subaverage} to $\nabla u$ with $q=2$ and with $p^-<p_0<2$.
\end{proof}

We now consider some bounds valid in the special case of $t$-independent coefficients.
\begin{lem}[{\cite[Proposition~2.1]{AlfAAHK11}}]\label{lem:slabs}
Suppose that $A$ is elliptic and $t$-inde\-pen\-dent, and suppose that $Q\subset\RR^n$ is a cube. If $u$ satisfies $\Div A\nabla u=0$ in $2Q\times(t-\ell(Q),t+\ell(Q))$, then whenever $0< p< p^+$ we have that
\begin{align*}
\biggl(\fint_Q \abs{\nabla u(x, t)}^p\,dx\biggr)^{1/p}
&\leq
C
\biggl(\fint_{2Q}\fint_{t-\ell(Q)/4}^{t+\ell(Q)/4} \abs{\nabla u(x, s)}^2\,ds\,dx\biggr)^{1/2}
.\end{align*}
\end{lem}

We remark that if $A$ is $t$-independent and $\Div A\nabla u=0$ in $2Q\times(t-\ell(Q)/4,t+\ell(Q)/4)$, then $\Div A\nabla \partial_t u=0$ in the same region. Thus, if $s$ is near $t$ then $\partial_t \nabla u(\,\cdot\,,s)$ is in $L^p(Q)$, and so $s\mapsto \nabla u(\,\cdot\,,s)$ is continuous $(t-\ell(Q)/4,t+\ell(Q)/4)\mapsto L^p(Q)$.
Thus, $\nabla u(\,\cdot\,,t)$ is well-defined as a $L^p(Q)$ function.

\begin{cor}\label{cor:NTM:slabs} Suppose that $\Div A\nabla u=0$ in $\R^{n+1}_+$ and that $\widetilde N_+(\nabla u)\in L^p(\R^n)$ for some $0<p<p^+$. Then 
\begin{equation*}\sup_{t>0}\doublebar{\nabla u(\,\cdot\,,t)}_{L^p(\R^n)}\leq C\doublebar{\widetilde N_+(\nabla u)}_{L^p(\R^n)}.\end{equation*}
\end{cor}

\begin{proof}
Choose some $t>0$. If $Q\subset\R^n$ is a cube of side-length $t/C$, for a large constant~$C$, then by Lemma~\ref{lem:slabs}, we have that
\begin{equation*}\int_Q \abs{\nabla u(x,t)}^p\,dx \leq C \abs{Q}\inf_{y\in Q} \widetilde N_+(\nabla u)(y)^p\leq C \int_Q \widetilde N_+(\nabla u)(y)^p\,dy.\end{equation*}
Summing over a grid of such cubes completes the proof.
\end{proof}

Finally, we remark that if $p^+$ is large enough, then Morrey's inequality guarantees validity of the De Giorgi-Nash-Moser condition; furthermore, we may put a lower bound on the number $\alpha_0$ in terms of~$p^+$. This is a straightforward generalization of the results in \cite[Appendix~B]{AlfAAHK11}, where it was established that in dimension $n+1=3$, all $t$-independent elliptic matrices satisfy the De Giorgi-Nash-Moser condition.

\begin{lem}\label{lem:morrey} Let $A$ be elliptic and let $p^+$ be as in Section~\ref{sec:dfn:elliptic}. 

If $p^+>n+1$, then $A$ and $A^*$ satisfy the De Giorgi-Nash-Moser condition with exponent $\alpha$ for any $\alpha$ with $0<\alpha<1-(n+1)/p^+$.

If in addition $A$ is $t$-independent, then $A$ and $A^*$ satisfy the De Giorgi-Nash-Moser condition with exponent $\alpha$ for any $\alpha$ with $0<\alpha<1-n/p^+$.
\end{lem}

\begin{proof}
If $0<\alpha<1-(n+1)/p^+$ then the bound \eqref{eqn:holder} follows from Lemma~\ref{lem:PDE2} by Morrey's inequality; see, for example, \cite[Section~5.6.2]{Eva98}.

Suppose $A$ is $t$-independent.
If $\Div A\nabla u=0$ in some ball $B((x_0,t),2r)$ and $n<p_0<p^+$, then by Lemma~\ref{lem:slabs}, we have that $\nabla_\parallel u(\,\cdot\,,t)\in L^{p_0}(\Delta(x_0,r))$. But then by Morrey's inequality, we have that $u(\,\cdot\,,t)\in\dot C^\alpha(\Delta(x_0,r))$, where $\alpha=1-n/p_0$. Because $\Div A\nabla \partial_t u=0$, we have that $\partial_t u$ is locally bounded; thus, we may show that $u\in \dot C^\alpha(B(X_0,r))$ whenever $\Div A\nabla u=0$ in $B(X_0,2r)$. This completes the proof.
\end{proof}

\subsection{The fundamental solution}
\label{sec:fundamental}

We now discuss the fundamental solution. Let $A$ be an elliptic matrix such that both $A$ and $A^*$ satisfy the De Giorgi-Nash-Moser condition. We say that $\Gamma_{(x,t)}^A$ is a fundamental solution for $\Div A\nabla u=0$ if 
\begin{equation}
\label{eqn:FS:point}
\int_{\R^{n+1}} A\nabla\Gamma^A_{(x,t)}\cdot\nabla\varphi
=\varphi(x,t)\quad\text{for all }\varphi\in C^\infty_0(\R^{n+1})
\end{equation}
and if $\Gamma_{(x,t)}^A$ satisfies the bound
\begin{equation}
\label{eqn:FS:far:space}
\int_{B((x,t),2r)\setminus B((x,t),r)} \abs{\nabla \Gamma_{(x,t)}^A(y,s)}^2\,dy\,ds\leq Cr^{1-n}\quad\text{for all }r>0
.\end{equation}
Such a fundamental solution exists whenever $n+1\geq 2$; it was constructed in \cite{HofK07} (in the case $n+1\geq 3$), and another construction (valid for $n+1\geq2$) was provided in \cite{Ros12A}. In the case $n+1=2$, see also \cite{AusMT98}, \cite{KenR09} and~\cite{Bar13}. We remark that the fundamental solution is unique up to an additive constant.

We will need some additional properties of the fundamental solution. First, by H\"older's inequality, if $1\leq p<(n+1)/n$, then
\begin{equation}
\label{eqn:FS:near:1}
\int_{B((x,t),r)} \abs{\nabla \Gamma_{(x,t)}^A(y,s)}^p\,dy\,ds\leq C(p) r^{1+n-pn}
.\end{equation}
This was observed in \cite[Proposition~5.2]{Ros12A} and \cite[Theorem~3.1]{HofK07}.
Next, by Lem\-ma~\ref{lem:slabs} we have that if $A$ is $t$-independent and $1<p<p^+$ then
\begin{equation}
\label{eqn:FS:far:slices}
\int_{\abs{x-y}>r}\abs{\nabla\Gamma_{(x,t)}^A(y,s)}^p \,dy \leq \frac{C}{(r+\abs{s-t})^{n(p-1)}}
\quad\text{for all }r+\abs{s-t}>0
.\end{equation}

The fundamental solution is only defined up to an additive constant. If $n+1\geq 3$, then by the Poincar\'e inequality, the bound \eqref{eqn:FS:far:space}, and H\"older continuity of solutions, we have that $\lim_{\abs{(y,s)}\to\infty}\Gamma_{(x,t)}^A(y,s)$ exists; we can normalize by taking this limit to be zero. If we do so, then the formula~\eqref{eqn:FS:switch} is valid; that is, $\Gamma_{(x,t)}^A(y,s)=\overline{\Gamma_{(y,s)}^{A^*}(x,t)}$.
If $n+1=2$ then this limit does not exist. However, there is still a possible normalization such that formula~\eqref{eqn:FS:switch} is valid.

Notice that because $\Gamma_{(x,t)}^A(y,s)$ is a solution to an elliptic equation in the variables $(x,t)$, we may apply the De Giorgi-Nash-Moser condition; by applying the Poincar\'e inequality to the bound~\eqref{eqn:FS:far:space}, we see that if $\abs{(x,t)-(x',t')}<\frac{1}{2}\abs{(x,t)-(y,s)}$, then 
\begin{equation}
\label{eqn:FS:holder}
\abs{\Gamma_{(x,t)}^A(y,s)-\Gamma_{(x',t')}^A(y,s)}
\leq \frac{C\abs{(x,t)-(x',t')}^\alpha}{\abs{(x,t)-(y,s)}^{n-1+\alpha}}
.\end{equation}

Moreover, by uniqueness of the fundamental solution we have that
\begin{align}
\label{eqn:FS:vertical}
\nabla\Gamma_{(x,t)}^A(y,s)&=\nabla\Gamma_{(x,t+\delta)}^A(y,s+\delta)
,\\
\label{eqn:FS:conjugate}
\overline{\nabla\Gamma_{(x,t)}^A(y,s)}&=\nabla\Gamma_{(x,t)}^{\overline A}(y,s)
\end{align}
for any $x$, $y\in\R^n$ and any $t$, $s$, $\delta\in\R$.

We remark that formulas~\eqref{eqn:Div:plus} and \eqref{eqn:Div:minus} are straightforward to derive from the properties of the fundamental solution given above.

Finally, we have the following lemma.
\begin{lem}\label{lem:FS:lax-milgram}
Suppose that $\F:\R^{n+1}\mapsto\C^{n+1}$ is smooth and compactly supported. Let 
\begin{equation*}
u(x,t) = -\int_{\R^{n+1}} \Gamma^A_{(y,s)}(x,t) \Div\F(y,s)\,dy\,ds
 = \int_{\R^{n+1}} \overline{\nabla\Gamma^{A^*}_{(x,t)} }\cdot\F(y,s)\,dy\,ds
.\end{equation*}
Then $\Div A\nabla u=\Div \F$ in the sense of formula~\eqref{eqn:solution}, and $u\in\dot W^2_1(\R^{n+1})$.

Furthermore, we have the bound
\begin{equation}\label{eqn:FS:lax-milgram}
\doublebar{\nabla u}_{L^2(\R^{n+1})}\leq \frac{1}{\lambda}\doublebar{\vec F}_{L^2(\R^{n+1})}.\end{equation}
\end{lem}

\begin{proof} If $n+1\geq 3$, we take our construction of the fundamental solution from \cite{HofK07}; their construction gives that $\Div A\nabla u=\Div \F$ and $u\in\dot W^2_1(\R^{n+1})$.

If $n+1=2$, the fact that $\Div A\nabla u=\Div \F$ may be seen by interchanging the order of integration in formula~\eqref{eqn:solution}. By the bound~\eqref{eqn:FS:near:1}, we have that $\nabla u\in L^p_{loc}(\R^{n+1})$ for any $p<2$; in particular, $\nabla u\in L^p_{loc}(\R^{n+1})$ for some $p>p^-$ and so by Lemma~\ref{lem:PDE2} we have that $\nabla u\in L^2_{loc}(\R^{n+1})$.

We wish to show that $\nabla u$ is \emph{globally} in $L^2_{loc}(\R^{n+1})$.
Suppose that $\F$ is supported in some ball $B(0,r)$. If
$\abs{(x,t)}=R>8r$, then
\begin{align*}
\abs{u(x,t)} 
&= 
	\abs[bigg]{\int_{B(0,r)} \overline{\nabla \Gamma_{(x,t)}^{A^*}(y,s)}\cdot \F(y,s)\,dy\,ds}
\\&\leq
	\doublebar{\F}_{L^2(\R^{1+1})}
	\biggl(\int_{B(0,r)} \abs{\nabla \Gamma_{(x,t)}^{A^*}(y,s)}^2\,dy\,ds\biggr)^{1/2}
.\end{align*}
Because $\Gamma_{(x,t)}^{A^*}$ is a solution to $\Div {A^*}\nabla \Gamma_{(x,t)}^{A^*}=0$ away from the point $(x,t)$, we may apply the Caccioppoli inequality to see that
\begin{align*}
\abs{u(x,t)} 
&\leq
	\frac{C}{r}\doublebar{\F}_{L^2(\R^{1+1})}
	\biggl(\int_{B(0,2r)} \abs{ \Gamma_{(x,t)}^{A^*}(y,s)-\Gamma_{(x,t)}^{A^*}(0,0)}^2\,dy\,ds\biggr)^{1/2}
.\end{align*}
By the bound~\eqref{eqn:FS:holder} and formula~\eqref{eqn:FS:switch},
\begin{align*}
\abs{u(x,t)} 
&\leq
	C\frac{r^\alpha}{R^\alpha}\doublebar{\F}_{L^2(\R^{1+1})}
.\end{align*}
Let $A_j=B(0,2^{j+1}r)\setminus B(0,2^jr)$; then if $j\geq 4$, by the Caccioppoli inequality we have that
\begin{equation*}\int_{A_j}\abs{\nabla u}^2 \leq C2^{-2j} r^{-2} \int_{\widetilde A_j} \abs{u}^2 \leq C 2^{-2j\alpha}\doublebar{\F}_{L^2(\R^{1+1})}^2
\end{equation*}
and so $\nabla u\in L^2(\R^{1+1}\setminus B(0,8r))$. Since $\nabla u\in L^2_{loc}(\R^{1+1})$, we have that $u\in \dot W^2_1(\R^{1+1})$.

We must now establish the quantitative bound \eqref{eqn:FS:lax-milgram}. We will use the following generalization of the Lax-Milgram lemma to complex spaces.
\begin{thm}[{\cite[Theorem~2.1]{Bab70}}]
\label{thm:lax-milgram}
Let $H_1$ and $H_2$ be two Hilbert spaces, and let $B$ be a bounded bilinear form on $H_1\times H_2$ that is coercive in the sense that
	\begin{equation*}\sup_{w\in H_1\setminus\{0\}} \frac{\abs{B(w,v)}}{\doublebar{w}_{H_1}}\geq \lambda \doublebar{v}_{H_2},\quad
	\sup_{w\in H_2\setminus\{0\}} \frac{\abs{B(u,w)}}{\doublebar{w}_{H_2}}\geq \lambda \doublebar{u}_{H_1}\end{equation*}
	for every $u\in {H_1}$, $v\in {H_2}$, for some fixed $\lambda>0$. Then for every linear functional $T$ defined on ${H_2}$ there is a unique $u_T\in {H_1}$ such that $B(u_T,v)=\overline{T(v)}$. Furthermore, $\doublebar{u_T}_{H_1}\leq \frac{1}{\lambda}\doublebar{T}_{H_1\mapsto H_2}$.
\end{thm}
Apply Theorem~\ref{thm:lax-milgram} to the bilinear form and linear operator
\begin{equation*}B(\xi,\eta) = \int_{\R^{n+1}} \nabla\bar\eta\cdot A\nabla \xi,
\qquad \overline{T(\eta)}=\int_{\R^{n+1}} \nabla\bar\eta\cdot \F
\end{equation*}
on the spaces $H_1=H_2=\dot W^2_1(\R^{n+1})$. The function $u_T$ produced by Theorem~\ref{thm:lax-milgram} satisfies the formula $\Div A\nabla u_T=\Div \F$ and the bound
\begin{equation*}\doublebar{\nabla u_T}_{L^2(\R^{n+1})}\leq \frac{1}{\lambda}\doublebar{\vec F}_{L^2(\R^{n+1})},\end{equation*}
and is unique among functions in $\dot W^2_1(\R^{n+1})$. But $u\in \dot W^2_1(\R^{n+1})$ with $\Div A\nabla u=\Div \F$ as well; thus, $u=u_T$ and the bound \eqref{eqn:FS:lax-milgram} is valid.
\end{proof}

%% file: sec-5-bounded.tex
\chapter{Boundedness of Integral Operators}
\label{chap:bounded}

In this chapter, we will bound the integral operators $\Pi^A$, $\D^A$ and $\s^A$ used to construct solutions to boundary-value problems; that is, we will prove Theorem~\ref{thm:bounded}. We will begin with the Newton potential $\Pi^A$ in Section~\ref{sec:newton:bounded} and will move to the layer potentials $\D^A$ and $\s^A$ in Section~\ref{sec:layers:bounded}.

\section{Boundedness of the Newton potential}
\label{sec:newton:bounded}
In this section we bound the Newton potential $\Pi^A:L(p,\theta,q)\mapsto \dot W(p,\theta,q)$; that is, we establish the bound \eqref{eqn:space:newton} for $p$, $q$, $\theta$ as in Theorem~\ref{thm:bounded}. We begin with smooth, compactly supported functions.

\begin{figure}[t]
\begin{tikzpicture}[scale=2]
\figureaxes

\fill [bounded] 
(1,1+\alph/\enn) node [black] {$\circ$} node [black, right] {$(1,1+\alpha/n)$}--
(1-\alph,1) node [black] {$\circ$} node [black, above, at = {(1-0.2-\alph,1)}] {$(1-\alpha,1)$}--
(1-\alph,1/2-\eps) --
(1,1/2-\eps) node [black] {$\circ$} node [black, right] {$(1,1/p^+)$}-- (1,1)--cycle;

\draw [line width=0.5pt, dotted] (0,1/2)--(1,1/2);
\end{tikzpicture}
\begin{tikzpicture}[scale=2]
\figureaxes
\fill [bounded] 
(0,0) -- (\alph,0) -- (\alph,1/2+\eps) -- 
(0,1/2+\eps) node [black] {$\circ$} node [black, left] {$1/p^-$}-- cycle;

\fill [bounded] 
(1,1+\alph/\enn)--
(1-\alph,1)--
(1-\alph,1/2-\eps) --
(1,1/2-\eps)-- (1,1)--cycle;

\draw [boundary bounded] (0,0)--(\alph,0)  node [black] {$\circ$} node [black, below] {$\alpha$};

\draw [line width=0.5pt, dotted] (0,1/2)--(1,1/2);
\end{tikzpicture}
\begin{tikzpicture}[scale=2]
\figureaxes
\fill [bounded] (0,0)--
(\alph,0) --
(1,1/2-\eps)--
(1,1+\alph/\enn)--
(1-\alph,1)--
(0,1/2+\eps)--
cycle;
\draw [boundary bounded] (0,0)--(\alph,0);
\draw [line width=0.5pt, dotted] (0,1/2)--(1,1/2);
\end{tikzpicture}
\caption{Values of $p$, $\theta$ such that $\Pi^A$ is bounded $L(p,\theta,q)\mapsto\dot W(p,\theta,q)$ given by Theorem~\ref{thm:newton:preliminary} (left), duality (middle) and interpolation (right).}
\label{fig:newton:1}
\end{figure}
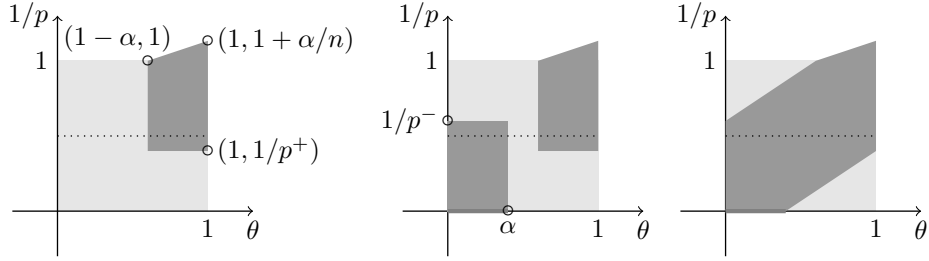

\begin{thm} \label{thm:newton:preliminary}
Let $\F$ be smooth and compactly supported in the half-space $\R^{n+1}_+$. Suppose that $A$ and $A^*$ are elliptic, $t$-independent matrices that satisfy the De Giorgi-Nash-Moser condition. Suppose that $p^-<q<p^+$, that $1-\alpha<\theta<1$, and that $1/p^+<1/p<1+(\alpha+\theta-1)/n$.

Then the bound \eqref{eqn:space:newton} is valid; that is,
\begin{multline*}
\int_{\R^{n+1}_+} 
\biggl(\fint_{\Omega(x,t)}\abs{\nabla \Pi^A\F}^q\biggr)^{p/q}\,t^{p-1-p\theta}\,dx\,dt
\\\leq 
C(p,\theta,q)\int_{\R^{n+1}_+} 
\biggl(\fint_{\Omega(x,t)}\abs{\F}^q\biggr)^{p/q}\,t^{p-1-p\theta}\,dx\,dt
.\end{multline*}
\end{thm}
The acceptable values of $p$ and $\theta$ are shown on the right of Figure~\ref{fig:newton:1}. After proving Theorem~\ref{thm:newton:preliminary} we will extend to the full range of $p$, $\theta$ allowed by Theorem~\ref{thm:bounded} using duality and interpolation.

\begin{proof}
Let $\G$ be a grid of dyadic Whitney cubes in $\R^{n+1}_+$.
We remark that 
\begin{equation*}\sum_{Q\in\G}\ell(Q)^{n+p-p\theta}
\biggl(\fint_{Q} \abs{\F}^q\biggr)^{p/q}
\approx 
\int_{\R^{n+1}_+} 
\biggl(\fint_{\Omega(x,t)}\abs{\F}^q\biggr)^{p/q}\,t^{p-1-p\theta}\,dx\,dt.\end{equation*}
The formulation in terms of dyadic cubes will be more convenient.

Let $\sum_{Q\in\G} \varphi_Q$ be a smooth partition of unity with $\varphi_Q$ supported in $Q'=(5/4)Q$.
We then write
\begin{equation*}
\Pi^A\F(x,t)=u(x,t)=\sum_{Q\in\G} u_{Q}(x,t)
\end{equation*}
where
\begin{align*}
u_{Q}(x,t)
&=
	\int_{Q'} \overline{\nabla\Gamma_{(x,t)}^{A^*}(y,s) } \cdot \F(y,s)\, \varphi_Q(y,s)\,ds\,dy
.\end{align*}
By Lemma~\ref{lem:FS:lax-milgram}, we have that $\Div A\nabla u_Q = \Div (\varphi_Q\F)$. By Lemmas~\ref{lem:FS:lax-milgram} and~\ref{lem:PDE2}, we have that
\begin{equation}\label{eqn:inhomogeneous:local}
\int_{N(Q)} \abs{\nabla u_Q(x,t)}^q\,dx\,dt \leq C(q)\int_{Q'} \abs{\F(x,t)}^q\,dx\,dt\end{equation}
for any $q$ with $2\leq q<p^+$. 
Here $N(Q)$ is the union of dyadic cubes $R\in \G$ with $\dist(Q,R)=0$; observe that $2Q\subset N(Q)\subset 5Q$.

Let $p^-<q<2$ so that $2<q'<p^+$; we claim that the bound \eqref{eqn:inhomogeneous:local} is still valid. To see this, recall that
$\Gamma_{(y,s)}^A(x,t)=\overline{\Gamma_{(x,t)}^{A^*}(y,s)}$. Therefore, the adjoint $(\nabla\Pi^A)^*$ to the operator $\nabla\Pi^A$ is $\nabla \Pi^{A^*}$. So since $\nabla\Pi^{A^*}$ is bounded $L^{q'}(N(Q))\mapsto L^{q'}(N(Q))$, we have that $\nabla\Pi^A$ is bounded $L^q(N(Q))\mapsto L^q(N(Q))$.

We now consider $\R^{n+1}_+\setminus N(Q)$. We divide $\R^{n+1}_+$ into thin slices as follows. Loosely, we would like $A_{j,k}(Q)$ to be the union of the dyadic cubes $R\in\G$ with $\ell(R)=\dist(R,\partial\R^{n+1}_+)=2^j\ell(Q)$ and with 
\begin{equation*}\dist(R,Q)\approx \diam(A_{j,k}(Q))\approx 2^k\max(1,2^j)\ell(Q).\end{equation*}
Estimating the volume of this region, we see that there are $C(2^k\max(1,2^{-j}))^n$ such cubes.

More precisely, we define $A_{j,k}(Q)$ as follows.
For any $R\in\G$ let $\pi(R)$ be the projection of $R$ onto $\R^{n}$. If $j\leq 0$, let $P_j(Q)=\pi(Q)$, and if $j\geq 1$ then let $P_j(Q)=\pi(R)$ where $R$ is the unique cube in $\G$ with $\ell(R)=2^j\ell(Q)$ and $\pi(Q)\subset\pi(R)$.

Then let
\begin{equation*}A_{j,0}(Q)=17P_j(Q)\times (2^j\ell(Q),2^{j+1}\ell(Q)) \setminus N(Q),\end{equation*}
and for $k\geq 1$ let \begin{equation*}A_{j,k}(Q)=((2^{k+4}+1)P_j(Q)\setminus(2^{k+3}+1)P_j(Q))\times (2^j\ell(Q),2^{j+1}\ell(Q)).\end{equation*}
Notice that $N(Q)\subset 5Q$ and so $N(Q)\cap A_{j,k}(Q)=\emptyset$ for any $k\geq 1$. We use odd numbers so that $A_{j,k}$ will be a union of dyadic Whitney cubes.

We will need two larger versions of $A_{j,k}$. First, let $A_{j,k}'(Q)=\cup_{R\subset A_{j,k}(Q)} (5/4)R$. 
Next, let
\begin{equation*}\widetilde A_{j,0}=32P_j(Q)\times (-\max(2,2^{j+2})\ell(Q),\max(2,2^{j+2})\ell(Q)) \setminus (3/2)Q\end{equation*}
and
\begin{align*}
\widetilde A_{j,k}&=\bigl(2^{k+5}P_j(Q)
\setminus 2^{k+2}P_j(Q)\bigr)
\\&\quad\times (-2^k\max(2,2^{j+2})\ell(Q),2^k\max(2,2^{j+2})\ell(Q)).
\end{align*}

By Lemma~\ref{lem:local-bound:gradient}, we have that 
\begin{equation*}
\sum_{R\subset A_{j,k}(Q)} \biggl(\fint_R \abs{\nabla u_Q}^q\biggr)^{p/q} \ell(R)^{n+p-p\theta}
\leq
	C(2^{j}\ell(Q))^{p-1-p\theta}
	\int_{A_{j,k}'(Q)} \abs{\nabla u_Q}^p
.\end{equation*}
If $p<p^+$ we may use Lemma~\ref{lem:slabs} to bound the right-hand side; for simplicity of exposition and because we will later use duality and interpolation to generalize to a broader range of~$p$, we will not consider the case $p\geq p^+$. Because $\Div A\nabla u_Q=0$ in $\widetilde A_{j,k}$, we have that
\begin{equation*}\int_{A_{j,k}'} \abs{\nabla u_Q}^{p}\leq C 2^j\ell(Q) (2^{k}\max(1,2^{j})\ell(Q))^{n-p}
\biggl(\fint_{\widetilde{A}_{j,k}} \abs{u_Q}^2\biggr)^{p/2}
.\end{equation*}
But if $(x,t)\in \widetilde A_{j,k}(Q)$, then for any $q$ with $q'<p^+$ we have that
\begin{align*}
\abs{u_Q(x,t)}
&\leq
	\int_{Q'}
	\abs{\nabla\Gamma_{(x,t)}^{A^*}(y,s)}
	\,\abs{\F(y,s)}\,ds\,dy
\\&\leq
	\biggl(\int_{Q'}\abs{\nabla\Gamma_{(x,t)}^{A^*}(y,s)}^{q'} \,dy\,ds\biggr)^{1/q'}
	\biggl(\int_{Q'}\abs{\F}^q\biggr)^{1/q}
\\&\leq
	\frac{C\ell(Q)} {(2^k\max(1,2^j))^{n-1+\alpha}}
	\biggl(\fint_{Q'}\abs{\F}^q\biggr)^{1/q}
\end{align*}
and so
\begin{multline}
\label{eqn:newton:proof:1}
\sum_{R\subset A_{j,k}(Q)} \biggl(\fint_R \abs{\nabla u_Q}^q\biggr)^{p/q} \ell(R)^{n+p-p\theta}
\\\leq
	C\ell(Q)^{n+p-p\theta}
		\frac{2^{j(p-p\theta)}} {(2^k\max(1,2^j))^{np+p\alpha-n}}
		\biggl(\fint_{Q'}\abs{\F}^q\biggr)^{p/q}
.\end{multline}

Suppose $p\leq 1$. Then
\begin{align*}
\sum_{R\in\G} \biggl(\fint_R \abs{\nabla u}^q\biggr)^{p/q}\ell(R)^{n+p-p\theta}
&\leq
	\sum_{R\in\G} 
	\biggl(\sum_{Q\in\G}\biggl(\fint_R \abs{\nabla u_Q}^q\biggr)^{1/q}\biggr)^p
	\ell(R)^{n+p-p\theta}
\\&\leq
	\sum_{Q\in\G}\sum_{R\in\G} 
	\biggl(\fint_R \abs{\nabla u_Q}^q\biggr)^{p/q}
	\ell(R)^{n+p-p\theta}
.\end{align*}

Let $A_{j,-1}(Q)=N(Q)\cap (\R^n\times(2^j,2^{j+1}))$; notice that $A_{j,-1}=\emptyset$ unless $j=-1$, $0$ or~$1$. By the bound~\eqref{eqn:inhomogeneous:local}, formula~\eqref{eqn:newton:proof:1} is still valid for $k=-1$, and $\R^{n+1}_+=\cup_{j=-\infty}^\infty \cup_{k=-1}^\infty A_{j,k}(Q)$.
Thus,
\begin{multline*}
\sum_{R\in\G} \biggl(\fint_R \abs{\nabla u}^q\biggr)^{p/q} \ell(R)^{n+p-p\theta}
\\\leq
	C\sum_{Q\in\G}
	\sum_{j=-\infty}^\infty \sum_{k=-1}^\infty
	\frac{\ell(Q)^{n+p-p\theta}
	2^{j(p-p\theta)}
	} {(2^k\max(1,2^j))^{np-n+p\alpha}}
	\biggl(\fint_{Q'}\abs{\F}^q\biggr)^{p/q}
.\end{multline*}
The sums in $j$, $k$ converge provided $\theta<1$ and $np-n+p\alpha>p-p\theta$; the second condition is equivalent to $1+(\alpha+\theta-1)/n>1/p$. Thus the bound~\eqref{eqn:space:newton} is valid for $p$ with $1\leq 1/p<1+(\alpha+\theta-1)/n$.

Now, suppose that $1\leq p<p^+$. Again we have that
\begin{equation*}
\sum_{R\in\G} \biggl(\fint_R \abs{\nabla u}^q\biggr)^{p/q}\ell(R)^{n+p-p\theta}
\leq
	\sum_{R\in\G} 
	\biggl(\sum_{Q\in\G}\biggl(\fint_R \abs{\nabla u_Q}^q\biggr)^{1/q}\biggr)^p
	\ell(R)^{n+p-p\theta}
.\end{equation*}
To exploit the bounds established above, we write the inner sum (for fixed~$R$) as
\[\biggl(\sum_{Q\in\G}\biggl(\fint_R \abs{\nabla u_Q}^q\biggr)^{1/q}\biggr)^{p}
=\biggl(\sum_{j=-\infty}^\infty \sum_{k=-1}^\infty \sum_{Q\in A_{j,k}^{-1}(R)}
\biggl(\fint_R \abs{\nabla u_Q}^q\biggr)^{1/q}\biggr)^{p}
\]
where $A_{j,k}^{-1}(R)$ is the set of all cubes $Q$ that satisfy $A_{j,k}(Q)\supset R$. We wish to move the summation outside the parentheses.
Because geometric series are convergent, we may use H\"older's inequality to show that if $\varepsilon>0$, then
\begin{multline*}
	\biggl(\sum_{j=-\infty}^\infty \sum_{k=-1}^\infty \sum_{Q\in A_{j,k}^{-1}(R)}
	\biggl(\fint_R \abs{\nabla u_Q}^q\biggr)^{1/q}
	\biggr)^p
\\\leq
	\sum_{j=-\infty}^\infty \sum_{k=-1}^\infty 
	\frac{C(\varepsilon) 2^{-j\varepsilon}} {(2^k\max(1,2^j))^{-2\varepsilon}}
	\biggl(\sum_{Q\in A_{j,k}^{-1}(R)}
	\biggl(\fint_R \abs{\nabla u_Q}^q\biggr)^{1/q}
	\biggr)^p	
.\end{multline*}
Observe that there are $C(2^k\max(1,2^j))^n$ cubes $Q$ with $R\subset A_{j,k}(Q)$, and so again by H\"older's inequality we have that
\begin{align*}
\biggl(\sum_{Q\in A_{j,k}^{-1}(R)}
	\biggl(\fint_R \abs{\nabla u_Q}^q\biggr)^{1/q}\biggl)^p
&\leq
	\frac{C} {(2^k\max(1,2^j))^{n-np}}
	\sum_{Q\in A_{j,k}^{-1}(R)}
	\biggl(\fint_R \abs{\nabla u_Q}^q\biggr)^{p/q}	
.\end{align*}
Combining the above bounds yields that
\begin{multline*}
\sum_{R\in\G} \biggl(\fint_R \abs{\nabla u}^q\biggr)^{p/q}\ell(R)^{n+p-p\theta}
\\\begin{aligned}
&\leq
	\sum_{R\in\G}
	\sum_{j=-\infty}^\infty \sum_{k=-1}^\infty 
	\frac{C(\varepsilon) 2^{-j\varepsilon} \ell(R)^{n+p-p\theta}} {(2^k\max(1,2^j))^{n-np-2\varepsilon}}
	\sum_{Q\in A_{j,k}^{-1}(R)}
	\biggl(\fint_R \abs{\nabla u_Q}^q\biggr)^{p/q}	
.\end{aligned}\end{multline*}
We may now interchange the order of summation. By the bound~\eqref{eqn:newton:proof:1}, we have that
\begin{multline*}
\sum_{R\in\G} \biggl(\fint_R \abs{\nabla u}^q\biggr)^{p/q}\ell(R)^{n+p-p\theta}
\\\begin{aligned}
&\leq
	\sum_{j=-\infty}^\infty \sum_{k=-1}^\infty 	
	\frac{C(\varepsilon)2^{j(p-p\theta-\varepsilon)}} {(2^k\max(1,2^j))^{p\alpha-2\varepsilon}}
	\sum_{Q\in\G}
	\ell(Q)^{n+p-p\theta}
		\biggl(\fint_{Q'}\abs{\F}^q\biggr)^{p/q}
.\end{aligned}\end{multline*}
If we can choose $\varepsilon>0$  such that $p\alpha-2\varepsilon>p-p\theta-\varepsilon>0$, then the sums in $j$ and~$k$ will converge and we will have completed the proof. Solving these inequalities, we see that it suffices to require $1-\alpha<\theta<1$, as desired.
\end{proof}

Smooth, compactly supported functions $\F$ are dense in $L(p,\theta,q)$ for any $p$, $q<\infty$. Thus, we may extend $\Pi^A$ to a bounded operator $L(p,\theta,q)\mapsto \dot W(p,\theta,q)$ for $p$, $q$, $\theta$ as in Theorem~\ref{thm:newton:preliminary}; the acceptable values of $(\theta,1/p)$ are shown on the left-hand side of Figure~\ref{fig:newton:1}.

Now, recall that $\Gamma_X^A(Y)=\overline{\Gamma_Y^{A^*}(X)}$ and so the adjoint to $\nabla\Pi^A$ is $\nabla\Pi^{A^*}$. Thus by duality, we have that the bound~\eqref{eqn:space:newton} is valid if $0<\theta<\alpha$, $p^-<p\leq\infty$, and if $p^-<q<p^+$.

Thus, $\Pi^A$ is bounded $L(p,\theta,q)\mapsto\dot W(p,\theta,q)$ provided $p^-<q<p^+$ and provided the point $(\theta,1/p)$ lies in one of the two regions indicated in the middle of Figure~\ref{fig:newton:1}.

We may use interpolation to fill in the gaps; see 
Theorem~\ref{thm:interpolation:L} and the interpolation properties \eqref{eqn:operator:interpolation:real} and~\eqref{eqn:operator:interpolation:complex}. Thus,  $\nabla\Pi^A$ is bounded $L(p,\theta,q)\mapsto L(p,\theta,q)$ for all $p$, $\theta$ in the convex region  shown on the right side of Figure~\ref{fig:newton:1}. (This is the same region as in Figure~\ref{fig:bounded}.) Thus, the bound~\eqref{eqn:space:newton} in Theorem~\ref{thm:bounded} is valid.

\section{Boundedness of the double and single layer potentials}
\label{sec:layers:bounded}

In this section, we will complete the proof of Theorem~\ref{thm:bounded} by establishing the bounds \eqref{eqn:space:D} and \eqref{eqn:space:S} on $\D$ and~$\s$; we will do this by bounding the combined operator~$\T$ of formula~\eqref{eqn:T}.

We begin by discussing the ways in which~$\T$ is a well-defined operator on $\dot B^{p,p}_\theta(\R^n)$.
Suppose that $A$ is $t$-independent and that $A$ and $A^*$ satisfy the De Giorgi-Nash-Moser condition. By the bound \eqref{eqn:FS:far:slices} we have that if $1<p<p^+$ then $\nabla\Gamma_{(x,t)}^{A^*}(\,\cdot\,,0)\in L^p(\R^n)$ with $\doublebar{\nabla\Gamma_{(x,t)}^{A^*}(\,\cdot\,,0)}_{L^p(\R^n)}\leq Ct^{-n(1-1/p)}$.
Thus, the integrals in the definition of $\T$ converge absolutely whenever $\vec f\in L^p(\R^n)$, $p^-<p<\infty$. Furthermore,
\begin{equation}
\label{eqn:pointwise}
\abs{\T\vec f(x,t)} \leq C t^{-n/p} \doublebar{\vec f}_{L^p(\R^n)}
\quad\text{provided }p^-<p<\infty
.\end{equation}
Because smooth, compactly supported functions are dense in $\dot B^{p,p}_\theta(\R^n)$ whenever $0<p<\infty$ and $0<\theta<1$, once we have established boundedness we will be able to extend $\T$ to a well-defined linear operator on~$\dot B^{p,p}_\theta(\R^n)$.

If $\vec f\in L^\infty(\R^n)$, we may define $\T\vec f$ up to an additive constant via the formula
\begin{equation}\label{eqn:T:difference}
\T\vec f(x,t)-\T\vec f(x',t')
=\int_{\R^n} \overline{\nabla_{A^*}\Gamma_{(x,t)}^{A^*}(y,0)
-\nabla_{A^*}\Gamma_{(x',t')}^{A^*}(y,0)} \cdot \vec f(y)\,dy
.\end{equation}
Compare to the definition~\eqref{eqn:T} of $\T f(x,t)$.
Recall that $\nabla_{A} u=(\nabla_\|u, \e_{n+1}\cdot A\nabla u)$. By the bound~\eqref{eqn:FS:holder} and Lemma~\ref{lem:slabs}, the integral in formula~\eqref{eqn:T:difference} converges absolutely whenever $(x,t)\in\R^{n+1}_+$ and $(x',t')\in\R^{n+1}_+$. 

We will need the following useful property of~$\T$.
\begin{lem}\label{lem:T-constant}
If $\vec f(y)=\vec v$ for any constant vector $\vec v$, then $\T\vec f$ is also a constant.
\end{lem}

\begin{proof}
Let $\varphi_R$ be smooth, supported in $\Delta(0,2R)$ and identically $1$ in $\Delta(0,R)$ with $\abs{\nabla\varphi_R}\leq C/R$. Again by the bound~\eqref{eqn:FS:holder} and Lemma~\ref{lem:slabs}, we have that $\T \vec v (x,t)-\T \vec v(x',t')=\lim_{R\to\infty}\T (\varphi_R\vec v) (x,t)-\T (\varphi_R\vec v)(x',t')$. 
Recall that $\T\vec f=(\s\nabla_\parallel)\cdot \vec f_\parallel -\D  f_\perp$. If $\vec f_\parallel=\vec v_\parallel$ is a constant, then integrating by parts we have that
\begin{multline*}(\s\nabla_\parallel)\cdot (\varphi_R\vec v_\parallel)(x,t)-(\s\nabla_\parallel)\cdot (\varphi_R\vec v_\parallel)(x',t')
\\=
\int_{\R^n} \bigl(\Gamma_{(y,0)}^{A}(x,t)
-\Gamma_{(y,0)}^{A}(x',t')\bigr)\, \nabla\varphi_R(y)\cdot \vec v_\parallel\,dy
.\end{multline*}
By the bound~\eqref{eqn:FS:holder}, this integral goes to zero as $R\to \infty$; thus $(\s\nabla_\parallel)\cdot \vec v_\parallel$ is a constant.

If $f_\perp$ is a constant, then so is $\D f_\perp$; see \cite[Remark~4.2]{HofMitMor}, or apply the definition~\eqref{eqn:conormal:1} of conormal derivative to compute $\D(\varphi_R f_\perp)$.
\end{proof}

We remark that the condition that $\T\vec v$ be constant for any constant~$\vec v$ is a necessary condition for $\T$ to be a well-defined operator on $\dot B^{p,p}_\theta(\R^n)$ for any $\theta>0$, because elements of $\dot B^{p,p}_\theta(\R^n)$ are only defined up to additive constants. Furthermore, it is straightforward to establish that if $\vec f(x_0)=0$ for some $x_0$ and if $\vec f\in \dot C^\theta(\R^n)=\dot B^{\infty,\infty}_\theta(\R^n)$ for some $0<\theta<\alpha$, then the integral \eqref{eqn:T:difference} converges absolutely.
Thus, the operator $\T$ may be defined on such spaces via the formula~\eqref{eqn:T:difference} and the condition $\T\vec v(x,t)-\T\vec v(x',t')=0$ for all constant vectors~$\vec v$.

\subsection{Boundedness in the range $p^-<p<p^+$} \label{sec:bounded:midrange}
In this section we will establish the bounds \eqref{eqn:space:D} and \eqref{eqn:space:S} under the additional assumption that $p^-<p<p^+$. We begin by proving the following theorem.
\begin{thm} \label{thm:space:vertical}
Suppose that $A$ is elliptic, $t$-independent and that both $A$ and $A^*$ satisfy the De Giorgi-Nash-Moser condition. If $0<\theta<1$ and $p^-<p<\infty$, then
\begin{equation}\label{eqn:space:vertical}
\int_{\R^{n+1}_+} \abs{\partial_t \T \vec f(x,t)}^p t^{p-1-p\theta} \,dt\,dx
\leq
	C(p,\theta)\doublebar{\vec f}_{\dot B^{p,p}_{\theta}(\R^n)}^p
.\end{equation}
\end{thm}

\begin{proof} 
By Lemma~\ref{lem:T-constant}, $\partial_t \T \vec f(x,t)=\partial_t\T  (\vec f-\vec f(x))(x,t)$. So by the definition \eqref{eqn:T} of $\T\vec f$,
\begin{multline*}
	\int_{\R^{n+1}_+} \abs{\partial_t \T \vec f(x,t)}^p t^{p-1-p\theta}\,dx\,dt
\\\begin{aligned}
&\leq
	C\int_0^\infty t^{p-1-p\theta}
	\int_{\R^n} 
	\biggl(\int_{\R^n} \abs{\vec f(y)-\vec f(x)}
	\,\abs{\nabla\partial_t\Gamma_{(x,t)}^{A^*}(y,0)} \,dy\biggr)^p
	\,dx\,dt
.\end{aligned}\end{multline*}
Recall that $\nabla \Gamma_{(x,t)}^{A^*}(y,0)$ denotes the gradient in $y$ and $s$ of $\Gamma_{(x,t)}^{A^*}(y,s)$ evaluated at $s=0$.

Applying H\"older's inequality, we see that if $\beta$ is a constant then
\begin{multline*}
\biggl(\int_{\R^n} \abs{\vec f(y)-\vec f(x)}
	\,\abs{\nabla\partial_t\Gamma_{(x,t)}^{A^*}(y,0)} \,dy\biggr)^p
\\\begin{aligned}
&\leq
	\int_{\R^n} \biggl( \frac{\abs{\vec f(y)-\vec f(x)}} {(t+\abs{x-y})^\beta} \biggr)^p dy
	\biggl(\int_{\R^n} \biggl(
	\frac{\abs{\nabla\partial_t\Gamma_{(x,t)}^{A^*}(y,0) }} {(t+\abs{x-y})^{-\beta }} \biggr)^{p'}\,dy\biggr)^{p/p'}
	.
\end{aligned}\end{multline*}
Let $A_0=\Delta(x,4t)$, and if $j\geq 1$ then let $A_j=\Delta(x,2^{j+2}t)\setminus \Delta(x,2^{j+1}t)$. Then
\begin{multline*}
\int_{\R^n} \Bigl((t+\abs{x-y})^\beta \abs{\nabla\partial_t\Gamma_{(x,t)}^{A^*}(y,0) }\Bigr)^{p'}\,dy
\\\begin{aligned}
&=
	\sum_{j=0}^\infty
	\int_{A_j} \Bigl((t+\abs{x-y})^\beta \abs{\nabla\partial_t\Gamma_{(x,t)}^{A^*}(y,0) }\Bigr)^{p'}\,dy
\\&\leq
	C\sum_{j=0}^\infty (2^j t)^{n+\beta p'}
	\fint_{A_j} \abs{\nabla\partial_t\Gamma_{(x,t)}^{A^*}(y,0) }^{p'}\,dy
.\end{aligned}\end{multline*}
If $p>p^-$, then $p'<p^+$, and so by Lemma~\ref{lem:slabs} and the bound~\eqref{eqn:FS:far:space},
\[\fint_{A_j} \abs{\nabla\partial_t\Gamma_{(x,t)}^{A^*}(y,0) }^{p'}\,dy
\leq
C(2^j t)^{-p'(n+1)}
.\]
Therefore,
\begin{multline*}
\int_{\R^n} \Bigl((t+\abs{x-y})^\beta \abs{\nabla\partial_t\Gamma_{(x,t)}^{A^*}(y,0) }\Bigr)^{p'}\,dy
\\\begin{aligned}
&\leq
	C\sum_{j=0}^\infty (2^j t)^{n+\beta p'-p'(n+1)}
	=Ct^{(p'/p)(\beta p-p-n)}
	\sum_{j=0}^\infty 2^{j(p'/p)(\beta p-p-n)}
.\end{aligned}\end{multline*}
If $\beta p-p-n<0$, then the sum converges and so
\begin{multline*}
\int_{\R^{n+1}_+} \abs{\partial_t \T\vec f(x,t)}^p t^{p-1-p\theta}\,dx\,dt
\\\begin{aligned}
&\leq
	C\int_0^\infty t^{p-1-p\theta}
	\int_{\R^n} 
	\int_{\R^n} \biggl(\frac{\abs{f(y)-f(x)}}{(t+\abs{x-y})^\beta}\biggr)^p \,dy
	\,t^{\beta p-p-n}
	\,dx\,dt
\\&=
	C\int_{\R^n} 
	\int_{\R^n} 	\abs{f(y)-f(x)}^p
	\int_0^\infty
	\frac{t^{-1-p\theta+p\beta -n}}{(t+\abs{x-y})^{p\beta}}
	 \,dt\,dy\,dx
\\&=
	C\int_{\R^n} 
	\int_{\R^n} 	
	\frac{\abs{f(y)-f(x)}^p}{\abs{x-y}^{n+p\theta}}
	\int_0^\infty
	\frac{s^{-1-p\theta+p\beta -n}}{(s+1)^{p\beta}}
	 \,ds\,dy\,dx
\end{aligned}\end{multline*}
where we have made the change of variables $s= \abs{x-y}t$. If we choose $\beta$ such that $-1-p\theta+\beta p-n>-1$, then the integral in $s$ converges. We recognize the remaining term from formula~\eqref{eqn:besov:norm:differences} as being comparable to $\doublebar{f}_{\dot B^{p,p}_\theta(\R^n)}$; thus, formula~\eqref{eqn:space:vertical} is valid 
provided we can find a real number $\beta$ such that 
\begin{equation*}-1<-1-n+\beta p-\theta p
\quad\text{and}\quad
\beta p-p-n<0.
\end{equation*}
Solving, we see that we need only require $n+\theta p<\beta p < n+p$, and so \eqref{eqn:space:vertical} is valid provided $p^-<p<\infty$ and $0<\theta<1$.
\end{proof}

We have established that if $p^-< p<\infty$, then
\begin{equation*}
\int_{\R^{n+1}_+} \abs{\partial_t \T\vec f(x,t)}^p t^{p-1-p\theta} \,dt\,dx
\leq
	C\doublebar{f}_{\dot B^{p,p}_{\theta}(\R^n)}^p
.\end{equation*}
To improve to an estimate on the full gradient, we first observe that if $p< p^+$, then by the Caccioppoli inequality and Lemma~\ref{lem:PDE2} applied in Whitney cubes, we have
\begin{gather}
\label{eqn:space:intermediate}
\int_{\R^{n+1}_+} \abs{\partial_t \nabla \T \vec f(x,t)}^p t^{2p-1-p\theta}\,dx\,dt
\leq C \doublebar{f}_{\dot B^{p,p}_\theta(\R^n)}^p
.\end{gather}
We will use the following lemma to remove the derivative $\partial_t$.
\begin{lem}\label{lem:space:down}
Let $E\subseteq\R^n$.
Suppose that $\beta>-1$, that $1<p<\infty$, and that $u\in W^2_{1,loc}(\R^{n+1}_+)$. Then for some constant $C(p,\beta)$ we have that
\begin{align}\label{eqn:space:down}
\int_0^\infty \int_{E}\abs{u(x,t)}^p t^{\beta}\,dx\,dt
&\leq
	C(p,\beta)\int_0^\infty \int_{E}\abs{\partial_t u(x,t)}^p t^{p+\beta}\,dx\,dt
	\\\nonumber&\qquad
	+C(p,\beta) \liminf_{t\to\infty} t^{1+\beta} \int_{E} \abs{u(x,t)}^p\,dx
.\end{align}
\end{lem}
We intend to apply Lemma~\ref{lem:space:down} with $u(x,t)=\nabla \T\vec f(x,t)$ and with $\beta=p-1-p\theta$. We will first prove Lemma~\ref{lem:space:down} {and then show how to deal with the second term on the right-hand side of \eqref{eqn:space:down}.}

\begin{proof}[Proof of Lemma~\ref{lem:space:down}]
This has been investigated for $p=2$ and $\beta$ an integer and can be proven using the same techniques; see  \cite[formula~(5.5)]{AlfAAHK11}, \cite{Gra12}, and \cite[Lemma~7.2]{BarM13B}.
Let $0<\varepsilon<R$. Then
\begin{align*}
\int_{\varepsilon}^R \int_{E}\abs{u(x,t)}^p t^{\beta}\,dx\,dt
&\leq
	\int_{\varepsilon}^R \int_{E}\abs{u(x,R)}^p t^{\beta}\,dx\,dt
	\\&\qquad
	+\int_{\varepsilon}^R \int_t^R \frac{d}{ds}\int_{E}\abs{u(x,s)}^p \,dx\,ds\,t^{\beta}\,dt.
\end{align*}
But 
\begin{align*}\frac{d}{ds}\int_{E}\abs{u(x,s)}^p \,dx
&\leq \int_{E}p\abs{u_s(x,s)}\abs{u(x,s)}^{p-1} \,dx
\\&\leq
C_p\int_{E} \eta^{1-p} s^{p-1} \abs{\partial_s u(x,s)}^p+\eta s^{-1}\abs{ u(x,s)}^p \,dx
\end{align*} for any $\eta$, $s>0$, and so
\begin{align*}
\int_{\varepsilon}^R \int_{E}\abs{u(x,t)}^p t^{\beta}\,dx\,dt
&\leq
	\int_{\varepsilon}^R \int_{E}\abs{u(x,R)}^p t^{\beta}\,dx\,dt
	\\&\qquad
	+C_p\eta^{1-p} \int_{\varepsilon}^R \int_t^R 
	\int_{E} s^{p-1} \abs{\partial_s u(x,s)}^p\,dx\,ds\,t^{\beta}\,dt
	\\&\qquad
	+C_p\eta \int_{\varepsilon}^R \int_t^R 
	\int_{E} s^{-1} \abs{ u(x,s)}^p \,dx\,ds\,t^{\beta}\,dt
.\end{align*}
Interchanging the order of integration and evaluating integrals in~$t$, we have that
\begin{align*}
\int_{\varepsilon}^R \int_{E}\abs{u(x,t)}^p t^{\beta}\,dx\,dt
&\leq
	\frac{R^{1+\beta}-\varepsilon^{1+\beta}}{\beta+1} \int_{E}\abs{u(x,R)}^p \,dx
	\\&\qquad
	+C_p\eta^{1-p} \int_{\varepsilon}^R \int_{E} 
	\frac{s^{p+\beta}-\varepsilon^{1+\beta}s^{p-1}}{1+\beta}
	 \abs{\partial_s u(x,s)}^p\,dx\,ds.
	\\&\qquad
	+C_p\eta \int_{\varepsilon}^R \int_{E} 
	\frac{s^{\beta}-\varepsilon^{1+\beta}/s}{1+\beta}
	\abs{ u(x,s)}^p \,dx\,ds.
\end{align*}
Because $\beta+1>0$ we may drop the terms involving powers of~$\varepsilon$. Observe that if the first two terms on the right-hand side are finite, then so is the last term. Therefore, if we choose $\eta=(1+\beta)/2C_p$, we may hide this term and see that
\begin{align*}
\int_{\varepsilon}^R \int_{E}\abs{u(x,t)}^p t^{\beta}\,dx\,dt
&\leq
	C_{p,\beta} R^{1+\beta}\int_{E}\abs{u(x,R)}^p \,dx
	\\&\qquad
	+C_{p,\beta}\int_{\varepsilon}^R  \int_{E} s^{p+\beta}\abs{\partial_s u(x,s)}^p \,dx\,ds.
\end{align*}
Taking the limit as $\varepsilon\to 0^+$ and the limit inferior as $R\to\infty$ we recover \eqref{eqn:space:down}.
\end{proof}

Let $R>0$ be a large number.
By Lemma~\ref{lem:space:down} and formula~\eqref{eqn:space:intermediate}, we have that if $p^-<p<p^+$ and $0<\theta<1$ then
\begin{multline*}
\int_0^\infty \int_{\Delta(0,R)} \abs{\nabla \T\vec f(x,t)}^p t^{p-1-p\theta}\,dx\,dt
\\\leq C(p,\theta) \doublebar{f}_{\dot B^{p,p}_\theta(\R^n)}^p
+C(p,\theta) \liminf_{t\to\infty} t^{p-p\theta} \int_{\Delta(0,R)} \abs{\nabla \T\vec  f(x,t)}^p\,dx
.\end{multline*}
By Lemma~\ref{lem:slabs}, the Caccioppoli inequality, and the bound~\eqref{eqn:pointwise}, we have that if $p^-<p<p^+$ and $t>R$, then
\begin{align*}
t^{p-p\theta} \int_{\Delta(0,R)} \abs{\nabla \T\vec  f(x,t)}^p\,dx
&\leq 
	t^{p-p\theta} \int_{\Delta(0,t)} \abs{\nabla \T\vec  f(x,t)}^p\,dx
\\&\leq 
	C(p) t^{p+n-p\theta} \biggl(\fint_{\Delta(0,2t)}\fint_{t/2}^{3t/2} \abs{\nabla \T\vec  f(x,s)}^2\,ds\,dx\biggr)^{p/2}
\\&\leq 
	C(p) t^{n-p\theta} 
\biggl(\fint_{\Delta(0,3t)}\fint_{t/4}^{7t/4} \abs{ \T\vec  f(x,s)}^2\,ds\,dx\biggr)^{p/2}
\\&\leq 
	C(p) t^{-p\theta} \doublebar{\vec f}_{L^p(\R^n)}^p
\end{align*}
which approaches zero as $t\to\infty$. Thus, if $f\in L^p(\R^n)\cap\dot B^{p,p}_\theta(\R^n)$, we may take the limit as $R\to \infty$ to see that if $0<\theta<1$ and $p^-<p<p^+$ then
\begin{equation}
\label{eqn:space:T:p}
\int_{\R^{n+1}_+} \abs{\nabla \T\vec f(x,t)}^p t^{p-1-p\theta}\,dx\,dt
\leq C(p,\theta) \doublebar{f}_{\dot B^{p,p}_\theta(\R^n)}^p
.\end{equation}
Applying Lemma~\ref{lem:local-bound:gradient}, we see that if $0<p\leq\infty$ and $0<q<p^+$, then the bound \eqref{eqn:space:T:p} implies that
\begin{equation}
\label{eqn:space:T}
\int_{\R^{n+1}_+} \biggl(\fint_{\Omega(x,t)}\abs{\nabla \T\vec f}^q\biggr)^{p/q} t^{p-1-p\theta}\,dx\,dt
\leq C(p,\theta,q) \doublebar{f}_{\dot B^{p,p}_\theta(\R^n)}^p
.\end{equation}

Because $L^p(\R^n)\cap\dot B^{p,p}_\theta(\R^n)$
is dense in $\dot B^{p,p}_\theta(\R^n)$ for any $0<\theta<1$ and any~$p<\infty$, we may extend $\T$ to a bounded operator $\dot B^{p,p}_\theta(\R^n)\mapsto \dot W(p,\theta,q)$.

\subsection{Atomic estimates\texorpdfstring{: boundedness in the range $p\leq 1$}{}}
\label{sec:bounded:atomic}

In this section we will establish the bounds \eqref{eqn:space:D} and \eqref{eqn:space:S} for $p\leq 1$ and satisfying the assumptions of Theorem~\ref{thm:bounded}. (That is, we will consider $p$ and $\theta$ that satisfy formula~\eqref{eqn:hexagon:1}.)
\begin{thm} \label{thm:atomic}
Let $A$  be elliptic, $t$-independent and such that both $A$ and $A^*$ satisfy the De Giorgi-Nash-Moser condition. Let $\vec a_Q:\R^n\mapsto\C^{n+1}$ be a vector-valued $(\theta, p)$ atom for some $1-\alpha<\theta<1$, $n/(n+\alpha)<p\leq 1$ with $1/p < 1+(\theta+\alpha-1)/n$. Then
\begin{gather}
\label{eqn:space:atomic}
\int_{\R^{n+1}_+} \abs{\nabla \T\vec a_Q(x,t)}^p t^{p-1-p\theta}\,dx\,dt
\leq C(p,\theta).
\end{gather}
\end{thm}
By Lemma~\ref{lem:local-bound:gradient}, and because $L(p,\theta,q)$ and $\dot B^{p,p}_\theta(\R^n)$ have $p$-norms when $p\leq 1$ (see Remark~\ref{rmk:atoms}), the bound \eqref{eqn:space:atomic} implies that the bounds~\eqref{eqn:space:T:p} and hence~\eqref{eqn:space:T} are valid for $p$, $\theta$ in this range as well.

\begin{proof}[Proof of Theorem~\ref{thm:atomic}]
First, observe that if $0<p<2$, then by H\"older's inequality
\begin{multline*}
\int_{8Q} \int_0^{16\ell(Q)} \abs{\nabla \T\vec a_Q(x,t)}^p t^{p-1-p\theta}\,dt\,dx
\\\begin{aligned}
&\leq 
	C\abs{Q} \int_0^{16\ell(Q)} \biggl(\fint_{8Q} \abs{\nabla \T\vec a_Q(x,t)}^2\,dx\biggr)^{p/2} t^{p-1-p\theta}\,dt
.\end{aligned}\end{multline*}
If $\beta>0$ then
\begin{multline*}
\int_{8Q} \int_0^{16\ell(Q)} \abs{\nabla \T\vec a_Q(x,t)}^p t^{p-1-p\theta}\,dt\,dx
\\\begin{aligned}
&\leq 
	C\abs{Q} 
	\biggl(\int_0^{16\ell(Q)} \fint_{8Q} \abs{\nabla \T\vec a_Q(x,t)}^2\,dx\,
	t^{(p-p/2-p\theta-\beta)(2/p)}\,dt\biggr)^{p/2}
	\\&\qquad\times
	\biggl(\int_0^{16\ell(Q)} t^{(p/2-1+\beta)/(1-p/2)} \,dt\biggr)^{1-p/2}
\\&\leq 
	C(\beta) \abs{Q}^{1+\beta/n} 
	\biggl(\int_0^{16\ell(Q)} \fint_{8Q} \abs{\nabla \T\vec a_Q(x,t)}^2\,dx\,
	t^{2-1-2(\theta+\beta/p)}\,dt\biggr)^{p/2}
.\end{aligned}\end{multline*}
Let $\beta/p=(1-\theta)/2$ so that $0<\theta+\beta/p<1$; then by \eqref{eqn:space:T} with $p$ and~$q$ replaced by~$2$ and $\theta$ replaced by $\theta+\beta/p$, we have that
\begin{equation*}
\int_{8Q} \int_0^{16\ell(Q)} \abs{\nabla \T\vec a_Q(x,t)}^p t^{p-1-p\theta}\,dt\,dx
\leq 
	C \abs{Q}^{1+\beta/n-p/2} 
	\doublebar{\vec a_Q}_{\dot B^{2,2}_{\theta+\beta/p}(\R^n)}^p
.\end{equation*}
By Remark~\ref{rmk:atom:strange-norm}, this reduces to
\begin{equation*}
\int_{8Q} \int_0^{16\ell(Q)} \abs{\nabla \T\vec a_Q(x,t)}^p t^{p-1-p\theta}\,dt\,dx
\leq 
	C
.\end{equation*}

Recall that $Q\subset\R^n$ is a dyadic cube.
Let $\G$ be a grid of dyadic Whitney cubes; then $\widetilde Q=Q	\times(\ell(Q),2\ell(Q))$ is an element of~$\G$.

If $j\leq 0$, let $P_j(Q)=Q$, and if $j\geq 1$ then let $P_j(Q)$ be the unique dyadic cube contained in $\R^n$ with side-length $2^j\ell(Q)$ and with $Q\subset P_j(Q)$. We let $N(Q)$ be the union of all dyadic cubes $R\in\G$ with $R\subset 8Q\times(0,16\ell(Q))$; observe that $5Q\times(0,4\ell(Q))\subset N(Q)$.

As in the proof of Theorem~\ref{thm:newton:preliminary}, we let
\begin{equation*}A_{j,0}(Q)=17P_j(Q)\times (2^j\ell(Q),2^{j+1}\ell(Q)) \setminus N(Q),\end{equation*}
and for $k\geq 1$ let \begin{equation*}A_{j,k}(Q)=((2^{k+4}+1)P_j(Q)\setminus(2^{k+3}+1)P_j(Q))\times (2^j\ell(Q),2^{j+1}\ell(Q)).\end{equation*}
We define the solid versions 
\begin{equation*}\widetilde A_{j,0}=32P_j(Q)\times (-\max(2,2^{j+2})\ell(Q),\max(2,2^{j+2})\ell(Q)) \setminus (3/2)Q\end{equation*}
and
\begin{align*}
\widetilde A_{j,k}&=\bigl(2^{k+5}P_j(Q)
\setminus 2^{k+2}P_j(Q)\bigr)
\\&\quad\times (-\max(2^{k+1},2^{k+j+1})\ell(Q),\max(2^{k+1},2^{k+j+1})\ell(Q)).
\end{align*}

Now, suppose that $(x,t)\in\R^{n+1}$ but $(x,t)\notin 4Q\times(0,4\ell(Q))$. Then, for any fixed $x_Q\in Q$, by by Lemma~\ref{lem:slabs}, the bound~\eqref{eqn:FS:far:space} and H\"older's inequality we have that
\begin{align*}
\abs{\T\vec a_Q(x,t)} 
&\leq
	C\int_{3Q} \abs{\nabla\Gamma_{(x,t)}^{A^*}(y,0)}
	\abs{\vec a_Q}\,dy
\\&\leq
	C\abs{Q}^{\theta/n-1/p}
	\int_{3Q}  \abs{\nabla\Gamma_{(x,t)}^{A^*}(y,0)}\,dy
\\&\leq
	C\abs{Q}^{\theta/n-1/p+1-1/n}
	\biggl(\fint_{4Q}\fint_{-\ell(Q)}^{\ell(Q)}  \abs{\Gamma_{(x,t)}^{A^*}(y,0)-\Gamma_{(x,t)}^{A^*}(x_Q,0)}^2\,dy\biggr)^{1/2}
\\&\leq
	C
	\frac{\abs{Q}^{\theta/n-1/p+1-1/n+\alpha/n}} {\abs{x-x_Q}^{n-1+\alpha}+t^{n-1+\alpha}}
.\end{align*}

So, applying Lemma~\ref{lem:slabs},
\begin{multline*}
\int_{\R^{n+1}_+\setminus N(Q)}
\abs{\nabla \T\vec a_Q(x,t)}^p\,t^{p-1-p\theta}\,dx\,dt
\\\begin{aligned}
&\leq C
	\sum_{j=-\infty}^\infty \sum_{k=0}^\infty
	2^{j(p-p\theta)}\ell(Q)^{n+p-p\theta}
	(2^k\max(1,2^j))^n
	\fint_{A_{j,k}} 
	\abs{\nabla \T\vec a_Q}^p
\\&\leq C
	\sum_{j=-\infty}^\infty \sum_{k=0}^\infty
	2^{j(p-p\theta)}\ell(Q)^{n-p\theta}
	(2^k\max(1,2^j))^{n-p}
	\biggl(\fint_{\widetilde A_{j,k}} 
	\abs{ \T\vec a_Q}^2\biggr)^{p/2}
\\&\leq C
	\sum_{j=-\infty}^\infty \sum_{k=0}^\infty	
	\frac{2^{j(p-p\theta)}} {(2^k\max(1,2^j))^{np-n+p\alpha}}	
.\end{aligned}\end{multline*}
The sum converges provided $\theta<1$ and provided $np-n+p\alpha>p-p\theta$, that is, provided
$1/p<1+(\theta+\alpha-1)/n$.
\end{proof}

\subsection{H\"older spaces\texorpdfstring{: boundedness results for $p=\infty$}{}}
\label{sec:bounded:holder}

In this section we will establish the bounds \eqref{eqn:space:D} and \eqref{eqn:space:S} for $p=\infty$. Recall that 
\begin{align*}
\doublebar{f}_{\dot B^{\infty,\infty}_\theta(\R^n)}
&\approx \esssup_{(x,y)\in \R^n\times\R^n} \frac{\abs{f(x)-f(y)}}{\abs{x-y}^\theta}
=\doublebar{f}_{\dot C^\theta(\R^n)}.
\end{align*}
Thus, we are interested in the H\"older spaces~$\dot C^\theta$.
We have the following theorem.
\begin{thm}[{\cite[formula (1.25)]{HofMitMor}}] \label{thm:bounded:holder:D}
Suppose that $A$ is elliptic and $t$-inde\-pen\-dent and that both $A$ and $A^*$ satisfy the De Giorgi-Nash-Moser condition.
If $f\in \dot C^\theta(\R^n)$ for some $0<\theta<\alpha$, then 
\begin{equation*}
\doublebar{\D f}_{\dot C^\theta(\R^{n+1}_+)}
\leq C(\theta) \doublebar{f}_{\dot C^\theta(\R^n)}.\end{equation*}
\end{thm}
Examining the proof of this theorem, we see that it is also valid for the single layer potential $\s\nabla_\parallel$; for the sake of completeness we include a proof.

\begin{thm}\label{thm:bounded:holder}
Suppose that $A$ is elliptic, $t$-independent and that both $A$ and $A^*$ satisfy the De Giorgi-Nash-Moser condition.
If $0<\theta<\alpha$ and $p^-<q<p^+$, then there is some $C(\theta,q)$ such that, if $\vec f\in \dot C^\theta(\R^n\mapsto\C^n)$, then 
\begin{align*}
\frac{1}{C(\theta,q)}\sup_{(x,t)\in\R^{n+1}_+}
t^{1-\theta} \biggl(\fint_{\Omega(x,t)}\abs{\nabla \T\vec f}^q\biggr)^{1/q}
&\leq
\doublebar{\T\vec f}_{\dot C^\theta(\R^{n+1}_+)}
\leq C(\theta,q) \doublebar{\vec f}_{\dot B^{\infty,\infty}_\theta(\R^n)}
.\end{align*}
\end{thm}

\begin{proof}
Let $Q\subset\R^n$ be a cube, and let $\widetilde Q=Q\times[a,a+\ell(Q)]$ for some $a\geq 0$. By Meyer's criterion in \cite{Mey64}, to show that
\begin{equation}\label{eqn:holder-to-Holder}
\doublebar{\T\vec f}_{\dot C^\theta(\R^{n+1}_+)}
\leq C \doublebar{\vec f}_{\dot B^{\infty,\infty}_\theta(\R^n)}\end{equation}
it suffices to show that for each such $\widetilde Q$ there is a constant $c_{\widetilde Q}$ (depending on~$\widetilde Q$) and a constant $C_\theta$ (independent of~$\widetilde Q$) such that
\begin{equation*}\fint_{\widetilde Q} \frac{\abs{\T\vec f(x,t)-c_{\widetilde Q}}} {\ell(Q)^\theta}\,dx\,dt
\leq C_\theta\doublebar{\vec f}_{\dot B^{\infty,\infty}_\theta(\R^n)}.\end{equation*}

By Lemma~\ref{lem:T-constant}, $\T (\fint_{4Q} \vec f)(x,t)$ is constant for $(x,t)\in\R^{n+1}_+$; thus, we may assume without loss of generality that $\fint_{4Q} \vec f=0$.

Let $A_0=4Q$, let $A_j=2^{j+2}Q\setminus 2^{j+1}Q$ for all $j\geq 1$, and let $\vec f_j = \vec f \1_{A_j}$. We define 
\begin{equation*} c_j = \fint_Q \T\vec f_j(y,a+\ell(Q))\,dy 
.\end{equation*}
We choose $c_{\widetilde Q}=\sum_{j=0}^\infty c_j$, so that $\T\vec f(x,t)-c_{\widetilde Q}$ may be given by formula~\eqref{eqn:T:difference}.

We may then write
\begin{align*}
\fint_{\widetilde Q}  \abs{\T\vec f(x,t)-c_{\widetilde Q}} \,dx\,dt
&\leq
	\sum_{j=0}^\infty \fint_{\widetilde Q}\abs{\T\vec f_j(x,t)-c_j}\,dx\,dt
\\&\leq
	\sum_{j=0}^\infty \fint_{\widetilde Q}\fint_Q \abs[big]{\T\vec f_j(x,t)-\T \vec f_j(y,a+\ell(Q))}\,dy\,dx\,dt
.\end{align*}
Now, observe that if $\fint_{4Q} \vec f=0$ and $z\in A_j$, then \begin{equation*}\abs{\vec f(z)}
=\abs[bigg]{ \fint_{4Q} {\vec f(z)-\vec f(x)}\,dx}
\leq \fint_{4Q} \abs{\vec f(z)-\vec f(x)}\,dx
\leq C\doublebar{\vec f}_{\dot C^\theta(\R^n)} 2^{j\theta} \ell(Q)^\theta.\end{equation*}

Certain bounds on layer potentials are well known. Notice that $\vec f\mapsto\T\vec f(\,\cdot\,,t)$ is the adjoint to the operator $g\mapsto\nabla_{A^*}\s^{A^*}g(\,\cdot\,,-t)$, and so by duality with \cite[Theorem~7.1]{Ros12A}, we have that 
\begin{equation*}\sup_{t>0} \doublebar{\T\vec f_0(\,\cdot\,,t)}_{L^2(\R^n)}
\leq C \doublebar{\vec f_0}_{L^2(\R^n)} \leq C\doublebar{\vec f}_{\dot C^\theta(\R^n)} \ell(Q)^{n/2+\theta}
\end{equation*}
and in particular 
\[c_0
\leq \ell(Q)^{-n/2}\doublebar{\T\vec f_0(\,\cdot\,,a+\ell(Q))}_{L^2(Q)}
\leq C\doublebar{\vec f}_{\dot C^\theta(\R^n)} \ell(Q)^{\theta}.\]
Thus,
\begin{align*}
\fint_{\widetilde Q}\abs{\T\vec f_0(x,t)-c_0}\,dx\,dt
&\leq
	\fint_a^{a+\ell(Q)} \ell(Q)^{-n/2}
	\biggl(\int_{ Q} \abs{\T\vec f_0(x,t)}^2\,dx\biggr)^{1/2}\,dt
	+\abs{c_0}
\\&\leq
	C\ell(Q)^{\theta}
	\doublebar{\vec f}_{\dot C^\theta(\R^n)} 
.\end{align*}

Now, if $j\geq 1$, and if $(x,t)\in \widetilde Q$ and $(y,s)\in {\widetilde Q}$, then
\begin{align*}
\abs[big]{\T\vec f_j(x,t)-\T \vec f_j(y,s)}
\leq C \int_{A_j} \abs{\vec f(z)} \,\abs{\nabla \Gamma_{(x,t)}(z,0)-\nabla\Gamma_{(y,s)}(z,0)}\,dz.
\end{align*}
By the bound~\eqref{eqn:FS:holder} and Lemma~\ref{lem:slabs},
\begin{equation*}
\int_{A_j} \abs{\nabla \Gamma_{(x,t)}(z,0)-\nabla \Gamma_{(y,s)}(z,0)} \,dz\leq C 2^{-j\alpha}
.\end{equation*}

Thus, we conclude that
\begin{align*}
\fint_{\widetilde Q}  \abs{\T\vec f(x,t)-c_{\widetilde Q}} \,dx\,dt
&\leq
	C\ell(Q)^\theta \doublebar{\vec f}_{\dot C^\theta(\R^n)} 
	\sum_{j=0}^\infty 2^{j(\theta-\alpha)} 
\end{align*}
which converges provided $\theta<\alpha$. This proves the bound~\eqref{eqn:holder-to-Holder}.

To complete the proof, we simply observe that by Caccioppoli's inequality and Lemma~\ref{lem:PDE2},
\begin{align*}
\biggl(\fint_{\Omega(x,t)}\abs{\nabla\T\vec f}^q\biggr)^{1/q}
&\leq
	Ct^{-1} \biggl(\fint_{B((x,t),2t/3)} \abs{\T \vec f(y,s)-\T\vec f(x,t)}^2 \,dy\,ds\biggr)^{1/2}
\\&\leq
	Ct^{\theta-1}\doublebar{\T \vec f}_{\dot C^\theta(\R^{n+1}_+)}.
	\qedhere
\end{align*}
\end{proof}

\subsection{Interpolation}

We have now shown that formula~\eqref{eqn:space:T} is valid, that is, that $\T:\dot B^{p,p}_\theta(\R^n)\mapsto \dot W(p,\theta,q)$ is valid
in the following three cases:
\begin{itemize}
\item $0<\theta<1$, $p^-<p<p^+$ (Section~\ref{sec:bounded:midrange}),
\item $1-\alpha<\theta<1$, $1\leq 1/p < 1+(\theta+\alpha-1)/n$ (Section~\ref{sec:bounded:atomic}),
\item $0<\theta<\alpha$, $p=\infty$ (Section~\ref{sec:bounded:holder}).
\end{itemize}

\begin{figure}[tbp]
\begin{tikzpicture}[scale=2]
\figureaxes
\draw [boundary bounded] (1,1)--(1-\alph,1);
\fill [bounded] 
(1,1/2-\eps) node [black] {$\circ$} node [black, right] {$(1,1/p^+)$}--
(1,1/2+\eps)--
(0,1/2+\eps) node [black] {$\circ$} node [black, left] {$(0,1/p^-)$}--
(0,1/2-\eps)--cycle;
\fill [bounded] (1,1)--
(1,1+\alph/\enn) node [black] {$\circ$} node [black, right] {$(1,1+\alpha/n)$}--
(1-\alph,1) node [black] {$\circ$} node [black, above, at = {(1-0.2-\alph,1)}] {$(1-\alpha,1)$}--
cycle;
\draw [boundary bounded] (0,0)--(\alph,0)  node [black] {$\circ$} node [black, below] {$(\alpha,0)$};
\draw [line width=0.5pt, dotted] (0,1/2)--(1,1/2);
\end{tikzpicture}
\qquad
\begin{tikzpicture}[scale=2]
\figureaxes
\boundedhexagon
\draw [line width=0.5pt, dotted] (0,1/2)--(1,1/2);
\end{tikzpicture}
\caption{Values of $p$, $\theta$ such that $\T$ is bounded $\dot B^{p,p}_\theta(\R^n)\mapsto L(p,\theta,q)$.}
\label{fig:bounded:1}
\end{figure}
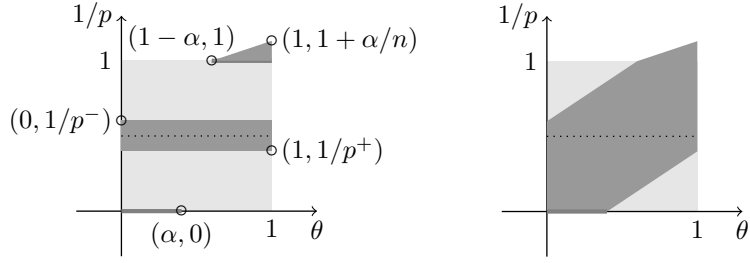

That is, we have that $\T$ extends to an operator that is bounded $\dot B^{p,p}_\theta(\R^n)\mapsto \dot W(p,\theta,q)$
for all $p^-<q<p^+$ and for the values of $\theta$ and $p$ indicated on the left in  Figure~\ref{fig:bounded:1}.

We may use interpolation to fill in the gaps; see 
Theorems~\ref{thm:B-B:interpolation:complex} and~\ref{thm:interpolation:L} and the interpolation property~\eqref{eqn:operator:interpolation:complex}. Thus, $\T$ is bounded $\dot B^{p,p}_\theta(\R^n)\mapsto \dot W(p,\theta,q)$ for all $p^-<q<p^+$ and for the values of $p$ and $\theta$ indicated on the right  in  Figure~\ref{fig:bounded:1}; this completes the proof of Theorem~\ref{thm:bounded}.

%% file: sec-6-trace.tex
\chapter{Trace Theorems}
\label{chap:trace}

In this section we will prove Theorem~\ref{thm:trace}.

This theorem concerns the trace and conormal derivatives of functions $u$ in the weighted, averaged Sobolev spaces $\dot W(p,\theta,q)$, that is, functions $u$ that satisfy
\[\int_{\R^{n+1}_+} \biggl(\fint_{\Omega(x,t)}\abs{\nabla u}^q\biggr)^{p/q} t^{p-1-p\theta}\,dx\,dt<\infty.\]
The main innovation in this section is to consider \emph{averaged} Sobolev spaces. Boundedness of the trace operator $\dot W(p,\theta)\mapsto \dot B^{p,p}_\theta(\R^n)$ was established in \cite{Liz60} and independently in \cite{Usp61}, where the space $\dot W(p,\theta)$ is given by 
\[\doublebar{u}_{\dot W(p,\theta)}^p=\int_{\R^{n+1}_+}\abs{\nabla u(x,t)}^p \, t^{p-1-p\theta}\,dx\,dt.\] See also \cite[Section~2.9.2]{Tri78}.
Similar results are valid in bounded $C^{m,\delta}$ domains; see \cite{Nik77b, Sha85, NikLM88, Kim07}. See also \cite{JerK95} and~\cite{FabMM98}, where the trace and the normal derivative of functions in weighted Sobolev spaces on Lipschitz domains are considered in the context of the Dirichlet problem for the Laplacian.





We begin with the following embedding theorem; compare to Lemma~\ref{lem:besov:embedding}.
\begin{thm}
\label{thm:whitney:embedding}
Let $0<p_0\leq p_1\leq \infty$, $0<q_1\leq q_0\leq\infty$, and suppose that the real numbers ${\theta_1}$, ${\theta_0}$ satisfy ${\theta_1}-n/p_1={\theta_0}-n/p_0$.

Then $L(p_0,\theta_0,q_0)\subset L(p_1,\theta_1,q_1)$.

Furthermore, recall that if $\F\in L(p,\theta,q)$ then
\begin{equation*}\doublebar{\F}_{L(p,\theta,q)}^{p}\approx \sum_{Q\in\G} \ell(Q)^{n+p-p{\theta}}\biggl(\fint_Q \abs{F}^{q}\biggr)^{p/q}\end{equation*}
where $\G$ is the standard grid of dyadic Whitney cubes.
For each $Q\in\G$ let $\pi(Q)$ be the projection of $Q$ onto $\R^n$, and let $T(Q)=\pi(Q)\times(0,2\ell(Q))$.

If $p_1\leq q_1$ and $n+p_1-p_1{\theta_1}-(n+1)p_1/q_1<0$, or if $p_1\geq q_1$ and $p_1-p_1\theta_1-p_1/q_1<0$, then
\begin{equation}\label{eqn:weighted:local}
\sum_{Q\in\G} \ell(Q)^{n+p_1-p_1{\theta_1}}\biggl(\fint_{T(Q)} \abs{\F}^{q_1}\biggr)^{p_1/q_1}
\leq C_1\doublebar{\F}_{L(p_1,{\theta_1},q_1)}^{p_1}
\leq C_0\doublebar{\F}_{L(p_0,{\theta_0},q_0)}^{p_1}
\end{equation}
for some constants $C_j$ depending only on the numbers $p_j$, $\theta_j$, $q_j$.
\end{thm}

Before proving this theorem we observe some immediate consequences. 
If $q_1=1$ then the bound \eqref{eqn:weighted:local} is valid whenever either $p_1\geq 1$ and $\theta_1>0$ or $p_1<1$ and $1/p_1-\theta_1/n=1/p_0-\theta_0/n<1$. Thus, if $u\in \dot W(p,\theta,q)$ for any $q\geq 1$ and any $p$, $\theta$ with $1/p<1+\theta/n$, then $\nabla u$ is in $L^1_{loc}(\overline{\R^{n+1}_+})$. 

In particular, if $\Div A\nabla u$ is smooth and compactly supported then the integral \eqref{eqn:conormal:2} in the definition of $\nu\cdot A\nabla u$ converges absolutely, and by the boundedness of the trace operator on unweighted Sobolev spaces, we have that $\Tr u$ exists and is locally in $L^1(\R^n)$. Thus, the trace and the conormal derivative exist.

\begin{proof}[Proof of Theorem~\ref{thm:whitney:embedding}]
Observe that

\begin{align*} 
	\sum_{Q\in\G} \ell(Q)^{n+p_1-p_1{\theta_1}}\biggl(\fint_Q \abs{\F}^{q_1}\biggr)^{p_1/q_1}
&= 
	\sum_{Q\in\G} 
	\biggl(\ell(Q)^{n+p_0-p_0{\theta_0}}\biggl(\fint_Q \abs{\F}^{q_1}\biggr)^{p_0/q_1}\biggr)^{p_1/p_0}
.\end{align*}
By applying H\"older's inequality twice, once in the sum and once in the integral, we see that $\doublebar{\F}_{L(p_1,\theta_1,q_1)} \leq C\doublebar{\F}_{L(p_0,\theta_0,q_0)}$, and so $L(p_0,\theta_0,q_0)\subset L(p_1,\theta_1,q_1)$.

Next, observe that
\begin{multline*}
\sum_{Q\in\G} \ell(Q)^{n+p_1-p_1{\theta_1}}\biggl(\fint_{T(Q)} \abs{\F}^{q_1}\biggr)^{p_1/q_1}
\\\begin{aligned}
&=
	\frac{1}{2^{p_1/q_1}}
	\sum_{Q\in\G} \ell(Q)^{n+p_1-p_1{\theta_1}-(n+1)p_1/q_1}
	\biggl(\sum_{R\subset T(Q)}\int_{R} \abs{\F}^{q_1}\biggr)^{p_1/q_1}
.\end{aligned}\end{multline*}
If $p_1\leq q_1$, then by H\"older's inequality
\begin{multline*}
\sum_{Q\in\G} \ell(Q)^{n+p_1-p_1{\theta_1}}\biggl(\fint_{T(Q)} \abs{\F}^{q_1}\biggr)^{p_1/q_1}
\\\begin{aligned}
&\leq
	\frac{1}{2^{p_1/q_1}}
	\sum_{Q\in\G} \ell(Q)^{n+p_1-p_1{\theta_1}-(n+1)p_1/q_1}
	\sum_{R\subset T(Q)}\biggl(\int_{R} \abs{\F}^{q_1}\biggr)^{p_1/q_1}
\\&=
	\frac{1}{2^{p_1/q_1}}
	\sum_{R\in\G} \biggl(\int_{R} \abs{\F}^{q_1}\biggr)^{p_1/q_1}
	\sum_{Q\in\G:R\subset T(Q)}
	\ell(Q)^{n+p_1-p_1{\theta_1}-(n+1)p_1/q_1}
.\end{aligned}\end{multline*}
If $Q$, $R$ are dyadic cubes and $R\subset T(Q)$, then $\ell(Q)=2^j\ell(R)$ for some integer $j\geq 0$; moreover, for each fixed $R\in\G$ and each $j\geq 0$ there is exactly one such cube~$Q$.
Thus
\begin{multline*}
\sum_{Q\in\G} \ell(Q)^{n+p_1-p_1{\theta_1}}\biggl(\fint_{T(Q)} \abs{\F}^{q_1}\biggr)^{p_1/q_1}
\\\begin{aligned}
&\leq
	\frac{1}{2^{p_1/q_1}}
	\sum_{R\in\G} 
	\ell(R)^{n+p_1-p_1{\theta_1}}\biggl(\fint_{R} \abs{\F}^{q_1}\biggr)^{p_1/q_1}
	\sum_{j=0}^\infty
	2^{j(n+p_1-p_1{\theta_1}-(n+1)p_1/q_1)}
.\end{aligned}\end{multline*}
But if $n+p_1-p_1{\theta_1}-(n+1)p_1/q_1<0$, then the sum in $j$ converges and so the bound~\eqref{eqn:weighted:local} is valid.

Conversely, if $p_1\geq q_1$, then 
\begin{align*}
\biggl(\sum_{R\subset T(Q)}\int_{R} \abs{\F}^{q_1}\biggr)^{p_1/q_1}
&=
	\biggl(\sum_{j=0}^\infty \sum_{\substack{R\subset T(Q)\\\ell(R)=2^{-j}\ell(Q)}} \int_{R} \abs{\F}^{q_1}\biggr)^{p_1/q_1}
\end{align*}
and for any $\varepsilon>0$,
\begin{align*}
\biggl(\sum_{R\subset T(Q)}\int_{R} \abs{\F}^{q_1}\biggr)^{p_1/q_1}
&\leq
	C(\varepsilon)\sum_{j=0}^\infty 2^{j\varepsilon}\biggl(\sum_{\substack{R\subset T(Q)\\\ell(R)=2^{-j}\ell(Q)}} \int_{R} \abs{\F}^{q_1}\biggr)^{p_1/q_1}
.\end{align*}
Observe that there are $2^{jn}$ cubes $R\in \G$ with $R\subset T(Q)$ and $\ell(R)=2^{-j}\ell(Q)$. Thus by H\"older's inequality,
\begin{align*}
\biggl(\sum_{R\subset T(Q)}\int_{R} \abs{\F}^{q_1}\biggr)^{p_1/q_1}
&\leq
	C(\varepsilon)\sum_{j=0}^\infty
	2^{j(\varepsilon+np_1/q_1-n)}\!\!\!\!
	\sum_{\substack{R\subset T(Q)\\\ell(R)=2^{-j}\ell(Q)}} \biggl(\int_{R} \abs{\F}^{q_1}\biggr)^{p_1/q_1}
\\&=
	C(\varepsilon)
	\sum_{R\subset T(Q)}
	\biggl(\frac{\ell(Q)}{\ell(R)}\biggr)^{\varepsilon+np_1/q_1-n}
	\biggl(\int_{R} \abs{\F}^{q_1}\biggr)^{p_1/q_1}
.\end{align*}
Thus,
\begin{multline*}
\sum_{Q\in\G} \ell(Q)^{n+p_1-p_1{\theta_1}}\biggl(\fint_{T(Q)} \abs{\F}^{q_1}\biggr)^{p_1/q_1}
\\\begin{aligned}
&=
	\frac{1}{2^{p_1/q_1}}
	\sum_{Q\in\G} \ell(Q)^{n+p_1-p_1{\theta_1}-(n+1)p_1/q_1}
	\biggl(\sum_{R\subset T(Q)}\int_{R} \abs{\F}^{q_1}\biggr)^{p_1/q_1}
\\&\leq
	C(\varepsilon)
	\sum_{Q\in\G} 
	\sum_{R\subset T(Q)}
	\frac{\ell(Q)^{p_1-p_1{\theta_1}-p_1/q_1+\varepsilon}} {\ell(R)^{\varepsilon+np_1/q_1-n}}
	\biggl(\int_{R} \abs{\F}^{q_1}\biggr)^{p_1/q_1}
.\end{aligned}\end{multline*}
Rearranging the sum, we have that
\begin{multline*}
\sum_{Q\in\G} \ell(Q)^{n+p_1-p_1{\theta_1}}\biggl(\fint_{T(Q)} \abs{\F}^{q_1}\biggr)^{p_1/q_1}
\\\begin{aligned}
&\leq
	C(\varepsilon)
	\sum_{R\in\G}
	\ell(R)^{n+p_1-p_1\theta_1}
	\biggl(\fint_{R} \abs{\F}^{q_1}\biggr)^{p_1/q_1}
	\sum_{Q:R\subset T(Q)}
	\biggl(\frac{\ell(Q)} {\ell(R)}\biggr)^{\varepsilon-p_1/q_1+p_1-p_1\theta_1}
.\end{aligned}\end{multline*}
If  $-p_1/q_1+p_1-p_1\theta_1<0$, choose $\varepsilon$ with $0<\varepsilon<p_1/q_1-p_1+p_1\theta_1$. As before, the second sum converges and the bound~\eqref{eqn:weighted:local} is valid.
\end{proof}

In order to bound the traces and conormal derivatives of a function $u\in\dot W(p,\theta,q)$ for $p>1$, we will make use the formulas
\[\langle \Tr u,\varphi\rangle = \int_{\R^{n+1}_+} \overline{\nabla\Phi} \cdot \nabla u,
\qquad
\langle \nu\cdot A\nabla u,\psi\rangle = \int_{\R^{n+1}_+} \overline{\nabla\Psi} \cdot A\nabla u\]
where $\Phi$ is harmonic and $\nu\cdot\nabla\Phi=\varphi$, and $\Psi$ is an arbitrary well-behaved function with $\Tr \Psi=\psi$. (We have that $\Phi=2\s^I\varphi$; it is convenient to take $\Psi=-2\D^I\psi$. The second formula must be suitably modified if $\Div A\nabla u\neq 0$; see formula~\eqref{eqn:conormal:2}.) Aside from the issue of defining $\nu\cdot A\nabla u$, we will consider all functions $u$ in $\dot W(p,\theta,q)$, not only solutions to elliptic equations, and so we will not use any special properties of the function~$u$. We will, however, need certain bounds on the extensions~$\Phi$ and~$\Psi$. For the sake of being self-contained and because our theorems are formulated in terms of the spaces $\dot W(p,\theta,q)$ we will continue to refer to our Theorem~\ref{thm:bounded} for these bounds; however, we observe that the necessary bounds on the layer potentials $\D^I$ and $\s^I$ for the Laplacian may also be found in \cite{FabMM98}. Because $-2\D^I$ coincides with the Poisson extension in the upper half-space, the necessary properties of $\D^I$ may also be found in~\cite{JerK95}; because $\Psi$ need not be harmonic it would suffice to instead use the extension theorem of \cite{Liz60} and the mollifier theorem of \cite{Sha85}.

\begin{thm}\label{thm:trace:banach}
Suppose $1< p < \infty$, $1\leq q<\infty$ and $0<\theta<1$. Then the trace operator $\Tr$ extends to an operator that is bounded 
\begin{equation*}\Tr:\dot W(p,\theta,q)\mapsto \dot B^{p,p}_\theta(\R^n).\end{equation*}

If  $u\in \dot W_+(\infty,\theta,q)$, then the function $\tilde u(x,t)=\fint_{\Omega(x,t)} u$ is H\"older continuous in $\R^{n+1}_+$ with exponent~$\theta$. Thus, if we extend the trace $\Tr$ by letting
\begin{equation*}\Tr u(x)=\lim_{t\to 0^+}\fint_{\Omega(x,t)} u
\quad\text{for all }u\in\dot W_+(\infty,\theta,q)\end{equation*}
then $\Tr$ is bounded $\dot W_+(\infty,\theta,q)\mapsto \dot B^{\infty,\infty}_\theta(\R^n)$.
\end{thm}

\begin{proof}
If $p<\infty$ and $q<\infty$, then smooth, compactly supported functions are dense in $\dot W(p,\theta,q)$, and so $\Tr$ is densely defined.

Choose some smooth, compactly supported $u\in \dot W(p,\theta,q)$, and some $\varphi\in (\dot B^{p,p}_\theta(\R^n))'=\dot B^{p',p'}_{-\theta}(\R^n)$.
We wish to show that 
\begin{equation*}\abs{\langle \Tr u,\varphi\rangle}\leq C\doublebar{u}_{\dot W(p,\theta,q)} \doublebar{\varphi}_{\dot B^{p',p'}_{-\theta}(\R^n)}.\end{equation*}

Let $\Phi=2\s^I \varphi$, where $I$ is the identity matrix; then $\Phi$ is harmonic and $-\partial_{n+1}\Phi$ is the Poisson extension of~$\varphi$, and so $\nu\cdot \nabla\Phi=\varphi$ on $\partial\R^{n+1}_+$. Then
\begin{align*}
{\langle \Tr u,\varphi\rangle}
&={\int_{\R^n} \Tr u(x)\, \overline{\nu\cdot \nabla \Phi(x)}\,dx }
\\&={\int_{\R^{n+1}_+} \nabla u(x,t)\cdot  \overline{\nabla \Phi(x,t)}\,dx \,dt}.
\end{align*}
Therefore,
\begin{align*}
\abs{\langle \Tr u,\varphi\rangle}
&\leq
	C\int_{\R^{n+1}_+}\fint_{\Omega(x,t)}\abs{\nabla \Phi(y,s)}\,\abs{\nabla u(y,s)}\,dy\,ds\,dx\,dt
\\&\leq
	C
	\biggl(\int_{\R^{n+1}_+}\biggl(\fint_{\Omega(x,t)}\abs{\nabla u}\biggr)^{p} t^{p-1-p\theta}\,dx\,dt\biggr)^{1/p}
	\\&\qquad\times
	\biggl(\int_{\R^{n+1}_+} 
	\sup_{(y,s)\in\Omega(x,t)} \abs{\nabla\Phi}^{p'} t^{p'-1-p'(1-\theta)}\,dx\,dt\biggr)^{1/p'}.
\end{align*}
We may bound the first term using H\"older's inequality and the second term using local boundedness of the harmonic function~$\nabla\Phi$. Thus,
\begin{align*}
\abs{\langle \Tr u,\varphi\rangle}
&\leq
	C
	\doublebar{u}_{\dot W(p,\theta,q)}
	\doublebar{\Phi}_{\dot W(p',1-\theta,2)}
.\end{align*}
Recall that $\Phi=2\s^I\varphi$; if $A\equiv I$ then $\alpha=1$ and $p^+=\infty$, and so the bound ~\eqref{eqn:space:S} is valid. Applying this bound, we see that
\begin{align*}
\abs{\langle \Tr u,\varphi\rangle}
\leq 
	C
	\doublebar{u}_{\dot W(p,\theta,q)}
	\doublebar{\varphi}_{\dot B^{p',p'}_{-\theta}(\R^n)}
\end{align*}
as desired.

We are left with the case $p=\infty$.
Observe that if $q\geq 1$, then 
\begin{equation*}\abs{\nabla \tilde u(x,t)} 
= \abs[bigg]{\nabla \fint_{B(0,1/2)} u(x+ty,s+ts)}
\leq C\fint_{\Omega(x,t)} \abs{\nabla u}\leq C t^{\theta-1}\doublebar{u}_{\dot W(\infty,\theta,q)}
.\end{equation*}
Let $(x,t)$ and $(z,r)\in \R^{n+1}_+$, and let $s=\abs{(x,t)-(z,r)}$. Then
\begin{align*}
\abs{\tilde u(z,r)-\tilde u(x,t)}
&\leq 
	\abs{\tilde u(z,r)-\tilde u(z,s+r)} 
	+ \abs{\tilde u(z,s+r)-\tilde u(x,s+r)} 
	\\&\qquad
	+ \abs{\tilde u(x,s+r)-\tilde u(x,t)}
\end{align*}
and using our bound on $\nabla \tilde u$, we see that each term is at most $Cs^\theta \doublebar{u}_{\dot W(\infty,\theta,q)}$.
Thus, $\tilde u$ is H\"older continuous and so $u\vert_{\partial\R^{n+1}_+}$ exists and is also H\"older continuous.
\end{proof}

\begin{cor}\label{cor:holder:weighted}
Suppose that $\Div A\nabla u=0$ in $\R^{n+1}_+$ for some elliptic matrix~$A$ that satisfies the De Giorgi-Nash-Moser condition. Let $0<\theta\leq\alpha$.

Then $u\in \dot W(p,\theta,2)$ if and only if $u\in  \dot C^\theta(\R^{n+1}_+)$ and $\nabla u\in L^2_{loc}(\R^{n+1}_+)$, and
$\doublebar{u}_{\dot C^\theta(\R^{n+1}_+)}\approx\doublebar{u}_{\dot W(\infty,\theta,2)}$.
\end{cor}

\begin{proof}
The bound $\doublebar{u}_{\dot W(\infty,\theta,2)} \leq C\doublebar{u}_{\dot C^\theta(\R^{n+1}_+)}$ follows immediately from the Caccioppoli inequality. For the reverse inequality, choose some $(x,t)$ and $(y,s)\in\R^{n+1}_+$. Let $r=\abs{(x,t)-(y,s)}$. If $(y,s)\in B((x,t),t/4)$, then by the De Giorgi-Nash-Moser condition
\[\abs{u(x,t)-u(y,s)}\leq C\biggl(\frac{r}{t}\biggr)^\alpha
\biggl(\fint_{\Omega(x,t)} \abs{u-\textstyle\fint_\Omega u}^2\biggr)^{1/2}\]
and by the Poincar\'e inequality and because $\theta\leq\alpha$, we have that
\[\abs{u(x,t)-u(y,s)}\leq Cr^\theta
t^{1-\theta}
\biggl(\fint_{\Omega(x,t)} \abs{\nabla u}^2\biggr)^{1/2}
\leq Cr^\theta \doublebar{u}_{\dot W(\infty,\theta,2)}.\]
Otherwise, $r>t/4$ and so $s/5\leq t/5 + r/5 < r$. Let $\tilde u$ be as in Theorem~\ref{thm:trace:banach}. Then
\[\abs{u(x,t)-u(y,s)}\leq \abs{u(x,t)-\tilde u(x,t)} + \abs{\tilde u(x,t)-\tilde u(y,s)}+\abs{\tilde u(y,s)-u(y,s)}.\]
We bound the middle term using Theorem~\ref{thm:trace:banach} and the first and last terms using the De Giorgi-Nash-Moser condition, as before; this completes the proof.
\end{proof}

\begin{thm}
\label{thm:conormal:banach}
Suppose $1< p\leq  \infty$, $1\leq q\leq\infty$ and $0<\theta<1$. 
If $u\in \dot W(p,\theta,q)$, and if $\Div A\nabla u=\Div \F$ in $\R^{n+1}_+$, for some $\F\in L(p,\theta,q)$ smooth and compactly supported in $\R^{n+1}_+$, then the normal derivative $\nu\cdot A\nabla u$ exists in the sense of formula~\eqref{eqn:conormal:2} and lies in $\dot B^{p,p}_{\theta-1}(\R^n)$.
\end{thm}
\begin{proof}
Recall that the formula~\eqref{eqn:conormal:2} for the normal derivative is given by
\begin{equation*}\langle \varphi,\nu\cdot A\nabla u\rangle
=\int_{\R^{n+1}_+} \nabla\Phi\cdot A\nabla u
-\int_{\R^{n+1}_+} \nabla\Phi\cdot \F
\end{equation*}
whenever $\Phi\in C^\infty_0(\R^{n+1})$ and $\varphi=\Phi\bigr\vert_{\partial\R^{n+1}_+}$. Recall that by definition of $\Div A\nabla u$, the value of the right-hand side does not depend on the choice of extension~$\Phi$.

We proceed as in the proof of Theorem~\ref{thm:trace:banach}. Let $\varphi\in \dot B^{p',p'}_{1-\theta}(\R^n)$, and choose $\Phi$ to be the Poisson extension of $\varphi$, that is, $\Phi=-2\D^I \varphi$. Then
\begin{align*}
\abs{\langle \varphi,\nu\cdot A\nabla u\rangle}
&\leq
	\int_{\R^{n+1}_+}\abs{\nabla \Phi(x,t)}\,
	(\Lambda\abs{\nabla u(x,t)}+\abs{\F(x,t)})\,dx\,dt
\\&\leq
	C
	\biggl(\int_{\R^{n+1}_+}\biggl(\fint_{\Omega(x,t)} \Lambda\abs{\nabla u}+\abs{\F}\biggr)^{p} t^{p-1-p\theta}\,dx\,dt\biggr)^{1/p}
	\\&\qquad\times
	\biggl(\int_{\R^{n+1}_+} 
	\sup_{(y,s)\in\Omega(x,t)} \abs{\nabla\Phi}^{p'} t^{p'-1-p'(1-\theta)}\,dx\,dt\biggr)^{1/p'}
.\end{align*}
Applying H\"older's inequality to $\nabla u$ and $\F$ and local boundedness to the harmonic function $\nabla\Phi$, we see that
\begin{align*}
\abs{\langle \varphi,\nu\cdot A\nabla u\rangle}
&\leq
	C
	(\doublebar{u}_{\dot W(p,\theta,q)}+\doublebar{\F}_{L(p,\theta,q)})
	\doublebar{\Phi}_{\dot W(p',\theta,2)}
.\end{align*}
Applying the bound~\eqref{eqn:space:D} to $\Phi=-2\D^I\varphi$, we have that 
\begin{align*}
\abs{\langle \varphi,\nu\cdot A\nabla u\rangle}
&\leq
	C
	(\doublebar{u}_{\dot W(p,\theta,q)}+\doublebar{\F}_{L(p,\theta,q)})
	\doublebar{\varphi}_{\dot B^{p',p'}_{1-\theta}(\R^n)}
\end{align*}
and so
\begin{equation*}\doublebar{\nu\cdot A\nabla u}_{\dot B^{p,p}_{\theta-1}(\R^n)}
\leq C \doublebar{u}_{\dot W(p,\theta,q)}+C\doublebar{\F}_{L(p,\theta,q)}
\end{equation*}
as desired.
\end{proof}

We now move to the case $p\leq 1$. Observe that if $p\leq 1$ then $\dot B^{p,p}_\theta(\R^n)$ is not the dual to $\dot B^{p',p'}_{1-\theta}(\R^n)$, and so the method of proof of Theorems~\ref{thm:trace:banach} and~\ref{thm:conormal:banach} cannot be used. Instead, we will use a wavelet decomposition of the Besov spaces. 

The homogeneous Daubechies wavelets were constructed in \cite[Section~4]{Dau88}. We will need the following properties.
\begin{lem}
For any integer $N>0$ there exist real functions $\psi$ and $\varphi$ defined on $\R$ that satisfy the following properties.
\begin{itemize}
\item $\abs{\frac{d^k}{dx^k}\psi(x)}\leq C(N)$, $\abs{\frac{d^k}{dx^k}\psi(x)}\leq C(N)$ for all $k<N$,
\item $\psi$ and $\varphi$ are supported in the interval $(-C(N),1+C(N))$,
\item $\int_{\R} \varphi(x)\,dx\neq 0$, $\int_{\R} \psi(x)\,dx = \int_{\R} x^k\,\psi(x)\,dx=0$ for all $0\leq k<N$.
\end{itemize}
Furthermore, suppose we let $\psi_{j,m}(x)=2^{jn/2}\psi(2^j x-m)$ and $\varphi_{j,m}(x)=2^{jn/2}\varphi(2^j x-m)$. Then $\{\psi_{j,m}:j,m\in\Z\}$ is an orthonormal basis for $L^2(\R)$, and  if $j_0$ is an integer then $\{\varphi_{j_0,m}:m\in\Z\}\cup\{\psi_{j,m}:m\in\Z,j\geq j_0\}$ is also an orthonormal basis for $L^2(\R)$.
\end{lem}
The functions $\varphi$ and $\psi$ are often referred to as a scaling function and a wavelet, or as a father wavelet and a mother wavelet.

We may produce an orthonormal basis of $L^2(\R^n)$ for $n\geq 1$ from these wavelets by considering the $2^n-1$ functions  $\Psi^\ell(x)=\eta_1(x_1)\eta_2(x_2)\dots\eta_n(x_n)$, where for each $j$ we have that either $\eta_j(x)=\varphi(x)$ or $\eta_j(x)=\psi(x)$ and $\eta_k(x)=\psi(x)$ for at least one~$k$. Let $\Psi^\ell_{j,m}=2^{j/2}\Psi^\ell(2^j x-m)$; then $\{\Psi^\ell_{j,m}:j\in\Z,\>m\in\Z^n\}$ is an orthonormal basis for~$L^2(\R^n)$. Notice that we may instead index the wavelets $\Psi^\ell_{j,m}$ by dyadic cubes~$Q$, with $\Psi^\ell_{j,m}=\Psi^\ell_Q$ if $Q=\{2^j(y+m):y\in[0,1]^n\}$. We then have that $\Psi^\ell_Q$
has the following properties:
\begin{itemize}
\item $\Psi_Q^\ell$ is supported in $CQ$,
\item $\abs{\partial^\beta\Psi^\ell_Q(x)}\leq C(N)\ell(Q)^{-n/2-\abs{\beta}}$ whenever $\abs{\beta}<N$,
\item $\int_{\R^n} x^\beta \Psi^\ell_Q(x)\,dx=0$ whenever  $\abs{\beta}<N$.
\end{itemize}
Here $\beta\in \NN^n$ is a multiindex.

Now, notice that if $N$ is large enough then $(1/C(N))\ell(Q)^{\theta+n/2-n/p}\Psi_Q^\ell$ is a $(\theta,p)$ atom in the sense of \cite{FraJ85} (see Section~\ref{sec:function} above). Because $\{\Psi_Q^\ell\}$ is an orthonormal basis of $L^2(\R^n)$, we have that if $f\in L^2(\R^n)$ then  
\begin{equation}
\label{eqn:besov:wavelet}
f(x)=\sum_{Q}\sum_{\ell=1}^{2^n-1} \langle f, \Psi_Q^\ell\rangle \Psi_Q^\ell(x)\end{equation}
and so by the atomic decomposition of $\dot B^{p,p}_\theta(\R^n)$, we have that
\begin{equation}
\label{eqn:besov:wavelet:norm}
\doublebar{f}_{\dot B^{p,p}_\theta(\R^n)}^p\leq C\sum_Q\sum_{\ell=1}^{2^n-1}
\abs{\langle f, \Psi_Q^\ell\rangle }^p \ell(Q)^{n-np/2-p\theta}.\end{equation}
The reverse inequality is also true for all $f\in \dot B^{p,p}_\theta(\R^n)$; see \cite[Theorem~4.2]{Kyr03}. Furthermore, by the same theorem, if $f\in \dot B^{p_1,r_1}_{\theta_1}(\R^n)$ for some $0<p_1,r_1\leq\infty$ and some $\theta_1\in\R$, then the decomposition~\eqref{eqn:besov:wavelet} is valid.

We now use this wavelet decomposition to investigate traces and conormal derivatives of $\dot W(p,\theta,q)$ functions for $p\leq 1$.

\begin{thm}\label{thm:trace:quasi-banach}
Suppose $1\leq q\leq \infty$, $0<p\leq 1$, and $0<\theta<1$ with $1/p<1+\theta/n$. Then the trace operator $\Tr$ extends to a bounded operator $\dot W(p,\theta,q)\mapsto \dot B^{p,p}_\theta(\R^n)$.
\end{thm}

\begin{proof}
By Theorems~\ref{thm:whitney:embedding} and~\ref{thm:trace:banach}, $f=\Tr u$ exists and lies in $\dot B^{p_1,p_1}_{\theta_1}(\R^n)$ for some $p_1>1$ and some $\theta_1>0$.
Let $\Psi_Q^\ell$ be as above. Then because $\Psi_Q^\ell$ is bounded,
\begin{equation*}\abs{\langle \Psi_Q^\ell,\Tr u\rangle}
\leq
C\ell(Q)^{-n/2} \int_{CQ} \abs{\Tr u}.
\end{equation*}
Let $T(Q)=CQ\times (0,C\ell(Q))$; observe that Theorem~\ref{thm:whitney:embedding} is still valid with this definition of $T(Q)$. It is well known (see, for example, \cite[Section~5.5]{Eva98}) that the trace operator is bounded $\dot W^1_1(T(Q))\mapsto L^1(\partial T(Q))$ with a constant independent of the size of~$Q$. Thus
\begin{equation*}\int_{3Q} \abs{\Tr u} \leq \int_{\partial T(Q)} \abs{\Tr u}
\leq C\int_{T(Q)} \abs{\nabla u}.\end{equation*}
So
\begin{align*}
\sum_Q  \sum_{\ell=1}^{2^n-1}\abs{\langle \Tr u,\Psi_Q^\ell\rangle}^p \ell(Q)^{n-np/2-p\theta}
&\leq
	C(2^n-1) \sum_Q  
	\biggl(\fint_{T(Q)} \abs{\nabla u}\biggr)^p
	\ell(Q)^{n+p-p\theta}.
\end{align*}
By Theorem~\ref{thm:whitney:embedding}, if $p\leq 1\leq q$ then the right-hand side is at most $C\doublebar{u}_{\dot W(p,\theta,q)}$; by the bound~\eqref{eqn:besov:wavelet:norm} the quantity $\doublebar{\Tr u}_{\dot B^{p,p}_\theta(\R^n)}$ is controlled by the left-hand side, as desired.
\end{proof}

\begin{thm}\label{thm:conormal:quasi-banach}
Suppose $1\leq q\leq \infty$, $0<p\leq 1$ and $0<\theta<1$ with $1/p<1+\theta/n$.
If $u\in \dot W(p,\theta,q)$, and if $\Div A\nabla u=\Div \F$ in $\R^{n+1}_+$, for some $\F\in L(p,\theta,q)$ smooth and compactly supported in $\R^{n+1}_+$, then the normal derivative $\nu\cdot A\nabla u$ exists in the sense of formula~\eqref{eqn:conormal:2} and lies in $\dot B^{p,p}_{\theta-1}(\R^n)$.
\end{thm}

\begin{proof}
Again by Theorems~\ref{thm:whitney:embedding} and~\ref{thm:conormal:banach} we have that $\nu\cdot A\nabla u$ exists and lies in a Besov space. Let $\Psi_Q^\ell$ be as above.
We may extend $\Psi_Q^\ell$ to a smooth function supported in $CQ\times (-\ell(Q),\ell(Q))$ and retain the property $\abs{\nabla\Psi_Q^\ell}\leq C\ell(Q)^{-n/2-1}$. We then have that
\begin{align*}
\abs{\langle \Psi_Q^\ell,\nu\cdot A\nabla u\rangle}
&=\abs[bigg]{\int_{\R^n} \Psi_Q^\ell \,\nu\cdot A\nabla u}
\\&=
	\abs[bigg]{\int_{\R^{n+1}_+} \nabla\Psi_Q^\ell \cdot A\nabla u
	-\int_{\R^{n+1}_+} \nabla\Psi_Q^\ell \cdot \F
	}
\\&\leq
	C\ell(Q)^{-n/2-1}\int_{3Q}\int_0^{\ell(Q)} 
	\abs{\nabla u(x,t)}+\abs{\F(x,t)}\,dt\,dx.
\end{align*}
Thus, if $T(Q)$ is as in the proof of Theorem~\ref{thm:trace:quasi-banach},
\begin{align*}
\doublebar{\nu\cdot A\nabla u}_{\dot B^{p,p}_{\theta-1}(\R^n)}^p
&\leq C\sum_Q  \sum_{\ell=1}^{2^n-1}
	\abs{\langle \Psi_Q^\ell, \nu\cdot A\nabla u\rangle}^p \ell(Q)^{n-np/2+p-p\theta}
\\&\leq
	C (2^n-1)\sum_Q 
	\biggl(\fint_{T(Q)} 
	\abs{\nabla u}+\abs{\F}\biggr)^p \ell(Q)^{n+p-p\theta}	
\\&\leq
	C \doublebar{u}_{\dot W(p,\theta,q)}+C\doublebar{\F}_{L(p,\theta,q)}
\end{align*}
where the last inequality follows by Theorem~\ref{thm:whitney:embedding}. This completes the proof.
\end{proof}

\begin{rmk}
\label{rmk:newton:trace}
Suppose that $p$, $q$ and $\theta$ satisfy the conditions of Theorem~\ref{thm:bounded}. 
We claim that the operators $\Pi^A_+$ and $(\partial_\nu^A\Pi^A)_+$ can be extended in a natural fashion to bounded  operators $\Pi^A_+:L(p,\theta,q)\mapsto\dot B^{p,p}_\theta(\R^n)$ and $(\partial_\nu^A\Pi^A)_+: L(p,\theta,q)\mapsto\dot B^{p,p}_{\theta-1}(\R^n)$.

If $p<\infty$, then smooth, compactly supported functions are dense in $L(p,\theta,q)$. Thus, $\Pi^A_+$ and $(\partial_\nu^A\Pi^A)_+$ may be extended to $L(p,\theta,q)$ by continuity, and by 
Theorem~\ref{thm:bounded} and Theorems~\ref{thm:trace:banach}, \ref{thm:conormal:banach}, \ref{thm:trace:quasi-banach} or~\ref{thm:conormal:quasi-banach}, these operators are bounded.

Consider the case $p=\infty$; in this case we are only concerned with $\theta$ such that $0<\theta<\alpha$. We may extend $\Pi^A_+$ as in Theorem~\ref{thm:trace:banach}; by that theorem $\Pi^A_+$ is bounded $L(\infty,\theta,q)\mapsto \dot B^{\infty,\infty}_{\theta}(\R^n)$.
If $\F$ is smooth and compactly supported, by inspecting the definitions \eqref{eqn:newton} and~\eqref{eqn:D}, we see that 
\[\int_{\R^n}\varphi\,\nu\cdot A\nabla\Pi^A\F=\int_{\R^{n+1}_+} \overline{\nabla\D^{A^*}\varphi}\cdot \F\]
in the sense of formula~\eqref{eqn:conormal:2}. But $\nabla\D^{A^*}$ is bounded $\dot B^{1,1}_{1-\theta}(\R^n)\mapsto L(1,1-\theta,q')$, and so its adjoint $(\nabla \D^{A^*})^*$ is bounded $L(\infty,\theta,q)\mapsto \dot B^{\infty,\infty}_{\theta-1}(\R^n)$; we may extend $(\partial_\nu^A\Pi^A)_+$ to an operator defined on $L(\infty,\theta,q)$ by taking $(\partial_\nu^A\Pi^A)_+\F=(\nabla \D^{A^*})^*\F$.
\end{rmk}

We have now completed the proof of Theorem~\ref{thm:trace}. In the remainder of this chapter, we will explore the sense in which $u$ and $\e_{n+1}\cdot A\nabla u$ approach their respective boundary values.

\begin{thm}\label{thm:interior-trace}
Let $A$ be elliptic and $t$-independent, and suppose that both $A$ and $A^*$ satisfy the De Giorgi-Nash-Moser condition. Fix some $\varepsilon>0$ and suppose that  $\Div A\nabla u=0$ in $\R^n \times(0,2\varepsilon)$. Let $u_\varepsilon(x,t)=u(x,t+\varepsilon)$, so $\Div A\nabla u_\varepsilon=0$ in $\R^n\times(-\varepsilon,\varepsilon)$. Suppose that $0<p\leq\infty$ and $0<\theta<1$. If $\theta<\alpha+(1-\alpha)p^+/p$, then
\begin{equation*}\doublebar{u_\varepsilon}_{\dot W(p,\theta,2)}
\leq
C(p,\theta)\doublebar{u}_{\dot W(p,\theta,2)}.\end{equation*}
\end{thm}

Combined with Theorem~\ref{thm:trace}, we have in particular that for appropriate $p$ and~$\theta$, if $f_\varepsilon(x)=u(x,\varepsilon)$ and $g_\varepsilon(x)=-\e_{n+1}\cdot A(x)\nabla u(x,\varepsilon)$, then $f_\varepsilon\in \dot B^{p,p}_\theta(\R^n)$ and $g_\varepsilon\in \dot B^{p,p}_{\theta-1}(\R^n).$

\begin{proof}[Proof of Theorem~\ref{thm:interior-trace}]
If $p=\infty$ then $\theta<\alpha$; by Corollary~\ref{cor:holder:weighted},
\[\doublebar{u}_{\dot C^\theta(\R^{n+1}_+)} \approx C \doublebar{u}_{\dot W(\infty,\theta,2)}\]
and so the conclusion follows immediately.

If $p<\infty$, then let $\Omega(x,t,\varepsilon)={B((x,t),\varepsilon/2)}$. We write 
\begin{multline*}
\int_{\R^{n+1}_+} \biggl(\fint_{\Omega(x,t)} \abs{\nabla u_\varepsilon}^2\biggr)^{p/2}\,t^{p-1-p\theta}\,dt\,dx
\\=
\int_{\R^n}\int_{\varepsilon}^\infty
\biggl(\fint_{\Omega(x,t,t-\varepsilon)} \abs{\nabla u}^2\biggr)^{p/2}\,(t-\varepsilon)^{p-1-p\theta}\,dt\,dx
.\end{multline*}
Observe that 
\begin{multline*}
\int_{\R^n}\int_{(3/2)\varepsilon}^\infty
\biggl(\fint_{\Omega(x,t,t-\varepsilon)} \abs{\nabla u}^2\biggr)^{p/2}\,(t-\varepsilon)^{p-1-p\theta}\,dt\,dx
\\
\leq
	C\int_{\R^n}\int_{(3/2)\varepsilon}^\infty
	\biggl(\fint_{\Omega(x,t)} \abs{\nabla u}^2\biggr)^{p/2}\,t^{p-1-p\theta}\,dt\,dx
.\end{multline*}
So we need only consider $\varepsilon<t<(3/2)\varepsilon$.

Let $\mathcal{Q}$ be a grid of cubes in $\R^n$ of side-length $\varepsilon/2$. Let $2\leq q_0<p^+$. Then by H\"older's inequality
\begin{multline*}
\int_{\R^n}\int_\varepsilon^{(3/2)\varepsilon}
\biggl(\fint_{\Omega(x,t,t-\varepsilon)} \abs{\nabla u}^2\biggr)^{p/2}\,(t-\varepsilon)^{p-1-p\theta}\,dt\,dx
\\\begin{aligned}
&\leq
	\int_\varepsilon^{(3/2)\varepsilon}
	\sum_{Q\in\mathcal{Q}} \int_Q
	\biggl(\fint_{\Omega(x,t,t-\varepsilon)} \abs{\nabla u}^{q_0}\biggr)^{p/q_0}\,dx
	\,(t-\varepsilon)^{p-1-p\theta}\,dt
.\end{aligned}\end{multline*}
Choose some $q_1$ with $q_1<\min(p,q_0)$. Again by H\"older's inequality
\begin{multline*}
\int_Q
	\biggl(\fint_{\Omega(x,t,t-\varepsilon)} \abs{\nabla u}^{q_0}\biggr)^{p/q_0}\,dx
\\\begin{aligned}
&\leq
	\biggl(\int_Q
	\fint_{\Omega(x,t,t-\varepsilon)} \abs{\nabla u}^{q_0}\,dx
	\biggr)^{q_1/q_0}
	\\&\qquad\qquad\times
	\biggl(\int_Q
	\biggl(\fint_{\Omega(x,t,t-\varepsilon)} \abs{\nabla u}^{q_0}\biggr)^{(p-q_1)/(q_0-q_1)}\,dx
	\biggr)^{1-q_1/q_0}
.\end{aligned}\end{multline*}
Now by Lemma~\ref{lem:slabs}, if $\varepsilon<t<(3/2)\varepsilon$ and $\ell(Q)=\varepsilon/2$ then
\begin{equation*}
	\biggl(\int_Q
	\fint_{\Omega(x,t,t-\varepsilon)} \abs{\nabla u}^{q_0}\,dx
	\biggr)^{q_1/q_0}
\leq
	C\abs{Q}^{q_1/q_0}
	\biggl(\fint_{2Q}\fint_{\varepsilon/2}^{2\varepsilon}
	\abs{\nabla u}^2
	\biggr)^{q_1/2}
.\end{equation*}
By Lemma~\ref{lem:PDE2} and the Caccioppoli inequality, 
\begin{multline*}
\biggl(\int_Q
	\biggl(\fint_{\Omega(x,t,t-\varepsilon)} \abs{\nabla u}^{q_0}\biggr)^{(p-q_1)/(q_0-q_1)}\,dx
	\biggr)^{1-q_1/q_0}
\\\begin{aligned}
&\leq
\biggl(\int_Q
	\biggl(
	(t-\varepsilon)^{-2}
	\fint_{\Omega(x,t,(3/2)(t-\varepsilon))} 
	\abs{u-\textstyle\fint_\Omega u}^{2}
	\biggr)^{q_0(p-q_1)/2(q_0-q_1)}\,dx
	\biggr)^{1-q_1/q_0}
\end{aligned}\end{multline*}
and using the De Giorgi-Nash-Moser condition and the Poincar\'e inequality to bound $\abs{u-\textstyle\fint_\Omega u}$, we have that
\begin{multline*}
\biggl(\int_Q
	\biggl(\fint_{\Omega(x,t,t-\varepsilon)} \abs{\nabla u}^{q_0}\biggr)^{(p-q_1)/(q_0-q_1)}\,dx
	\biggr)^{1-q_1/q_0}
\\\begin{aligned}
&\leq
(t-\varepsilon)^{-(1-\alpha)(p-q_1)} 
\varepsilon^{(1-\alpha)(p-q_1)} 
\abs{Q}^{1-q_1/q_0}
	\biggl(	
	\fint_{2Q}\fint_{\varepsilon/2}^{2\varepsilon}
	\abs{\nabla u}^2
	\biggr)^{(p-q_1)/2}
.\end{aligned}\end{multline*}
So
\begin{multline*}
\int_{\R^n}\int_\varepsilon^{(3/2)\varepsilon}
\biggl(\fint_{\Omega(x,t,t-\varepsilon)} \abs{\nabla u}^2\biggr)^{p/2}\,(t-\varepsilon)^{p-1-p\theta}\,dt\,dx
\\\begin{aligned}
&\leq
	\sum_{Q\in\mathcal{Q}} 
	\abs{Q}
	\biggl(\fint_{2Q}\fint_{\varepsilon/2}^{2\varepsilon}
	\abs{\nabla u}^2
	\biggr)^{p/2}
	\\&\qquad\qquad\times
	\varepsilon^{(1-\alpha)(p-q_1)} 
	\int_\varepsilon^{(3/2)\varepsilon}
	(t-\varepsilon)^{p-1-p\theta-(1-\alpha)(p-q_1)}\,dt
.\end{aligned}\end{multline*}
If we can choose $q_0<p^+$ and $q_1<\min(p,q_0)$ such that \begin{equation*}{p-1-p\theta-(1-\alpha)(p-q_1)}>-1,\end{equation*}
then the integral will converge and the right-hand side will be at most $C\doublebar{\nabla u}_{L(p,\theta,2)}$. But this condition is equivalent to the condition $q_1>p(\theta-\alpha)/(1-\alpha)$. To ensure that an appropriate $q_1$ exists we need only require that $p(\theta-\alpha)/(1-\alpha)<\min(p,p^+)$; since $\theta<1$ we only need $\theta<\alpha+(1-\alpha)p^+/p$, as specified in the theorem statement.
\end{proof}

\begin{thm}\label{thm:trace:converge}
Let $A$ be elliptic and $t$-independent, and suppose that both $A$ and $A^*$ satisfy the De Giorgi-Nash-Moser condition.
Suppose that  $\Div A\nabla u=0$ in $ \R^{n+1}_+$ and that $u\in\dot W(p,\theta,2)$ for some $\theta$, $p$ satisfying the conditions of Theorems~\ref{thm:trace} (and hence also of Theorem~\ref{thm:interior-trace}). 

Let $f=\Tr u$, $g=\nu\cdot A\nabla u\big\vert_{\partial\R^{n+1}_+}$, and let $f_\delta(x)=u(x,\delta)$, $g_\delta(x)=-\e_{n+1}\cdot A(x)\nabla u(x,\delta)$.

If $1<p<\infty$ then $f_\delta\to f$ in $\dot B^{p,p}_{\theta}(\R^n)$ and $g_\delta \to g$ in $\dot B^{p,p}_{\theta-1}(\R^n)$. If $p=\infty$ then $f_\delta\rightharpoonup f$ in the weak-$*$ topology of  $\dot B^{\infty,\infty}_{\theta}(\R^n)=(\dot B^{1,1}_{-\theta}(\R^n))'$, and $g_\delta\rightharpoonup g$ in the weak-$*$ topology of $\dot B^{\infty,\infty}_{\theta-1}(\R^n)=(\dot B^{1,1}_{1-\theta}(\R^n))'$.

Finally, if $1<p<\infty$ then $f_\delta\to f$ in $L^p(\R^n)$ with \[\doublebar{f_\delta-f}_{L^p(\R^n)}\leq C(p,\theta)\delta^\theta\doublebar{u}_{\dot W(p,\theta,2)}.\]
\end{thm}

\begin{proof}
By Theorems~\ref{thm:trace} and~\ref{thm:interior-trace}, we have that $f_\delta$, $f\in \dot B^{p,p}_\theta(\R^n)$ and that $g_\delta$, $g\in \dot B^{p,p}_{\theta-1}(\R^n)$.

We begin with $g_\delta$. Fix some $\varphi\in \dot B^{p',p'}_{1-\theta}(\R^n)$ with $\doublebar{\varphi}_{\dot B^{p',p'}_{1-\theta}(\R^n)}=1$. It suffices to prove that 
\begin{equation*}
\lim_{\delta\to 0^+}\abs{\langle\varphi,g_\delta-g\rangle}=0 
\end{equation*}
and that, if $p<\infty$, then this limit is uniform in functions $\varphi$ with $\doublebar{\varphi}_{\dot B^{p',p'}_{1-\theta}(\R^n)}=1$.

Let $\Phi=-2\D^I\varphi$ be the Poisson extension of $\varphi$, so $\Phi\in \dot W(p',1-\theta,2)$. If $\eta>\delta$ then
\begin{align*}
\abs{\langle\varphi,g_\delta-g\rangle}
&\leq
	\abs[bigg]{
	\int_{\R^n}\int_\eta^\infty (\nabla\Phi(x,t-\delta)-\nabla\Phi(x,t))\cdot A(x)\nabla u(x,t)\,dx\,dt
	}
	\\&\qquad+
	\abs[bigg]{
	\int_{\R^n}\int_\delta^\eta (\nabla\Phi(x,t-\delta)-\nabla\Phi(x,t))\cdot A(x)\nabla u(x,t)\,dx\,dt
	}
	\\&\qquad+
	\abs[bigg]{\int_{\R^n}\int_0^\delta \nabla\Phi(x,t)\cdot A(x)\nabla u(x,t)\,dx\,dt
	}.
\end{align*}
By Theorem~\ref{thm:interior-trace}, the second and third terms are at most 
\begin{equation*}
C\doublebar{\1_{t<\eta}\nabla \Phi}_{L(p',1-\theta,2)} \doublebar{\1_{t<\eta}\nabla u}_{L(p,\theta,2)}
\leq 
C\doublebar{\varphi}_{\dot B^{p',p'}_{1-\theta}(\R^n)} \doublebar{\1_{t<\eta}\nabla u}_{L(p,\theta,2)}
.\end{equation*}
If we choose $\eta$ so that $\eta\to 0$ as $\delta\to 0$, then the left-hand side will always go to zero; furthermore, if $p<\infty$ then the right-hand side will approach zero uniformly in~$\varphi$. We need to consider the first term. For the sake of definiteness we will choose $\eta=\sqrt{\delta}$ and deal only with $\delta\ll 1$.

Let $\G$ be the standard grid of dyadic Whitney cubes. Then
\begin{multline*}
\abs[bigg]{
	\int_{\R^n}\int_{\sqrt{\delta}}^\infty (\nabla\Phi(x,t-\delta)-\nabla\Phi(x,t))\cdot A(x)\nabla u(x,t)\,dx\,dt
	}
\\\begin{aligned}
&\leq
	C
	\sum_{Q\in\G,\>\ell(Q)>\sqrt{\delta}/2}
	\int_Q
	\abs{\nabla\Phi(x,t-\delta)-\nabla\Phi(x,t)}
	\abs{\nabla u(x,t)}\,dx\,dt
.\end{aligned}\end{multline*}
Since $\Phi$ is harmonic, we have that if $\ell(Q)/\delta$ is large enough then
\begin{equation*}\max_{(x,t)\in \overline Q} \abs{\nabla\Phi(x,t-\delta)-\nabla\Phi(x,t)}
\leq 
\frac{C\delta}{\ell(Q)} \fint_{(5/4)Q} \abs{\nabla \Phi}
\end{equation*}
and so
\begin{multline*}
\abs[bigg]{
	\int_{\R^n}\int_{\sqrt{\delta}}^\infty (\nabla\Phi(x,t-\delta)-\nabla\Phi(x,t))\cdot A(x)\nabla u(x,t)\,dx\,dt
	}
\\\begin{aligned}
&\leq
	C\sqrt{\delta}
	\sum_{Q\in\G,\>\ell(Q)>\sqrt{\delta}/2}
	\fint_{(5/4)Q} \abs{\nabla \Phi}
	\int_Q
	\abs{\nabla u(x,t)}\,dx\,dt
\\&\leq
	C\sqrt{\delta}
	\doublebar{\Phi}_{\dot W(p',1-\theta,1)}
	\doublebar{u}_{\dot W(p,\theta,1)}
.\end{aligned}\end{multline*}
Thus, $\abs{\langle\varphi,g_\delta-g\rangle}\to 0$ as $\delta\to0^+$. If $p<\infty$ then $g_\delta\to g$ in $\dot B^{p,p}_{\theta-1}(\R^n)$, as desired. If $p=\infty$ then $\langle \varphi, g_\delta\rangle \to \langle\varphi,g\rangle$ for all $\varphi\in \dot B^{1,1}_{1-\theta}(\R^n)$, and so $g_\delta\rightharpoonup g$ in the weak-$*$ topology on~$\dot B^{\infty,\infty}_{\theta-1}(\R^n)$.

Next, consider $f_\delta$. In this case we choose $\varphi\in \dot B^{p',p'}_{-\theta}(\R^n)$ and let $\Phi=-2\s^I\varphi$, so that $\nu_+\cdot \nabla \Phi=\varphi$ on $\partial\R^{n+1}_+$.
Then
\begin{align*}
\abs{\langle\varphi,f_\delta-f\rangle}
&\leq
	\abs[bigg]{
	\int_{\R^n}\int_\eta^\infty (\nabla\Phi(x,t-\delta)-\nabla\Phi(x,t))\cdot \nabla u(x,t)\,dx\,dt
	}
	\\&\qquad+
	\abs[bigg]{
	\int_{\R^n}\int_\delta^\eta (\nabla\Phi(x,t-\delta)-\nabla\Phi(x,t))\cdot \nabla u(x,t)\,dx\,dt
	}
	\\&\qquad+
	\abs[bigg]{\int_{\R^n}\int_0^\delta \nabla\Phi(x,t)\cdot \nabla u(x,t)\,dx\,dt
	}
\end{align*}
and arguing as before we see that the right-hand side goes to zero as $\delta\to 0$, and if $p<\infty$ the convergence is uniform in functions $\varphi$ with $\doublebar{\varphi}_{\dot B^{p',p'}_{-\theta}(\R^n)}=1$.

Finally, we consider limits of $f_\delta$ in $L^p(\R^n)$. Observe that if $0\leq\varepsilon<\delta$, then by H\"older's inequality and the fact that $\partial_t u(x,t)$ satisfies $\Div A\nabla (\partial_t u)=0$ and so is locally bounded, we have that
\begin{align}
\label{eqn:besov-trace-convergence}
\int_{\R^n} \abs{f_\delta-f_\varepsilon}^p 
&
	\leq \int_{\R^n} \abs[bigg]{\int_\varepsilon^\delta \partial_t u(x,t)\,dt}^p\,dx
\\\nonumber&
	\leq \int_{\R^n} \int_\varepsilon^\delta\abs{\partial_t u(x,t)}^p t^{p-1-p\theta} \,dt
	\,\biggl(\int_\varepsilon^\delta  t^{\theta p'-1} dt\biggr)^{p/p'}
	dx
\\\nonumber&
	\leq C \delta^{\theta p}\int_{\R^n} \int_\varepsilon^\delta \abs{\partial_t u(x,t)}^p t^{p-1-p\theta}dx \,dt
\\\nonumber&
	\leq C \delta^{\theta p}\int_{\R^n} \int_\varepsilon^\delta 
	\biggl(\fint_{\Omega(x,t)} \abs{\nabla u}^2\biggr)^{p/2} t^{p-1-p\theta}dx \,dt
\end{align}
and so $f_\delta\to f$ in $L^p(\R^n)$; furthermore, notice that $u(x,t)\to f(x)$ as $t\to 0^+$ almost everywhere.
\end{proof}

%% file: sec-7-sobolev.tex
\chapter[Lebesgue and Sobolev Spaces]{Results for Lebesgue and Sobolev Spaces: Historic Account and some Extensions}
\label{chap:sobolev}

We wish to establish well-posedness of the boundary-value problems $(D)^A_{p,\theta}$ and $(N)^A_{p,\theta}$ under certain assumptions on $A$, $p$ and~$\theta$. To do this, we will make use of the extensive literature concerning the Neumann and regularity problems $(N)^A_{p,1}$ and~$(D)^A_{p,1}$ with boundary data in Lebesgue or Sobolev spaces; recall that Corollary~\ref{cor:well-posed:1} assumes well-posedness of $(D)^A_{p,1}$ or $(N)^A_{p,1}$. Thus, in this section we will review some known results concerning such problems. We will review only the results we will use in Chapters~\ref{chap:green} and \ref{chap:invertibility}; this chapter is by no means a complete history of the topic.

We will begin by considering layer potentials acting on Lebesgue and Sobolev spaces, that is, with the analogues to Theorems~\ref{thm:bounded} and~\ref{thm:trace}. Suppose that $A$ is elliptic and $t$-independent and that $A$ and $A^*$ satisfy the De Giorgi-Nash-Moser condition.
By \cite{HofMitMor}, there is some $\varepsilon>0$ such that the following estimates hold. We remark that it is not clear whether the exponent $p$ given by $1/p=1/2-\varepsilon$ is necessarily equal to the $p^+$ of Lemma~\ref{lem:PDE2}.
\begin{align}
\label{eqn:S:NTM}
\doublebar{\widetilde N_+(\nabla\s^A f)}_{L^p(\R^n)}
	&\leq C(p)\doublebar{f}_{H^p(\R^n)} 
	&& \text{if } \frac{1}{2}-\varepsilon < \frac{1}{p} < 1+\frac{\alpha}{n}
,\\
\label{eqn:S:slices}
\sup_{t\neq 0}\doublebar{\nabla\s^A f(\,\cdot\,,t)}_{L^p(\R^n)}
	&\leq C(p)\doublebar{f}_{H^p(\R^n)} 
	&& \text{if } \frac{1}{2}-\varepsilon < \frac{1}{p} < 1+\frac{\alpha}{n}
,\\
\label{eqn:D:NTM:0}
\doublebar{ N(\D^A f)}_{L^p(\R^n)}
	&\leq C(p)\doublebar{f}_{L^p(\R^n)} 
	&& \text{if } 0<\frac{1}{p}<\frac{1}{2}+\varepsilon
.\end{align}
We observe that Theorem~\ref{thm:bounded:holder} also comes from \cite{HofMitMor}; we included a proof above for the sake of completeness and also because the theorem in \cite{HofMitMor} does not formally state the result in the case of the single layer potential.

By \cite[Proposition~5.19]{HofKMP13},
\begin{align}
\label{eqn:D:NTM:1:banach}
\doublebar{\widetilde N_+(\nabla\D f)}_{L^p(\R^n)}
	&\leq C(p)\doublebar{\nabla_\parallel f}_{L^p(\R^n)} 
	&& \text{if } 1/2-\varepsilon < 1/p < 1
.\end{align}

Although we will not make use of these bounds, we remark that by \cite[Corollary~3.3]{HofMayMou}, we have the square-function estimates
\begin{align}
\label{eqn:square:0}
\int_{\R^n}\biggl(\int_{\abs{x-y}<t} \abs{\nabla\s f(y,t)}^2 \frac{dy\,dt}{t^{n-1}}\biggr)^{p/2}dx \leq C(p) \doublebar{f}_{\dot W^p_{-1}(\R^n)}^p
	&& \text{if } 0<\frac{1}{p}\leq\frac{1}{2}
,\\
\label{eqn:square:S:1}
\int_{\R^n}\biggl(\int_{\abs{x-y}<t} \abs{\nabla\partial_t \s f(y,t)}^2 \frac{dy\,dt}{t^{n-1}}\biggr)^{p/2}dx \leq C(p) \doublebar{f}_{L^p(\R^n)}^p
	&& \text{if } 0<\frac{1}{p}<1
\end{align}
and the Carleson-measure estimate
\begin{align}
\label{eqn:carleson:0}
\fint_Q\int_0^{\ell(Q)} \abs{\nabla (\s\nabla)\cdot \vec f(x,t)}^2\,t\,dt\,dx 
	&\leq C \doublebar{\vec f}_{L^\infty(\R^n)}^2 && \text{for all cubes } Q\subset\R^n
.\end{align}

\begin{rmk}\label{rmk:not-interpolation}
The boundedness results of Section~\ref{sec:layers:bounded} do \emph{not} follow by interpolating the bounds \eqref{eqn:S:NTM}--\eqref{eqn:carleson:0}.

If $p=2$ or $p=\infty$, then we may obtain the estimate $\doublebar{\T\vec f}_{\dot W(p,0,2)}\leq C\doublebar{\vec f}_{L^p(\R^n)}$ from the bounds \eqref{eqn:square:0} or~\eqref{eqn:carleson:0}; interpolating gives validity of these bounds for $2<p<\infty$. This provides the $\theta=0$ endpoint of the bound \eqref{eqn:space:T} for such~$p$. However, the bound $\doublebar{\T\vec f}_{\dot W(p,\theta,1)}\leq C\doublebar{\vec f}_{\dot H^p_1(\R^n)}$ is \emph{not} valid; notice that if $u\in \dot W(p,1,2)$ then $\Tr u$ is necessarily constant!

Nontangential maximal estimates inherently involve $L^\infty$ norms, and most interpolation methods require that at least one of the spaces considered be separable; see, for example, the remark on p.~519 of \cite{MenM00}. Thus, interpolating nontangential estimates is highly problematic. 

It is possible to interpolate between the bounds \eqref{eqn:square:0} and \eqref{eqn:square:S:1} (using the Caccioppoli inequality); this yields the strange bound
\[\int_{\R^n}\biggl(\int_{\abs{x-y}<t} \abs{\nabla\partial_t \s f(y,t)}^2 \frac{dy\,dt}{t^{n-3+2\theta}}\biggr)^{p/2}dx \leq C(p) \doublebar{f}_{\dot B^{p,p}_{\theta-1}(\R^n)}^p
.\]
If $p=2$ then we may use Lemma~\ref{lem:space:down} to eliminate the $\partial_t$ derivative (and recover the familiar formula $\doublebar{\nabla \s f}_{\dot W(2,\theta,2)}\leq C\doublebar{f}_{\dot B^{2,2}_{\theta-1}(\R^n)}$); however, it is not clear how this formula may be rewritten for general~$p$, and it is also unclear how to interpolate between this space and~$\dot W(\infty,\theta,2)$.
\end{rmk}

The bound \eqref{eqn:D:NTM:1:banach} is also valid in a wider range.
\begin{lem} \label{lem:D:NTM:1} Suppose that $\max(1/p^+,1/2-\varepsilon)< 1/p\leq 1+\alpha/n$. Then
\begin{align}
\label{eqn:D:NTM:1}
\doublebar{\widetilde N_+(\nabla\D f)}_{L^p(\R^n)}
	&\leq C(p)\doublebar{\nabla_\parallel f}_{L^p(\R^n)} 
.\end{align}
\end{lem}

\begin{proof}
The only new result is the case $1/p\geq 1$. It suffices to establish the bound \eqref{eqn:D:NTM:1} for $\dot H^p_1(\R^n)$ atoms.

Let $a_Q$ be an atom supported in a cube~$Q$. Because the lemma is valid for $p=2$ we have that 
\begin{align*}\int_{16Q} \abs{\widetilde N_+(\nabla \D a_Q)(x)}^p\,dx
&\leq 
	C\ell(Q)^{n-np/2}\biggl(\int_{16Q} \abs{\widetilde N_+(\nabla \D a_Q)(x)}^2\,dx\biggr)^{p/2}
\\&\leq 
	C\ell(Q)^{n-np/2}\doublebar{\nabla a_Q}_{L^2(\R^n)}^{p/2}
\leq C.
\end{align*}

To bound $\int_{\R^n\setminus 16Q} \abs{\widetilde N_+(\nabla \D a_Q)(x)}^p\,dx$, we follow \cite{HofMitMor} and \cite{KenP93}.

If $j\geq 1$ then let $A_j = 2^{j+4}Q\setminus 2^{j+3}Q$. We let 
\begin{align*}
\widetilde N_1 F(x)&=\sup_{\abs{x-z}<t<2^{j-3}\ell(Q)} \biggl(\fint_{\Omega(z,t)} \abs{F}^2\biggr)^{1/2}
,\\
\widetilde N_2 F(x)&=\sup_{\abs{x-z}<t, \>t\geq 2^{j-3}\ell(Q)} \biggl(\fint_{\Omega(z,t)} \abs{F}^2\biggr)^{1/2}
\end{align*}
By the Caccioppoli inequality, if $x\in A_j$ then
\begin{align*}
\widetilde N_2 (\nabla \D a_Q)(x)
&\leq \frac{C}{2^j\ell(Q)} \sup_{t>\abs{x-z}/2}\abs{\D a_Q(z,t)}.
\end{align*}
But if ${t>\abs{x-z}/2}$, then
\begin{align*}
\abs{\D a_Q(z,t)} 
&=
	\abs[bigg]{\int_Q \overline{\nu\cdot A(y)\nabla\Gamma_{(z,s)}^{A^*}(y,0)} \,a_Q(y)\,dy}
\\&
\leq
	C \ell(Q)^{1-n/p}
	\int_Q \abs{\nabla\Gamma_{(z,s)}^{A^*}(y,0)} \,dy
.\end{align*}
Applying Lemma~\ref{lem:slabs}, the Caccioppoli inequality, and the bound~\eqref{eqn:FS:holder}, we see that  if $\dist(Q,(z,t))\approx 2^j\ell(Q)$ then
\begin{align}
\label{eqn:D:atom:decay}
\abs{\D a_Q(z,t)} 
&\leq
	C \ell(Q)^{1-n/p} 2^{-j(n-1+\alpha)}
\end{align}
and so if $x\in A_j$ then
\begin{align*}
\widetilde N_2 (\nabla \D a_Q)(x)
&\leq \ell(Q)^{-n/p} 2^{-j(n+\alpha)}
.\end{align*}

Now we consider $N_1(\nabla\D a_Q)(x)$ for  $x\in A_j$, $j\geq 1$. If $\abs{x-z}<t<2^{j-3}\ell(Q)$ and $c$ is a constant, then by the Caccioppoli inequality,
\begin{align*}
\biggl(\fint_{B((z,t),t/2)} \abs{\nabla\D a_Q}^2\biggr)^{1/2}
&\leq
	\frac{C}{t}\biggl(\fint_{B((z,t),3t/4)} \abs{\D a_Q-c}^2\biggr)^{1/2}
.\end{align*}
Following \cite{KenP93}, we then write
\begin{align*}
\biggl(\fint_{B((z,t),t/2)} \abs{\nabla\D a_Q}^2\biggr)^{1/2}
&\leq
	\frac{C}{t}\biggl(\fint_{B((z,t),3t/4)} \abs{\D a_Q(y,s)-\D a_Q(y,0)}^2\,dy\biggr)^{1/2}
	\\&\qquad+
	\frac{C}{t}\biggl(\fint_{\Delta(z,3t/4)} \abs{\D a_Q(y,0)-c}^2\,dy\biggr)^{1/2}
.\end{align*}
Because $\Div A\nabla \partial_t \D a_Q=0$ away from~$Q$, we have that $\partial_s \D a_Q(y,s)$ is locally bounded; thus, by the bounds~\eqref{eqn:D:atom:decay} and~\eqref{eqn:local-bound:2}, \[\abs{\partial_s \D a_Q(y,s)}\leq C \ell(Q)^{-n/p} 2^{-j(n+\alpha)}\] and so
\begin{equation*}\frac{C}{t}\biggl(\fint_{B((z,t),3t/4)} \abs{\D a_Q(y,s)-\D a_Q(y,0)}^2\,dy\biggr)^{1/2}
\leq
C \ell(Q)^{-n/p} 2^{-j(n+\alpha)}.
\end{equation*}
By Sobolev's inequality, if we choose $c$ appropriately,
\begin{align*}
	\frac{C}{t}\biggl(\fint_{\Delta(z,3t/4)} \abs{\D a_Q(y,0)-c}^2\,dy\biggr)^{1/2}
&\leq
	\frac{C}{t}\biggl(\fint_{\Delta(x,2t)} \abs{\D a_Q(y,0)-c}^2\,dy\biggr)^{1/2}
\\&\leq
	C\biggl(\fint_{\Delta(x,2t)} \abs{\nabla_\|\D a_Q(y,0)}^{2^*}\,dy\biggr)^{1/{2^*}}
\\&\leq
	C M(\1_{\widetilde A_j}\abs{\nabla_\|\D a_Q}^{2^*})(x)^{1/{2^*}}
\end{align*}
where ${2^*}=2n/(n+2)<2$, $M$ is the Hardy-Littlewood maximal operator and $\widetilde A_j = 2^{j+5}Q\setminus 2^{j+2}Q$.

Thus, we have that
\begin{align*}
\fint_{A_j} \widetilde N_+(\nabla\D a_Q)(z)\,dz
&\leq
	\fint_{A_j} \widetilde N_1(\nabla\D a_Q)(z)\,dz
	+	\fint_{A_j} \widetilde N_2(\nabla\D a_Q)(z)\,dz
\\&\leq
	C \ell(Q)^{-n/p} 2^{-j(n+\alpha)}
	+	C\fint_{A_j} M(\1_{\widetilde A_j}\abs{\nabla_\|\D a_Q}^{2^*})(z)^{1/{2^*}}\,dz.
\end{align*}
The Hardy-Littlewood maximal operator is bounded on $L^{2/2^*}(\R^n)$ because $2/2^*>1$, and so by H\"older's inequality
\begin{align*}
\fint_{A_j} \widetilde N_+(\nabla\D a_Q)(z)\,dz
&\leq
	C \ell(Q)^{-n/p} 2^{-j(n+\alpha)}
	+	C\biggl(\fint_{\widetilde A_j} \abs{\nabla_\|\D a_Q}^{2}\biggr)^{1/2}.
\end{align*}
We may bound the second term using Lemma~\ref{lem:slabs}, and so by H\"older's inequality and because $p\leq 1$ we have that
\[\int_{A_j} \widetilde N_+(\nabla\D a_Q)(z)^p\,dz
\leq
C2^{-j(np-n+p\alpha)}
.\]
If $1/p<1+\alpha/n$, then $np-n+p\alpha>0$ and so we may sum in $j$ to see that $\widetilde N_+(\nabla\D a_Q)\in L^p(\R^n)$, as desired.
\end{proof}

The analogue to the trace theorem~\ref{thm:trace} for solutions to $(D)^A_{p,1}$ and $(N)^A_{p,1}$ is as follows.

\begin{thm}[{\cite[Lemmas~6.1 and~6.2]{HofMitMor}}]
\label{thm:trace:NTM}
Suppose that $A$ is elliptic, $t$-in\-de\-pen\-dent and that $A$ and $A^*$ satisfy the De Giorgi-Nash-Moser condition. Suppose that $u$ satisfies
$\Div A\nabla u=0$ in $\R^{n+1}_+$ and that $\widetilde N_+(\nabla u)\in L^p(\R^n)$ for some $p$ with $n/(n+1)<p<\infty$. Then there is a function $f\in \dot H^p_1(\R^n)$ and a function $g\in H^p(\R^n)$ such that
\begin{equation*}\doublebar{f}_{\dot H^p_1(\R^n)}\leq C(p)\doublebar{\widetilde N_+(\nabla u)}_{L^p(\R^n)},
\quad
\doublebar{g}_{H^p(\R^n)}\leq C(p)\doublebar{\widetilde N_+(\nabla u)}_{L^p(\R^n)}\end{equation*}
and such that $\nu\cdot A\nabla u\big\vert_{\partial\R^{n+1}_+}=g$ in the sense of formula~\eqref{eqn:conormal:1} and such that $u\to f$ in the sense of nontangential limits, that is, in the sense that
\begin{equation*}\lim_{(y,s)\to (x,0),\>(y,s)\in\gamma_+(x)} u(y,s)=f(x) \text{ for almost every $x\in\R^n$}\end{equation*}
where the limit is taken only over points $(y,s)$ in the cone $\gamma_+(x)$ of formula~\eqref{eqn:cone}.

Furthermore, if $1/2-\varepsilon<1/p<1$, then
\begin{equation*}\sup_{t>0}\doublebar{\nabla u}_{L^p(\R^n)}\leq C(p)\doublebar{\widetilde N_+(\nabla u)}_{L^p(\R^n)}\end{equation*}
and  $-\e_{n+1}\cdot A\nabla u(\,\cdot\,,t)\rightharpoonup g$,
$\nabla_\parallel u(\,\cdot\,,t)\rightharpoonup \nabla_\parallel f$
weakly in $L^p(\R^n)$ as $t\to 0^+$.
\end{thm}

\begin{rmk} Let $p_0$ satisfy $n/(n+1)<p_0<\infty$ and let $p_1$, $\theta_1$ satisfy the conditions of Theorem~\ref{thm:trace}. Suppose that $\Div A\nabla u=0$ in $\R^{n+1}_+$ and that both 
$\widetilde N_+(\nabla u)\in L^{p_0}(\R^n)$ and $u\in \dot W(p_1,\theta,2)$. Then by the bound~\eqref{eqn:besov-trace-convergence}, \[f(x)=\lim_{t\to 0^+} u(x,t)=\Tr u(x)\] for almost every $x\in\R^n$, where $f$ is as in Theorem~\ref{thm:trace:NTM} and~$\Tr u$ is as in Theorem~\ref{thm:trace}; thus, for such~$u$ the traces given by Theorems~\ref{thm:trace:NTM} and~\ref{thm:trace} are the same.
\end{rmk}

Furthermore, if $n/(n+1)<p_0<\infty$ and $\widetilde N_+(\nabla u)\in L^{p_0}(\R^n)$, then by the following theorem, there is always some $p_1$ and $\theta$ satisfying the conditions of Theorem~\ref{thm:trace} such that $u\in \dot W(p_1,\theta,2)$. 
\begin{thm}\label{thm:NTM:embedding}
Let $0<p_0<p_1\leq \infty$, and suppose that ${\theta_1}-n/p_1=1-n/p_0$. Then $u\in \dot W(p_1,\theta_1,2)$ whenever $\widetilde N_+(\nabla u)\in L^{p_0}(\R^n)$.
\end{thm}

Recalling the definitions of $\widetilde N_+$ and $\dot W(p,\theta,2)$, 
we see that Theorem~\ref{thm:NTM:embedding} follows by applying the following lemma to the function $F(x,t)=\bigl(\fint_{\Omega(x,t)} \abs{\nabla u}^2\bigr)^{p/2}$.

\begin{lem}\label{lem:NTM:embedding}
Suppose that $N_+F \in L^1(\R^n)$. Let $\sigma>-1$. Then 
\[\int_{\R^{n+1}_+}\abs{F(x,t)}^{1+1/n+\sigma/n} \,t^\sigma\,dx\,dt\leq C(\sigma)\doublebar{NF}_{L^1(\R^n)}^{1+1/n+\sigma/n} \]
where the constant $C(\sigma)$ depends only on the dimension $n+1$ and the number~$\sigma$.
\end{lem}

The $\sigma=0$ case of this lemma is by now fairly well known;  see, for example, \cite[Lemma~2.2]{HofMitMor}.

\begin{proof}[Proof of Lemma~\ref{lem:NTM:embedding}]
If $\beta>0$, then let $e(\beta)={\{x\in\R^n:N_+F(x)>\beta\}}$. Notice that by definition of the Lebesgue integral,
\[\int_{\R^n} N_+F(x)\,dx = \int_0^\infty \abs{e(\beta)}\,d\beta\]
where $\abs{e(\beta)}$ denotes Lebesgue measure.

Now, let $E(\beta,t)={\{x\in\R^n:\abs{F(x,t)}>\beta\}}$.
We again apply the definition of the Lebesgue integral to see that
\[\int_0^\infty\int_{\R^n}\abs{F(x,t)}^\rho \,dx\,t^\sigma\,dt = \int_{0}^\infty \int_{0}^\infty \rho \beta^{\rho -1}\abs{E(\beta,t)}\,d\beta\,t^\sigma\,dt\]
where $\rho=1+1/n+\sigma/n>1$.
Observe that $E(\beta,t)\subset e(\beta)$ for all $t>0$. Furthermore, if there is even single point $x$ such that $\abs{F(x,t)}>\beta$, then $B(x,t)\subset e(\beta)$; thus, if $t>c(n)\abs{e(\beta)}^{1/n}$, then $E(\beta,t)=\emptyset$. Thus,
\[\int_0^\infty\int_{\R^n}\abs{F(x,t)}^\rho \,dx\,t^\sigma\,dt 
\leq
	\int_0^\infty \int_{0}^{c(n)\abs{e(\beta)}^{1/n}} \rho \beta^{\rho -1} \abs{e(\beta)}\,t^\sigma\,dt\,d\beta 
.\]
But if $\sigma>-1$, then we may evaluate this integral to see that
\begin{align*}
\int_{\R^{n+1}_+}\abs{F(x,t)}^\rho  \,t^\sigma\,dx\,dt
&\leq
	C(\sigma)
	\int_0^\infty  \beta^{\rho -1}
	\abs{e(\beta)}^{1+(\sigma+1)/n}\,d\beta 
\\&=
	C(\sigma)
	\int_0^\infty  \beta^{\rho -1}
	\abs{e(\beta)}^{\rho }\,d\beta 
.\end{align*}
Applying the bound $\beta\abs{e(\beta)}\leq {\doublebar{N_+F}_{L^1(\R^n)}}$, we see that
\begin{align*}
\int_{\R^{n+1}_+}\abs{F(x,t)}^\rho  \,t^\sigma\,dx\,dt
&\leq
	C(\sigma) \doublebar{N_+F}_{L^1(\R^n)}^{\rho -1}
	\int_0^\infty 
	\abs{e(\beta)}\,d\beta 
=
	C(\sigma)\doublebar{N_+F}_{L^1(\R^n)}^{\rho }
\end{align*}
as desired.
\end{proof}

\begin{rmk} Theorems~\ref{thm:whitney:embedding} and~\ref{thm:NTM:embedding} allow us to resolve the issue of compatibility for certain values of~$p$. Specifically, suppose that $(D)^A_{2n/(n+1),1}$ is solvable. But if $u$ is a $(D)^A_{2n/(n+1),1}$ solution then $\widetilde N_+(\nabla u)\in L^{2n/(n+1)}(\R^n)$, and so $\nabla u \in L^2(\R^{n+1}_+)=L(2,1/2,2)$. Thus, all $(D)^A_{2n/(n+1),1}$-solutions are in fact $(D)^A_{2,1/2}$ solutions, and so if $(D)^A_{2n/(n+1),1}$ is solvable then it must be compatibly solvable. The same is true of $(N)^A_{2n/(n+1),1}$.

A similar result is valid for special values of~$p$ in the range $0<\theta<1$. Suppose that $1/2-n/2=\theta-n/p$. If $1/2<\theta\leq1$ and $(D)^A_{p,\theta}$ or $(N)^A_{p,\theta}$ is solvable, then by the argument above, that problem must be compatibly solvable. Conversely, if $0<\theta<1/2$ then any $(D)^A_{2,1/2}$ or $(N)^A_{2,1/2}$ solution is a $(D)^A_{p,\theta}$ or $(N)^A_{p,\theta}$ solution, and so if $(D)^A_{p,\theta}$ or $(N)^A_{p,\theta}$ is well-posed then it must be compatibly well-posed.
\end{rmk}

We now consider the De Giorgi-Nash-Moser condition and well-posedness of boundary-value problems. That is, we will remind the reader of some known sufficient conditions for the matrix $A$ to satisfy the De Giorgi-Nash-Moser condition or for the Neumann and regularity problems $(N)^A_{p,1}$ and $(D)^A_{p,1}$ to be well-posed. 

We begin with the De Giorgi-Nash-Moser condition. 
\begin{thm}\label{thm:DGNM}
Let $A$ be an elliptic matrix. Suppose that one of the following conditions holds.
\begin{itemize}
\item $A$ is constant,
\item $A$ has real coefficients,
\item The ambient dimension $n+1=2$, or
\item $A$ is $t$-independent and the ambient dimension $n+1=3$.
\end{itemize}
Then $A$ satisfies the De Giorgi-Nash-Moser condition, with constants $H$ and $\alpha$ depending only on the dimension $n+1$ and the ellipticity constants~$\lambda$ and~$\Lambda$.
\end{thm}
The theorem was proven for real symmetric coefficients $A$ by De Giorgi and Nash in \cite{DeG57,Nas58} and extended to real nonsymmetric coefficients by Morrey in \cite{Mor66}. In dimension $n+1=2$ the theorem follows from Lemma~\ref{lem:PDE2} and Lemma~\ref{lem:morrey}; the case $n+1=3$ was established in \cite[Appendix~B]{AlfAAHK11}.
We observe that in dimension $n+1\geq 4$, or in dimension $n+1=3$ for $t$-dependent coefficients, the De Giorgi-Nash-Moser condition may fail; see \cite{Fre08} for an example.

We now summarize some sufficient conditions for well-posedness of boundary-value problems in the case $p=2$.
\begin{thm}\label{thm:well-posed:L2}
Let $A$ be an elliptic, $t$-independent matrix. 


Suppose that $A$ is either constant or self-adjoint.
Then $(N)^A_{2,1}$ and $(D)^A_{2,1}$ are both well-posed in $\R^{n+1}_+$ and $\R^{n+1}_-$.

If $A$ is block upper triangular then $(N)^A_{2,1}$ is solvable, and if $A$ is block lower triangular then $(D)^A_{2,1}$ is solvable. Thus, if $A$ is a block matrix then both problems are solvable.
\end{thm}

Well-posedness results for such coefficients are also valid in the case of the Dirichlet problem, although because we do not require $A$ to satisfy the De Giorgi-Nash-Moser condition, the Dirichlet problem must be formulated slightly differently from our $(D)^A_{2,0}$ (for example, by taking the square-function estimate~\eqref{eqn:square-function} or the estimate $\widetilde N_+ u\in L^p(\R^n)$ in place of the bound $N_+u\in L^p(\R^n)$).

Well-posedness of $(N)^A_{2,1}$ and $(D)^A_{2,1}$ for real symmetric coefficients were established in \cite{KenP93}. (Well-posedness of $(D)^A_{2,0}$ for such $A$ was established in \cite{JerK81a}.)
In the case of complex self-adjoint matrices, the result was established in \cite{AusAM10}, as was solvability of $(D)^A_{2,0}$ with square-function or modified nontangential estimates.


Block and block triangular matrices are defined as follows. We may write the $t$-independent matrix $A$ in the form
\begin{equation*}
A(x)=\begin{pmatrix}A_\parallel(x) &\vec a_{\parallel \perp}
\\\noalign{\smallskip} \vec a_{ \perp \parallel}^T & a_\perp(x)\end{pmatrix}
\end{equation*}
for some $n\times n$ matrix $A_\parallel$, some $1\times n$ column vectors $\vec a_{\perp\parallel}$ and $\vec a_{\parallel\perp}$, and some complex-valued function $a_\perp$. 
A block upper or lower triangular matrix is a matrix for which the vectors $\vec a_{\perp\parallel}$ or $\vec a_{\parallel\perp}$, respectively, are identically equal to zero; a block matrix is one for which both vectors are zero.
Well-posedness in the block case follows from validity of the Kato conjecture, as explained in \cite[Remark~2.5.6]{Ken94}; see \cite{AusHLMT02} for the proof of the Kato conjecture and \cite[Consequence~3.8]{AxeKM06} for the case $a_\perp\not\equiv 1$. Solvability in the block triangular case was proven in \cite{AusMM13}. (Solvability of an appropriately modified $(D)^A_{2,0}$ for upper block triangular~$A$ was also established.)

Both well-posedness and the De Giorgi-Nash-Moser condition are stable under $t$-independent perturbation. That is, we have the following theorems.
\begin{thm}\label{thm:perturb:A:DGNM}
Let $A_0$ be an elliptic matrix. Suppose that $A_0$ and $A_0^*$ both satisfy the De Giorgi-Nash-Moser condition. Then there is some constant~$\varepsilon>0$, depending only on the dimension $n+1$ and the constants $\lambda$, $\Lambda$ in formula~\eqref{eqn:elliptic} and $H$, $\alpha$ in formula~\eqref{eqn:holder}, such that if $\doublebar{A-A_0}_{L^\infty(\R^{n+1})}<\varepsilon$, then $A$ satisfies the De Giorgi-Nash-Moser condition.
\end{thm}
This result is from 
\cite{Aus96}; see also \cite[Chapter~1, Theorems~6 and~10]{AusT98}.
\begin{thm}
\label{thm:perturb:A:well-posed}
Let $A_0$ be an elliptic, $t$-independent matrix. Suppose that $(N)^{A_0}_{2,1}$ or $(D)^{A_0}_{2,1}$ is well-posed in $\R^{n+1}_+$.

Then there is some constant~$\varepsilon>0$ such that if $\doublebar{A-A_0}_{L^\infty(\R^{n+1})}<\varepsilon$, then $(N)^{A}_{2,1}$ or $(D)^{A}_{2,1}$, respectively, is also well-posed in $\R^{n+1}_+$.
\end{thm}
This was proven for arbitrary elliptic, $t$-independent matrices $A$ in \cite{AusAM10}; the Dirichlet problem was also considered but their formulation is again somewhat different from our $(D)^A_{2,0}$. If $A_0$ and $A_0^*$ satisfy the De Giorgi-Nash-Moser condition and all three problems are well-posed, then this perturbation result was proven using layer potentials in \cite{AlfAAHK11}. $L^\infty$ perturbation results have also been investigated in \cite{FabJK84, AusAH08, Bar13}.

We also mention that stability of well-posedness under $t$-dependent perturbation in the sense of Carleson measures has also been investigated extensively; see \cite{Dah86a, Fef89, FefKP91, Fef93, KenP93, KenP95, DinPP07, DinR10, AusA11, AusR12, HofMayMou}. However, as the present monograph is concerned exclusively with $t$-independent coefficients, we will not discuss such perturbations further.

Some further statements may be made in the case of real nonsymmetric coefficients.
The following theorem was established in \cite{KenKPT00} in dimension $n+1=2$, and in higher dimensions in \cite{HofKMP12}.
\begin{thm}\label{thm:real:dirichlet}
Let $A$ be elliptic, $t$-independent and have real coefficients. Then there is some $p'<\infty$, depending only on the dimension $n+1$ and the ellipticity constants $\lambda$ and~$\Lambda$, such that $(D)^A_{p',0}$ is solvable in $\R^{n+1}_+$ and in $\R^{n+1}_-$. Furthermore, the square-function estimate \eqref{eqn:square-function} is valid for solutions to $(D)^A_{p',0}$.
\end{thm}

There are also a number of theorems which allow us to extrapolate well-posedness results.
\begin{thm}[{\cite[Theorem~1.11]{HofKMP13}}]\label{thm:extrapolation:HKMP}
Let $A$ and $A^*$ be elliptic, $t$-independent and satisfy the De Giorgi-Nash-Moser condition. Then there exists a number $\varepsilon>0$ such that the following is true.

Suppose that $p$ satisfies $1/2-\varepsilon<p<1$, and that $1/p+1/p'=1$. Then the following two statements are equivalent.
\begin{itemize}
\item $(R)^A_p$ is solvable in $\R^{n+1}_+$ and $\R^{n+1}_-$.

\item  $(D)^{A^*}_{p'}$ is solvable in $\R^{n+1}_+$ and $\R^{n+1}_-$, and solutions to $(D)^{A^*}_{p'}$ satisfy the square-function estimate~\eqref{eqn:square-function}.
\end{itemize}
Furthermore, if either of these two equivalent conditions is true, then the operator $\s^{A}_+$ is bounded and invertible $L^{p}(\R^n)\mapsto\dot W^{p}_{1}(\R^n)$, and the operator  $\s^{A^*}_+$ is bounded and invertible $\dot W^{p'}_{-1}(\R^n)\mapsto L^{p'}(\R^n)$.
\end{thm}

Combined with Theorem~\ref{thm:real:dirichlet}, we see that if $A$ is real then there is some $p>1$ such that $(D)^A_{p,1}$ is solvable. In Chapter~\ref{chap:invertibility}, we will see that $(D)^A_{p,1}$ is \emph{compatibly} solvable.
In dimension~2, by \cite{KenR09} there is also some $p>1$ such that $(N)^A_{p,1}$ is compatibly solvable; it is not known if this result holds in higher dimensions.

Given some additional good behavior of solutions near $\partial\R^{n+1}_+$, a further extrapolation theorem is available.
\begin{thm}[\cite{AusM}]\label{thm:extrapolation:AusM}
Let $A$ be elliptic and $t$-independent, and let $A^\sharp$ be as in formula~\eqref{eqn:AusM}. Suppose that $A$, $A^*$, $A^\sharp$ and $(A^\sharp)^*$ all satisfy the De Giorgi-Nash-Moser condition with exponent~$\alpha^\sharp$.

If $1<p_0\leq 2$, and if $(D)^A_{p_0,1}$ is compatibly well-posed in $\R^{n+1}_+$, then $(D)^A_{p,1}$ is compatibly well-posed for all $p$ with $n/(n+\alpha^\sharp)<p<p_0$. 

If $1<p_0\leq 2$, and if $(N)^A_{p_0,1}$ is compatibly well-posed in $\R^{n+1}_+$, then $(N)^A_{p,1}$ is compatibly well-posed for all $p$ with $n/(n+\alpha^\sharp)<p<p_0$.

If $\theta=n(1/p-1)$ and $0<\theta<1$ (equivalently $n/(n+\alpha)<p<1$),  then compatible well-posedness of $(D)^A_{p,1}$ implies compatible well-posedness of $(D)^{A^*}_{\infty,\theta}$, and compatible well-posedness of $(N)^A_{p,1}$ implies compatible well-posedness of $(N)^{A^*}_{\infty,\theta}$.

Finally, if $1<p<\infty$, then $(D)^A_{p,1}$ is compatibly well-posed if and only if $(D)^{A^*}_{p',0}$, with the square-function estimate~\eqref{eqn:square-function} in place of the bound $N_+u\in L^p(\R^n)$, is compatibly well-posed. Similarly, if $1<p<\infty$ then $(N)^A_{p,1}$ is compatibly well-posed if and only if $(N)^{A^*}_{p',0}$ with square-function estimates is compatibly well-posed.
\end{thm}

Notice that if $A$ is real-valued or if $n+1=2$, then the conditions of Theorem~\ref{thm:extrapolation:AusM} are valid.

The condition that $A^\sharp$ satisfy the De Giorgi-Nash-Moser condition is connected to the ``boundary De Giorgi-Nash-Moser condition''. It is not clear whether the interior De Giorgi-Nash-Moser condition of Definition~\ref{dfn:DGNM} implies its boundary analogue.
Notice that there is necessarily some $\alpha\geq \alpha^\sharp$ such that $A$ and $A^*$ satisfy the De Giorgi-Nash-Moser condition with exponent $\alpha$. In some cases $\alpha$ is greater than~$\alpha^\sharp$; in the example given in Section~\ref{sec:sharp}, for example, we may take $\alpha=1$, but $\alpha^\sharp>0$ can be made arbitrarily small.

%% file: sec-8-green.tex
\chapter{The Green's Formula Representation for a Solution}
\label{chap:green}

In this section, we will prove a Green's formula representation for solutions to $\Div A\nabla u=0$; combined with the results of Section~\ref{chap:invertibility}, this representation formula will let us prove uniqueness of solutions.

\begin{thm} \label{thm:green}
Let $A$ be as in Theorem~\ref{thm:bounded}.
Suppose that $\Div A\nabla u=0$ in $\R^{n+1}_+$, and that $u\in\dot W(p,\theta,2)$
for some $\theta$, $p$ that satisfy the conditions of Theorem~\ref{thm:bounded}. In addition suppose $\theta<\alpha+n/p$.
By Theorem~\ref{thm:trace}, $f=\Tr u$ and $g=\nu_+\cdot A\nabla u\big\vert_{\partial\R^{n+1}_+}$ exist, and $f\in \dot B^{p,p}_\theta(\R^n)$, $g\in \dot B^{p,p}_{\theta-1}(\R^n)$.

Then $u=-\D f+\s g$ up to an additive constant.

If $\widetilde N_+(\nabla u)\in L^p(\R^n)$ for some $p$ with $n/(n+\alpha)<p<p^+$, then again $f=\Tr u$ and $g=\nu_+\cdot A\nabla u\big\vert_{\partial\R^{n+1}_+}$ exist and $u=-\D f+\s g$ up to an additive constant. We think of this as the $\theta=1$ endpoint of our representation formula.
\end{thm}

Notice that we may strengthen the lemma to include all $\theta$ with $\theta<\alpha_0+n/p$. Recall that by Lemma~\ref{lem:morrey}, if $A$ is $t$-independent then $\alpha_0\geq 1-n/p^+$, and so the condition $\theta<\alpha_0+n/p$ poses no restraints beyond those imposed by the conditions of Theorem~\ref{thm:bounded}.

\begin{rmk}\label{rmk:green:lower} If $\Div A\nabla u=0$ in the lower half-plane $\R^{n+1}_-$ instead of~$\R^{n+1}_+$, and if $u$ satisfies appropriate estimates, then
\begin{equation*}u(x,t)=\D f-\s g\end{equation*}
where $f=\Tr u=u\big\vert_{\partial\R^{n+1}_-}$ and $g=\nu_+\cdot A\nabla u\big\vert_{\partial\R^{n+1}_-}=-\nu_-\cdot A\nabla u\big\vert_{\partial\R^{n+1}_-}$.
\end{rmk}

\begin{proof}[Proof of Theorem~\ref{thm:green}]
By Theorem~\ref{thm:whitney:embedding} or Theorem~\ref{thm:NTM:embedding}, we may assume without loss of generality that $p>1$ and that $0<\theta<1$.

Fix some $(x,t)\in\R^{n+1}_+$, and choose some $R$, $\delta$ with $0<\delta<t/6$, $\abs{t}+\abs{x}<R$. 
Define the following cutoff functions. Let $\xi_R:\R^n\mapsto\R$ be smooth, supported in $\Delta(0,2R)$ and 1 in $\Delta(0,R)$. Let $\eta_\delta$ be smooth with $\eta_\delta(s)=0$ for $s<\delta/6$ and $\eta_\delta(s)=1$ for $s>\delta/3$, and let $\zeta_R$ be smooth with $\zeta_R(s)=1$ for $s<R$ and $\zeta_R(s)=0$ for $s>2R$. Let $\varphi(y,s)=\varphi_{R,\delta}(y,s)=\xi_R(y)\zeta_R(s)\eta_\delta(s)$. We impose the obvious bounds on $\nabla\varphi$.

We will want to use formula \eqref{eqn:FS:point}; however, this formula is only known to be valid for smooth functions. Thus, we will need a smooth approximation to~$u$. Choose some (very small) numbers $\beta$, $\gamma>0$. Let $\theta_\beta$ be a smooth mollifier, supported in $B(0,\beta)$ and integrating to~$1$, and let $V_\beta=u*\theta_\beta$. Because $\Gamma_{(x,t)}^{A^*}$ has a pole at $(x,t)$ it will be convenient if $U_{\beta,\gamma}$ is constant in a neighborhood of~$(x,t)$. Let $\psi_\gamma$ be a smooth cutoff function supported in $B((x,t),\gamma)$ and identically equal to 1 in $B((x,t),\gamma/2)$. Let
\begin{equation*}U_{\beta,\gamma}(y,s)=V_\beta(y,s)(1-\psi_\gamma(y,s))+u(x,t)\psi_\gamma(y,s).\end{equation*}

Then because $U_{\beta,\gamma}$ is smooth, we have that by the formula~\eqref{eqn:FS:point},
\begin{align*}
u(x,t) &=U_{\beta,\gamma}(x,t)\varphi(x,t)
	=\int_{\R^{n+1}_+} \nabla (U_{\beta,\gamma}\varphi)\cdot \overline{A^*\nabla \Gamma_{(x,t)}^{A^*}}
\\&=
	\int_{\R^n}\int_0^\delta \nabla (U_{\beta,\gamma}\varphi)\cdot \overline{A^*\nabla \Gamma_{(x,t)}^{A^*}}
	+
	\int_{\R^n}\int_\delta^\infty\nabla (U_{\beta,\gamma}\varphi)\cdot \overline{A^*\nabla \Gamma_{(x,t)}^{A^*}}
.\end{align*}
Because $\Div A^*\nabla \Gamma_{(x,t)}^{A^*}=0$ in $\R^n\times(0,\delta)$, by the formula~\eqref{eqn:conormal:1} for the conormal derivative, the first integral is equal to $-\D (U_{\beta,\gamma}(\,\cdot\,,\delta)\xi_R)(x,t-\delta)$; we observe that if $\gamma$ is small enough then $U_{\beta,\gamma}(\,\cdot\,,\delta)\to u(\,\cdot\,,\delta)$ as $\beta\to 0^+$, uniformly on compact sets. Let $f_\delta(x)=u(x,\delta)$. Then
\begin{align*}
u(x,t) 
&=
	-\D (f_\delta\xi_R)(x,t-\delta)
	+\lim_{\beta\to 0^+}
	\int_{\R^n}\int_\delta^\infty\nabla (U_{\beta,\gamma}\varphi)\cdot \overline{A^*\nabla \Gamma_{(x,t)}^{A^*}}
.\end{align*}
Let $\kappa_\delta$ be smooth with $\kappa_\delta(y,s)=\kappa_\delta(s)=1$ for $0<s<2\delta$, $\kappa_\delta(s)=0$ for $s>3\delta$. Then
\begin{align}\label{eqn:green:proof:3}
u(x,t) 
&=
	-\D (f_\delta\xi_R)(x,t-\delta)
	\\&\qquad\nonumber
	+\lim_{\beta\to 0^+}
	\int_{\R^n}\int_\delta^\infty\nabla (U_{\beta,\gamma}\kappa_\delta\varphi)\cdot \overline{A^*\nabla \Gamma_{(x,t)}^{A^*}}
	\\&\qquad\nonumber
	+\lim_{\beta\to 0^+}
	\int_{\R^n}\int_\delta^\infty\nabla (U_{\beta,\gamma}(1-\kappa_\delta)\varphi)\cdot \overline{A^*\nabla \Gamma_{(x,t)}^{A^*}}
.\end{align}
Of course the third term is also equal to $u(x,t)$; we will expand the second and third terms and exploit some cancellation features.

Consider the second term of formula~\eqref{eqn:green:proof:3}. Notice that we are integrating only over the region $S=\Delta(0,2R)\times(\delta,3\delta)$, and that both $\nabla u$ and $\nabla \Gamma_{(x,t)}^{A^*}$ lie in $L^2(S)$. Furthermore, $\nabla U_{\beta,\gamma}\to \nabla u$ as $\beta\to 0$ in $L^2(S)$. Thus, the second term satisfies
\begin{multline*}
\lim_{\beta\to 0^+}
	\int_{\R^n}\int_\delta^\infty\nabla (U_{\beta,\gamma}\kappa_\delta\varphi)\cdot \overline{A^*\nabla \Gamma_{(x,t)}^{A^*}}
\\\begin{aligned}
&=
	\int_{\R^n}\int_\delta^\infty\nabla (u\kappa_\delta\varphi)\cdot \overline{A^*\nabla \Gamma_{(x,t)}^{A^*}}
\\&=
	\int_{\R^n}\int_\delta^\infty	
	A\nabla u\cdot \nabla (\kappa_\delta\varphi\overline{\Gamma_{(x,t)}^{A^*}})
-
	\int_{\R^n}\int_\delta^\infty	\overline{\Gamma_{(x,t)}^{A^*}}\,
	A\nabla u\cdot \nabla (\kappa_\delta\varphi)
\\&\qquad+
	\int_{\R^n}\int_\delta^\infty
	u
	\nabla (\kappa_\delta\varphi)\cdot \overline{A^*\nabla \Gamma_{(x,t)}^{A^*}}
.\end{aligned}\end{multline*}
Approximating $\Gamma_{(x,t)}^{A^*}$ by smooth functions, we see that the first term is equal to
\begin{equation*}\int_{\R^n} -\e_{n+1}\cdot A(y)\nabla u(y,\delta)\, \xi_R(y) \, \overline{\Gamma_{(x,t)}^{A^*}(y,\delta)}\,dy
=\s (\xi_R g_\delta)(x,t-\delta)
\end{equation*}
where $g_\delta(y)=-\e_{n+1}\cdot A\nabla u(y,\delta)$. 

Now, consider the third term of formula~\eqref{eqn:green:proof:3}. We have that
\begin{multline*}
	\int_{\R^n}\int_\delta^\infty\nabla (U_{\beta,\gamma}(1-\kappa_\delta)\varphi)\cdot \overline{A^*\nabla \Gamma_{(x,t)}^{A^*}}
\\\begin{aligned}
&=
	\int_{\R^n}\int_\delta^\infty A\nabla U_{\beta,\gamma}\cdot \nabla ((1-\kappa_\delta)\varphi\overline{ \Gamma_{(x,t)}^{A^*}})
\\&\qquad+
	\int_{\R^n}\int_\delta^\infty
	U_{\beta,\gamma}\,
	\nabla ((1-\kappa_\delta)\varphi)\cdot \overline{A^*\nabla \Gamma_{(x,t)}^{A^*}}
	-\overline{ \Gamma_{(x,t)}^{A^*}}
	\,
	A\nabla U_{\beta,\gamma}\cdot \nabla ((1-\kappa_\delta)\varphi)
\end{aligned}\end{multline*}
and $U_{\beta,\gamma}\to u$ uniformly and in $\dot W^2_1$ in $\supp \nabla ((1-\kappa_\delta)\varphi)$; thus,
\begin{multline*}
\lim_{\beta\to 0^+}
	\int_{\R^n}\int_\delta^\infty\nabla (U_{\beta,\gamma}(1-\kappa_\delta)\varphi)\cdot \overline{A^*\nabla \Gamma_{(x,t)}^{A^*}}
\\\begin{aligned}
&=
	\lim_{\beta\to 0^+}
	\int_{\R^n}\int_\delta^\infty A\nabla U_{\beta,\gamma}\cdot \nabla ((1-\kappa_\delta)\varphi\overline{ \Gamma_{(x,t)}^{A^*}})
\\&\qquad
+
	\int_{\R^n}\int_\delta^\infty
	u\,
	\nabla ((1-\kappa_\delta)\varphi)\cdot \overline{A^*\nabla \Gamma_{(x,t)}^{A^*}}
-
	\overline{ \Gamma_{(x,t)}^{A^*}}
	\,
	A\nabla u\cdot \nabla ((1-\kappa_\delta)\varphi)
.\end{aligned}\end{multline*}

Thus, we may rewrite formula~\eqref{eqn:green:proof:3} as
\begin{align*}
u(x,t) 
&=
	-\D (f_\delta\xi_R)(x,t-\delta)
	+
	\s (\xi_R g_\delta)(x,t-\delta)
\\&\qquad-
	\int_{\R^n}\int_\delta^\infty	\overline{\Gamma_{(x,t)}^{A^*}}\,
	A\nabla u\cdot \nabla \varphi
	+
	\int_{\R^n}\int_\delta^\infty
	u
	\nabla \varphi\cdot \overline{A^*\nabla \Gamma_{(x,t)}^{A^*}}
\\&\qquad
	+	\lim_{\beta\to 0^+}
	\int_{\R^n}\int_\delta^\infty A\nabla U_{\beta,\gamma}\cdot \nabla ((1-\kappa_\delta)\varphi\overline{ \Gamma_{(x,t)}^{A^*}})
.\end{align*}
We now deal with the final term. Essentially, we wish to exploit the fact that $\Div A\nabla u=0$ and thus $\int A\nabla u\cdot \nabla\Psi=0$ for any smooth, compactly supported function~$\Psi$.

Recall that we have a second parameter, $\gamma$, at our disposal, and that $U_{\beta,\gamma}$ is constant in $B=B((x,t),\gamma/2)$. We define a smooth approximation to $\overline{\Gamma_{(x,t)}^{A^*}}$, analogous to~$U_{\beta,\gamma}$. Specifically, let $H_\beta=\overline{\Gamma_{(x,t)}^{A^*}}*\theta_\beta$, and let $G_{\beta,\gamma}=(1-\widetilde\psi_\gamma) H_\beta + \widetilde\psi_\gamma \overline{\Gamma_{(x,t)}^{A^*}(y_0,s_0)}$ for some $(y_0,s_0)$ with $\abs{(x,t)-(y_0,s_0)}=\gamma/4$ and some $\widetilde\psi$ supported in $B((x,t),\gamma/2)$. Let $\widetilde B=B((x,t),\gamma)$ and let $B=B((x,t),\gamma/2)$. Notice that $\varphi$, $1-\kappa_\delta$ are identically $1$ in~$\widetilde B$.

Then 
\begin{multline*}
\lim_{\beta \to 0^+}
\int_{\R^n}\int_\delta^\infty A\nabla U_{\beta,\gamma}\cdot \nabla ((1-\kappa_\delta)\varphi\overline{ \Gamma_{(x,t)}^{A^*}})
\\\begin{aligned}
&=
	\lim_{\beta \to 0^+}
	\int_{\R^{n+1}_+\setminus \widetilde B} A\nabla U_{\beta,\gamma}\cdot \overline{ \nabla ((1-\kappa_\delta)\varphi\Gamma_{(x,t)}^{A^*})}
+
	\int_{\widetilde B\setminus B}A\nabla U_{\beta,\gamma}\cdot  \nabla \overline{\Gamma_{(x,t)}^{A^*}}
.\end{aligned}\end{multline*}
Let $S=\supp (1-\kappa_\delta)\varphi\setminus \widetilde B$. Notice that $\nabla U_{\beta,\gamma}\to \nabla u$ and $\nabla G_{\beta,\gamma}\to \overline{\Gamma^{A^*}_{(x,t)}}$ as $\beta\to 0$ in $L^2(S)$; thus
\begin{multline*}
	\lim_{\beta \to 0^+}
	\int_{\R^{n+1}_+\setminus \widetilde B} A\nabla U_{\beta,\gamma}\cdot \overline{ \nabla ((1-\kappa_\delta)\varphi\Gamma_{(x,t)}^{A^*})}
\\\begin{aligned}
&=
	\lim_{\beta \to 0^+}
	\int_{\R^{n+1}_+\setminus\widetilde B} A\nabla u\cdot { \nabla ((1-\kappa_\delta)\varphi G_{\beta,\gamma})}
.\end{aligned}\end{multline*}
Therefore, 
\begin{multline*}
	\lim_{\beta \to 0^+}
	\int_{\R^{n+1}_+\setminus \widetilde B} A\nabla U_{\beta,\gamma}\cdot \overline{ \nabla ((1-\kappa_\delta)\varphi\Gamma_{(x,t)}^{A^*})}
\\\begin{aligned}
&=
	\lim_{\beta \to 0^+}
	\int_{\R^{n+1}_+} A\nabla u\cdot { \nabla ((1-\kappa_\delta)\varphi G_{\beta,\gamma})}
-
	\lim_{\beta \to 0^+}
	\int_{\widetilde B} A\nabla u\cdot { \nabla  G_{\beta,\gamma}}
\end{aligned}\end{multline*}
and because $\Div A\nabla u=0$ in $\R^{n+1}_+$ the first integral is zero for all~$\beta$. Thus,
\begin{multline*}
\lim_{\beta \to 0^+}
\int_{\R^n}\int_\delta^\infty A\nabla U_{\beta,\gamma}\cdot \nabla ((1-\kappa_\delta)\varphi\overline{ \Gamma_{(x,t)}^{A^*}})
\\\begin{aligned}
&=
	\lim_{\beta \to 0^+}
	\int_{\widetilde B\setminus B} A\nabla U_{\beta,\gamma}\cdot \overline{ \nabla \Gamma_{(x,t)}^{A^*}}
-
	\lim_{\beta \to 0^+}
	\int_{\widetilde B} A\nabla u\cdot { \nabla G_{\beta,\gamma}}
.\end{aligned}\end{multline*}
Notice that
\begin{equation*}\nabla U_{\beta,\gamma}=(1-\psi_\gamma) \nabla V_\beta + (u(x,t)-V_\beta)\nabla\psi_\gamma.\end{equation*}
If $\gamma$ is sufficiently small (compared to~$t$) and if $\beta$ is sufficiently small (compared to~$\gamma$), then
\begin{equation*}\doublebar{\nabla U_{\beta,\gamma}}_{L^2(\widetilde B)} 
\leq C \doublebar{\nabla u}_{L^2(B((x,t),2\gamma)} + C\gamma^{n/2-1/2} \sup_{(y,s)\in\widetilde B} \abs{u(x,t)-V_\beta(y,s)}. \end{equation*}
We bound the first term using the Caccioppoli inequality and the second term using the definition of~$V_\beta$; thus
\begin{equation*}\doublebar{\nabla U_{\beta,\gamma}}_{L^2(\widetilde B)}  \leq C\gamma^{n/2-1/2} \sup_{(y,s)\in B((x,t),3\gamma)} \abs{u(x,t)-u(y,s)}.\end{equation*}
So by the De Giorgi-Nash-Moser condition and the bound~\eqref{eqn:FS:far:space},
\begin{equation*}\int_{\widetilde B\setminus B} \abs{\nabla U_{\beta,\gamma}}\,\abs{\nabla \Gamma_{(x,t)}^{A^*}}
\leq C\biggl(\frac{\gamma}{t}\biggr)^{\alpha} \biggl(\fint_{B((x,t),t/2)}\abs{u}^2\biggr)^{1/2}
.\end{equation*}

Similarly, by the bounds \eqref{eqn:FS:far:space} and~\eqref{eqn:FS:holder}, $\doublebar{\nabla G_{\beta,\gamma}}_{L^2(\widetilde B)}\leq C\gamma^{1/2-n/2}$, and by the Caccioppoli inequality and the De Giorgi-Nash-Moser condition,
\begin{equation*}\int_{\widetilde B} \abs{\nabla u}\,\abs{ \nabla G_{\beta,\gamma}}
\leq C \biggl(\frac{\gamma}{t}\biggr)^{\alpha} \biggl(\fint_{B((x,t),t/2)}\abs{u}^2\biggr)^{1/2}.\end{equation*}
Thus, taking the limit as $\beta\to 0$ and then taking the limit as $\gamma\to 0$, we have that
\begin{align*}
u(x,t) 
&=
	-\D (f_\delta\xi_R)(x,t-\delta)
	+
	\s (\xi_R g_\delta)(x,t-\delta)
\\&\qquad-
	\int_{\R^n}\int_\delta^\infty	\overline{\Gamma_{(x,t)}^{A^*}}\,
	A\nabla u\cdot \nabla \varphi
	+
	\int_{\R^n}\int_\delta^\infty
	u
	\nabla \varphi\cdot \overline{A^*\nabla \Gamma_{(x,t)}^{A^*}}
.\end{align*}
We now take the limits as $R\to\infty$ and $\delta\to 0^+$.

We remark that if $c_{R,\delta}$ is a constant, we may consider $u-c_{R,\delta}$ instead of~$u$; thus
\begin{align*}
u(x,t)-c_{R,\delta}
&=
	\int_{\R^n}\int_\delta^\infty (u-c_{R,\delta})\,\nabla \varphi\cdot \overline{A^*\nabla \Gamma_{(x,t)}^{A^*}}
	-
	\int_{\R^n}\int_\delta^\infty A\nabla u\cdot \nabla \varphi \,\overline{\Gamma_{(x,t)}^{A^*}}
	\\&\qquad
	+
	\s (\xi_R g_\delta)(x,t-\delta)
	-\D (\xi_R (f_\delta-c_{R,\delta}))(x,t-\delta)
\\&=
	I_{R,\delta}(x,t)
	-
	II_{R,\delta}(x,t)
	\\&\qquad
	+
	\s (\xi_R g_\delta)(x,t-\delta)
	-\D (\xi_R (f_\delta-c_{R,\delta}))(x,t-\delta)
.\end{align*}

Let $\Omega_{R,\delta}=\supp\nabla\varphi\cap\{(x,t):t>\delta\}$.
Observe that if $(x',t')\in B((x,t),t/4)$, then 
\begin{align*}
\abs{II_{R,\delta}(x,t)-II_{R,\delta}(x',t')}
&\leq
	C\int_{\Omega_{R,\delta}}
		\abs{\Gamma_{(x,t)}^{A^*}-\Gamma_{(x',t')}^{A^*}}
		\abs{\nabla u}
		\abs{\nabla\varphi}
.\end{align*}
We may bound $\abs{\Gamma_{(x,t)}^{A^*}-\Gamma_{(x',t')}^{A^*}}$ in $\Omega_{R,\delta}$ using the bound \eqref{eqn:FS:far:space}, the Poincar\'e inequality and the De Giorgi-Nash-Moser condition. We may bound $\abs{\nabla u}$ using 
Theorem~\ref{thm:whitney:embedding}. Then
\begin{align*}
\abs{II_{R,\delta}(x,t)-II_{R,\delta}(x',t')}
&\leq
	\frac{Ct^\alpha}{R^{n+\alpha}}
	\int_{\Omega_{R,\delta}} 
		\abs{\nabla u}
\leq
	\frac{Ct^\alpha}{R^{\alpha-\theta+n/p}}
	\doublebar{u}_{\dot W(p,\theta,2)}
.\end{align*}
So by assumption on $\theta$, $p$, we have that $\lim_{R\to\infty} {II_{R,\delta}(x,t)-II_{R,\delta}(x',t')}=0$, and this limit is uniform in~$\delta$.

Now,
\begin{align*}
\abs{I_{R,\delta}(x,t)-I_{R,\delta}(x',t')}
&\leq \frac{C}{R}
	\int_{\Omega_{R,\delta}}
		\abs{u-c_{R,\delta}}
		\, \abs{\nabla\Gamma_{(x,t)}^{A^*}-\nabla\Gamma_{(x',t')}^{A^*}} 
\\&\leq 
	\frac{C}{R}
	\biggl(\int_{\Omega_{R,\delta}}
		\abs{u-c_{R,\delta}}^2\biggr)^{1/2}
	\biggl(\int_{\Omega_{R,\delta}}
		\abs{\nabla(\Gamma_{(x,t)}^{A^*} -\Gamma_{(x',t')}^{A^*})}^2\biggr)^{1/2}
\end{align*}
and by the bound \eqref{eqn:FS:holder} and the Caccioppoli inequality,		
\begin{align*}
\abs{I_{R,\delta}(x,t)-I_{R,\delta}(x',t')}
&\leq 
	CR^{-n/2-1/2-\alpha} t^\alpha
	\biggl(\int_{\Omega_{R,\delta}}
		\abs{u-c_{R,\delta}}^2\biggr)^{1/2}
.\end{align*}
Recall that by Theorem~\ref{thm:whitney:embedding}, we may assume that $p>1$.
We will later need $c_{R,\delta}=\fint_{\Delta(0,2R)\setminus \Delta(0,R)} f_\delta$. Recall that by the bound~\eqref{eqn:besov-poincare},
\begin{equation*}
\doublebar[big]{f_\delta-\textstyle\fint_{\Delta(0,2R)} f_\delta}_{L^p(\Delta(0,2R))} 
\leq CR^\theta \doublebar{f_\delta}_{\dot B^{p,p}_\theta(\R^n)}
\end{equation*}
and it is straightforward to establish that
\begin{equation*}
\doublebar{f_\delta-c_{R,\delta}}_{L^p(\Delta(0,2R))} 
\leq CR^\theta \doublebar{f_\delta}_{\dot B^{p,p}_\theta(\R^n)}
\end{equation*}
as well; by Theorem~\ref{thm:interior-trace}, we have that
\begin{equation*}
\doublebar{f_\delta-c_{R,\delta}}_{L^p(\Delta(0,2R))} 
\leq CR^\theta \doublebar{u}_{\dot W(p,\theta,2)}.
\end{equation*}

If $p\geq 2$, 
then by H\"older's inequality, we have that 
\begin{align*}
\abs{I_{R,\delta}(x,t)-I_{R,\delta}(x',t')}^p
&\leq 
	\frac{Ct^{p\alpha}}{R^{1+n+p\alpha}} 
	\int_{\Omega_{R,\delta}}
		\abs{u(y,s)-c_{R,\delta}}^p\,dy\,ds
\\&\leq 
	\frac{Ct^{p\alpha}}{R^{1+n+p\alpha}} 
	\int_\delta^{2R} 	\doublebar{u(\,\cdot\,,s)-f_\delta}_{L^p(\Delta(0,2R))}^p \,ds
	\\&\qquad
	+\frac{Ct^{p\alpha}}{R^{1+n+p\alpha}} 
	\int_\delta^{2R} \doublebar{f_\delta-c_{R,\delta}}_{L^p(\Delta(0,2R))}^p \,ds
.\end{align*}
By the bound~\eqref{eqn:besov-trace-convergence}, we have that
\begin{align*}
\abs{I_{R,\delta}(x,t)-I_{R,\delta}(x',t')}^p
&\leq 
	\frac{Ct^{p\alpha}}{R^{n+p\alpha-p\theta}} 
	\doublebar{u}_{\dot W(p,\theta,2)}^p
\end{align*}
which again goes to zero as $R\to\infty$, uniformly in~$\delta$.

If $p<2$, choose some $s$ with $\delta<s<2R$. We may cover $\Delta(0,2R)$ by at most $C(R/s)^n$ disjoint cubes $Q_j$ of side-length~$s$. Then by Corollary~\ref{cor:local-bound},
\begin{align*}
\int_{\Delta(0,2R)} \abs{u(y,s)-c_{R,\delta}}^2\,dy
&\leq 
	C s^n\sum_{j}
	\biggl( \fint_{2Q_j} \fint_{s/2}^{3s/2} \abs{u(y,r)-c_{R,\delta}}^p\,dr\,dy\biggr)^{2/p}
\\&\leq 
	\frac{Cs^n}{s^{2n/p+2/p}}
	\biggl( \int_{\Delta(0,4R)} \int_{s/2}^{3s/2} \abs{u(y,r)-c_{R,\delta}}^p\,dr\,dy
	\biggr)^{2/p}
.\end{align*}
Again by formulas~\eqref{eqn:besov-trace-convergence} and~\eqref{eqn:besov-poincare}, we have that
\begin{align*}
\int_{\Delta(0,2R)} \abs{u(y,s)-c_{R,\delta}}^2\,dy
&\leq 
	C R^{2\theta}s^{n-2n/p} \doublebar{u}_{\dot W(p,\theta,2)}^2 	
.\end{align*}
Therefore,
\begin{align*}
\abs{I_{R,\delta}(x,t)-I_{R,\delta}(x',t')}
&\leq  
	CR^{-n/2-1/2-\alpha} t^\alpha
	\biggl(\int_{\Omega_{R,\delta}}
		\abs{u(y,s)-c_{R,\delta}}^2\,dy\,ds\biggr)^{1/2}
\\&\leq  
	C t^\alpha
	\biggl(R^{-n-1-2\alpha+2\theta}\int_{\delta}^{2R} s^{n-2n/p}
		\,ds\biggr)^{1/2}\doublebar{u}_{\dot W(p,\theta,2)}
.\end{align*}
If $n-2n/p>-1$, then the term inside the parentheses evaluates to $CR^{-2(\alpha-\theta+n/p)}$, which goes to zero as $R\to\infty$, uniformly in~$\delta$. If $n-2n/p<-1$ or $n-2n/p=-1$, then this term is $CR^{-n-1-2\alpha+2\theta}\delta^{n-2n/p+1}$ or $R^{-n-1-2\alpha+2\theta}(\ln R-\ln \delta)$, respectively; in each case this term goes to zero as $R\to \infty$ (although in this case the convergence is \emph{not} uniform in~$\delta$).

Thus, whether $p<2$ or $p\geq 2$ we have that $\lim_{R\to\infty}{I_{R,\delta}(x,t)-I_{R,\delta}(x',t')}=0$.

We have now shown that for all $(x',t')\in B((x,t),t/4)$,
\begin{multline}
\label{eqn:representation:1}
u(x,t+\delta)-u(x',t'+\delta)
\\\begin{aligned}
 =\lim_{R\to\infty}\Bigl(& 
	 \s (g_\delta\xi_R)(x,t)
	- \s (g_\delta\xi_R)(x',t')
\\ &
	-\D ((f_\delta-c_{R,\delta})\xi_R)(x,t)
	+\D ((f_\delta-c_{R,\delta})\xi_R)(x',t') 
	\Bigr)
.\end{aligned}\end{multline}
Furthermore, if $p\geq 2$ then the limit is uniform in~$\delta$. We want to evaluate the limit as $\delta\to 0^+$.

We begin by considering the case where $p<\infty$.
By Theorem~\ref{thm:whitney:embedding} we may assume that $p>1$ and so Lemma~\ref{lem:besov:localize} and Theorem~\ref{thm:trace:converge} are valid.
By Lemma~\ref{lem:besov:localize}, we have that $(f_\delta-c_{R,\delta})\xi_R\to f_\delta$ in $\dot B^{p,p}_\theta(\R^n)$.

By Theorem~\ref{thm:bounded} and the De Giorgi-Nash-Moser condition, if $(x,t)$ and $(x',t')\in\R^{n+1}_+$ with $\abs{(x,t)-(x',t')}<t/4$, then 
\begin{align*}\abs{\D h(x,t)-\D h(x',t')} 
&\leq Ct^{\theta-n/p}\doublebar{h}_{\dot B^{p,p}_\theta(\R^n)}
\end{align*}
and so
\begin{equation*}\lim_{R\to\infty}\Bigl(
\D ((f_\delta-c_{R,\delta})\xi_R)(x',t') 
-\D ((f_\delta-c_{R,\delta})\xi_R)(x,t)
\Bigr)
=\D f_\delta(x',t')-\D f_\delta(x,t).\end{equation*}

Recall that we assumed that $u\in \dot W(p,\theta,2)$ for some $\theta<1$, and so $g_\delta$ lies in a negative smoothness space. Notice that 
\begin{multline}
\label{eqn:representation:2}
{\s (g_\delta(1-\xi_R))(x,t)
	- \s (g_\delta(1-\xi_R))(x',t')}
\\\begin{aligned}
&=
	\int_{\R^n} (1-\xi_R(y))
	\overline{(\Gamma_{(x,t)}^{A^*}(y,0) -\Gamma_{(x',t')}^{A^*}(y,0))}
	g_\delta(y)\,dy
.\end{aligned}\end{multline}
We intend to apply Lemma~\ref{lem:besov:localize} not to $g_\delta$ but to $\gamma(y)=\Gamma_{(x,t)}^{A^*}(y,0) -\Gamma_{(x',t')}^{A^*}(y,0)$. First, as before we have that
\begin{equation*}
\abs{\s h(x,t)-\s h(x',t')} 
\leq Ct^{\theta-n/p}\doublebar{h}_{\dot B^{p,p}_{\theta-1}(\R^n)}
\end{equation*}
for any $h\in {\dot B^{p,p}_{\theta-1}(\R^n)}$, and because $\s h(x,t)-\s h(x',t')=\langle \gamma,h\rangle$, this implies that
\begin{equation*}
\doublebar{\gamma}_{\dot B^{p',p'}_{1-\theta}(\R^n)}=
\doublebar{\Gamma_{(x,t)}^{A^*}(\,\cdot\,,0) -\Gamma_{(x',t')}^{A^*}(\,\cdot\,,0)}_{\dot B^{p',p'}_{1-\theta}(\R^n)}
\leq C t^{\theta-n/p}.
\end{equation*}
By Lemma~\ref{lem:besov:localize}, 
\begin{equation*}\lim_{R\to\infty} \doublebar{(\gamma-\widetilde\gamma_R)(1-\xi_R)}_{\dot B^{p',p'}_{1-\theta}(\R^n)} = 0\end{equation*}
where 
\begin{equation*}\widetilde\gamma_R=\fint_{\Delta(0,2R)\setminus\Delta(0,R)}
{\Gamma_{(x,t)}^{A^*}(y,0) -\Gamma_{(x',t')}^{A^*}(y,0)}\,dy.\end{equation*}
But by the bound~\eqref{eqn:FS:holder}, we have that 
$\abs{\widetilde\gamma_R}\leq CR^{1-\alpha-n} t^\alpha$, and so \begin{equation*}\doublebar{\widetilde\gamma_R(1-\xi_R)}_{\dot B^{p',p'}_{1-\theta}(\R^n)} = \doublebar{\widetilde\gamma_R\xi_R}_{\dot B^{p',p'}_{1-\theta}(\R^n)}
\leq CR^{n/p'-1+\theta}\abs{\widetilde\gamma_R}
\leq CR^{\theta-\alpha-n/p} t^\alpha\end{equation*}
which also approaches zero as $R\to \infty$.
Thus, 
\[\lim_{R\to\infty}\s (g_\delta(1-\xi_R))(x,t)
	- \s (g_\delta(1-\xi_R))(x',t')
=
	\lim_{R\to\infty} \langle (1-\xi_R)\gamma,g_\delta\rangle
	=0
\]
and so
\begin{align*}
u(x,t+\delta)-u(x',t'+\delta)
&=
	\D f_\delta(x',t') 
	-\D f_\delta(x,t)
	- \s g_\delta(x',t')
	+ \s g_\delta(x,t)
.\end{align*}
Taking the limit as $\delta\to 0$ and applying Theorems~\ref{thm:trace:converge} and~\ref{thm:bounded} completes the proof.

Now, consider the case $p=\infty$. Because  the limit \eqref{eqn:representation:1} is uniform in $\delta$ we may take the limit as $\delta\to 0$ before taking the limit as $R\to\infty$. We first consider the terms involving~$\D$.

By Corollary~\ref{cor:holder:weighted},
\begin{equation*}\int_{\Delta(0,2R)}\abs{f_\delta-f}^2
\leq
C R^{n}\delta^{2\theta} \doublebar{u}_{\dot W(\infty,\theta,2)}^2\end{equation*}
and so by the bound~\eqref{eqn:pointwise}, \begin{multline*}\lim_{\delta\to 0}\Bigl(
\D ((f_\delta-c_{R,\delta})\xi_R)(x',t') 
-\D ((f_\delta-c_{R,\delta})\xi_R)(x,t)
\Bigr)
\\
=\D ((f-c_{R})\xi_R)(x',t') 
-\D ((f-c_{R})\xi_R)(x,t)
.\end{multline*}

To deal with the terms involving~$\s$, recall that by Theorem~\ref{thm:trace:converge}, we have that $g_\delta\rightharpoonup g$ in the weak-$*$ topology. Notice that formula~\eqref{eqn:representation:2} is still valid. Observe that  by the bound~\eqref{eqn:FS:holder} and Lemma~\ref{lem:slabs}, the function $\gamma(y)=\xi_R(y)(\Gamma_{(x,t)}^{A^*}(y,0) -\Gamma_{(x',t')}^{A^*}(y,0))$ lies in $L^{1}(\R^n)\cap \dot W^{1}_1(\R^n)\subset\dot B^{1,1}_{1-\theta}(\R^n)$. So
\begin{equation*}
\lim_{\delta\to 0}\Bigl(
\s (g_\delta\xi_R)(x',t') 
-\s (g_\delta\xi_R)(x,t)
\Bigr)
=\s (g\xi_R)(x',t') 
-\s (g\xi_R)(x,t)
.\end{equation*}

We have now established that
\begin{align*}
u(x,t)-u(x',t')
& =\lim_{R\to\infty}\Bigl(
	\s (g\xi_R)(x,t)
	-\s (g\xi_R)(x',t') 
\\ &\qquad
	+\D ((f-c_{R})\xi_R)(x',t') 
	-\D ((f-c_{R})\xi_R)(x,t)
	\Bigr)
.\end{align*}

Observe that by Lemma~\ref{lem:slabs} and the bound \eqref{eqn:FS:far:space}, $\Gamma_{(x,t)}^{A^*}-\Gamma_{(x',t')}^{A^*}$ is a (constant multiple of a) $\dot B^{1,1}_{1-\theta}(\R^n)$ molecule in the sense of Lemma~\ref{lem:molecule:plus}. Thus, we may apply Lemma~\ref{lem:besov:localize} to $\Gamma_{(x,t)}^{A^*}-\Gamma_{(x',t')}^{A^*}$, as before, to see that
\begin{equation*}\lim_{R\to\infty}\bigl(
	\s (g\xi_R)(x,t)
	-\s (g\xi_R)(x',t')\bigr)=\s g(x,t)
	-\s g(x',t').\end{equation*}
To prove an equivalent result for~$\D$,
let $A_0=\Delta(x,2R)$, and for each $j\geq 1$ let $A_j=\Delta(x,2^{j+1}R)\setminus \Delta(x,2^{j}R)$. Observe that 
\begin{equation*}\doublebar{f-c_R}_{L^\infty(A_j)}
\leq 
C 2^{j\theta}R^\theta\doublebar{f}_{\dot B^{\infty,\infty}_\theta(\R^n)}\end{equation*}
and that if $j\geq 1$, then by Lemma~\ref{lem:slabs} and the bound~\eqref{eqn:FS:holder},
\begin{equation*}\doublebar{\nabla\Gamma_{(x,t)}^{A^*}(\,\cdot\,,0)-\nabla\Gamma_{(x',t')}^{A^*}(\,\cdot\,,0)}_{L^{1}(A_j)}\leq C2^{-j\alpha} \biggl(\frac{t}{R}\biggr)^\alpha.\end{equation*}
Thus, recalling the definition \eqref{eqn:D} of the double layer potential, we have that
\begin{multline*}
\abs{\D ((f-c_R)(1- \xi_R))(x,t)-\D ((f-c_R)(1- \xi_R))(x',t')}
\\\leq
	\sum_{j=1}^\infty 
	C 2^{-j(\alpha-\theta)}R^{\theta-\alpha} t^\alpha
	\doublebar{f}_{\dot B^{\infty,\infty}_\theta(\R^n)}
.\end{multline*}
Because $\theta<\alpha$ the geometric series converges, and so we still have that
\begin{equation*}\lim_{R\to\infty}\Bigl(
\D ((f-c_R)\xi_R)(x',t') 
-\D ((f-c_R)\xi_R)(x,t)
\Bigr)
=\D f(x',t')-\D f(x,t)\end{equation*}
as desired.
\end{proof}

%% file: sec-9-invertible.tex
\chapter[Invertibility of Layer Potentials and Well-Posedness]{Invertibility of Layer Potentials and Well-Posedness of Boundary-Value Problems}
\label{chap:invertibility}

In this chapter we will prove our remaining main theorems, that is,  Theorems~\ref{thm:well-posed:invertible}, \ref{thm:well-posed:invertible:1}, and~\ref{thm:well-posed:compatible}, and Corollaries~\ref{cor:well-posed:open}, \ref{cor:well-posed:1}, \ref{cor:well-posed:duality}, \ref{cor:well-posed:AusM} and~\ref{cor:well-posed:real}. That is, we will show that well-posedness of boundary-value problems is equivalent to invertibility of layer potentials, and will explore some consequences of this equivalence, with particular attention to the case of real coefficients.

We remark that a necessary condition for most of our arguments is \emph{boundedness} of layer potentials; we will also need the Green's formula representation of Theorem~\ref{thm:green}.
Thus, throughout this chapter, we will let $A$ be an elliptic, $t$-independent matrix such that both $A$ and $A^*$ satisfy the De Giorgi-Nash-Moser condition. We will let $p$, $\theta$ be numbers such that either $0<\theta<1$ and $p$, $\theta$ satisfy the conditions of Theorem~\ref{thm:bounded}, or $\theta=1$ and $p$ is such that the bounds \eqref{eqn:S:NTM} and~\eqref{eqn:D:NTM:1} are valid, and furthermore such that $p$, $\theta=1$ satisfy the conditions of Theorem~\ref{thm:green}. Note that there is some $\varepsilon>0$ such that, if $1/2-\varepsilon<1/p<1+\alpha/n$, then $p$ satisfies the conditions specified in the case $\theta=1$.

\section[Invertibility and well-posedness]{Invertibility and well-posedness: Theorems~\ref{thm:well-posed:invertible}, \ref{thm:well-posed:invertible:1} and~\ref{thm:well-posed:compatible}}
\label{sec:invertible:well-posed}

In this section we will show that invertibility of layer potentials is equivalent to well-posedness of certain boundary-value problems.

We remind the reader of the classic method of layer potentials. By the bound \eqref{eqn:trace:S}, $\s^A_+$ is bounded $\dot B^{p,p}_{\theta-1}(\R^n)\mapsto \dot B^{p,p}_\theta(\R^n)$ for all $p$, $\theta$ satisfying the conditions of Theorem~\ref{thm:bounded}. 
Suppose that $\s^A_+$ is surjective $\dot B^{p,p}_{\theta-1}(\R^n)\mapsto \dot B^{p,p}_\theta(\R^n)$. By formula~\eqref{eqn:Div:plus}, the definition of~$\s^A_+$, and the bound~\eqref{eqn:space:S}, if $\s^A_+ g=f$ then $u=\s^A g$ is a solution to $(D)^A_{p,\theta}$ with boundary data~$f$, and so surjectivity of $\s^A_+$ implies the existence property. (If $\s^A_+$ is invertible with bounded inverse then we have solvability, that is, the existence property and the estimate $\doublebar{u}_{\dot W(p,\theta,2)}\leq C\doublebar{f}_{\dot B^{p,p}_\theta(\R^n)}$.)

Similarly, if $(\partial_\nu^A\D^A)_+$ is invertible $\dot B^{p,p}_{\theta}(\R^n)\mapsto \dot B^{p,p}_{\theta-1}(\R^n)$, then $(N)^A_{p,\theta}$ is solvable; if $\s^A_+:H^p(\R^n)\mapsto \dot H^p_1(\R^n)$ or $(\partial_\nu^A\D^A)_+:\dot H^p_1(\R^n)\mapsto H^p(\R^n)$ is bounded and invertible, then $(D)^A_{p,1}$ or $(N)^A_{p,1}$ is solvable. We remark that many of the theorems of this chapter follow this pattern; that is, the arguments involving $(\partial_\nu^A\D^A)_+$ and the Neumann problem, or for the case $\theta=1$, closely parallel those for $\s^A_+$ and the Dirichlet problem in the case $0<\theta<1$. We will generally present complete arguments for $(D)^A_{p,\theta}$, $0<\theta<1$ and leave the details of the corresponding arguments for $(N)^A_{p,\theta}$ or $\theta=1$ to the reader.

We have seen that invertibility of $\s^A_+$ and $(\partial_\nu^A\D^A)_+$ implies solvability of boundary-value problems. Given Theorem~\ref{thm:green}, uniqueness may also be reduced to invertibility of layer potentials.
\begin{thm} \label{thm:unique:D:N} 
Let $A$, $p$ and $\theta$ satisfy the conditions specified at the beginning of this chapter. Suppose that $\s^A_+$ is one-to-one $\dot B^{p,p}_{\theta-1}(\R^n)\mapsto \dot B^{p,p}_{\theta}(\R^n)$, or (if $\theta=1$) is bounded and one-to-one $H^p(\R^n)\mapsto \dot H^p_1(\R^n)$. 

Then $(D)^A_{p,\theta}$ has the uniqueness property; that is, if 
$\Div A\nabla u=0$ in $\R^{n+1}_+$, $\Tr u=0$, and either $0<\theta<1$ and $u\in \dot W(p,\theta,2)$ or $\theta=1$ and $\widetilde N_+(\nabla u)\in L^p(\R^n)$, then $u\equiv 0$. 

If instead $\nu\cdot A\nabla u=0$ on $\partial\R^{n+1}_+$ and $(\partial_\nu^A\D^A)_+$ is bounded and one-to-one $\dot B^{p,p}_{\theta}(\R^n)\mapsto \dot B^{p,p}_{\theta-1}(\R^n)$ or $\dot W^p_1(\R^n)\mapsto L^p(\R^n)$, then $(N)^A_{p,\theta}$ has the uniqueness property.
\end{thm}

\begin{proof} Let $\Div A\nabla u=0$ in $\R^{n+1}_+$ and either $u\in \dot W(p,\theta,2)$ or $\widetilde N_+(\nabla u)\in L^p(\R^n)$.
By Theorem~\ref{thm:green}, we have that (up to an additive constant) $u=-\D f+\s g$, where $f=\Tr u$ and $g=\nu\cdot A\nabla u$. By Theorem~\ref{thm:trace} or Theorem~\ref{thm:trace:NTM}, $g\in \dot B^{p,p}_{\theta-1}(\R^n)$ or $g\in H^p(\R^n)$.

If $\Tr u=0$ then $f=0$ and so $u=\s g$. Therefore $\s^A_+ g$ is constant. Recall that the space $\dot B^{p,p}_\theta(\R^n)$ (or $\dot H^p_1(\R^n)$) is only defined modulo constants, and so $\s^A_+ g=0$  in this space. But by injectivity of $\s^A_+$, this implies that $g=0$, and so $u$ is a constant. Because $\Tr u=0$ we must have that this constant is zero.

If $\nu\cdot A\nabla u=0$ on $\partial\R^{n+1}_+$ then the same argument is valid with the roles of $\D$ and $\s$ reversed, except that we cannot conclude that $u\equiv 0$, only that $u$ is constant.
\end{proof}

We have established the following facts: if layer potentials are onto then we have the existence property; if layer potentials are one-to-one then we have the uniqueness property; and if layer potentials are invertible then we have well-posedness. To complete the proofs of  Theorems~\ref{thm:well-posed:invertible} and~\ref{thm:well-posed:invertible:1}, we need only show the converses.

We will use a classic argument of Verchota \cite{Ver84} to establish that layer potentials have bounded inverses; surjectivity comes from an argument in \cite{BarM13B}, and injectivity is straightforward. All three arguments involve jump relations, that is, the interaction between the values of $\D f$ and $\s f$ on $\partial\R^{n+1}_+$ and~$\partial\R^{n+1}_-$. Specifically, we will need the following result.

\begin{prp}\label{prp:jump}
Suppose that $A$ and $A^*$ are elliptic, $t$-independent and satisfy the De Giorgi-Nash-Moser condition. 
If $f$ is smooth and compactly supported, then
\begin{align}
\label{eqn:D:jump}
\D^A_+ f-\D^A_- f &= - f
,\\
\label{eqn:S:cts}
\s^A_+ f - \s^A_- f &= 0
,\\
\label{eqn:S:jump}
\nu_+\cdot A\nabla \s^A f\vert_{\partial\R^{n+1}_+} - \nu_+\cdot A\nabla \s^Ag\vert_{\partial\R^{n+1}_-}f&= f
,\\
\label{eqn:D:cts}
\nu_+\cdot A\nabla \D^A f\big\vert_{\partial\R^{n+1}_+}-\nu_+\cdot A\nabla \D^A f\big\vert_{\partial\R^{n+1}_-} &= 0
.\end{align}
If $\varphi$ is also smooth and compactly supported, then
\begin{equation}
\label{eqn:D:conormal:adjoint}
\Bigl\langle\varphi,\nu_+\cdot A\nabla \D^A f\big\vert_{\partial\R^{n+1}_+}\Bigr\rangle
=
-\Bigl\langle \nu_-\cdot A^*\nabla \D^{A^*} \varphi\big\vert_{\partial\R^{n+1}_-}, f\Bigr\rangle
.\end{equation}
\end{prp}

Notice that by our boundedness and trace results (Theorems~\ref{thm:bounded} and~\ref{thm:trace}, and the bounds~\eqref{eqn:S:NTM} and~\eqref{eqn:D:NTM:1} and Theorem~\ref{thm:trace:NTM}), if $p$, $\theta$ satisfy the conditions at the beginning of this chapter and if $p<\infty$, then the formulas~\eqref{eqn:D:jump}, \eqref{eqn:S:cts}, \eqref{eqn:S:jump} and~\eqref{eqn:D:cts} extend to all $f$ in (as appropriate) $\dot B^{p,p}_\theta(\R^n)$, $\dot B^{p,p}_{\theta-1}(\R^n)$, $H^p(\R^n)$ or $\dot H^p_1(\R^n)$. If $p=\infty$, we will see that these formulas may be extended by duality. (See formulas~\eqref{eqn:S:adjoint} and~\eqref{eqn:S-D:adjoint} below as well as~\eqref{eqn:D:conormal:adjoint}.)

\begin{proof}[Proof of Proposition~\ref{prp:jump}]
For $f\in L^2(\R^n)$, the formulas \eqref{eqn:D:jump}, \eqref{eqn:S:cts} and \eqref{eqn:S:jump} were established in \cite[Lemma~4.18]{AlfAAHK11} under certain assumptions on $\s^A$; that these assumptions are always valid was proven in \cite{Ros12A}. See also \cite[Propositions~3.3 and 3.4]{HofMitMor}, where these results are gathered.

We are left with formulas~\eqref{eqn:D:cts} and~\eqref{eqn:D:conormal:adjoint}.
Let $F$ be a smooth, compactly supported extension of~$f$ to $\R^{n+1}$ and let $\Phi$ be a smooth, compactly supported extension of~$\varphi$. 
Observe that if $X\in \R^{n+1}_+$ then 
\[\D^A f(X) = 
-\int_{\R^{n+1}_-}\nabla F(Y)\cdot \overline{A^*(Y)\nabla\Gamma^{A^*}_{X}(Y)} \,dY\]
and so
\begin{multline}
\label{eqn:D:cts:1}
\Bigl\langle\varphi,\nu_+\cdot A\nabla \D^A f\big\vert_{\partial\R^{n+1}_+}\Bigr\rangle
\\\begin{aligned}
&=
	\int_{\R^{n+1}_+} \overline{\nabla\Phi(X)} \cdot A(X)\nabla  \D^A f(X)\,dX
\\&=
	-\int_{\R^{n+1}_+} \overline{\nabla\Phi(X)} \cdot A(X)\nabla  \int_{\R^{n+1}_-} \nabla F(Y)\cdot \overline{A^*(Y)\nabla\Gamma^{A^*}_X(Y)} \,dY\,dX
.\end{aligned}\end{multline}
Similarly,
\begin{multline}
\label{eqn:D:cts:2}
\Bigl\langle\varphi,\nu_-\cdot A\nabla \D^A f\big\vert_{\partial\R^{n+1}_-}\Bigr\rangle
\\\begin{aligned}
&=
	\int_{\R^{n+1}_-} \overline{\nabla\Phi(X)} \cdot A(X)\nabla  \int_{\R^{n+1}_+} \nabla F(Y)\cdot \overline{A^*(Y)\nabla\Gamma^{A^*}_X(Y)} \,dY\,dX
.\end{aligned}\end{multline}

Using the Caccioppoli inequality and the weak definition of the gradient, it is straightforward to establish that $\nabla_{X}\nabla_{Y}\Gamma_X^{A^*}(Y)$ is locally in $L^2(\R^{n+1}\times\R^{n+1}\setminus D)$ where $D$ is the diagonal $\{(X,X):X\in\R^{n+1}\}$. Furthermore, using Lemma~\ref{lem:slabs} and the bound~\eqref{eqn:FS:far:space}, we may compute that if $t\neq s$ and $Q$ is a cube of side-length $\abs{t-s}$, then
\begin{equation*}\int_{Q}\int_{A_j(Q)}\abs{\nabla_{x,t}\nabla_{y,s} \Gamma_{(x,t)}(y,s)}^2\,dx\,dy\leq \frac{C}{2^{j(n+2\alpha)}\abs{t-s}^{2}}\end{equation*}
where $A_0(Q)=Q$ and $A_j(Q)=2^{j}Q\setminus 2^{j-1}Q$ for $j\geq 1$. Applying H\"older's inequality and then summing over~$j$, we have that
\begin{equation*}\int_{Q}\int_{\R^n}\abs{\nabla_{x,t}\nabla_{y,s} \Gamma_{(x,t)}(y,s)}\,dx\,dy\leq \frac{C\abs{Q}}{\abs{t-s}}.\end{equation*}
Thus, if $R\geq \abs{t-s}$, then
\begin{equation*}\int_{\Delta(0,R)}\int_{\R^n}\abs{\nabla_{x,t}\nabla_{y,s} \Gamma_{(x,t)}(y,s)}\,dx\,dy\leq \frac{CR^n}{\abs{t-s}}.\end{equation*}

This implies that if $F$, $\Phi$ are smooth and compactly supported, then the integrals in formulas~\eqref{eqn:D:cts:1} and~\eqref{eqn:D:cts:2} converge absolutely. Therefore, we may interchange the order of integration. Using formulas~\eqref{eqn:FS:switch} and~\eqref{eqn:FS:conjugate}, we may easily derive the adjoint formula~\eqref{eqn:D:conormal:adjoint}.

We are left with the continuity relation~\eqref{eqn:D:cts}.
Let \[u(X)=\int_{\R^{n+1}_-} \nabla F(Y)\cdot \overline{A^*(Y)\nabla\Gamma^{A^*}_X(Y)} \,dY.\]
Recall from Lemma~\ref{lem:FS:lax-milgram} that $\doublebar{\nabla u}_{L^2(\R^{n+1})}\leq C \doublebar{\nabla F}_{L^2(\R^{n+1})}$. So
\begin{multline*}
\Bigl\langle\varphi,\nu_+\cdot A\nabla \D^A f\big\vert_{\partial\R^{n+1}_+}\Bigr\rangle
\\\begin{aligned}
&=
	-\int_{\R^{n+1}_+} \overline{\nabla\Phi(X)} \cdot A(X)\nabla  u(X)\,dX
\\&=
	-\int_{\R^{n+1}} \overline{\nabla\Phi(X)} \cdot A(X)\nabla  u(X)\,dX
	+\int_{\R^{n+1}_-} \overline{\nabla\Phi(X)} \cdot A(X)\nabla  u(X)\,dX
.\end{aligned}\end{multline*}
But $\Div A\nabla u = \Div (\1_{\R^{n+1}_-}A\nabla F)$, and so by the weak definition~\eqref{eqn:solution} of $\Div A\nabla u$,
\begin{multline*}
\Bigl\langle\varphi,\nu_+\cdot A\nabla \D^A f\big\vert_{\partial\R^{n+1}_+}\Bigr\rangle
\\\begin{aligned}
&=
	-\int_{\R^{n+1}_-} \overline{\nabla\Phi(X)} \cdot A(X)\nabla  F(X)\,dX
	\\&\qquad
	+\int_{\R^{n+1}_-} \overline{\nabla\Phi(X)} \cdot A(X)\nabla \int_{\R^{n+1}_-} \nabla F(Y)\cdot \overline{A^*(Y)\nabla\Gamma^{A^*}_X(Y)} \,dY\,dX
.\end{aligned}\end{multline*}
By formula~\eqref{eqn:FS:point}, $F(X)=\int_{\R^{n+1}} \nabla F(Y)\cdot \overline{A^*(Y)\nabla\Gamma^{A^*}_X(Y)} \,dY$ and so
\begin{multline*}
\Bigl\langle\varphi,\nu_+\cdot A\nabla \D^A f\big\vert_{\partial\R^{n+1}_+}\Bigr\rangle
\\\begin{aligned}
&=
	-\int_{\R^{n+1}_-} \overline{\nabla\Phi(X)} \cdot A(X)\nabla  \int_{\R^{n+1}} \nabla F(Y)\cdot \overline{A^*(Y)\nabla\Gamma^{A^*}_X(Y)} \,dY\,dX
	\\&\qquad
	+\int_{\R^{n+1}_-} \overline{\nabla\Phi(X)} \cdot A(X)\nabla \int_{\R^{n+1}_-} \nabla F(Y)\cdot \overline{A^*(Y)\nabla\Gamma^{A^*}_X(Y)} \,dY\,dX
\\&=
	-\int_{\R^{n+1}_-} \overline{\nabla\Phi(X)} \cdot A(X)\nabla  \int_{\R^{n+1}_+} \nabla F(Y)\cdot \overline{A^*(Y)\nabla\Gamma^{A^*}_X(Y)} \,dY\,dX
.\end{aligned}\end{multline*}
But the right-hand side is precisely equal to our value of  $-\Bigl\langle \varphi, \nu_-\cdot A\nabla \D^A f\big\vert_{\partial\R^{n+1}_-}\Bigr\rangle $ from above, as desired.
\end{proof}

We remark that it is straightforward to establish the following adjoint formulas for the operators~$\D^A_\pm$, $\s^A_\pm$ and $(\partial_\nu^A\s^A)_\pm$:
\begin{align}
\label{eqn:S:adjoint}
\langle \varphi, \s^A_+ f\rangle &= \langle \s^{A^*}_- \varphi,f\rangle
,\\
\label{eqn:S-D:adjoint}
\langle \varphi, \D^A_+ f\rangle &= -\langle (\partial_\nu^{A^*}\s^{A^*})_-\varphi,f\rangle
,&
\langle \varphi, \D^A_- f\rangle &= \langle (\partial_\nu^{A^*}\s^{A^*})_+\varphi,f\rangle
.\end{align}

We now use the above jump relations to complete the proofs of Theorems~\ref{thm:well-posed:invertible} and~\ref{thm:well-posed:invertible:1}. We begin with the converse to Theorem~\ref{thm:unique:D:N}.
\begin{thm} \label{thm:unique:converse}
Let $A$, $p$ and $\theta$ satisfy the conditions specified at the beginning of this chapter. Suppose that $(D)^A_{p,\theta}$ has the uniqueness property in both $\R^{n+1}_+$ and $\R^{n+1}_-$.

Then $\s^A_+$ is one-to-one $\dot B^{p,p}_{\theta-1}(\R^n)\mapsto \dot B^{p,p}_{\theta}(\R^n)$ or (if $\theta=1$) $H^p(\R^n)\mapsto \dot H^p_1(\R^n)$.

Similarly, if $(N)^A_{p,\theta}$ has the uniqueness property in both $\R^{n+1}_+$ and $\R^{n+1}_-$, then $(\partial_\nu^A\D^A)_+$ is one-to-one $\dot B^{p,p}_{\theta}(\R^n)\mapsto \dot B^{p,p}_{\theta-1}(\R^n)$ or $\dot H^p_1(\R^n)\mapsto H^p(\R^n)$.
\end{thm}

\begin{proof}
Suppose that $(D)^A_{p,\theta}$ has the uniqueness property.
Let $g$ satisfy $\s^A_+ g=0$. By the bound \eqref{eqn:space:S} or~\eqref{eqn:S:NTM} and the uniqueness property, $\s^A g\equiv 0$ in $\R^{n+1}_+$ and $\R^{n+1}_-$. But then by the jump relation~\eqref{eqn:S:jump}, $g\equiv 0$, as desired. 

A similar argument concerning $\D^A g$ is valid if solutions to $(N)^A_{p,\theta}$ has the uniqueness property; we need only replace formulas~\eqref{eqn:space:S}, \eqref{eqn:S:NTM} and~\eqref{eqn:S:jump} by formulas~\eqref{eqn:space:D}, \eqref{eqn:D:NTM:1} and~\eqref{eqn:D:jump}.
\end{proof}

We now consider surjectivity. We use an argument from \cite{BarM13B}.

\begin{thm}\label{thm:onto}
Let $A$, $p$ and $\theta$ satisfy the conditions specified at the beginning of this chapter.

If $(D)^A_{p,\theta}$ has the existence property in $\RR^{n+1}_+$ and $\RR^{n+1}_-$, then the operator $\s^A_+$ is surjective $\dot B^{p,p}_{\theta-1}(\RR^n)\mapsto  \dot B^{p,p}_\theta(\RR^n)$ or (if $\theta=1$) $H^p(\R^n)\mapsto \dot H^p_1(\R^n)$). If $(D)^A_{p,\theta}$ is compatibly solvable then $\s^A_+$ is also surjective $\dot B^{p,p}_{\theta-1}(\RR^n)\cap \dot B^{2,2}_{-1/2}(\R^n)\mapsto  \dot B^{p,p}_\theta(\RR^n)\cap \dot B^{2,2}_{1/2}(\R^n)$ or $H^p(\RR^n)\cap \dot B^{2,2}_{-1/2}(\R^n)\mapsto  \dot H^p_1(\RR^n)\cap \dot B^{2,2}_{1/2}(\R^n)$.

Similarly, if $(N)^A_{p,\theta}$ has the existence property in $\RR^{n+1}_+$ and~$\RR^{n+1}_-$, then the operator $(\partial_\nu^A \D^A)_+$ is surjective $\dot B^{p,p}_{\theta}(\RR^n)\mapsto  \dot B^{p,p}_{\theta-1}(\RR^n)$ or $\dot H^p_1(\R^n)\mapsto H^p(\R^n)$, and if $(N)^A_{p,\theta}$ is compatibly solvable then $(\partial_\nu^A \D^A)_+$ is surjective $\dot B^{p,p}_{\theta}(\RR^n)\cap \dot B^{2,2}_{1/2}(\R^n)\mapsto  \dot B^{p,p}_{\theta-1}(\RR^n)\cap \dot B^{2,2}_{-1/2}(\R^n)$ or $\dot H^p_1(\RR^n)\cap \dot B^{2,2}_{-1/2}(\R^n)\mapsto  H^p(\RR^n)\cap \dot B^{2,2}_{1/2}(\R^n)$.
\end{thm}

\begin{proof}
Suppose that $(D)^A_{p,\theta}$ has the existence property for some $0<\theta<1$.
Choose some $f\in  \dot B^{p,p}_\theta(\R^n)$. Let $u_\pm\in \dot W_\pm(p,\theta,2)$ be the solutions to $(D)^A_{p,\theta}$ with boundary data~$f$ in $\RR^{n+1}_\pm$. Here $\dot W_+(p,\theta,q)=\dot W(p,\theta,q)$, and $\dot W_-(p,\theta,q)$ is the space of functions $u$ defined on $\R^{n+1}_-$ such that $u(x,t)=v(x,-t)$ for some $v\in\dot W(p,\theta,q)$.

By Theorem~\ref{thm:trace}, we have that the conormal derivatives $g_\pm=\nu_+ \cdot A\nabla u\vert_{\partial\R^{n+1}_\pm}$ exist in the weak sense and lie in $\dot B^{p,p}_{\theta-1}(\RR^n)$. By Theorem~\ref{thm:green} and Remark~\ref{rmk:green:lower}, we have that $u_+=-\D f+\s g_+ + c_+$ and $u_-=\D f -\s g_- + c_-$.

By the jump and continuity relations \eqref{eqn:D:jump} and \eqref{eqn:S:cts}, we have that
\begin{align*}
2\nabla_\parallel f 
&=\nabla_\parallel (\Tr u_+ + \Tr u_-)
=\nabla_\parallel (-\D_+ f + \s_+ g_+ + \D_- f -\s_- g_-)
\\&=\nabla_\parallel f + \nabla_\parallel (\s_+ g_+  -\s_+ g_-)
\end{align*}
and so $\s_+ (g_+-g_-)=f$ up to an additive constant, as desired.

If $(D)^A_{p,\theta}$ is compatibly solvable then we can choose $u_\pm $ in the intersection $\dot W(p,\theta,2)\cap \dot W(2,1/2,2)$, and thus $g_+-g_- \in \dot B^{p,p}_{\theta-1}(\RR^n)\cap \dot B^{2,2}_{-1/2}(\RR^n)$.

The same argument is valid in the case that $(D)^A_{p,1}$ has the existence property, with Theorem~\ref{thm:trace} replaced by Theorem~\ref{thm:trace:NTM}.
If $(N)^A_{p,\theta}$ or $(N)^A_{p,1}$ has the existence property, then a similar argument is valid, but instead we let $u_\pm$ be solutions to $(N)^A_{p,\theta}$ and use the jump and continuity relations \eqref{eqn:S:jump} and~\eqref{eqn:D:cts} rather than \eqref{eqn:D:jump} and~\eqref{eqn:S:cts}.
\end{proof}

If $(D)^A_{p,\theta}$ is well-posed, then it has both the existence and uniqueness properties, and so $\s^A_+$ is one-to-one and onto; we use a classic argument of Verchota \cite{Ver84} to bound the inverse.
\begin{thm} \label{thm:verchota}
Let $A$, $p$ and $\theta$ satisfy the conditions specified at the beginning of this chapter.

If $(D)^A_{p,\theta}$ is well-posed in $\RR^{n+1}_+$ and $\RR^{n+1}_-$ for some $p$ and $\theta$ satisfying the conditions of Theorem~\ref{thm:bounded}, then the operator $\s^A_+$ is invertible $\dot B^{p,p}_{\theta-1}(\RR^n)\mapsto  \dot B^{p,p}_\theta(\RR^n)$ or $H^p(\R^n)\mapsto \dot H^p_1(\R^n)$, and the inverse map $(\s^A_+)^{-1}$ is bounded with operator norm depending only on $p$, $\theta$, the standard constants, and the constant $C$ in the definition of solvability.

Similarly, if $(N)^A_{p,\theta}$ is well-posed in $\RR^{n+1}_+$ and $\RR^{n+1}_-$ for some $p$ and $\theta$ satisfying the conditions of Theorem~\ref{thm:bounded}, then the operator $(\partial_\nu^A \D^A)_+$ is invertible $\dot B^{p,p}_{\theta}(\RR^n)\mapsto  \dot B^{p,p}_{\theta-1}(\RR^n)$ or $\dot H^p_1(\R^n)\mapsto H^p(\R^n)$, and its inverse is bounded with norm depending only on the standard constants and on the constant $C$ in the definition of well-posedness.
\end{thm}

\begin{proof}
We begin with the case where $(D)^A_{p,\theta}$ is well-posed for some $0<\theta<1$. By formula~\eqref{eqn:S:jump}, if $p<\infty$ and $f\in \dot B^{p,p}_{\theta-1}(\R^n)$, then
\begin{equation*}f=\nu_+\cdot A\nabla \s^A f\big\vert_{\partial\R^{n+1}_+}-\nu_+\cdot A\nabla \s^A f\big\vert_{\partial\R^{n+1}_-} =(\partial_\nu^A \s^A)_+ f+(\partial_\nu^A \s^A)_- f.\end{equation*}
If $p=\infty$, then formula~\eqref{eqn:S:jump} is still valid by the adjoint formula~\eqref{eqn:S-D:adjoint} and duality with formula~\eqref{eqn:D:jump}.

Therefore, we have that
\begin{equation*}
\doublebar{f}_{\dot B^{p,p}_{\theta-1}(\RR^n)}
=\doublebar{(\partial_\nu^A \s^A)_+ f+(\partial_\nu^A \s^A)_- f}_{\dot B^{p,p}_{\theta-1}(\RR^n)}
.\end{equation*}
By Theorem~\ref{thm:trace}, we have that 
\begin{equation*}\doublebar{(\partial_\nu^A \s^A)_\pm f}_{\dot B^{p,p}_{\theta-1}(\RR^n)}
\leq C \biggl(\int_{\R^{n+1}_\pm} \biggl(\fint_{\Omega(x,t)} \abs{\nabla \s^A f}^2\biggr)^{p/2} t^{p-1-p\theta}\,dx\,dt\biggr)^{1/p}.\end{equation*}
By Theorem~\ref{thm:bounded}, we have that the right-hand side is finite. Therefore, by the uniqueness property for~$(D)^A_{p,\theta}$, we have that $\s^A f$ must equal the solution with boundary data~$\s^A_+ f$ in the definition of solvability; thus the right-hand side is at most $C\doublebar{\s^A_\pm g}_{\dot B^{p,p}_\theta(\R^n)}$. Thus, by formula~\eqref{eqn:S:cts},
\begin{equation*}\doublebar{f}_{\dot B^{p,p}_{\theta-1}(\RR^n)}
\leq C \doublebar{S_+ f }_{\dot B^{p,p}_\theta(\R^n)}+C \doublebar{S_- f }_{\dot B^{p,p}_\theta(\R^n)}
= 2C \doublebar{S_+f }_{\dot B^{p,p}_\theta(\R^n)}.\end{equation*}

This bounds the operator norm of~$(\s^A_+)^{-1}$.

We may easily modify the argument above to establish invertibility of $(\partial_\nu^A \D^A)_+$; we need only substitute the jump relations \eqref{eqn:D:jump} and~\eqref{eqn:D:cts} for the relations \eqref{eqn:S:jump} and~\eqref{eqn:S:cts}, and use well-posedness of $(N)^A_{p,\theta}$ instead of~$(D)^A_{p,\theta}$.
Similarly, we may modify the argument above to establish invertibility of the operators on Hardy spaces by using Theorem~\ref{thm:trace:NTM} instead of Theorem~\ref{thm:trace} and the bounds \eqref{eqn:S:NTM} and \eqref{eqn:D:NTM:1} instead of Theorem~\ref{thm:bounded}.
\end{proof}

We have now proven Theorems~\ref{thm:well-posed:invertible} and~\ref{thm:well-posed:invertible:1}.
We will now prove Theorem~\ref{thm:well-posed:compatible}.

\begin{proof}[Proof of Theorem~\ref{thm:well-posed:compatible}]
Suppose that $(D)^A_{p,\theta}$ is compatibly well-posed in both $\R^{n+1}_+$ and $\R^{n+1}_-$, and so $\s^A_+$ is invertible $\dot B^{p,p}_{\theta-1}(\R^n)\mapsto \dot B^{p,p}_\theta(\R^n)$. We want to show that the inverses of the operators $(\s^A_+)^{-1}:\dot B^{p,p}_{\theta-1}(\R^n)\mapsto \dot B^{p,p}_\theta(\R^n)$ and $(\s^A_+)^{-1}:\dot B^{2,2}_{-1/2}(\R^n)\mapsto \dot B^{2,2}_{1/2}(\R^n)$ coincide.  By Theorem~\ref{thm:2,1/2}, we have that $\s^A_+$ is also invertible $\dot B^{2,2}_{-1/2}(\R^n)\mapsto \dot B^{2,2}_{1/2}(\R^n)$. 

By Theorem~\ref{thm:onto}, we have that $\s^A_+$ is surjective $\dot B^{p,p}_{\theta-1}(\RR^n)\cap \dot B^{2,2}_{-1/2}(\R^n)\mapsto  \dot B^{p,p}_\theta(\RR^n)\cap \dot B^{2,2}_{1/2}(\R^n)$. Thus, if $f\in \dot B^{p,p}_\theta(\RR^n)\cap \dot B^{2,2}_{1/2}(\R^n)$, then the unique function $g\in \dot B^{p,p}_{\theta-1}(\R^n)$ with $\s^A_+g=f$, whose existence is guaranteed by well-posedness of $(D)^A_{p,\theta}$, lies in $\dot B^{2,2}_{-1/2}(\R^n)$, and so is necessarily equal to the unique function $h$ in $\dot B^{2,2}_{-1/2}(\R^n)$ with $\s^A_+ h=f$. 

Thus, the operator $(\s^A_+)^{-1}$ is compatible on the spaces $\dot B^{p,p}_\theta(\R^n)$ and~$\dot B^{2,2}_{1/2}(\R^n)$.

As usual, a similar argument is valid for the operator $(\partial_\nu^A\D^A)_+$ and for the endpoint cases $(D)^A_{p,1}$ and~$(N)^A_{p,1}$.
\end{proof}
We remark that if $(D)^A_{p_0,\theta_0}$ and $(D)^A_{p_1,\theta_1}$ are both compatibly well-posed for $p_j<\infty$, then the operator $(\s^A_+)^{-1}$ is compatible on the spaces $\dot B^{p_0,p_0}_{\theta_0}(\R^n)$ and $\dot B^{p_1,p_1}_{\theta_1}(\R^n)$.

This theorem is important in the proof of Corollaries~\ref{cor:well-posed:open} and~\ref{cor:well-posed:1}. Specifically, we wish to use interpolation to show that $\s^A_+$ is invertible by showing that $(\s^A_+)^{-1}$ is bounded. If $(\s^A_+)^{-1}$ is not compatible, then $(\s^A_+)^{-1}$ is essentially two \emph{different} operators $\dot B^{p,p}_{\theta}(\RR^n)\mapsto  \dot B^{p,p}_{\theta-1}(\RR^n)$ and $\dot B^{2,2}_{1/2}(\R^n)\mapsto \dot B^{2,2}_{-1/2}(\R^n)$. Thus, $(\s^A_+)^{-1}$ is not well-defined on $\dot B^{p,p}_{\theta}(\RR^n)+\dot B^{2,2}_{1/2}(\R^n)$, and so interpolation methods do not apply. 

\section[Invertibility and functional analysis]{Invertibility and functional analysis: Corollaries~\ref{cor:well-posed:open}, \ref{cor:well-posed:1} and~\ref{cor:well-posed:duality}}
\label{sec:invertible:function}

In this section, we use some results of functional analysis, together with the equivalence of well-posedness and invertibility of the previous section, to prove some useful corollaries.

\begin{proof}[Proof of Corollary~\ref{cor:well-posed:open}]
Recall that $WP(D)$ is the set of points $(\theta,1/p)$ with $p<\infty$ such that $(D)^A_{p,\theta}$ is compatibly well-posed and such that $\theta$ and~$p$ satisfy the conditions of Theorem~\ref{thm:bounded}. The set $WP(N)$ is defined similarly. By Theorem~\ref{thm:2,1/2}, $WP(D)$ and $WP(N)$ both contain the point $(1/2,1/2)$. We wish to show that $WP(D)$ and $WP(N)$ are open and convex. But we have established that well-posedness of boundary-value problems is equivalent to invertibility of layer potentials; thus, convexity follows from Theorem~\ref{thm:B-B:interpolation:complex}, and openness follows from 
Theorem~\ref{thm:interpolation:extrapolation}. 
\end{proof}

\begin{proof}[Proof of Corollary~\ref{cor:well-posed:1}]
We give the proof of Corollary~\ref{cor:well-posed:1} for the Dirichlet problems~$(D)^A_{p,\theta}$; a similar technique establishes the results for the Neumann problems~$(N)^A_{p,\theta}$. Let $\varepsilon$ be as at the start of the chapter. Recall that we assumed that $(D)^A_{p_0,1}$ is compatibly well-posed for some $p_0$ with $1/2-\varepsilon<1/p_0<1+\alpha/n$. We wish to show that $(D)^A_{p,\theta}$ is compatibly well-posed for all $(\theta,1/p)$ as in Figure~\ref{fig:well-posed:1}.

By Theorem~\ref{thm:well-posed:invertible:1} and formulas \eqref{eqn:H-F} and~\eqref{eqn:H1-F}, we have that $\s^A_+$ is invertible $\dot F^{p_0,2}_0(\R^n)\mapsto \dot F^{p_0,2}_1(\R^n)$.

We apply Theorem~\ref{thm:interpolation:extrapolation} to the spaces $\dot F^{p,2}_0(\R^n)=H^p(\R^n)$ and $\dot F^{p,2}_1(\R^n)=\dot H^p_1(\R^n)$; we may do this due to the boundedness results~\eqref{eqn:S:NTM} and Theorem~\ref{thm:trace:NTM} and the interpolation formula~\eqref{eqn:F-F:interpolation:complex}. We conclude that there is some $p_1$ and $p_3$, with $1/2-\varepsilon<1/p_1<1/p_0<1/p_3<1+\alpha/n$, such that $\s^A_+$ is invertible $\dot F^{p_2,2}_0(\R^n)\mapsto \dot F^{p_2,2}_1(\R^n)$ whenever $1/p_1<1/p_2<1/p_3$.

Recall from Theorem~\ref{thm:verchota} that $\s^A_+$ is invertible $\dot B^{2,2}_{-1/2}(\RR^n)\mapsto  \dot B^{2,2}_{1/2}(\RR^n)$. By the definition of Besov and Triebel-Lizorkin spaces, $\dot B^{2,2}_{\theta}(\RR^n)=\dot F^{2,2}_{\theta}(\RR^n)$ for any real number~$\theta$.

We apply complex interpolation, that is, formula~\eqref{eqn:F-F:interpolation:complex}. By Theorem~\ref{thm:well-posed:compatible} we may use interpolation to bound~$(\s^A_+)^{-1}$. We then have that $\s^A_+$ is invertible $\dot F^{p,2}_{\theta-1}(\R^n)\mapsto\dot F^{p,2}_{\theta}(\R^n)$ whenever $1/2<\theta<1$ and whenever $1/p=\theta/p_2+(1-\theta)/p_2'$ for some $p_2$ with $1/p_1<1/p_2<1/p_3$. See Figure~\ref{fig:interpolation:F}.

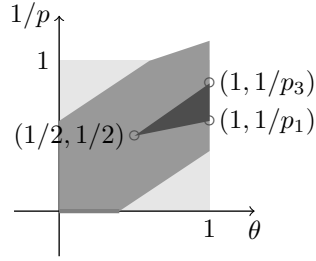
\begin{figure}[tbp]
\begin{tikzpicture}[scale=2]
\figureaxes
\boundedhexagon
\fill [well-posed] 
	(1,0.85) node [right,black] {$(1,1/p_3)$} node {$\circ$} --
	(1,0.6) node [right,black] {$(1,1/p_1)$} node {$\circ$} --
	(1/2,1/2) node [left,black] {$(1/2,1/2)$} node {$\circ$}-- cycle;	
\end{tikzpicture}
\caption{Under the assumptions of Corollary~\ref{cor:well-posed:1}, the operator $\s^A_+$ is invertible $\dot F^{p,2}_{\theta-1}(\R^n)\mapsto\dot F^{p,2}_{\theta}(\R^n)$ whenever the point $(\theta,1/p)$ lies in the indicated triangle. }\label{fig:interpolation:F}
\end{figure}

We wish to return to the spaces $\dot B^{p,p}_\theta(\R^n)$.
Notice that if $\theta$ and $p$ satisfy these conditions, then there is some $\theta_1$, $\theta_2$ with $\theta_1<\theta<\theta_2$ such that $\s^A_+$ is invertible $\dot F^{p,2}_{\theta_1-1}(\R^n)\mapsto\dot F^{p,2}_{\theta_1}(\R^n)$
and $\dot F^{p,2}_{\theta_2-1}(\R^n)\mapsto\dot F^{p,2}_{\theta_2}(\R^n)$. We may thus use real interpolation (that is, formula~\eqref{eqn:F-F:interpolation:real}) to see that $\s^A_+$ is invertible $\dot B^{p,p}_{\theta-1}(\R^n)\mapsto\dot B^{p,p}_{\theta}(\R^n)$, as desired.
\end{proof}

\begin{proof}[Proof of Corollary~\ref{cor:well-posed:duality}]
Suppose $(D)^A_{p,\theta}$ is well-posed in $\R^{n+1}_+$ and $\R^{n+1}_-$ for some $p$ and $\theta$ that satisfy the conditions of Theorem~\ref{thm:bounded} with $p<\infty$. We wish to show that $(D)^{A^*}_{q,\sigma}$ is also well-posed for appropriate $q$ and~$\sigma$.

We have that $\s^A_+$ is invertible $\dot B^{p,p}_{\theta-1}(\R^n)\mapsto \dot B^{p,p}_{\theta}(\R^n)$. 
The adjoint to $\s^A_+$ is $\s^{A^*}_+$, and so the adjoint to $(\s^A_+)^{-1}$ is $(\s^{A^*}_+)^{-1}$. Recall the dual spaces to the Besov spaces given by formulas~\eqref{eqn:dual:banach} and~\eqref{eqn:dual:quasi-banach}. We have that if $1\leq p<\infty$, then $\s^{A^*}_+$ is invertible $\dot B^{p',p'}_{-\theta}(\R^n)\mapsto \dot B^{p',p'}_{1-\theta}(\R^n)$, and so $(D)^{A^*}_{p',1-\theta}$ is well-posed. Furthermore, if $0<p<1$, then $\s^{A^*}_+$ is invertible $\dot B^{\infty,\infty}_{-\theta+n(1/p-1)}(\R^n)\mapsto \dot B^{\infty,\infty}_{1-\theta+n(1/p-1)}(\R^n)$ and so $(D)^A_{\infty,1-\theta+n(1/p-1)}$ is well-posed. If $(D)^A_{p,\theta}$ is compatibly well-posed, then the inverses to  $\s^A_+:\dot B^{p,p}_{\theta-1}(\R^n)\mapsto \dot B^{p,p}_{\theta}(\R^n)$ and $\s^A_+:\dot B^{2,2}_{-1/2}(\R^n)\mapsto \dot B^{2,2}_{1/2}(\R^n)$ coincide; thus, so do the two inverses to $\s^{A^*}_+$, and so $(D)^{A^*}_{p',1-\theta}$ or $(D)^{A^*}_{\infty,1-\theta+n(1/p-1)}$ is compatibly well-posed.

We recall some basic results of functional analysis. If $T:B_1\mapsto B_2$ is a linear operator and the $B_j$s are reflexive Banach spaces, then $T$ is one-to-one if and only if its adjoint $T^*:B_2'\mapsto B_1'$ is onto. If $1<p<\infty$ then $\dot B^{p,p}_\theta(\R^n)$ is reflexive; thus, surjectivity of $\s^A_+:\dot B^{p,p}_{\theta-1}(\R^n)\mapsto \dot B^{p,p}_{\theta}(\R^n)$ is equivalent to injectivity of $\s^{A^*}_+:\dot B^{p',p'}_{-\theta}(\R^n)\mapsto \dot B^{p',p'}_{1-\theta}(\R^n)$, and so the existence property for $(D)^A_{p,\theta}$ is equivalent to the uniqueness property for $(D)^{A^*}_{p',1-\theta}$, as desired.




Again, the same argument is valid for $(N)^A_{p,\theta}$.
\end{proof}

We conclude this section with the following remark.

\begin{rmk}\label{rmk:compatible} Suppose that $p$, $\theta$ satisfy the conditions specified at the beginning of this chapter and that $p<\infty$. We claim that compatible solvability implies compatible well-posed\-ness. 

Suppose that $(D)^A_{p,\theta}$ is compatibly solvable for $\theta<1$.
We need only establish the uniqueness property; by Theorem~\ref{thm:unique:D:N}, we need only show that $\s^A_+$ is one-to-one. We apply the method of Theorem~\ref{thm:verchota}. Let $f\in \dot B^{p,p}_{\theta-1}(\R^n)\cap \dot B^{2,2}_{-1/2}(\R^n)$. Then $\s^A f\in \dot W_\pm(p,\theta,2)\cap \dot W_\pm(2,1/2,2)$; because $(D)^A_{2,1/2}$ is well-posed, $\s^A f$ must be the solution to $(D)^A_{p,\theta}$ with boundary data~$\s^A_+ f$ in the definition of compatible solvability. Therefore, we have that $\doublebar{\s^A f}_{\dot W(p,\theta,2)}\leq C\doublebar{\s^A_+ f}_{\dot B^{p,p}_\theta(\R^n)}$. Proceeding as in the proof of Theorem~\ref{thm:verchota}, we see that 
$\doublebar{f}_{\dot B^{p,p}_{\theta-1}(\RR^n)} \leq C \doublebar{\s^A_+f }_{\dot B^{p,p}_\theta(\R^n)}$
whenever $f\in \dot B^{p,p}_{\theta-1}(\RR^n)\cap \dot B^{2,2}_{-1/2}(\R^n)$.
If $p<\infty$ then $\dot B^{p,p}_{\theta-1}(\RR^n)\cap \dot B^{2,2}_{-1/2}(\R^n)$ is dense in $\dot B^{p,p}_{\theta-1}(\RR^n)$, and so $\s^A_+$ extends to a one-to-one operator. 

A similar argument is valid if $(D)^A_{p,1}$ is compatibly solvable; an argument with $(\partial_\nu^A\D^A)_+$ in place of~$\s^A_+$ is valid if $(N)^A_{p,\theta}$ is compatibly solvable.
\end{rmk}

\section[Extrapolation of well-posedness and real coefficients]{Extrapolation of well-posedness and real coefficients: Corollaries~\ref{cor:well-posed:AusM} and~\ref{cor:well-posed:real}}
\label{sec:real}

In this section, we use the extrapolation theorem of Auscher and Mourgoglou (Theorem~\ref{thm:extrapolation:AusM} above) and the known well-posedness results for real coefficients to prove an extrapolation theorem and well-posedness in Besov spaces.

\begin{proof}[Proof of Corollary~\ref{cor:well-posed:AusM}]
Suppose that $A$ is elliptic and $t$-independent and that $A$, $A^*$, $A^\sharp$ and $(A^\sharp)^*$ all satisfy the De Giorgi-Nash-Moser condition with exponent~$\alpha^\sharp$, where $A^\sharp$ is as in formula~\eqref{eqn:AusM}. Notice that it might be the case that $A$ and  $A^*$ satisfy the De Giorgi-Nash-Moser condition with exponent~$\alpha$ for some $\alpha>\alpha^\sharp$; for example, this is true for the matrices $A_k$ in Section~\ref{sec:sharp}. 

Corollary~\ref{cor:well-posed:AusM} concerns matrices~$A$ such that $(D)^A_{p_0,1}$ or $(N)^A_{p_0,1}$ is compatibly well-posed; we will present the argument in the case where $(D)^A_{p_0,1}$ is compatibly well-posed and leave the case of the Neumann problem to the reader.

\def\figurelabel#1#2{
	\node at (-0.2,1.55) [right] {#2};
	\node at (-0.2,1.0) [left] {(\textsl{#1})};}
\def\phantoms{
	\node at (1,0) [right] {$\phantom{(1,1+\alpha^\sharp/n)}$};
	\node at (0,0) [left] {$\phantom{1/p_0'}$};
	\node at (-0.2,1.0) [left] {\phantom{(\textsl{b})}};
	}

\begin{figure}[tbp]
	
\begin{tikzpicture}[scale=2]
\drawunitsquare
\plainfigureaxes
\boundedhexagon
\phantoms
\fill [well-posed] 
	(1,1+\alphasharp/\enn) node [right,black] {$(1,1+\alpha^\sharp/n)$} node [black] {$\circ$} --
	(1,\pnought) node [right,black] {$(1,1/p_0)$} node [black] {$\circ$} --
	(1/2,1/2)-- cycle;	
\figurelabel{a}{Invertibility of $\s^A_+$}
\end{tikzpicture}
\qquad
\begin{tikzpicture}[scale=2]
\drawunitsquare
\plainfigureaxes
\boundedhexagon
\phantoms
\path[name path=x axis] (0,0)--(1,0);
\path[name path=extrapolation] (1/2,1/2)--(0,-\alphasharp/\enn);
\fill [name intersections={of=x axis and extrapolation, by=x}, well-posed] 
	(x) --
	(0,0)--
	(0,1-\pnought) node [left,black] {$1/p_0'$} node [black] {$\circ$} --
	(1/2,1/2)-- cycle;	
\draw [boundary well-posed] (0,0)--(\alphasharp,0) node [black] {$\circ$} node [below,black] {$\alpha^\sharp$};
\figurelabel{b}{Invertibility of $\s^{A^*}_+$}
\end{tikzpicture}

\begin{tikzpicture}[scale=2]
\drawunitsquare
\plainfigureaxes
\boundedhexagon
\phantoms
\fill [well-posed] 
	(1,1+\alphasharp/\enn) node [right,black] {$(1,1+\alpha^\sharp/n)$} node [black] {$\circ$} --
	(1,\pnought) node [right,black] {$(1,1/p_0)$} node [black] {$\circ$} --
	(1/2,1/2)-- cycle;	
\fill [well-posed] (1/2,1/2)--(1,\pnought)--(1,1)--
	(1-\alphasharp,1) 
	node [black] {$\circ$} 
	node [above left,black, at = {(1-0.3*\alphasharp,1)}] {$(1-\alpha^\sharp,1)$} 
	-- cycle;
\figurelabel{d}{}
\end{tikzpicture}
\qquad
\begin{tikzpicture}[scale=2]
\drawunitsquare
\plainfigureaxes
\boundedhexagon
\phantoms
\fill [well-posed] 
	(\alphasharp,0) node [below,black] {$\alpha^\sharp$} node [black] {$\circ$} --
	(0,0)--
	(0,1-\pnought) node [left,black] {$1/p_0'$} node [black] {$\circ$} --
	(1/2,1/2)-- cycle;	
\draw [boundary well-posed] (0,0)--(\alphasharp,0);
\figurelabel{c}{}
\end{tikzpicture}

\begin{tikzpicture}[scale=2]
\drawunitsquare
\plainfigureaxes
\boundedhexagon
\phantoms
\fill [well-posed] 
	(1,1+\alphasharp/\enn) node [right,black] {$(1,1+\alpha^\sharp/n)$} node [black] {$\circ$} --
	(1,\pnought) node [right,black] {$(1,1/p_0)$} node [black] {$\circ$} --
	(1/2,1/2)--
	(1-\alphasharp,1) 
	node [black] {$\circ$} 
	node [above left,black, at = {(1-0.3*\alphasharp,1)}] {$(1-\alpha^\sharp,1)$} 
	-- cycle;
\figurelabel{e}{}
\end{tikzpicture}
\qquad
\begin{tikzpicture}[scale=2]
\phantoms
\end{tikzpicture}

\caption{Invertibility of $\s^A_+$, and thus well-posedness of $(D)^A_{p,\theta}$, under the assumption that $A$ satisfies the conditions of \cite{AusM} and $(D)^A_{p_0,1}$ is compatibly well-posed. In each figure we plot the values of $(\theta,1/p)$ such that $\s^A_+$ or $\s^{A^*}_+$ is invertible on $\dot B^{p,p}_\theta(\R^n)$; as we move through the figure, more and more values of $(\theta,1/p)$ are allowed.} \label{fig:AusM:1}
\end{figure}
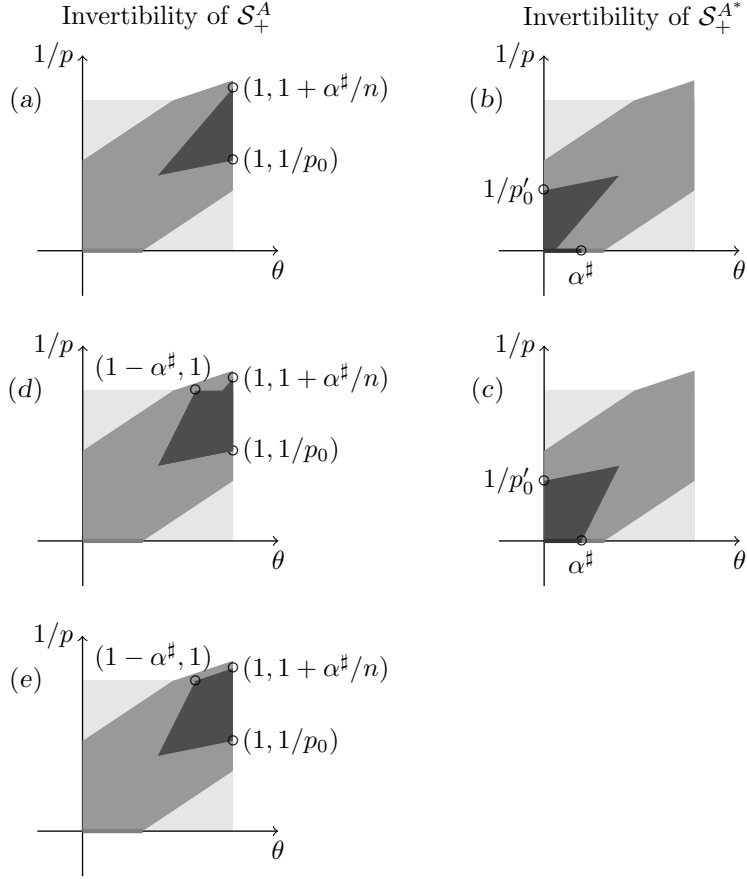

\begin{figure}

\begin{tikzpicture}[scale=2]
\drawunitsquare
\plainfigureaxes
\boundedhexagon
\phantoms
\fill [well-posed] 
	(1,1+\alphasharp/\enn) node [right,black] {$(1,1+\alpha^\sharp/n)$} node [black] {$\circ$} --
	(1,\pnought) node [right,black] {$(1,1/p_0)$} node [black] {$\circ$} --
	(1/2,1/2)--
	(1-\alphasharp,1) 
	node [black] {$\circ$} 
	node [above left,black, at = {(1-0.3*\alphasharp,1)}] {$(1-\alpha^\sharp,1)$} 
	-- cycle;
\figurelabel{a}{Invertibility of $\s^A_+$}
\end{tikzpicture}
\qquad
\begin{tikzpicture}[scale=2]
\drawunitsquare
\plainfigureaxes
\boundedhexagon
\phantoms
\fill [well-posed] 
	(\alphasharp,0) node [below,black] {$\alpha^\sharp$} node [black] {$\circ$} --
	(0,0)--
	(0,1-\pnought) node [left,black] {$1/p_0'$} node [black] {$\circ$} --
	(1/2,1/2)-- cycle;	
\draw [boundary well-posed] (0,0)--(\alphasharp,0);
\figurelabel{b}{Invertibility of $\s^{A^*}_+$}
\end{tikzpicture}

\begin{tikzpicture}[scale=2]
\drawunitsquare
\plainfigureaxes
\boundedhexagon
\phantoms
\fill [well-posed] 
	(1,1+\alphasharp/\enn) node [right,black] {$(1,1+\alpha^\sharp/n)$} node [black] {$\circ$} --
	(1,\pnought) node [right,black] {$(1,1/p_0)$} node [black] {$\circ$} --
	(1/2,1/2)--
	(1-\alphasharp,1) 
	node [black] {$\circ$} 
	node [above left,black, at = {(1-0.3*\alphasharp,1)}] {$(1-\alpha^\sharp,1)$} 
	-- cycle;
\fill [well-posed] 
	(\alphasharp,0) node [below,black] {$\alpha^\sharp$} node [black] {$\circ$} --
	(0,0)--
	(0,1-\pone) node [left,black] {$1/p_1'$} node [black] {$\circ$} --
	(1/2,1/2)-- cycle;	
\draw [boundary well-posed] (0,0)--(\alphasharp,0);
\figurelabel{c}{}
\end{tikzpicture}
\qquad
\begin{tikzpicture}[scale=2]
\drawunitsquare
\plainfigureaxes
\boundedhexagon
\phantoms
\fill [well-posed] 
	(\alphasharp,0) node [below,black] {$\alpha^\sharp$} node [black] {$\circ$} --
	(0,0)--
	(0,1-\pnought) node [left,black] {$1/p_0'$} node [black] {$\circ$} --
	(1/2,1/2)-- cycle;	
\fill [well-posed] 
	(1,1+\alphasharp/\enn) node [right,black] {$(1,1+\alpha^\sharp/n)$} node [black] {$\circ$} --
	(1,\pone) node [right,black] {$(1,1/p_1)$} node [black] {$\circ$} --
	(1/2,1/2)-- 
	(1-\alphasharp,1) 
	node [black] {$\circ$} 
	node [above left,black, at = {(1-0.3*\alphasharp,1)}] {$(1-\alpha^\sharp,1)$} 
	--
	cycle;	
\draw [boundary well-posed] (0,0)--(\alphasharp,0);
\figurelabel{d}{}
\end{tikzpicture}

\begin{tikzpicture}[scale=2]
\drawunitsquare
\plainfigureaxes
\boundedhexagon
\fill [well-posed] 
	(\alphasharp,0) node [below,black] {$\alpha^\sharp$} node [black] {$\circ$} --
	(0,0)--
	(0,1-\pone) node [left,black] {$1/p_1'$} node [black] {$\circ$} --
		(1-\alphasharp,1) 
		node [black] {$\circ$} 
		node [above left,black, at = {(1-0.3*\alphasharp,1)}] {$(1-\alpha^\sharp,1)$} 
		-- 
	(1,1+\alphasharp/\enn) node [right,black] {$(1,1+\alpha^\sharp/n)$} node [black] {$\circ$} --
	(1,\pnought) node [right,black] {$(1,1/p_0)$} node [black] {$\circ$} --
	cycle;
\draw [boundary well-posed] (0,0)--(\alphasharp,0);
\figurelabel{e}{}
\end{tikzpicture}
\qquad
\begin{tikzpicture}[scale=2]
\phantoms
\end{tikzpicture}

\caption{The continuation of Figure~\ref{fig:AusM:1} under the assumption that $(D)^{A^*}_{p_1,1}$ is also well-posed for some $1<p_1\leq 2$.} \label{fig:AusM:2}
\end{figure}
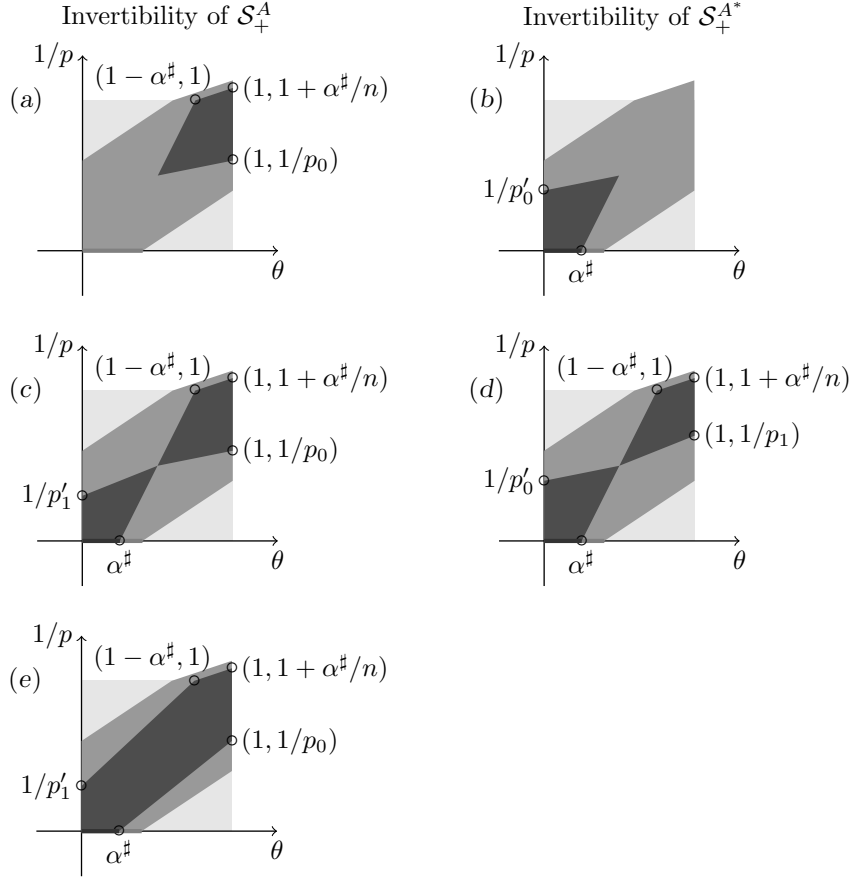

Suppose that $(D)^A_{p_0,1}$ is compatibly well-posed for some $1<p_0\leq 2$. By Theorem~\ref{thm:extrapolation:AusM}, we have that $(D)^A_{p,1}$ is compatibly well-posed for all $p$ with $n/(n+\alpha)<p<p_0$.

By Corollary~\ref{cor:well-posed:1} and Theorem~\ref{thm:well-posed:invertible},  $\s^A_+$ is invertible $\dot B^{p,p}_{\theta-1}(\R^n)\mapsto \dot B^{p,p}_\theta(\R^n)$ whenever $(\theta,1/p)$ lies in the triangular region in Figure~\ref{fig:AusM:1}(\textsl{a}). 
By Corollary~\ref{cor:well-posed:duality} and Theorem~\ref{thm:well-posed:invertible}, $\s^{A^*}_+$ is invertible whenever $(\theta,1/p)$ lies in the trapezoidal region in Figure~\ref{fig:AusM:1}(\textsl{b}), and is also invertible $\dot B^{\infty,\infty}_{\theta-1}(\R^n)\mapsto \dot B^{\infty,\infty}_\theta(\R^n)$ whenever $0<\theta<\alpha$.

Recall that Corollary~\ref{cor:well-posed:1} was proved by interpolation. Notice that all of the steps above preserve compatibility with $(D)^A_{2,1/2}$.
We may interpolate again to see that 
$\s^{A^*}_+$ is invertible $\dot B^{p,p}_{\theta-1}(\R^n)\mapsto \dot B^{p,p}_\theta(\R^n)$ for all $(\theta,1/p)$ in the four-sided region in  Figure~\ref{fig:AusM:1}(\textsl{c}).

By duality, we have that $\s^A_+$ is invertible $\dot B^{p,p}_{\theta-1}(\R^n)\mapsto \dot B^{p,p}_\theta(\R^n)$ whenever $(\theta,1/p)$ lies in the oddly shaped region in Figure~\ref{fig:AusM:1}(\textsl{d}). We may immediately use interpolation to fill in the missing triangle in the region $1/p>1$ to yield a convex region; see Figure~\ref{fig:AusM:1}(\textsl{e}).

Now we pass to the case where $(D)^{A^*}_{p_1,1}$ is also well-posed. We have already established that $\s^A_+$ and $\s^{A^*}_+$ are invertible $\dot B^{p,p}_{\theta-1}(\R^n)\mapsto \dot B^{p,p}_\theta(\R^n)$ whenever 
$(\theta,1/p)$ lies in the four-sided regions in Figure~\ref{fig:AusM:2}(\textsl{a}) or Figure~\ref{fig:AusM:2}(\textsl{b}), respectively. By the same argument, we have well-posedness in the new four-sided regions in Figure~\ref{fig:AusM:2}(\textsl{c}) or Figure~\ref{fig:AusM:2}(\textsl{d}).
We may then  interpolate to see that $\s^A_+$ is invertible $\dot B^{\infty,\infty}_{\theta-1}(\R^n)\mapsto \dot B^{\infty,\infty}_\theta(\R^n)$, and thus that $(D)^A_{p,\theta}$ is well-posed, whenever 
$(\theta,1/p)$ lies in the hexagonal region in Figure~\ref{fig:AusM:1}(\textsl{e}).
\end{proof}

\begin{proof}[Proof of Corollary~\ref{cor:well-posed:real}] We wish to show that if $A$ is real, then $A$ satisfies the conditions of Corollary~\ref{cor:well-posed:AusM}, at least in the case of the Dirichlet problem. (In this theorem the corresponding arguments for the Neumann problem are \emph{not} known to be valid!)
Recall from Theorem~\ref{thm:DGNM} that all real elliptic matrices satisfy the De Giorgi-Nash-Moser condition in the interior and at the boundary. Thus, $A$ satisfies the conditions of our Theorems~\ref{thm:bounded} and~\ref{thm:trace}. Furthermore, since the matrix $A^\sharp$ is also real, $A$ satisfies the conditions of Theorem~\ref{thm:extrapolation:AusM}.

By Theorems~\ref{thm:real:dirichlet} and~\ref{thm:extrapolation:HKMP}, we have that $(D)^A_{p_0,1}$ and $(D)^{A^*}_{p_1,1}$ are solvable for some $p_0$, $p_1$ with $1<p_j\leq 2$. Furthermore, by Theorem~\ref{thm:well-posed:L2} if $A$ is real and symmetric then we may take $p_j=2$, and furthermore $(N)^A_{2,1}$ and $(N)^{A^*}_{2,1}$ are also solvable. Corollary~\ref{cor:well-posed:real} will then follow from Corollary~\ref{cor:well-posed:AusM} if we can show that these boundary-value problems are \emph{compatibly} solvable.

We intend to use the method of continuity, as used in \cite[Section~5]{HofMitMor} or \cite[Section~9]{AlfAAHK11}. As a starting point, observe that in the case of harmonic functions we may write solutions explicitly using the Poisson kernel, and so if $I$ is the identity matrix then $(D)^I_{p,\theta}$ is compatibly well-posed for any $0\leq\theta\leq 1$ and any $p\leq \infty$ with $1/p<1+\theta/n$. Thus, by Theorems~\ref{thm:well-posed:invertible:1} and~\ref{thm:well-posed:compatible}, $\s^I_+$ is invertible with bounded inverse $H^{p_0}(\R^n)\cap \dot B^{2,2}_{-1/2}(\R^n)\mapsto \dot H^{p_0}_1(\R^n)\cap \dot B^{2,2}_{1/2}(\R^n)$. 

Let $A_z=(1-z) I+z A$. Note that if $\sigma$ is real with $0\leq \sigma\leq 1$, then $A_\sigma$ is real-valued, elliptic and $t$-independent. Moreover, the ellipticity constants of $A_\sigma$ may be taken to depend only on~$A$, not on~$\sigma$. By Theorems~\ref{thm:real:dirichlet} and~\ref{thm:extrapolation:HKMP}, we have that $(D)^{A_\sigma}_{p_0,1}$ is solvable for some $p_0$ with $1<p_0\leq 2$, and that both $p_0$ and the constant~$C$ in the definition of solvability are independent of~$\sigma$.

By Theorem~\ref{thm:perturb:A:DGNM}, the bound \eqref{eqn:trace:S}, and the bound~\eqref{eqn:S:NTM} and Theorem~\ref{thm:trace:NTM}, we have that $\s^{A_z}_+:H^{p_0}(\R^n)\cap \dot B^{2,2}_{-1/2}(\R^n)\mapsto \dot H^{p_0}_1(\R^n)\cap \dot B^{2,2}_{1/2}(\R^n)$ is bounded for all $z$ in a complex neighborhood of the interval $[0,1]$. Furthermore, $\s^{A_0}_+=\s^I_+$ is invertible with bounded inverse. Thus by analytic perturbation theory, $\s^{A_\sigma}_+$ is invertible and $(D)^{A_\sigma}_{p_0,1}$ is not merely solvable but compatibly well-posed for all $\sigma$ in a small neighborhood of~$0$. 

But by Theorem~\ref{thm:verchota}, we have that the operator norm of $(\s^{A_\sigma}_+)^{-1}$ depends only on the constant~$C$ in the definition of solvability. Thus, if $0\leq \sigma\leq 1$ and $(D)^{A_{\sigma}}_{p_0,1}$ is compatibly well-posed, then $(D)^{A_z}_{p_0,1}$ is also compatibly well-posed for all $z$ in a small neighborhood of~$\sigma$, and the size of this neighborhood may be taken to be independent of~$\sigma$. Thus, we may go from $A_0=I$ to $A_1=A$ in small steps, and see that $\s^A_+$ is also invertible with bounded inverse $H^{p_0}(\R^n)\cap \dot B^{2,2}_{-1/2}(\R^n)\mapsto \dot H^{p_0}_1(\R^n)\cap \dot B^{2,2}_{1/2}(\R^n)$, as desired.

If $A$ is symmetric then the above argument is valid with $p_0=2$, and a similar argument involving the Neumann problem and the potential $(\partial_\nu^A\D^A)_+$ is valid.
\end{proof}

%% file: sec-10-besov.tex
\chapter{Besov Spaces and Weighted Sobolev Spaces}
\label{chap:besov}

Recall that we constructed solutions to the Dirichlet problem
\begin{equation}\label{eqn:dirichlet:late}
\left\{\begin{aligned} \Div A\nabla u&=0 &&\text{in }\R^{n+1}_+,\\
\Tr u &=f && \text{on } \R^n=\partial\R^{n+1}_+,\\
\end{aligned}\right.
\end{equation}
that satisfy the estimate
\begin{equation}\label{eqn:sobolev:bound}
\doublebar{u}_{\dot W(p,\theta,q)}
\leq C \doublebar{f}_{\dot B^{p,p}_\theta(\R^n)}
.\end{equation}
Recall also that in \cite{JerK95,May05,MayMit04a}, the authors investigated the Dirichlet problem
\begin{equation}\label{eqn:harmonic:dirichlet:late}
\left\{\begin{aligned} \Delta u&=0 &&\text{in }\Omega,\\
\Tr u &=f && \text{on } \partial\Omega,\\
\doublebar{u}_{\dot B^{p,p}_{\theta+1/p}(\Omega)}
&\leq C \doublebar{f}_{\dot B^{p,p}_\theta(\partial\Omega)}
\hskip-20pt
\end{aligned}\right.
\end{equation}
for $\Omega$ a Lipschitz domain. (In \cite{FabMM98,Zan00,May05,MayMit04a} a similar Neumann problem was investigated.) 

This suggests investigation of the Dirichlet problem \eqref{eqn:dirichlet:late} subject to the bound
\begin{equation}\doublebar{u}_{\dot B^{p,p}_{\theta+1/p}(\R^{n+1}_+)}
\leq C \doublebar{f}_{\dot B^{p,p}_\theta(\R^n)}
\end{equation}
instead of the bound~\eqref{eqn:sobolev:bound}. As explained in the introduction, if $\theta+1/p>1$, then this bound cannot hold for general bounded measurable coefficients.

However, for many values of $\theta$, $p$ with  $\theta+1/p\leq 1$, we can relate the two formulations of the Dirichlet problem. We will generalize the following theorem of Jerison and Kenig.
\begin{thm}[{\cite[Theorem~4.1]{JerK95}}]
Let $u$ be harmonic in a bounded Lipschitz domain $\Omega$ and let $\delta(X)=\dist(X,\partial\Omega)$. Let $k\geq 0$ be an integer, let $0<\theta<1$, and let $1\leq p\leq \infty$. 

Then $u\in \dot B^{p,p}_{k+\theta}(\Omega)\cap L^p(\Omega)$ if and only if $\delta^{1-\theta}\abs{\nabla^{k+1} u}+\abs{\nabla^k u}+\abs{u} \in L^p(\Omega)$.
\end{thm}
In fact, the inequality $\doublebar{u}_{\dot B^{p,p}_{k+\theta}(\Omega)\cap L^p(\Omega)}\leq C \doublebar{\delta^{1-\theta}\abs{\nabla^{k+1} u}+\abs{\nabla^k u}+\abs{u}}_{L^p(\Omega)}$ does not require that $u$ be harmonic; this is only needed for the reverse inequality.

We will prove the following theorem; this theorem implies that for appropriate $p$ and~$\theta$, we may require solutions to the homogeneous Dirichlet (or Neumann) problems to lie in either $\dot W(p,\theta,2)$ or $\dot B^{p,p}_{\theta+1/p}(\R^{n+1})$. Notice that if $\theta>0$ and $\theta+1/p<1$ then $p>1$, and so we will not concern ourselves with the case $p\leq 1$. (The boundary-value problem~\eqref{eqn:harmonic:dirichlet:late} with $p\leq 1$, and also the corresponding Neumann problem, were investigated in \cite{May05,MayMit04a}.) 

We remark that the proofs of Lemmas~\ref{lem:Besov-to-weighted} and~\ref{lem:weighted-to-Besov} follow closely the proof of {\cite[Theorem~4.1]{JerK95}}.

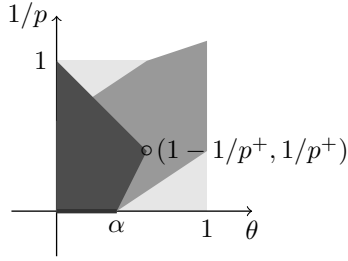
\begin{figure}
\begin{tikzpicture}[scale=2]
\figureaxes
\boundedhexagon

\fill[well-posed] (0,0)--(\alph,0) node [black, below] {$\alpha$}
--(1/2+\eps,1/2-\eps) node [black] {$\circ$} node [black, right] {$(1-1/p^+,1/p^+)$}
--(0,1)
;

\draw [boundary well-posed] (0,0)--(\alph,0);

\end{tikzpicture}
\caption{If $(\theta,1/p)$ lie in the indicated region, then solutions $u$ to $\Div A\nabla u=0$ lie in $\dot W(p,\theta,2)$ if and only if they lie in $\dot B^{p,p}_{\theta+1/p}(\R^{n+1}_+)$.}
\label{fig:besov}
\end{figure}

\begin{thm}\label{thm:weighted-and-besov}
Suppose that $A$ is an elliptic, $t$-independent matrix and that $A$ and $A^*$ satisfy the De Giorgi-Nash-Moser condition.
Let $u$ be a solution to $\Div A\nabla u=0$ in $\R^{n+1}_+$.

Suppose that the numbers $p$ and $\theta$ are such that the point $(\theta,1/p)$ lies in the region shown in Figure~\ref{fig:besov}. Then $u\in \dot W(p,\theta,2)$ if and only if $u\in \dot B^{p,p}_{\theta+1/p}(\R^{n+1}_+)$, and 
\[\frac{1}{C(p,\theta)}\doublebar{u}_{\dot W(p,\theta,2)}
\leq \doublebar{u}_{\dot B^{p,p}_{\theta+1/p}(\R^{n+1}_+)}
\leq C(p,\theta)\doublebar{u}_{\dot W(p,\theta,2)}.
\]
\end{thm}
Notice that at the endpoint $\theta+1/p=1$, we have that $\dot W(p,\theta,p)=\dot W^p_1(\R^{n+1}_+)$, the \emph{unweighted} Sobolev space of functions whose gradient lies in $L^p(\R^{n+1}_+)$. By formula~\eqref{eqn:L-F} this space is equal to $\dot F^{p,2}_1(\R^{n+1}_+)$. Recall that $\dot B^{p,p}_\theta=\dot F^{p,p}_\theta$ and that the Besov and Triebel-Lizorkin spaces are increasing in the second exponent; thus, if $p\leq 2$ then $\dot B^{p,p}_1(\R^{n+1}_+)\subset \dot W(p,1-1/p,p)$, and if $p\geq 2$ then $\dot W(p,1-1/p,p)\subset\dot B^{p,p}_1(\R^{n+1}_+)$.

We refer the reader to \cite{Tri83} for an extensive discussion of Besov spaces on subsets of $\R^{n+1}$; here we will only need that they obey the interpolation formula \eqref{eqn:L-W:interpolation} and the norm estimate~\eqref{eqn:besov:norm:differences}.

We begin with the following lemma.

\begin{lem}
\label{lem:Besov-to-weighted}
Let $A$ be as in Theorem~\ref{thm:weighted-and-besov} and suppose that $\Div A\nabla u=0$ in $\R^{n+1}_+$. Suppose that $0<\theta+1/p<1$ and that $1\leq p\leq\infty$. Then \[\doublebar{u}_{\dot W(p,\theta,2)}\leq C(p,\theta)\doublebar{u}_{\dot B^{p,p}_{\theta+1/p}(\R^{n+1}_+)}.\]
\end{lem}


\begin{proof} If $p=\infty$, then $\dot B^{\infty,\infty}_{\theta}(\R^{n+1}_+)=\dot C^{\theta}(\R^{n+1}_+)$; the bound on $\doublebar{u}_{\dot W(p,\theta,2)}$ follows from Caccioppoli's inequality.

Recall that if $1\leq p<\infty$  and $0<{\theta+1/p}<1$ then
\begin{equation*}\doublebar{u}_{\dot B^{p,p}_{\theta+1/p}(\R^{n+1}_+)}^p
\approx
\int_{\R^{n+1}_+}\int_{\R^{n+1}_+}
\frac{\abs{u(x,t)-u(y,r)}^p}{(\abs{x-y}+\abs{t-r})^{n+1+p(\theta+1/p)}}\,dx\,dt\,dy\,dr
.\end{equation*}
Furthermore, if $\G$ is a grid of dyadic Whitney cubes, then
\[\doublebar{u}_{\dot W(p,\theta,2)}^p
\approx 
	C\sum_{Q\in\G}\ell(Q)^{n+p-p\theta}\biggl(\fint_{Q} \abs{\nabla u}^2\biggr)^{p/2}.
\]
By Caccioppoli's inequality and Corollary~\ref{cor:local-bound}, if $Q\in\G$ then for any $(y,s)\in \R^{n+1}_+$,
\[
\ell(Q)^{n+p-p\theta}\biggl(\fint_{Q} 
\abs{\nabla u}^2\,dt\,dx\biggr)^{p/2}
\leq 
	C\ell(Q)^{-1-p\theta}\int_{(3/2)Q} \abs{u(x,t)-u(y,s)}^p\,dt\,dx
.\]
We may average over $(y,s)\in Q$ to see that
\begin{multline*}
\ell(Q)^{n+p-p\theta}\biggl(\fint_{Q} 
\abs{\nabla u}^2\,dt\,dx\biggr)^{p/2}
\\\begin{aligned}
&\leq 
	C\ell(Q)^{-1-p\theta-n-1}\int_{(3/2)Q} \int_Q \abs{u(x,t)-u(y,s)}^p\,ds\,dy\,dt\,dx
\\&\leq 
	C\int_{(3/2)Q} \int_Q \frac{\abs{u(x,t)-u(y,s)}^p}{(\abs{x-y}+\abs{t-r})^{n+1+p\theta+1}} \,ds\,dy\,dt\,dx
.\end{aligned}\end{multline*}
Adding, we see that by the bound~\eqref{eqn:besov:norm:differences},
\begin{align*}
\doublebar{u}_{\dot W(p,\theta,2)}^p
&\leq
	C\sum_{Q\in\G}\int_{(3/2)Q} \int_Q \frac{\abs{u(x,t)-u(y,r)}^p}{(\abs{x-y}+\abs{t-r})^{n+1+p\theta+1}} \,dr\,dy\,dt\,dx
\\&\leq C \doublebar{u}_{\dot B^{p,p}_{\theta+1/p}(\R^{n+1}_+)}^p
.\end{align*}
This completes the proof.
\end{proof}

The converse is more complicated. We will begin by considering functions in $\dot W(p,\theta,p)$. Notice that by Lemma~\ref{lem:local-bound:gradient}, if $\Div A\nabla u=0$ and $u\in \dot W(p,\theta,2)$, then $u\in \dot W(p,\theta,p)$ provided $p<p^+$.

\begin{lem} 
\label{lem:weighted-to-Besov}
Suppose that $0<\theta+1/p<1$, $1\leq p<\infty$ and that $u$ is any function in $\dot W(p,\theta,p)$. Then 
\[\doublebar{u}_{\dot B^{p,p}_{\theta+1/p}(\R^{n+1}_+)}\leq C(p,\theta)\doublebar{u}_{\dot W(p,\theta,p)}.\]
\end{lem}

\begin{proof}
Notice that 
\[\doublebar{u}_{\dot W(p,\theta,p)}
\approx
\int_{\R^{n+1}_+} \abs{\nabla u(x,t)}^p \, t^{p-1-p\theta}\,dx\,dt.\]

Let $\sigma=\theta+1/p$.
Recall that $\dot B^{p,p}_\sigma(\R^{n+1}_+)$ is an interpolation space between $L^p(\R^{n+1}_+)$ and $\dot W^p_1(\R^{n+1}_+)$; specifically, formula~\eqref{eqn:L-W:interpolation} tells us that
\[B^{p,p}_{\sigma}(\R^{n+1}_+)=(L^p(\R^{n+1}_+),\dot W^p_1(\R^{n+1}_+))_{\sigma,p}.\]
By the definition of real interpolation, if $u_s$ is a family of functions with $u_s\in \dot W^p_1(\R^{n+1}_+)$ and $u-u_s\in L^p(\R^{n+1}_+)$, then
\begin{align*}
\doublebar{u}_{\dot B^{p,p}_{\sigma}(\R^{n+1}_+)}
&\leq C
	\biggl(\int_0^\infty (s^{-\sigma} \{s\doublebar{u_s}_{\dot W^p_1(\R^{n+1}_+)} + \doublebar{u-u_s}_{L^p(\R^{n+1}_+)}\})^p \frac{ds}{s}\biggr)^{1/p}
.\end{align*}
We choose $u_s(x,t)=u(x,s+t)$. 
The first term may be controlled as follows:
\begin{align*}
\int_0^\infty  s^{p-1-p\sigma}\doublebar{u_s}_{\dot W^p_1(\R^{n+1}_+)}^p ds &=
	\int_0^\infty s^{p-1-p\sigma}\int_{\R^n}\int_s^\infty \abs{\nabla u(x,t)}^p\,dt\,dx\, ds
\\&=
	\frac{1}{p-p\sigma}\int_{\R^n}\int_0^\infty t^{p-p\sigma}\, \abs{\nabla u(x,t)}^p\,dt\,dx
\end{align*}
provided $\sigma<1$.
The second term may be controlled as follows:
\begin{multline*}
\int_0^\infty  s^{-1-p\sigma}\doublebar{u-u_s}_{L^p(\R^{n+1}_+)}^p ds 
\\\begin{aligned}&=
	\int_0^\infty s^{-1-p\sigma}\int_{\R^n}\int_0^\infty
	\abs{u(x,t+s)-u(x,t)}^p\,dt\,dx\, ds
\\&=
	\int_0^\infty \int_{\R^n}\int_0^\infty
	s^{-1-p\sigma}\abs[bigg]{\int_0^s \partial_t u(x,r+t)\,dr}^p\,dt\,dx\, ds.
\end{aligned}\end{multline*}
By H\"older's inequality, if $\beta$ is a real number and $p'\beta<1$, then
\begin{align*}
\abs[bigg]{\int_0^s \partial_t u(x,r+t)\,dr}^p
&\leq
	\int_0^s \abs{\partial_t u(x,r+t)}^p\,r^{p\beta}\,dr
	\biggl(\int_0^s r^{-p'\beta}\,dr\biggr)^{p/p'}
\\&=
	C(p,\beta)s^{p/p'-p\beta}\int_0^s \abs{\partial_t u(x,r+t)}^p\,r^{p\beta}\,dr
.\end{align*}

Thus, 
\begin{multline*}
\int_0^\infty  s^{-1-p\sigma}\doublebar{u-u_s}_{L^p(\R^{n+1}_+)}^p ds 
\\\begin{aligned}
&\leq
	C(p,\beta)\int_0^\infty \int_{\R^n}\int_0^\infty
	s^{p/p'-p\beta-1-p\sigma}\int_0^s \abs{\partial_t u(x,r+t)}^p\,r^{p\beta}\,dr\,dt\,dx\, ds
.\end{aligned}\end{multline*}
Interchanging the order of integration and evaluating the integral in~$s$, we see that if $p/p'-p\beta-p\sigma<0$, then
\begin{multline*}
\int_0^\infty  s^{-1-p\sigma}\doublebar{u-u_s}_{L^p(\R^{n+1}_+)}^p ds 
\\\begin{aligned}
&\leq
	C(p,\beta) \int_{\R^n}\int_0^\infty
	\int_0^\infty r^{p/p'-p\beta-p\sigma}\,r^{p\beta}\, \abs{\partial_t u(x,r+t)}^p\,dr\,dt\,dx
.\end{aligned}\end{multline*}
To ensure that $p'\beta<1$ and that $p/p'-p\beta-p\sigma<0$, we need only require $1-p'\sigma<p'\beta<1$, which is possible provided $\sigma>0$.

Making the change of variables $t\mapsto t-r  $, we have that
\begin{multline*}
\int_0^\infty  s^{-1-p\sigma}\doublebar{u-u_s}_{L^p(\R^{n+1}_+)}^p ds 
\\\begin{aligned}
&\leq
	C(p,\beta) \int_{\R^n}\int_0^\infty
	\int_0^\infty r^{p/p'-p\sigma}\, \abs{\partial_t u(x,r+t)}^p\,dr\,dt\,dx
\\&=
	C(p,\beta) \int_{\R^n}\int_0^\infty
	\int_r^\infty r^{p/p'-p\sigma}\, \abs{\partial_t u(x,t)}^p\,dt\,dr\,dx
.\end{aligned}\end{multline*}
and integrating in $r$ yields that
\[
\int_0^\infty  s^{-1-p\sigma}\doublebar{u-u_s}_{L^p(\R^{n+1}_+)}^p ds 
\leq
	C(p,\beta) \int_{\R^n}\int_0^\infty
	t^{p-p\sigma}\,\abs{\partial_t u(x,t)}^p\,dt\,dx
\]
provided $p-p\sigma>0$, that is, provided $\sigma<1$; combining with the bound on $\doublebar{u_s}_{\dot W^p_1(\R^{n+1}_+)}$, we have that
\[\doublebar{u}_{\dot B^{p,p}_{\sigma}(\R^{n+1}_+)}\leq C(p,\theta)\doublebar{u}_{\dot W(p,\sigma+1/p,p)}\]
as desired.
\end{proof}

By Lemmas~\ref{lem:Besov-to-weighted} and~\ref{lem:weighted-to-Besov}, we have that if $0<\theta<1$ and $0<1/p<1-\theta$ then
\[\frac{1}{C(p,\theta)}\doublebar{u}_{\dot W(p,\theta,2)}
\leq \doublebar{u}_{\dot B^{p,p}_{\theta+1/p}(\R^{n+1}_+)}
\leq C(p,\theta)\doublebar{u}_{\dot W(p,\theta,p)}.
\]
If $p<p^+$, then by H\"older's inequality or Lemma~\ref{lem:PDE2}, we have that $\doublebar{u}_{\dot W(p,\theta,p)}\leq C(p)\doublebar{u}_{\dot W(p,\theta,2)}$, and so the conclusion of
Theorem~\ref{thm:weighted-and-besov} is valid whenever $0<\theta<1$ and $1/p^+<1/p<1-\theta$.

We now complete the proof of Theorem~\ref{thm:weighted-and-besov}; we are left with the case $p\geq p^+$. We need only show that if $\Div A\nabla u=0$ in $\R^{n+1}_+$ and
$u\in \dot W(p,\theta,2)$ for some $(\theta,1/p)$ as in Figure~\ref{fig:besov}, then $u\in \dot B^{p,p}_{\theta+1/p}(\R^{n+1}_+)$.

Notice that such $(\theta,1/p)$ with $p\geq p^+$ satisfy the conditions of Theorem~\ref{thm:green}, and so we have that $u=-\D f+\s g$ for some $f\in \dot B^{p,p}_\theta(\R^n)$ and some $g\in \dot B^{p,p}_{\theta-1}(\R^n)$. Thus, to show that $u\in \dot B^{p,p}_{\theta+1/p}(\R^{n+1}_+)$, it suffices to show that the combined operator $\T$ of Section~\ref{sec:layers:bounded} is bounded $\dot B^{p,p}_\theta(\R^n)\mapsto \dot B^{p,p}_{\theta+1/p}(\R^{n+1}_+)$.

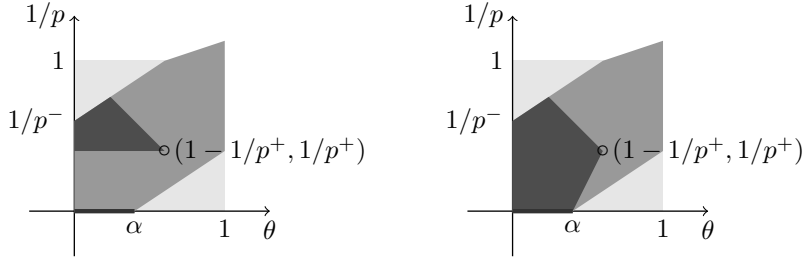
\begin{figure}
\begin{tikzpicture}[scale=2]
\figureaxes
\boundedhexagon

	\path [name path = main diagonal] (1,0)--(0,1);
	\path [name path = top boundary] (1-\alph,1) -- (0,1/2+\eps);
	\path [name intersections={of=main diagonal and top boundary, by={C}}];

\fill[well-posed] (0,1/2-\eps)
--(1/2+\eps,1/2-\eps) node [black] {$\circ$} node [black, right] {$(1-1/p^+,1/p^+)$}
--(C)
--(0,1/2+\eps) node [black,left] {$1/p^-$}
;

\draw[boundary well-posed] (0,0)--(\alph,0) node [below,black] {$\alpha$};

\end{tikzpicture}
\qquad
\begin{tikzpicture}[scale=2]
\figureaxes
\boundedhexagon

	\path [name path = main diagonal] (1,0)--(0,1);
	\path [name path = top boundary] (1-\alph,1) -- (0,1/2+\eps);
	\path [name intersections={of=main diagonal and top boundary, by={C}}];

\fill[well-posed] (0,0)--(\alph,0) node [black, below] {$\alpha$}
--(1/2+\eps,1/2-\eps) node [black] {$\circ$} node [black, right] {$(1-1/p^+,1/p^+)$}
--(C)
--(0,1/2+\eps) node [black,left] {$1/p^-$}
;

\draw[boundary well-posed] (0,0)--(\alph,0);

\end{tikzpicture}
\caption{If $(\theta,1/p)$ lies in the indicated region, then $\T$ is bounded $\dot B^{p,p}_\theta(\R^n)\mapsto \dot B^{p,p}_{\theta+1/p}(\R^{n+1}_+)$. The left-hand figure shows the result before interpolation; the right-hand figure shows the result after interpolation. The right-hand region is the intersection of the regions in Figures~\ref{fig:bounded} and~\ref{fig:besov}.}
\label{fig:besov:2}
\end{figure}

But by Theorem~\ref{thm:bounded:holder}, we have that $\T$ is bounded $\dot B^{\infty,\infty}_\theta(\R^n)\mapsto \dot B^{\infty,\infty}_{\theta}(\R^{n+1}_+)$ whenever $0<\theta<\alpha$. By Lemma~\ref{lem:weighted-to-Besov} and Theorem~\ref{thm:bounded}, we have that $\T$ is bounded 
$\dot B^{p,p}_\theta(\R^n)\mapsto \dot B^{p,p}_{\theta+1/p}(\R^{n+1}_+)$ whenever $(\theta,1/p)$ lies in the region shown on the left of Figure~\ref{fig:besov:2}; we may interpolate to extend this to the region on the right of Figure~\ref{fig:besov:2}.

\begin{rmk} In general, we cannot formulate the inhomogeneous problem in terms of Besov spaces.

The inhomogeneous Dirichlet problem for the Laplacian, studied in \cite{JerK95}, may be stated as follows. Given functions $f\in \dot B^{p,p}_{\theta}(\partial\Omega)$  and $\F\in \dot B^{p,p}_{\theta+1/p-1}(\Omega)$ for some bounded Lipschitz domain~$\Omega$, find a function $u\in \dot B^{p,p}_{\theta+1/p}(\Omega)$ with $\Delta u=\Div\F$ in $\Omega$ and $u=0$ on~$\partial\Omega$.

As described in the introduction, for general coefficients, we do not expect solutions $u$ to have more than one order of smoothness; that is, we do not expect to be able to solve problems with estimates of the form $u\in \dot B^{p,p}_{\theta+1/p}(\Omega)$ if $\theta+1/p>1$.

Consider the inhomogeneous problem in the whole space $\R^{n+1}$.
Recall that we may solve $\Div A\nabla u=\Div \F$ in $\R^{n+1}$ by setting $u=\Pi^A\F$; solvability of $\Div A\nabla u=\Div \F$ for $\F\in  \dot B^{p,p}_{\theta+1/p-1}(\R^{n+1})$ and $u\in \dot B^{p,p}_{\theta+1/p}(\R^{n+1})$ is equivalent to boundedness of the operator $\nabla\Pi^A:\dot B^{p,p}_{\theta+1/p-1}(\R^{n+1})\mapsto \dot B^{p,p}_{\theta+1/p-1}(\R^{n+1})$. We do not expect well-posedness if $\theta+1/p-1>0$.

However, note that the adjoint $(\nabla\Pi^A)^*$ to $\nabla\Pi^A$ is $\nabla\Pi^{A^*}$; thus, if we can solve $\Div A\nabla u=\Div \F$ for $\F\in  \dot B^{p,p}_{\theta+1/p-1}(\R^{n+1})$, then we can solve $\Div A^*\nabla v=\Div \vec G$ for $\vec G$ in the dual space $\dot B^{p',p'}_{-\theta-1/p+1}(\R^{n+1})$. We do not expect solvability of this problem for $-\theta-1/p+1>0$.

Thus, if $\theta+1/p\neq 1$, we do not expect solvability of $\Div A\nabla u=\Div \F$ for $\F\in  \dot B^{p,p}_{\theta+1/p-1}(\R^{n+1})$ and $u\in \dot B^{p,p}_{\theta+1/p}(\R^{n+1})$.
\end{rmk}

%% file: barton-mayboroda.bbl
\newcommand{\etalchar}[1]{$^{#1}$}
\providecommand{\bysame}{\leavevmode\hbox to3em{\hrulefill}\thinspace}
\providecommand{\MR}{\relax\ifhmode\unskip\space\fi MR }
\providecommand{\MRhref}[2]{%
  \href{http://www.ams.org/mathscinet-getitem?mr=#1}{#2}
}
\providecommand{\href}[2]{#2}
\begin{thebibliography}{GCRdF85}

\bibitem[AA11]{AusA11}
Pascal Auscher and Andreas Axelsson, \emph{Weighted maximal regularity
  estimates and solvability of non-smooth elliptic systems {I}}, Invent. Math.
  \textbf{184} (2011), no.~1, 47--115. \MR{2782252}

\bibitem[AAA{\etalchar{+}}11]{AlfAAHK11}
M.~Angeles Alfonseca, Pascal Auscher, Andreas Axelsson, Steve Hofmann, and
  Seick Kim, \emph{Analyticity of layer potentials and {$L^2$} solvability of
  boundary value problems for divergence form elliptic equations with complex
  {$L^\infty$} coefficients}, Adv. Math. \textbf{226} (2011), no.~5,
  4533--4606. \MR{2770458}

\bibitem[AAH08]{AusAH08}
Pascal Auscher, Andreas Axelsson, and Steve Hofmann, \emph{Functional calculus
  of {D}irac operators and complex perturbations of {N}eumann and {D}irichlet
  problems}, J. Funct. Anal. \textbf{255} (2008), no.~2, 374--448. \MR{2419965
  (2009h:35079)}

\bibitem[AAM10]{AusAM10}
Pascal Auscher, Andreas Axelsson, and Alan McIntosh, \emph{Solvability of
  elliptic systems with square integrable boundary data}, Ark. Mat. \textbf{48}
  (2010), no.~2, 253--287. \MR{2672609 (2011h:35070)}

\bibitem[Agr07]{Agr07}
M.~S. Agranovich, \emph{On the theory of {D}irichlet and {N}eumann problems for
  linear strongly elliptic systems with {L}ipschitz domains}, Funktsional.
  Anal. i Prilozhen. \textbf{41} (2007), no.~4, 1--21, 96, English translation:
  Funct. Anal. Appl. \textbf{41} (2007), no.~4, 247--263. \MR{2411602
  (2009b:35070)}

\bibitem[Agr09]{Agr09}
\bysame, \emph{Potential-type operators and conjugation problems for
  second-order strongly elliptic systems in domains with a {L}ipschitz
  boundary}, Funktsional. Anal. i Prilozhen. \textbf{43} (2009), no.~3, 3--25.
  \MR{2583636 (2011b:35362)}

\bibitem[AHL{\etalchar{+}}02]{AusHLMT02}
Pascal Auscher, Steve Hofmann, Michael Lacey, Alan McIntosh, and Ph.
  {Tcha\-mit\-chian}, \emph{The solution of the {K}ato square root problem for
  second order elliptic operators on {$\mathbb{R}^n$}}, Ann. of Math. (2)
  \textbf{156} (2002), no.~2, 633--654. \MR{1933726 (2004c:47096c)}

\bibitem[AKM06]{AxeKM06}
Andreas Axelsson, Stephen Keith, and Alan McIntosh, \emph{Quadratic estimates
  and functional calculi of perturbed {D}irac operators}, Invent. Math.
  \textbf{163} (2006), no.~3, 455--497. \MR{2207232 (2007k:58029)}

\bibitem[AM13]{AusM}
Pascal Auscher and Mihalis Mourgoglou, \emph{{Boundary layers, {R}ellich
  estimates and extrapolation of solvability for elliptic systems}}, ArXiv
  e-prints (2013).

\bibitem[AMM13]{AusMM13}
Pascal Auscher, Alan McIntosh, and Mihalis Mourgoglou, \emph{On {$L^2$}
  {S}olvability of {BVP}s for {E}lliptic {S}ystems}, J. Fourier Anal. Appl.
  \textbf{19} (2013), no.~3, 478--494. \MR{3048587}

\bibitem[AMT98]{AusMT98}
Pascal Auscher, Alan McIntosh, and Philippe Tchamitchian, \emph{Heat kernels of
  second order complex elliptic operators and applications}, J. Funct. Anal.
  \textbf{152} (1998), no.~1, 22--73. \MR{1600066 (99e:47062)}

\bibitem[AR12]{AusR12}
P.~{Auscher} and A.~{Ros\'en}, \emph{Weighted maximal regularity estimates and
  solvability of non-smooth elliptic systems {II}}, Analysis \& PDE \textbf{5}
  (2012), no.~5, 983--1061.

\bibitem[AT95]{AusT95}
Pascal Auscher and Philippe Tchamitchian, \emph{Calcul fontionnel pr\'ecis\'e
  pour des op\'er\-a\-teurs elliptiques complexes en dimension un (et
  applications \`a certaines \'equations elliptiques complexes en dimension
  deux)}, Ann. Inst. Fourier (Grenoble) \textbf{45} (1995), no.~3, 721--778.
  \MR{1340951 (96f:35036)}

\bibitem[AT98]{AusT98}
\bysame, \emph{Square root problem for divergence operators and related
  topics}, Ast\'erisque (1998), no.~249, viii+172. \MR{MR1651262 (2000c:47092)}

\bibitem[Aus96]{Aus96}
Pascal Auscher, \emph{Regularity theorems and heat kernel for elliptic
  operators}, J. London Math. Soc. (2) \textbf{54} (1996), no.~2, 284--296.
  \MR{1405056 (97f:35034)}

\bibitem[Axe10]{Axe10}
Andreas Axelsson, \emph{Non-unique solutions to boundary value problems for
  non-symmetric divergence form equations}, Trans. Amer. Math. Soc.
  \textbf{362} (2010), no.~2, 661--672. \MR{2551501 (2010j:35103)}

\bibitem[Bab71]{Bab70}
Ivo Babu{\v{s}}ka, \emph{Error-bounds for finite element method}, Numer. Math.
  \textbf{16} (1970/1971), 322--333. \MR{MR0288971 (44 \#6166)}

\bibitem[Bar13]{Bar13}
Ariel Barton, \emph{Elliptic partial differential equations with almost-real
  coefficients}, Mem. Amer. Math. Soc. \textbf{223} (2013), no.~1051, vii+108.

\bibitem[BL76]{BerL76}
J{\"o}ran Bergh and J{\"o}rgen L{\"o}fstr{\"o}m, \emph{Interpolation spaces.
  {A}n introduction}, Springer-Verlag, Berlin, 1976, Grundlehren der
  Mathematischen Wissenschaften, No. 223. \MR{0482275 (58 \#2349)}

\bibitem[BM13]{BarM13B}
Ariel Barton and Svitlana Mayboroda, \emph{The {D}irichlet problem for higher
  order equations in composition form}, J. Funct. Anal. \textbf{265} (2013),
  49--107.

\bibitem[CFK81]{CafFK81}
Luis~A. Caffarelli, Eugene~B. Fabes, and Carlos~E. Kenig, \emph{Completely
  singular elliptic-harmonic measures}, Indiana Univ. Math. J. \textbf{30}
  (1981), no.~6, 917--924. \MR{MR632860 (83a:35033)}

\bibitem[CGT81]{ChaGT81}
J.~A. Chao, J.~E. Gilbert, and P.~A. Tomas, \emph{Molecular decompositions in
  {$H^{p}$}-theory}, Proceedings of the {S}eminar on {H}armonic {A}nalysis
  ({P}isa, 1980), no. suppl. 1, 1981, pp.~115--119. \MR{639473 (83c:42016)}

\bibitem[Dah86]{Dah86a}
Bj{\"o}rn E.~J. Dahlberg, \emph{On the absolute continuity of elliptic
  measures}, Amer. J. Math. \textbf{108} (1986), no.~5, 1119--1138. \MR{859772
  (88i:35061)}

\bibitem[Dau88]{Dau88}
Ingrid Daubechies, \emph{Orthonormal bases of compactly supported wavelets},
  Comm. Pure Appl. Math. \textbf{41} (1988), no.~7, 909--996. \MR{951745
  (90m:42039)}

\bibitem[DG57]{DeG57}
Ennio De~Giorgi, \emph{Sulla differenziabilit\`a e l'analiticit\`a delle
  estremali degli integrali multipli regolari}, Mem. Accad. Sci. Torino. Cl.
  Sci. Fis. Mat. Nat. (3) \textbf{3} (1957), 25--43. \MR{0093649 (20 \#172)}

\bibitem[DJK84]{DahJK84}
Bj{\"o}rn E.~J. Dahlberg, David~S. Jerison, and Carlos~E. Kenig, \emph{Area
  integral estimates for elliptic differential operators with nonsmooth
  coefficients}, Ark. Mat. \textbf{22} (1984), no.~1, 97--108. \MR{735881
  (85h:35021)}

\bibitem[DK87]{DahK87}
Bj{\"o}rn E.~J. Dahlberg and Carlos~E. Kenig, \emph{Hardy spaces and the
  {N}eumann problem in {$L^p$} for {L}aplace's equation in {L}ipschitz
  domains}, Ann. of Math. (2) \textbf{125} (1987), no.~3, 437--465. \MR{890159
  (88d:35044)}

\bibitem[DKP11]{DinKP11}
Martin Dindos, Carlos Kenig, and Jill Pipher, \emph{B{MO} solvability and the
  {$A_\infty$} condition for elliptic operators}, J. Geom. Anal. \textbf{21}
  (2011), no.~1, 78--95. \MR{2755677}

\bibitem[DKPV97]{DahKPV97}
B.~E.~J. Dahlberg, C.~E. Kenig, J.~Pipher, and G.~C. Verchota, \emph{Area
  integral estimates for higher order elliptic equations and systems}, Ann.
  Inst. Fourier (Grenoble) \textbf{47} (1997), no.~5, 1425--1461. \MR{1600375
  (98m:35045)}

\bibitem[DPP07]{DinPP07}
Martin Dindos, Stefanie Petermichl, and Jill Pipher, \emph{The {$L^p$}
  {D}irichlet problem for second order elliptic operators and a {$p$}-adapted
  square function}, J. Funct. Anal. \textbf{249} (2007), no.~2, 372--392.
  \MR{2345337 (2008f:35108)}

\bibitem[DR10]{DinR10}
Martin Dindo{\v{s}} and David~J. Rule, \emph{Elliptic equations in the plane
  satisfying a {C}arleson measure condition}, Rev. Mat. Iberoam. \textbf{26}
  (2010), no.~3, 1013--1034. \MR{2789374}

\bibitem[Eva98]{Eva98}
Lawrence~C. Evans, \emph{Partial differential equations}, Graduate Studies in
  Mathematics, vol.~19, American Mathematical Society, Providence, RI, 1998.
  \MR{MR1625845 (99e:35001)}

\bibitem[Fef89]{Fef89}
R.~Fefferman, \emph{A criterion for the absolute continuity of the harmonic
  measure associated with an elliptic operator}, J. Amer. Math. Soc. \textbf{2}
  (1989), no.~1, 127--135. \MR{955604 (90b:35068)}

\bibitem[Fef93]{Fef93}
Robert~A. Fefferman, \emph{Large perturbations of elliptic operators and the
  solvability of the {$L^p$} {D}irichlet problem}, J. Funct. Anal. \textbf{118}
  (1993), no.~2, 477--510. \MR{1250271 (94k:35082)}

\bibitem[FJ85]{FraJ85}
Michael Frazier and Bj{\"o}rn Jawerth, \emph{Decomposition of {B}esov spaces},
  Indiana Univ. Math. J. \textbf{34} (1985), no.~4, 777--799. \MR{808825
  (87h:46083)}

\bibitem[FJ90]{FraJ90}
\bysame, \emph{A discrete transform and decompositions of distribution spaces},
  J. Funct. Anal. \textbf{93} (1990), no.~1, 34--170. \MR{1070037 (92a:46042)}

\bibitem[FJK84]{FabJK84}
Eugene~B. Fabes, David~S. Jerison, and Carlos~E. Kenig, \emph{Necessary and
  sufficient conditions for absolute continuity of elliptic-harmonic measure},
  Ann. of Math. (2) \textbf{119} (1984), no.~1, 121--141. \MR{736563
  (85h:35069)}

\bibitem[FKP91]{FefKP91}
R.~A. Fefferman, C.~E. Kenig, and J.~Pipher, \emph{The theory of weights and
  the {D}irichlet problem for elliptic equations}, Ann. of Math. (2)
  \textbf{134} (1991), no.~1, 65--124. \MR{MR1114608 (93h:31010)}

\bibitem[FMM98]{FabMM98}
Eugene Fabes, Osvaldo Mendez, and Marius Mitrea, \emph{Boundary layers on
  {S}obolev-{B}esov spaces and {P}oisson's equation for the {L}aplacian in
  {L}ipschitz domains}, J. Funct. Anal. \textbf{159} (1998), no.~2, 323--368.
  \MR{1658089 (99j:35036)}

\bibitem[Fre08]{Fre08}
Jens Frehse, \emph{An irregular complex valued solution to a scalar uniformly
  elliptic equation}, Calc. Var. Partial Differential Equations \textbf{33}
  (2008), no.~3, 263--266. \MR{2429531 (2009h:35084)}

\bibitem[FS72]{FefS72}
C.~Fefferman and E.~M. Stein, \emph{{$H^{p}$} spaces of several variables},
  Acta Math. \textbf{129} (1972), no.~3-4, 137--193. \MR{MR0447953 (56 \#6263)}

\bibitem[GCRdF85]{GarR85}
Jos{\'e} Garc{\'{\i}}a-Cuerva and Jos{\'e}~L. Rubio~de Francia, \emph{Weighted
  norm inequalities and related topics}, North-Holland Mathematics Studies,
  vol. 116, North-Holland Publishing Co., Amsterdam, 1985, Notas de
  Matem{\'a}tica [Mathematical Notes], 104. \MR{807149 (87d:42023)}

\bibitem[GdlH12]{Gra12}
Ana Grau de~la Herran, \emph{Local ${T}b$ {T}heorems for {S}quare functions and
  generalizations}, Ph.D. thesis, University of Missouri, Columbia, Missouri,
  2012.

\bibitem[Gri85]{Gri85}
P.~Grisvard, \emph{Elliptic problems in nonsmooth domains}, Monographs and
  Studies in Mathematics, vol.~24, Pitman (Advanced Publishing Program),
  Boston, MA, 1985. \MR{775683 (86m:35044)}

\bibitem[HK07]{HofK07}
Steve Hofmann and Seick Kim, \emph{The {G}reen function estimates for strongly
  elliptic systems of second order}, Manuscripta Math. \textbf{124} (2007),
  no.~2, 139--172. \MR{MR2341783 (2008k:35110)}

\bibitem[HKMP12]{HofKMP12}
Steve Hofmann, Carlos Kenig, Svitlana Mayboroda, and Jill Pipher, \emph{Square
  function\slash non-tangential maximal function estimates and the {D}irichlet
  problem for non-symmetric elliptic operators}, ArXiv e-prints (2012).

\bibitem[HKMP13]{HofKMP13}
\bysame, \emph{The regularity problem for second order elliptic operators with
  complex-valued bounded measurable coefficients}, ArXiv e-prints (2013).

\bibitem[HMMa]{HofMayMou}
Steve Hofmann, Svitlana Mayboroda, and Mihalis Mourgoglou, \emph{${L}^p$ and
  endpoint solvability results for divergence form elliptic equations with
  complex ${L}^\infty$ coefficients}, preprint.

\bibitem[HMMb]{HofMitMor}
Steve Hofmann, Marius Mitrea, and Andrew Morris, \emph{The method of layer
  potentials in ${L}^p$ and endpoint spaces for elliptic operators with
  ${L}^\infty$ coefficients}, preprint.

\bibitem[Jaw77]{Jaw77}
Bj{\"o}rn Jawerth, \emph{Some observations on {B}esov and {L}izorkin-{T}riebel
  spaces}, Math. Scand. \textbf{40} (1977), no.~1, 94--104. \MR{0454618 (56
  \#12867)}

\bibitem[JK81]{JerK81a}
David~S. Jerison and Carlos~E. Kenig, \emph{The {D}irichlet problem in
  nonsmooth domains}, Ann. of Math. (2) \textbf{113} (1981), no.~2, 367--382.
  \MR{MR607897 (84j:35076)}

\bibitem[JK95]{JerK95}
David Jerison and Carlos~E. Kenig, \emph{The inhomogeneous {D}irichlet problem
  in {L}ipschitz domains}, J. Funct. Anal. \textbf{130} (1995), no.~1,
  161--219. \MR{1331981 (96b:35042)}

\bibitem[JW84]{JonW84}
Alf Jonsson and Hans Wallin, \emph{Function spaces on subsets of {${\bf
  R}^n$}}, Math. Rep. \textbf{2} (1984), no.~1, xiv+221. \MR{820626
  (87f:46056)}

\bibitem[Ken94]{Ken94}
Carlos~E. Kenig, \emph{Harmonic analysis techniques for second order elliptic
  boundary value problems}, CBMS Regional Conference Series in Mathematics,
  vol.~83, Published for the Conference Board of the Mathematical Sciences,
  Washington, DC, 1994. \MR{MR1282720 (96a:35040)}

\bibitem[Kim07]{Kim07}
Doyoon Kim, \emph{Trace theorems for {S}obolev-{S}lobodeckij spaces with or
  without weights}, J. Funct. Spaces Appl. \textbf{5} (2007), no.~3, 243--268.
  \MR{2352844 (2008k:46100)}

\bibitem[KKPT00]{KenKPT00}
C.~Kenig, H.~Koch, J.~Pipher, and T.~Toro, \emph{A new approach to absolute
  continuity of elliptic measure, with applications to non-symmetric
  equations}, Adv. Math. \textbf{153} (2000), no.~2, 231--298. \MR{MR1770930
  (2002f:35071)}

\bibitem[KM98]{KalM98}
Nigel Kalton and Marius Mitrea, \emph{Stability results on interpolation scales
  of quasi-{B}anach spaces and applications}, Trans. Amer. Math. Soc.
  \textbf{350} (1998), no.~10, 3903--3922. \MR{1443193 (98m:46094)}

\bibitem[KMM07]{KalMM07}
Nigel Kalton, Svitlana Mayboroda, and Marius Mitrea, \emph{Interpolation of
  {H}ardy-{S}obolev-{B}esov-{T}riebel-{L}izorkin spaces and applications to
  problems in partial differential equations}, Interpolation theory and
  applications, Contemp. Math., vol. 445, Amer. Math. Soc., Providence, RI,
  2007, pp.~121--177. \MR{MR2381891 (2009g:46031)}

\bibitem[KMR01]{KozMR01}
V.~A. Kozlov, V.~G. Maz'ya, and J.~Rossmann, \emph{Spectral problems associated
  with corner singularities of solutions to elliptic equations}, Mathematical
  Surveys and Monographs, vol.~85, American Mathematical Society, Providence,
  RI, 2001. \MR{1788991 (2001i:35069)}

\bibitem[KO83]{KonO83}
V.A. Kondrat'ev and O.A. Olejnik, \emph{{Boundary-value problems for partial
  differential equations in non-smooth domains.}}, {Russ. Math. Surv.}
  \textbf{38} (1983), no.~2, 1--86 (English).

\bibitem[KP93]{KenP93}
Carlos~E. Kenig and Jill Pipher, \emph{The {N}eumann problem for elliptic
  equations with nonsmooth coefficients}, Invent. Math. \textbf{113} (1993),
  no.~3, 447--509. \MR{MR1231834 (95b:35046)}

\bibitem[KP95]{KenP95}
\bysame, \emph{The {N}eumann problem for elliptic equations with nonsmooth
  coefficients. {II}}, Duke Math. J. \textbf{81} (1995), no.~1, 227--250
  (1996), A celebration of John F. Nash, Jr. \MR{1381976 (97j:35021)}

\bibitem[KR09]{KenR09}
Carlos~E. Kenig and David~J. Rule, \emph{The regularity and {N}eumann problem
  for non-symmetric elliptic operators}, Trans. Amer. Math. Soc. \textbf{361}
  (2009), no.~1, 125--160. \MR{MR2439401 (2009k:35050)}

\bibitem[Kyr03]{Kyr03}
G.~Kyriazis, \emph{Decomposition systems for function spaces}, Studia Math.
  \textbf{157} (2003), no.~2, 133--169. \MR{1981430 (2004e:42036)}

\bibitem[Liz60]{Liz60}
P.~I. Lizorkin, \emph{Boundary properties of functions from ``weight''
  classes}, Soviet Math. Dokl. \textbf{1} (1960), 589--593. \MR{0123103 (23
  \#A434)}

\bibitem[May05]{May05}
Svitlana Mayboroda, \emph{The {P}oisson problem on {L}ipschitz domains},
  ProQuest LLC, Ann Arbor, MI, 2005, Thesis (Ph.D.)--University of Missouri -
  Columbia. \MR{2709638}

\bibitem[May10]{May10a}
\bysame, \emph{The connections between {D}irichlet, regularity and {N}eumann
  problems for second order elliptic operators with complex bounded measurable
  coefficients}, Adv. Math. \textbf{225} (2010), no.~4, 1786--1819.
  \MR{2680190}

\bibitem[Mey63]{Mey63}
Norman~G. Meyers, \emph{An {$L^{p}$}e-estimate for the gradient of solutions of
  second order elliptic divergence equations}, Ann. Scuola Norm. Sup. Pisa (3)
  \textbf{17} (1963), 189--206. \MR{0159110 (28 \#2328)}

\bibitem[Mey64]{Mey64}
\bysame, \emph{Mean oscillation over cubes and {H}\"older continuity}, Proc.
  Amer. Math. Soc. \textbf{15} (1964), 717--721. \MR{0168712 (29 \#5969)}

\bibitem[Mit08]{Mit08}
Dorina Mitrea, \emph{A generalization of {D}ahlberg's theorem concerning the
  regularity of harmonic {G}reen potentials}, Trans. Amer. Math. Soc.
  \textbf{360} (2008), no.~7, 3771--3793. \MR{2386245 (2009b:35046)}

\bibitem[MM00]{MenM00}
Osvaldo Mendez and Marius Mitrea, \emph{The {B}anach envelopes of {B}esov and
  {T}riebel-{L}izorkin spaces and applications to partial differential
  equations}, J. Fourier Anal. Appl. \textbf{6} (2000), no.~5, 503--531.
  \MR{1781091 (2001k:46056)}

\bibitem[MM04]{MayMit04a}
Svitlana Mayboroda and Marius Mitrea, \emph{Sharp estimates for {G}reen
  potentials on non-smooth domains}, Math. Res. Lett. \textbf{11} (2004),
  no.~4, 481--492. \MR{MR2092902 (2005i:35059)}

\bibitem[MM11]{MitM11}
Dorina Mitrea and Irina Mitrea, \emph{On the regularity of {G}reen functions in
  {L}ipschitz domains}, Comm. Partial Differential Equations \textbf{36}
  (2011), no.~2, 304--327. \MR{2763343}

\bibitem[MMS10]{MazMS10}
V.~Maz'ya, M.~Mitrea, and T.~Shaposhnikova, \emph{The {D}irichlet problem in
  {L}ipschitz domains for higher order elliptic systems with rough
  coefficients}, J. Anal. Math. \textbf{110} (2010), 167--239. \MR{2753293
  (2011m:35088)}

\bibitem[Mor66]{Mor66}
Charles~B. Morrey, Jr., \emph{Multiple integrals in the calculus of
  variations}, Die Grund\-leh\-ren der mathematischen Wissenschaften, Band 130,
  Springer-Verlag New York, Inc., New York, 1966. \MR{0202511 (34 \#2380)}

\bibitem[Nas58]{Nas58}
J.~Nash, \emph{Continuity of solutions of parabolic and elliptic equations},
  Amer. J. Math. \textbf{80} (1958), 931--954. \MR{0100158 (20 \#6592)}

\bibitem[Nik77]{Nik77b}
S.~M. Nikol'ski{\u\i}, \emph{Priblizhenie funktsii mnogikh peremennykh i
  teoremy vlozheniya [{A}pproximation of functions of several variables and
  imbedding theorems]}, ``Nauka'', Moscow, 1977, Second edition, revised and
  supplemented. \MR{506247 (81f:46046)}

\bibitem[NLM88]{NikLM88}
S.~M. Nikol'ski{\u\i}, P.~I. Lizorkin, and N.~V. Miroshin, \emph{Weighted
  function spaces and their applications to the investigation of boundary value
  problems for degenerate elliptic equations}, Izv. Vyssh. Uchebn. Zaved. Mat.
  (1988), no.~8, 4--30. \MR{971868 (90b:35104)}

\bibitem[Pee71]{Pee71}
Jaak Peetre, \emph{Interpolation functors and {B}anach couples}, Actes du
  {C}ongr\`es {I}nternational des {M}ath\'ematiciens ({N}ice, 1970), {T}ome 2,
  Gauthier-Villars, Paris, 1971, pp.~373--378. \MR{0425636 (54 \#13590)}

\bibitem[{Ros}12]{Ros12A}
Andreas {Ros\'en}, \emph{{Layer potentials beyond singular integral
  operators}}, ArXiv e-prints (2012).

\bibitem[RS96]{RunS96}
Thomas Runst and Winfried Sickel, \emph{Sobolev spaces of fractional order,
  {N}emytskij operators, and nonlinear partial differential equations}, de
  Gruyter Series in Nonlinear Analysis and Applications, vol.~3, Walter de
  Gruyter \& Co., Berlin, 1996. \MR{1419319 (98a:47071)}

\bibitem[Rul07]{Rul07}
David~J. Rule, \emph{Non-symmetric elliptic operators on bounded {L}ipschitz
  domains in the plane}, Electron. J. Differential Equations (2007), No. 144,
  8. \MR{MR2366037 (2008m:35070)}

\bibitem[Sha85]{Sha85}
V.~V. Shan'kov, \emph{The averaging operator with variable radius, and the
  inverse trace theorem}, Sibirsk. Mat. Zh. \textbf{26} (1985), no.~6,
  141--152, 191. \MR{816512 (87d:47044)}

\bibitem[{\v{S}}ne74]{Sne74}
I.~Ja. {\v{S}}ne{\u\i}berg, \emph{Spectral properties of linear operators in
  interpolation families of {B}anach spaces}, Mat. Issled. \textbf{9} (1974),
  no.~2(32), 214--229, 254--255. \MR{0634681 (58 \#30362)}

\bibitem[Ste93]{Ste93}
Elias~M. Stein, \emph{Harmonic analysis: real-variable methods, orthogonality,
  and oscillatory integrals}, Princeton Mathematical Series, vol.~43, Princeton
  University Press, Princeton, NJ, 1993, With the assistance of Timothy S.
  Murphy, Monographs in Harmonic Analysis, III. \MR{MR1232192 (95c:42002)}

\bibitem[Tri78]{Tri78}
Hans Triebel, \emph{Interpolation theory, function spaces, differential
  operators}, North-Holland Mathematical Library, vol.~18, North-Holland
  Publishing Co., Amsterdam, 1978. \MR{503903 (80i:46032b)}

\bibitem[Tri83]{Tri83}
\bysame, \emph{Theory of function spaces}, Monographs in Mathematics, vol.~78,
  Birkh\"auser Verlag, Basel, 1983. \MR{781540 (86j:46026)}

\bibitem[TW80]{TaiW80}
Mitchell~H. Taibleson and Guido Weiss, \emph{The molecular characterization of
  certain {H}ardy spaces}, Representation theorems for {H}ardy spaces,
  Ast\'erisque, vol.~77, Soc. Math. France, Paris, 1980, pp.~67--149.
  \MR{604370 (83g:42012)}

\bibitem[Usp61]{Usp61}
S.~V. Uspenski{\u\i}, \emph{Imbedding theorems for classes with weights}, Trudy
  Mat. Inst. Steklov. \textbf{60} (1961), 282--303. \MR{0136980 (25 \#440)}

\bibitem[Ver84]{Ver84}
Gregory Verchota, \emph{Layer potentials and regularity for the {D}irichlet
  problem for {L}aplace's equation in {L}ipschitz domains}, J. Funct. Anal.
  \textbf{59} (1984), no.~3, 572--611. \MR{MR769382 (86e:35038)}

\bibitem[Zan00]{Zan00}
Daniel~Z. Zanger, \emph{The inhomogeneous {N}eumann problem in {L}ipschitz
  domains}, Comm. Partial Differential Equations \textbf{25} (2000), no.~9-10,
  1771--1808. \MR{1778780 (2001g:35056)}

\end{thebibliography}
